دانشگاه تهران

دانشکده فنی

گروه مهندسی برق و کامپیوتر

عنوان:

# بررسی عملکرد موتور القائی سه فاز قفسه ای در حالت ناهم محوری بین روتور و استاتور و سایر خطاها

نگارش:

ایمان طباطبایی اردکانی

استاد راهنما:

دکتر جواد فیض

استاد مشاور:

دکتر حمید ابوالحسنی تولیت

دانشگاه تهران

دانشکده فنی

گروه مهندسی برق و کامپیوتر

عنوان:

# بررسی عملکرد موتور القائی سه فاز قفسه ای در حالت ناهم محوری بین روتور و استاتور و سایر خطاها

نگارش:

ایمان طباطبایی اردکانی

از پایاننامه حاضر در تاریخ 81/12/18 در مقابل هیات داوران دفاع بعمل آمد و مورد تصویب قرار گرفت.

سرپرست تحصیلات تکمیلی دانشکده فنی: دکتر جواد فیض

مدیر گروه آموزشی برق و کامپیوتر: دکتر پرویز جبه دار مارالانی

سرپرست تحصیلات تکمیلی گروه برق و کامپیوتر: دکتر حمیدرضا جمالی

استاد راهنما: دکتر جواد فیض

عضو هیات داوران: دکتر هاشم اورعی

عضو هیات داوران: دکتر حمید لسانی

عضو هیات داوران: دکتر سعید افشارنیا

تقدیم به همه آنانی که از آنان

## چگونه زیستن

و

## چگونه اندیشیدن

را آموختم.

تقدیم به پدر و مادر عزیزم

که شمع وجودشان همواره

روشنی بخش

و

گرمی بخش

زندگی ام بوده است.

# تشکر و قدردانی

بدینوسیله از رهنمودهای *جناب آقای پروفسور جواد فیض* استاد راهنما و *جناب آقای دکتر حمید ابوالحسنی تولیت* استاد مشاور این پایاننامه و *جناب آقای دکتر همایون مشکین کلک* که که در ایجاد محیطی پژوهشی و پویا برای انجام این پروژه نقش مهمی داشته‌اند کمال تشکر و قدردانی به عمل می‌آورم. همچنین از *جناب آقای مهندس عماد شریفی قزوینی* مدیر گروه ماشینهای الکتریکی پژوهشگاه نیرو، *جناب آقای مهندس مهدی میرباقری*، *جناب آقای مهندس غلامرضا سلطانی خسروشاهی*، *جناب آقای مهندس حسام میر عابدینی* و سرکار *خانم مهندس الهام خسروشاهلی* کارشناسان ارشد گروه ماشینهای الکتریکی پژوهشگاه نیرو و *خانم لیلی خلقی* تشکر و قدردانی می‌نمایم.

از *سرکار خانم گیتی جلالی* و *سرکار خانم مینو غفوری* متصدیان کتابخانه گروه مهندسی برق و کامپیوتر دانشکده فنی که نه تنها در طول مدت انجام این پروژه بلکه در طول مدت تحصیل اینجانب مرا در تهیه منابع و مراجع یاری نموده‌اند، کمال تشکر و قدردانی به عمل می‌آورم.

به این امید که در آینده بتوانم گوشه‌ای از محبتهای این بزرگواران را جبران کنم.

# چکیده

پایاننامه حاضر، در آغاز به معرفی عیبهای متداول موتور القائی سه فاز قفسه ای و بررسی روشهای عمده مدلسازی این موتور تحت شرایط مختلف ناهم‌محوری بین روتور و استاتور و شکستگی میله‌های روتور می‌پردازد. سپس ضمن بازنگری نظریه تابع سیم‌پیچی اشتباه موجود در این نظریه تصحیح و برای اولین بار نشان داده شده که در حالت ناهم‌محوری میان روتور و استاتور نیز برای هر دو مدار دلخواهی در موتور باز هم تساوی $L_{12} = L_{21}$ صادق است. همچنین به محاسبه دقیق شرایط هندسی روتور در حالتهای مختلف ناهم‌محوری میان روتور و استاتور پرداخته می‌شود. در ادامه برای نخستین بار به محاسبه تحلیلی کلیه اندوکتانسهای مغناطیسی خودی و متقابل موتور پرداخته می‌شود و روابط جبری بسته‌ای جهت محاسبه این اندوکتانسها ارائه می‌شود. این روابط جبری از ساده‌سازی عبارتهای انتگرالی پیچیده‌ای بدست آمده‌اند.

در ادامه نشان داده می‌شود که اندوکتانسهای استاتور در حالت ناهم‌محوری ایستا مستقل از زاویه مکانیکی روتور و اندوکتانسهای روتور تابعی از زاویه مکانیکی روتور است. در حالت ناهم‌محوری پویا این وضعیت برعکس است. اما در حالت ناهم‌محوری مرکب کلیه اندوکتانسهای موتور تابعی از موقعیت زاویه‌ای روتور است. بنابراین جهت شبیه‌سازی حالت ناهم‌محوری مرکب ذخیره‌سازی و درونیابی این اطلاعات به مراتب طولانی‌تر از حالت ناهم‌محوری ایستا و پویاست. روش ارائه شده در این پایاننامه منجر به ایجاد توابع تحلیلی ماتریسی جهت محاسبه ماتریسهای اندوکتانس موتور می‌شود که با داشتن زاویه مکانیکی روتور و درجه ناهم‌محوری ایستا و پویا مقدار اندوکتانس را محاسبه می‌کند بنابراین نیازی به ذخیره‌سازی و درونیابی داده‌ها نیست.

و در انتها به شبیه‌سازی کامپیوتری موتور تحت شرایط مختلف ناهم‌محوری بین روتور و استاتور و شکستگی میله‌های روتور  پرداخته می‌شود همچنین طیف فرکانسی جریان خط موتور در حالتهای مختلف تا فرکانسهای حدود 20 برابر فرکانس اصلی محاسبه و با هم مقایسه می‌شوند.

# فهرست

فصل اول

# عیب های متداول در

# موتورهای القائی سه‌فاز قفسه‌ای

امروزه موتورهای القایی به عنوان متداول‌ترین موتور سرعت متغیر در صنایع کوچک و بزرگ بکار می‌روند. عملکرد این موتورها در شرایط وقوع عیب علاوه بر اینکه موجب اختلال در کارآیی آنها می‌شود، باعث کاهش طول عمر آنها نیز خواهد شد. بنابراین تشخیص و آشکارسازی عیب در موتورهای القایی علاوه بر بهبود کارایی آنها به طول عمر مفید آنها نیز می‌افزاید. این مساله درمورد موتورهای الکتریکی بزرگ که درنیروگاهها، پالایشگاهها، صنایع پتروشیمی و غیره بکار می‌روند بسیار حائز اهمیت است.

موتورهای القایی در انواع مختلفی ساخته می‌شوند ولی بطورکلی ساختار همه آنها متشکل از سه قسمت اصلی استاتور، روتور و مکانیکی است. عیبهای احتمالی موتورهای القایی نیز براساس تقسیم‌بندی مزبور دسته‌بندی می‌شوند. در ادامه به معرفی عیوب متداول قسمتهای مختلف یک موتور القایی متوسط پرداخته می‌شود.

## 1-1- عیبهای متداول در استاتور

در موتورهای القایی استاتور مجموعه‌ای از سیم پیچهای چند فازی است که تولید میدان مغناطیسی گردان می‌کند. سیم پیچی هر فاز استاتور شامل چند گروه موازی از کلافهاست. کلاف‌های هر گروه نیز از اتصال موازی یا سری کلافهای متعددی تشکیل می‌شود. کلافها نیز از اتصال سری یا موازی حلقه های سیم پیچی ایجاد می‌شوند. عایق‌بندی این سیم پیچها عامل بروز عیوب متعددی هستند که در ادامه به مهمترین آنها اشاره خواهد شد.

### 1-1-1- اتصال کوتاه حلقه به حلقه

در این حالت دو یا چند حلقه از یک کلاف سیم پیچهای استاتور اتصال کوتاه شده و جریان بیش از حد مجاز از حلقه های اتصال کوتاه می‌گذرد. چنین جریانی می‌تواند که می‌تواند منجر به گرم شدن و حتی ذوب شدن تدریجی عایق سیم پیچهای استاتور شود. بنابراین عدم تشخیص چنین عیبی کوچکی می‌تواند در دراز مدت موجب اتصال کوتاه قسمت بیشتری از سیم پیچی استاتور و سرانجام سوختن آن شود. در فازی که دچار اتصال کوتاه حلقه به حلقه شده استMMF سیم پیچی فاز تغییر می‌کند و اندوکتانس متقابل آن فاز و همه مدارهای دیگر در ماشین تغییر می‌کند. بعلاوه یک فاز جدید که فاز اتصال کوتاه شده D نامیده می‌شود به مجموعه سیم بندی ماشین اضافه

می شود. اگرچه این فاز تماس رسانایی با فازهای دیگر ندارد اما با مدارهای روتور و استاتور تزویج مغناطیسی دارد. این فاز بصورت یک مدار در مدل ماشین ظاهر می شود.

### 1-1-2- اتصال کوتاه کلاف به کلاف

در صورتیکه دو یا چند کلاف از سیم پیچهای یک فاز استاتور اتصال کوتاه شوند، اتصال کوتاه کلاف به کلاف رخ می دهد. این عیب می تواند نتیجه پیشرفت تدریجی اتصال کوتاه حلقه به حلقه باشد. در صورت وقوع این عیب، شدت عبور جریان غیر مجاز سیم پیچی های استاتور بیشتر از جریان اتصال کوتاه حلقه به حلقه می‌شود. شرایط ایجاد شده توسط این عیب نیز مانند اتصال کوتاه حلقه به حلقه است.

### 1-1-3- اتصال کوتاه فاز به فاز

اتصال کوتاه فاز به فاز حالت کلی تری از اتصال کوتاه کلاف به کلاف است که در آن دو یا چند فاز استاتور اتصال کوتاه می‌شوند. این عیب به مراتب شدیدتر از اتصال کوتاه حلقه به حلقه و کلاف به کلاف است. در این حالت فازهای اتصال کوتاه شده هر کدام تبدیل به دو فاز می شوند. یعنی در صورات اتصال کوتاه دو فاز، مجموعا 5 مدار استاتور وجود خواهد داشت که برای هر کدام بایستی معادلات الکترومغناطیسی و الکتریکی جداگانه نوشت.

### 1-1-4- قطع سیم پیچهای استاتور

بهنگام قطع یک یا چند حلقه از سیم پیچهای استاتور توزیع متقارن سیم پیچهای استاتور به هم می‌خورد. بنابـراین هارمونیکهای اضافی در میدان مغناطیسی استاتور ایجاد می‌شود. این امر منجر به افزایش تموج[1] گشتاور و نیز سرعت ماشین می‌شود. در این حالت بر حسب محل و شدت این عیب ممکن است یک یا چند سیم‌پیچی از مجموعه موتور حذف شود و یا توزیع سیم پیچی ها تغییر کند. احتمال وقوع این عیب در استاتور از سایر عیبهای آن کمتر است.

---

[1] Ripple

### 5-1-1- اتصال کوتاه فاز به زمین

این عیب شدیدترین عیب ممکن در سیم پیچهای استاتور است که می‌تواند منجر به از بین رفتن سریع عایق سیم پیچهای استاتور و سوختن آنها شود. تعمیر سیم پیچی سوخته حتی در صورا امکان نیز بسیار دشوار است. این عیب می‌تواند نتیجه پیشرفت و عدم جلوگیری از عیبهای کوچک تری مانند اتصال کوتاه حلقه به حلقه و یا کلاف به کلاف باشد.

### 6-1-1- نارسائیهای عایقی

علت اصلی بروز مشکلات عایقی سهل انگاری در پیچیدن سیم پیچهای استاتور و یا وجود حبابها و حفره‌هایی است که معمولا در هنگام ساخت عایق پدید می‌آید. چون ساختمان عایقی سیم پیچهای استاتور لایه لایه است حذف این حفره‌ها امکان‌پذیر نیست. علاوه بر این با بهره برداری از موتور، عایق های آن تحت تاثیر عوامل مختلفی فرسوده می‌شوند.

درصورت عدم ممانعت بهنگام از این عیب احتمال وقوع سایر عیبهای استاتور بیشتر می‌شود. در بعضی از تقسیم بندیها کلیه عیبهای موجود در سیم پیچی استاتور تحت عنوان نارسائیهاهی عایقی سیم پیچی استاتور طبقه‌بندی می‌شوند. در شکل 1-1 انواع عیبهای متداول استاتور موتورهای الکتریکی و رابطه بین آنها نشان داده شده است. دهد. مطابق این شکل عدم جلوگیری از یک عیب ساده و پیشرفت آن ممکن است منجر به عیبهای شدیدتری مانند اتصال کوتاه فاز به زمین در سیم پیچی استاتور گردد.

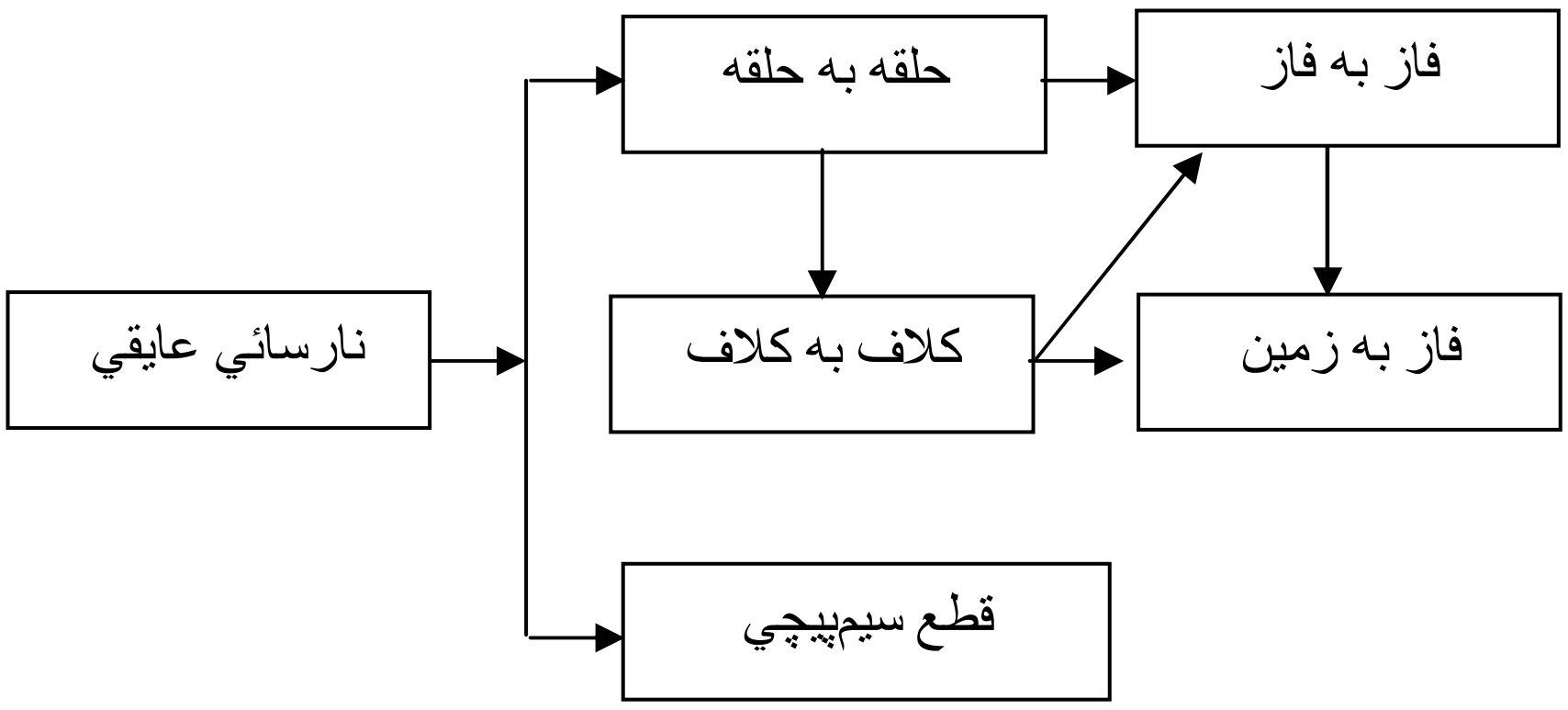

شکل 1-1: انواع عیبهای استاتور موتورهای القائی و رابطه بین آنها

## 1-2- شرایط عیب در روتور

در موتورهای الکتریکی روتور قسمت دواری است که بر هم کنش میدان مغناطیسی آن با میدان مغناطیسی استاتور موجب چرخش آن می‌شود. روتور موتورهای الکتریکی را در انواع مختلفی مانند روتور مغناطیس دائم سطحی، روتور مغناطیس دائم داخلی، روتور رلوکتانسی، روتور سیم پیچی شده و روتور قفسی می‌سازند. در این روتورهای قفسی موتورهای القایی گسترده‌ترین کاربرد را در صنعت دارند.

در ماشینهای بزرگ القایی، روتور قفسی را از میله‌های مسی می‌سازند اما در موتورهای القایی کوچک قفس آلومینیومی را درون یک استوانه آهنی تزریق می‌کنند تا به استحکام روتور بیافزایند. گاهی دو قفس درون روتور جا می‌دهند که به این روتورها روتور قفس مضاعف می‌گویند. همچنین به منظور کاهش هارمونیک‌های اضافی جریان، میله های روتور قفسی را بصورت مورب قرار می‌دهند.

### 1-2-1- شکستگی (ترک خوردگی) میله های روتور :

حدود 5 تا 10 درصد عیبهای متداول در موتورهای القایی روتور قفسی به دلیل شکستگی و یا ترک خوردگی میله‌های روتور و یا حلقه‌های انتهای آن است [1]. اغلب، دلایل زیر موجب شکستگی و یا ترک خوردگی میله‌های روتور می‌شوند [2]:

- ✓ تنش‌های حرارتی که در اثر اضافه بار، توزیع غیریکنواخت حرارت، نقاط داغ و جرقه (بیشتر در روتورهای شمشی) ایجاد می‌شود.
- ✓ تنش‌های مغناطیسی که به علت نیروهای الکترومغناطیسی، نیروهای نامتقارن مغناطیسی و نویزها و ارتعاشات الکترومغناطیسی ایجاد می‌شوند.
- ✓ تنش‌های پسماند  مرحله ساخت
- ✓ تنش‌های دینامیکی ناشی از گشتاور محور و نیروهای گریز از مرکز
- ✓ تنش‌های محیطی در اثر آلودگی و ساییدگی مواد روتور توسط مواد شیمیایی و رطوبت
- ✓ تنش‌های مکانیکی ناشی از خستگی مکانیکی قسمت‌های مختلف، خرابی بلبرینگ و غیره

اگر موتوری به طور مناسب طراحی، ساخته، نصب و نگهداری شود، این تنش‌ها در حد مجاز بوده و اثر تخریبی ندارد[3].

### 1-2-2- شکستگی (ترک خوردگی) حلقه انتهایی

در این حالت اتصال بین یک یا چند میله قفس روتور با حلقه انتهایی قفس یا خود میله انتهایی قفس ترک می‌خورد و یا می‌شکند [4]. همه دلایلی که منجر به ترک خوردگی و یا شکستگی میله های روتور می‌شوند، می‌توانند منجر به این نوع عیب نیز گردند. شکل 1-2 انواع عیبهای روتور در موتورهای الکتریکی و رابطه بروز آنها را نشان می‌دهد.

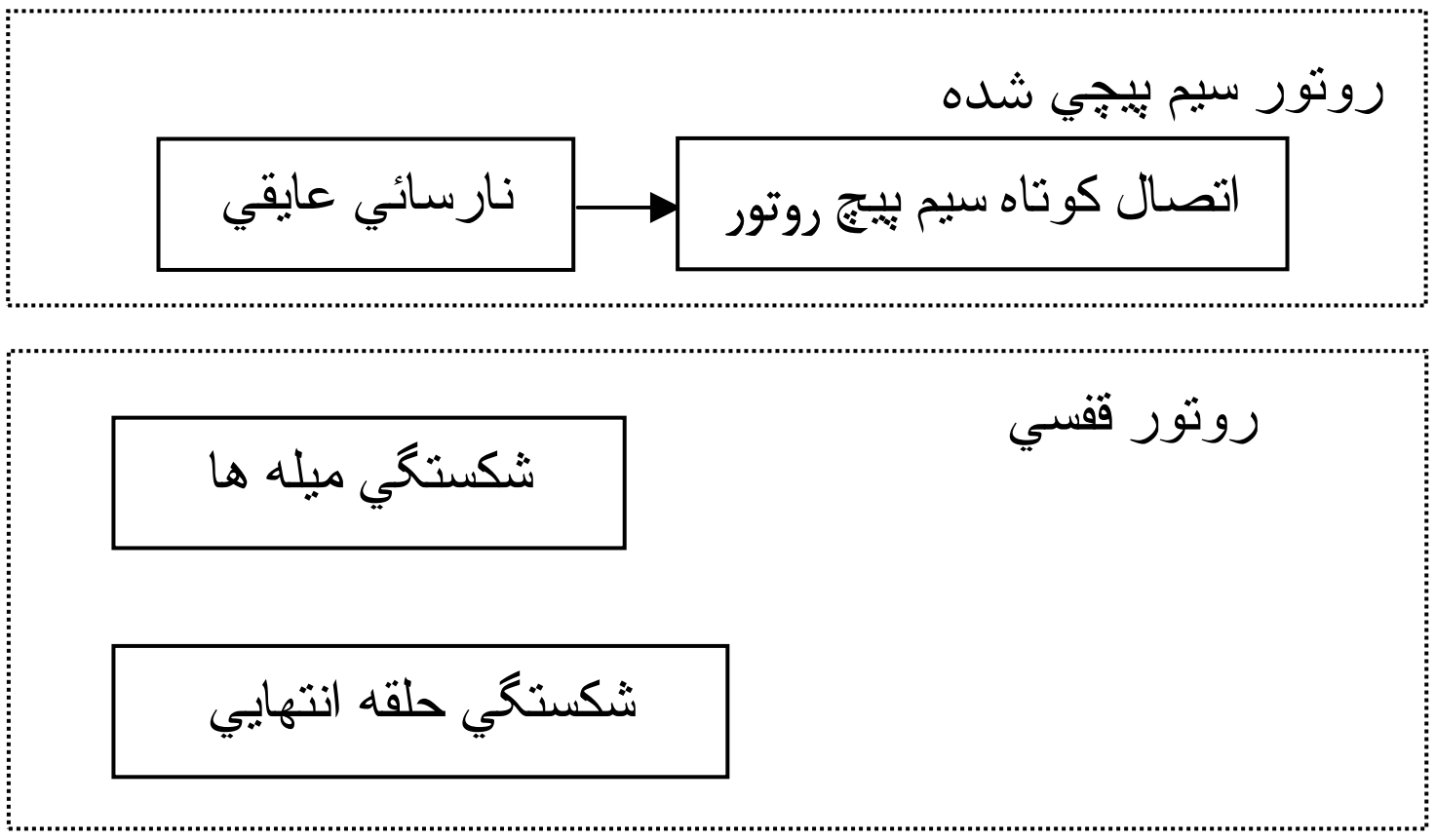

**شکل 1-2:** انواع عیبهایی روتور موتورهای القایی و رابطه بروز آنها

## 1-3- شرایط عیب در قسمتهای مکانیکی

به آن دسته از عیبهایی که به دلیل نقص در ملحقات مکانیکی موتورهای الکتریکی مانند یاتاقانها و چرخ دنده رخ می‌دهد، عیبهای مکانیکی می‌گویند. احتمال وقوع انواع این عیب در موتورهای الکتریکی از سایر عیبهای متداول بیشتر و حدود 50 تا 60 درصد است. مهمترین این عیب ها عبارتند از [5]:

### 1-3-1- عیبهای یاتاقان

حدود 40 تا 50 درصد از عیبهای موتورهای القایی مربوط به یاتاقانهای آنهاست [1]. در حالتی که از یاتاقانهای ساچمه‌ای استفاده می‌شود این عیب به سه دسته کلی زیر طبقه بندی می‌شود [6].

- ✓ عیبهای موجود در ساچمه ها (سائیدگی و یا ترک خوردگی)
- ✓ عیبهای موجود در حلقه بیرونی یاتاقان مانند سائیدگی و یا خوردگی
- ✓ عیبهای موجود در حلقه داخلی یاتاقان

### 2-3-1- ناهم محوری روتور و استاتور

در این حالت بنا به دلایلی تطابق محورهای تقارن روتور، تقارن استاتور و چرخش روتور به هم می‌خورد و فاصله هوایی بصورت غیر یکنواخت در می‌آید. تقریبا 80 درصد عیبهای مکانیکی منجر به ناهم محوری روتور و استاتور می‌شود [3] البته امکان وقوع مستقیم این عیب (بطور مثال در زمان ساخت و نصب رتور) نیز وجود دارد.

در بسیاری از تقسیم بندیها عیبهای ناهم محوری خود یک دسته جداگانه از عیبهای احتمالی موتور بشمار می‌روند. با این حال از آنجا که این عیب موجب نامتقارنی وضعیت مکانیکی روتور و استاتور نسبت به یکدیگر می‌شود، جزء عیبهای مکانیکی در نظر گرفته شده است. انواع این عیبها عبارتند از:

- ✓ ناهم محوری ایستا: در این حالت توزیع فاصله هوایی غیر یکنواخت است ولی موقعیت زاویه‌ای حداقل فاصله هوایی ثابت است. این حالت بیشتر به دلیل بی‌دقتی در نصب روتور اتفاق می‌افتد.
- ✓ ناهم محوری پویا: در این حالت حداقل فاصله هوایی بصورت تابعی از موقعیت زاویه‌ای روتور بوده و پیرامون روتور می‌چرخد. این حالت ممکن است ناشی از Misalignment و یا خمیدگی محور روتور باشد [10].

✓ ناهم محوری مرکب : در صورتیکه هر دو ناهم محوری ایستا و پویا اتفاق بیفتد، ناهم محوری را مرکب می‌گویند.

بسیاری از عیبهای مکانیکی مانند Misalignment، خمیدگی محور روتور و حتی ضعفهای یاتاقان منجر به حالتهای مختلف ناهم محوری می‌شود [10].

## 1-4- آمار بروز عیب در ماشینهای الکتریکی

امروزه ماشینهای الکتریکی کاربرد وسیعی در صنایع کوچک و بزرگ دارند، بگونه‌ای که بدون حضور آنها تمامی نیروگاهها، پالایشگاهها و کارخانه‌ها متوقف خواهند شد. بنابراین همواره محافظت مناسبی جهت پیشگیری از حالت عیب از آنها به عمل می‌آید. درصد وقوع انواع عیب در ماشینهای الکتریکی بصورت جدول (1-1) می‌باشد [6].

این آمار در خصوص موتورهای القایی روتور قفسی بصورت جدول (1-2) می‌باشد [6].

جدول **1-1** : درصد وقوع انواع عیب در ماشینهای الکتریکی

| نوع عیب | درصد وقوع |
|---|---|
| عیبهای یاتاقان | 51/6 |
| عیبهای استاتور | 24/8 |
| عیبهای روتور | 6 |
| عیبهای محور یا کوپلینگ موتور | 3/2 |
| سایر عیبها | 14/4 |

جدول **2-1**: درصد وقوع انواع عیب در موتورهای القایی روتور قفسی

| نوع عیب | درصد وقوع |
|---|---|
| عیبهای مکانیکی | 50 تا 60 |
| عیبهای استاتور | 30 تا 40 |
| عیبهای روتور | 5 تا 10 |

فصل دوم

# بررسی روشهای مدلسازی دینامیکی

# موتور القایی معیوب

روشهای مختلفی جهت تحلیل موتورهای القایی وجود دارد با اینحال بکارگیری این روشها در شرایطی که موتور دچار عیب شده است مستلزم بازنگری در فرضها و معادلات اساسی روش مورد نظر است. به عنوان مثال وقتی در مدل برداری موتور القایی از کلیه هارمونیکهای میدان فاصله هوایی صرفنظر می‌شود دیگر نمیتوان انتظار داشت که به راحتی بتوان شرایط خطای سیم پیچی استاتور و یا شکستگی میله های روتور را به مدل اعمال کرد. بطور کلی روشهایی که با بکارگیری فرضهایی معادلات پیچیده الکترومغناطیسی و الکترودینامیکی موتور القایی را ساده می کنند به آسانی قادر به لحاظ کردن شرایط عیب در موتور نیستند. بنابراین برای مدلسازی و تحلیل موتور القایی معیوب بایستی به بازنگری روشهای موجود پرداخت.

مدلسازی و تحلیل موتور القایی از دو جنبه استاتیکی (حالت دائم) و دینامیکی (حالت گذرا) قابل توجه است. در بررسی حاضر هر دو رفتار حالت دائم و گذرای موتور القایی تحت شرایط ناهم محوری بین روتور و استاتور، همچنین شکستگی و ترک خوردگی میله‌های روتور مدنظر میباشند.

در حالت کلی، عملکرد حالت دائم ماشینهای الکتریکی در سیستمهای الکتریکی و مکانیکی اهمیت دارد. تحلیل عملکرد این حالت از طریق حل مدار معادل ماشین صورت میگیرد. البته از حل معادلات دیفرانسیل حاکم بر عملکرد موتور نیز می توان به رفتار حالت دائم ماشین پی برد. بنابراین با تعیین رفتار حالت گذرای ماشین می توان رفتار حالت دائمی آن را نیز به آسانی پیش بینی کرد.

در یک موتور القایی قفس سنجابی سه فاز انرژی الکتریکی از طریق پایانه های سه فاز استاتور وارد موتور میشود، سپس این انرژی باعث ایجاد میدان مغناطیسی گردان در فضای فاصله هوایی میشود. این میدان نقش واسط در تبدیل انرژی الکتریکی به مکانیکی دارد. بنابر این دو دینامیک تبدیل انرژی الکتریکی به میدان مغناطیسی و تبدیل نیروی مغناطیسی حاصل از این میدان به نیروی مکانیکی در معادلات توصیف مدل موتور القایی وجود دارد. چگونگی تولید و توزیع میدان مغناطیسی ناشی از جریان الکتریکی فازهای استاتور وابستگی زیادی به نحوه سیم پیچی و هندسه فاصله هوایی و حتی هندسه هسته های مغناطیسی روتور و استاتور دارد. در اینجا منظور از هندسه روتور قطر، ابعاد دندانه ها و ابعاد شیارهای آن و منظور از هندسه استاتور قطر داخلی، قطر خارجی، ابعاد شیارها و ابعاد دندانه های آن است.

در ادامه به معرفی روشهای جامع تحلیل موتور القایی قفس سنجابی سه فاز که قادر به پیش بینی رفتار حالت گذرای این موتور در حالت ناهم محوری بین روتور و استاتور، همچنین شکستگی و یا ترک خوردگی میله‌های روتور می‌باشند، می‌پردازیم. البته، روشهای ساده شده دیگری نیز وجود دارند که قادر به تحلیل موتور القایی معیوب در شرایط خاص و با دقت نسبتا کمی هستند. اما در اینجا تنها معرفی روشهای جامع و دقیق مورد نظر می‌باشد.

## 2-1- روش اجزاء محدود

روش اجزاء محدود (FEM) قادر است توزیع میدان مغناطیسی درون موتور را از روی ابعاد هندسی و پارامترهای مغناطیسی موتور محاسبه کند. با اطلاع از توزیع میدان می توان سایر کمیتهای موتور، نظیر شکل موج ولتاژ القایی، چگالی شار مغناطیسی درون فاصله هوایی و اندوکتانس سیم پیچهای مختلف را بدست آورد[7],[6].

اگر $A_m$ نشان دهنده پتانسیل برداری مغناطیسی، $r_m$ مقاومت ویژه مغناطیسی و $J$ مولفه محور $z$ چگالی جریان باشد در هر نقطه از فضا معادله دیفرانسیل زیر برقرار است:

$$\frac{\partial}{\partial x}\left(r_m \frac{\partial \vec{A}_m}{\partial x}\right) + \frac{\partial}{\partial y}\left(r_m \frac{\partial \vec{A}_m}{\partial y}\right) = -J \qquad (1\text{-}2)$$

محور $z$ در امتداد محور طولی ماشین و $x$ و $y$ دو محور متعامد در صفحه سطح مقطع ماشین میباشند. معادله (1-2) معادله اساسی پتانسیل بردار مغناطیسی در ماشین های الکتریکی است. به منظور حل این معادله دیفرانسیل، عبارت تابعی زیر برای تفاضل بین انرژی کل ذخیره شده در میدان مغناطیسی و انرژی ورودی به ماشین تعریف می‌شود.

$$Fu(A_m) = \iint_s (0.5\ \vec{B}.\vec{H} - \vec{J}.\vec{A}_m)dxdy \qquad (2\text{-}2)$$

در معادله فوق $s$ سطح مقطع ماشین است. کمینه سازی معادله فوق به دستگاه معادلاتی منجر میشود که با حل عددی آن توزیع پتانسیل برداری مغناطیسی درون ناحیه $s$ بدست می‌آید. توزیع میدان مغناطیسی در این ناحیه از روی پتانسیل برداری مغناطیسی بصورت زیر قابل محاسبه است.

$$\vec{B} = \nabla \times \vec{A}_m \qquad (3\text{-}2)$$

کمیتهای دیگر مانند شکل موج ولتاژ القایی در روتور، چگالی شار مغناطیسی در فاصله هوایی و اندوکتانس سیم پیچهای مختلف با توجه به توزیع میدان مغناطیسی محاسبه می‌شوند[10],[9].

در گذشته بعلت طولانی شدن حل کامپیوتری معادلات این روش تنها برای تحلیل عملکرد دائمی ماشینهای الکتریکی از این روش استفاده می شد. امروزه با وجود کامپیوترهای قوی از FEM در تحلیل حالت گذرای ماشینهای الکتریکی نیز استفاده می شود. با اینحال هنوز معادلات ناشی از این روش با صرف زمان بسیار طولانی و با اعمال تکنیکهای ریاضی پیچیده ای قابل حل می باشند. یکی از روشهای کاهش زمان حل معادلات فوق استفاده از تقارن هندسی ماشین است. برای مثال تحت شرایط خاصی می توان با حل معادلات دیفرانسیل اجزاء محدود در یک چهارم ماشین و به دست آوردن توزیع میدان مغناطیسی در آن توزیع میدان مغناطیسی در کل ماشین را بدست آورد. واضح است که این ترفند در حالتی که خطای ناهم محوری میان روتور و استاتور و یا شکستگی یک میله در روتور رخ داده باشد کارساز نیست چون شرط تقارن هندسی موتور به هم‌خورده است.

یکی دیگر از دشواریهای استفاده از این روش در تحلیل موتورهای القایی سه فاز قفس سنجابی این است که بایستی در هر ناحیه چگالی جریان $(J)$ را به عنوان ورودی به معادلات وارد کرد. وقتی عملکرد موتور در حالت دائم مورد نظر باشد و یا بتوان از هارمونیکهای جریان خط صرفنظر کرد این کار امکان پذیر است، اما اگر هدف تحلیل گذرای موتور و یا به دست آوردن طیف فرکانسی جریان خط باشد این امر مستلزم محاسبه چگالی جریان سیمبندیهای مختلف ماشین با روشهای دیگری است که منجر به پیچیدگی بیشتر روش اجزاء محدود می شود.[11]در بیشتر تحلیلهای که بر روی موتور القایی در شرایط ناهم محوری و میله شکسته روتور انجام میگیرد هدف

محاسبه هارمونیکهای اضافی تولید شده توسط موتور معیوب است. بنابر این استفاده از روش اجزاء محدود در این حالت مستلزم محاسبه $J$ از روش دیگری است.

یکی از روشهایی که برای حل این مشکل ارائه شده است ترکیب روش اجزاء محدود و فضای حالت است که Coupled Finite Element State Space نام دارد. در این روش چگالی جریان از معادلات فضای حالت غیر خطی موتور به دست می آید و در معادلات اجزاء محدود قرار می گیرد، پس از حل معادلات اجزاء محدود، توزیع میدان مغناطیسی در فاصله هوایی به دست می آید. این توزیع به عنوان ورودی به یک روالنرم افزاری داده می شود تا از روی آن اندوکتانسهای مغناطیسی کلیه اجزاء مغناطیسی موتور محاسبه شود. این اندوکتانسهای دوباره بصورت ورودی معادلات فضای حالت قرار می گیرند تا از روی آنها و ولتاژ و بار مکانیکی اعمالی به موتور مجددا چگالی جریان محاسبه شود[12]. شکل 2-1 فلوچارت اجرای روالهای نرم افزاری این روش و پارامترهای ورودی-خروجی هر روال را نشان می دهد.

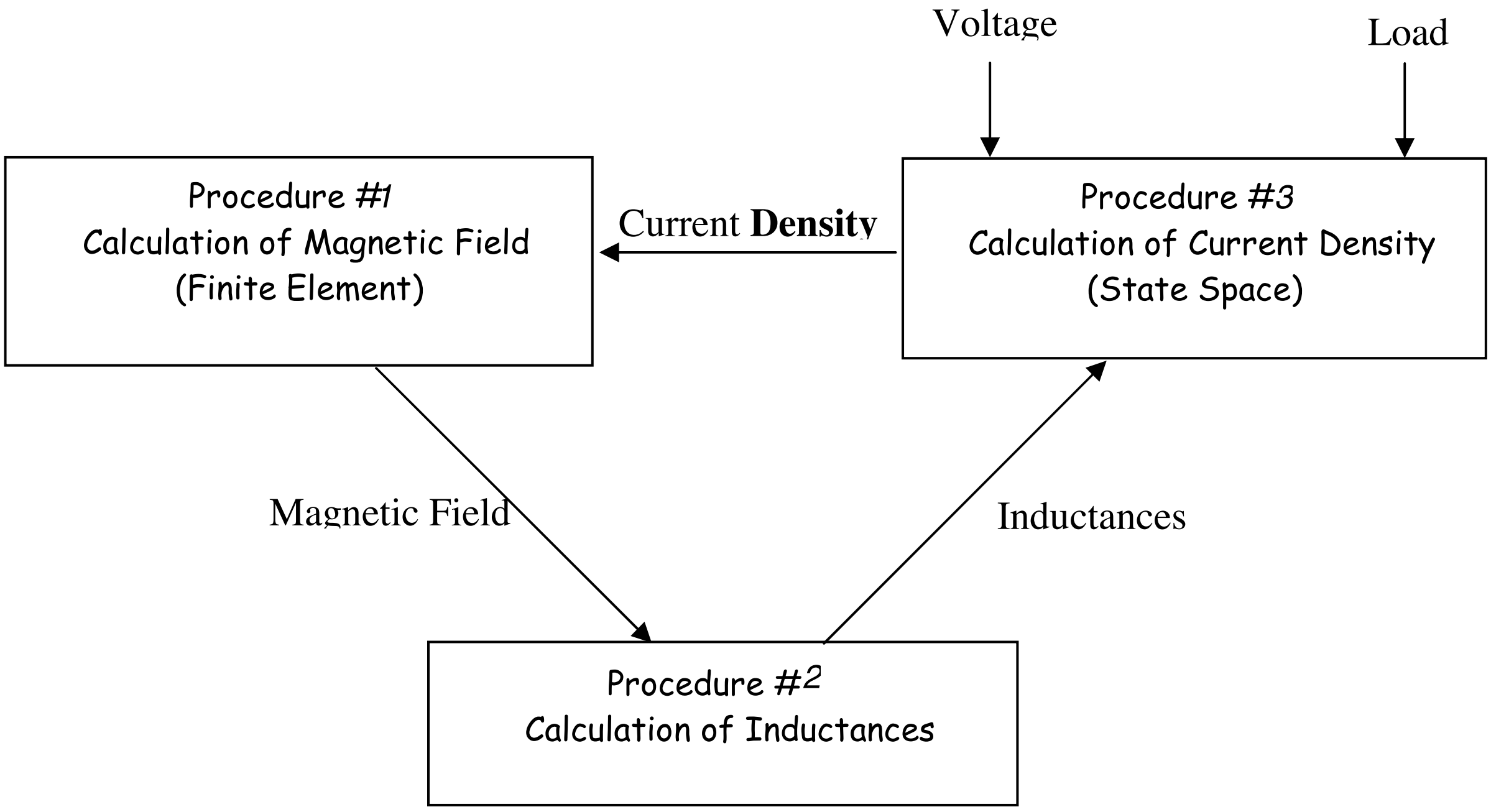

شکل**2-1** : اجزاء روش اجزاء محدود فضای حالت و ارتباطات بین آنها

## 2-2- روش تفاضلات متناهی

از روش تفاضلات متناهی (FDM) نیز می توان برای محاسبه توزیع میدان مغناطیسی موتور القایی استفاده نمود. در این روش معادلات دیفرانسیل توصیف مدل اجزاء محدود با معادلات تفاضلی جبری جایگزین میشود. همچنین از شبکه های گسسته کننده فضای مورد نظر نیز استفاده می شود. به شبکه های گسسته محدود به دست آمده از این روش نواحی گویند که بصورت چهار ضلعی می باشند و با یک مقاومت مغناطیسی مشخص می‌شوند[13]. پس از گرد آوری معادلات تفاضلی جبری برای تمامی نواحی به یک دستگاه معادلات میرسیم که از حل آن توزیع میدان در تمام نواحی به دست می آید. در مقایسه این روش با روش اجزاء محدود می توان به موارد زیر اشاره نمود[14]:

- ✓ دقت FEM به مراتب بیشتر از FDM است.
- ✓ اجرای کامپیوتری FDM کوتاهتر از FEM است.
- ✓ برنامه نویسی FDM آسانتر از FEM است.
- ✓ FDM در مقایسه با FEM احتیاج به حافظه کامپیوتری کمتری دارد.
- ✓ در FDM حد بالایی وجود دارد که دقت جوابها نمی تواند از آن بیشتر باشد در حالیکه در FEM هرچه شبکه بندی دقیقتر صورت گیرد دقت جوابها بصورت نمایی افزایش می یابد.

همچنین معادلات جبری این روش نیز مانند روش اجزاء محدود به چگالی جریان به عنوان ورودی نیاز دارند. بنابراین در کنار این معادلات بایستی معادلات دیگری جهت محاسبه اندوکتانسها و همچنین محاسبه چگالی جریان تعریف نمود.

با توجه به آنچه گفته شد می توان استنتاج کرد که در حالتهای مختلف عیب اعم از ناهم محوری بین روتور و استاتور و خطای میله شکسته استفاده از FEM نسبت به FDM منطقی تر و ارجح است.

## 2-3- روش مدار معادل مغناطیسی

یکی دیگر از تکنیکهای عددی تحلیل حالت گذرای موتورهای القایی روش مدار معادل مغناطیسی (MEC)

است [15]. میان این روش و روش اجزاء محدود دو تفاوت اساسی وجود دارد:

تعداد اجزاء بکار رفته در این روش به مراتب کمتر از روش اجزاء محدود است.

قبل از اعمال روش MEC باید جهت شارمغناطیسی معلوم باشد در حالیکه در FDM جهت شار مغناطیسی یکی از جوابهای معادلات است.

در این روش برای تمامی قسمتهای ماشین مدار معادل مغناطیسی ارائه می شود. دندانه، شیار و یوغهای استاتور و روتور به همراه فاصله هوایی و سیم پیچهای استاتور و میله های روتور دارای مدار معادل مغناطیسی می باشند. سپس با استفاده از جهتی که برای عبور شار مغناطیسی فرض می شود این مدارهای معادل به همدیگر متصل می‌شوند. مدار معادل حاصله دارای منابع نیروی محرکه مغناطیسی و هدایتهای مغناطیسی خواهد بود. گره های این مدار مشخص کننده پتانسیل مغناطیسی اسکالر موجود در نقاط مختلف ماشین و جریان گذرنده از هر پرمانس نشان دهنده شار عبوری از آن جزء خواهد بود[15].

مزیت اصلی این روش نسبت به FEM کوتاهی زمان محاسبه آن است. به این ترتیب که زمان لازم برای تحلیل حالات گذرای این روش به مراتب کمتر از FEM است[16]. اما دقت آن نیز کمتر از FEM است.

## 2-4- روش تابع سیم پیچی

روش دیگری که اخیراً مورد توجه قرار گرفته است تئوری تابع سیم پیچی(WFM) است. این تئوری با استفاده از مدارهای الکتریکی با تزویج مغناطیسی روابط بین شار حلقوی، جریانهای سیم پیچهای استاتور، جریان میله های روتور و ولتاژهای موتور را تدوین می‌کند. این روابط بر حسب پارامترهای مقاومت سیم پیچهای استاتور و میله های روتور و اندوکتانس خودی و متقابل آنهاست. محاسبه اندوکتانس متقابل سیم پیچهای استاتور و حلقه های روتور قدم بعدی در بکارگیری این نظریه است. تا کنون مطالعات بسیاری جهت محاسبه این اندوکتانسها با در نظر گرفتن اثر تمامی هارمونیکهای فضایی ناشی از توزیع سیم پیچهای استاتور و شیار گذاری روتور انجام گرفته است. یکی از مزایای این روش آن است که در بررسی رفتار هر ماشینی با هر توزیع سیم بندی و طول فاصله هوایی ضمن در نظر گرفتن اثر تمامی هارمونیکهای فضایی و زمانی قادر به محاسبه پاسخ حالت گذرا و پایای ماشین می‌باشد. این امر باعث می‌شود تا بتوان کلیه عیبهای ناشی سیم پیچهای استاتور، حلقه های روتور و ناهم محوری

روتور و استاتور را در مدل به دست آمده از این تئوری لحاظ کرد.

این تئوری نخستین بار بمنظور محاسبه اندوکتانسهای مغناطیسی موتور القایی تک فاز در سال 1965 ارائه شد[17]. در 1969 با استفاده از این تئوری مدارهای معادل حالت دائمی موتورهای القایی دوفاز نامتقارن با در نظر گرفتن هارمونیکهای زوج و فرد در میدان مغناطیسی ارائه شد[18]. در 1979 این روش در تحلیل موتورهای القایی خطی بکار گرفته شد[19]. در 1985 با تعمیم این تئوری و ترکیب آن با مفاهیم مدارهای الکتریکی تزویج شده مغناطیسی مدل معادل dqn برای ماشینهای سنکرون و القایی ارائه شد [20]و در 1991 در تحلیل ماشینهای‌القایی با سیم پیچهای متمرکز استفاده شد[21].

این تئوری نخستین بار در سال 1992 در بررسی رفتار حالت گذرای موتورهای القایی تحت عیبهای داخلی استفاده شد[22]. سپس با استفاده از آن مطالعاتی در زمینه بررسی رفتار دینامیکی موتور القایی و هارمونیکهای جریان استاتور در حالتهای مختلف عیب سیم پیچی استاتور، شکستگی میله های روتور و ناهم محوری روتور و استاتور انجام شد[23].

در بررسی حاضر استفاده از این روش برای تحلیل رفتار موتور القایی روتور قفسی سه فاز در شرایط ناهم محوری روتور و استاتور و همچنین در شرایط ترک خوردگی و شکستگی میله های روتور انجام شده است. به منظور تحلیل رفتار موتور در شرایط ناهم محوری روتور و استاتور نیاز به تعمیم نظریه تابع سیم پیچی در حالت غیر یکنواختی فاصله هوایی است که بدان پرداخته خواهد شد.

این نکته قابل ذکر است که با اینکه تاکنون تحقیقات بسیاری جهت تحلیل رفتار موتورهای القایی در حالت ناهم محوری روتور و استاتور با استفاده از این روش شده ولی همگی از نظریه تابع سیم پیچی استفاده کرده اند که از شرایط تقارن فاصله هوایی به دست آمده است[24,25]. در مطالعه حاضر با بازنگری فرضها و روابط اساسی این نظریه نشان داده می شود که روش فوق الذکر در حالتی که فاصله هوایی یکنواخت نیست نسبت به حالتی که فاصله هوایی یکنواخت است متفاوت است[26].

فصل سوم

# تعمیم نظریه تابع سیم‌پیچی در حالت ناهم‌محوری روتور و استاتور

نظریه تابع سیم پیچی در 1965 جهت محاسبه اندوکتانسهای مغناطیسی موتور القایی ارائه شد[17]. این نظریه روش حلی برای یک سری معادلات تزویج الکترومغناطیسی است که با در نظر گرفتن قوانین ساده تئوری مدارهای الکتریکی تشکیل می‌شوند. یک قسمت اساسی از این نظریه محاسبه اندوکتانسهای موتور است. این اندوکتانسها با تعریف توابع سیم پیچی و روابط تعریف شده در این نظریه قابل محاسبه هستند.

در این فصل ضمن بازنگری این نظریه و تعمیم آن در شرایط فاصله هوایی غیر یکنواخت، تفاوت رابطه محاسبه اندوکتانس در حالت فاصله هوایی غیر یکنواخت با حالتی که فاصله هوایی یکنواخت است مطرح می‌شود. این در حالی است که تاکنون کلیه مقالاتی که با استفاده از این نظریه به تحلیل موتور القایی در حالتهای غیریکنواختی فاصله هوایی پرداخته‌اند از رابطه محاسبه اندوکتانسی استفاده کرده‌اند که فقط در حالت فاصله هوایی یکنواخت معتبر است. این مساله باعث شده که در محاسبات آنان تساوی $M_{12} = M_{21} = M$ نقص شود. شکلهای نشان داده شده در [24] و [25] به خوبی بیانگر این دعا هستند.

با استفاده از تعمیم ارائه شده خواهیم دید همانطور که از همه مدارهای مغناطیسی خطی انتظار می‌رود، همواره با هر توزیع فاصله هوایی تساوی $M_{12} = M_{21} = M$بر قرار است. در اینصورت نسبت به آنچه در مقالات [24] و [25] ارائه شده تردید اساسی وجود خواهد داشت.

## 3-1- مفهوم نظریه تابع سیم پیچی

با استفاده از نظریه تابع سیم پیچی میتوان کلیه اندوکتانسهای مغناطیسی یک موتور القایی سه فاز روتور قفسی را محاسبه نمود. این امر از طریق محاسبه شدت میدان مغناطیسی ناشی از هر سیم‌پیچی توسط این نظریه ممکن می شود. در ابتدای بکارگیری این نظریه بایستی تابع دور متناظر با هریک از فازهای استاتور و یا حلقه های روتور تعیین گردد. این تابع بیانگر توزیع سیم پیچی مورد نظر در اطراف فاصله هوایی با در نظر گرفتن جهت عبور جریان در سیم پیچی است. نظریه تابع سیم پیچی با اطلاع از تابع دور متناظر با هر فاز استاتور و حلقه روتور و با اعمال قوانین الکترومغناطیسی در موتور القایی، متناظر با هر تابع دور یک تابع سیم پیچی تعریف میکند که شدت میدان

مغناطیسی ناشی از آن سیم پیچی در هر نقطه فاصله هوایی متناسب با حاصلضرب این تابع در جریان عبوری از آن است. ضریب تناسب این رابطه عکس طول فاصله هوایی در آن نقطه از فاصله هوایی است.

آنچه تا کنون در مطالعات منتشر شده وجود دارد استفاده از تعریفی از تابع سیم پیچی است که با فرض یکنواختی فاصله هوایی بدست آمده است. در حالیکه در حالت غیر یکنواختی فاصله هوایی این تعریف صادق نیست. بنابراین درحالتهای مختلف ناهم محوری نمی توان از تعریف موجود استفاده نمود. این در حالی است که تا کنون به اشتباه از این تعریف برای مدلسازی و شبیه سازی موتور القایی در حالتهای خطای ناهم محوری ایستا [24] و پویا [25] استفاده شده است.

## 3-2- مفهوم تابع دور

با توجه به قانون مداری آمپر در الکترو مغناطیس میتوان برای هر توزیع سیم بندی در موتور القایی یک تابع دور تعریف نمود. این تابع که برای سیم پیچی دلخواه $x$ بصورت $n_x(\varphi)$ نشان داده می شود وابسته به متغیر $\varphi$ در دستگاه مرجع استاتور است. در صورتیکه قانون مداری آمپر را به مسیر نشان داده شده در شکل 3-1 اعمال کنیم خواهیم داشت :

$$\oint \vec{H}_x(r,\varphi)d\vec{\ell} = \int_S \vec{J}_x(r,\varphi)d\vec{s} \qquad (1-3)$$

در این رابطه $H_x$ شدت میدان مغناطیسی ناشی از سیم پیچی $x$ و $J_x$ توزیع چگالی جریان این سیم پیچی است.

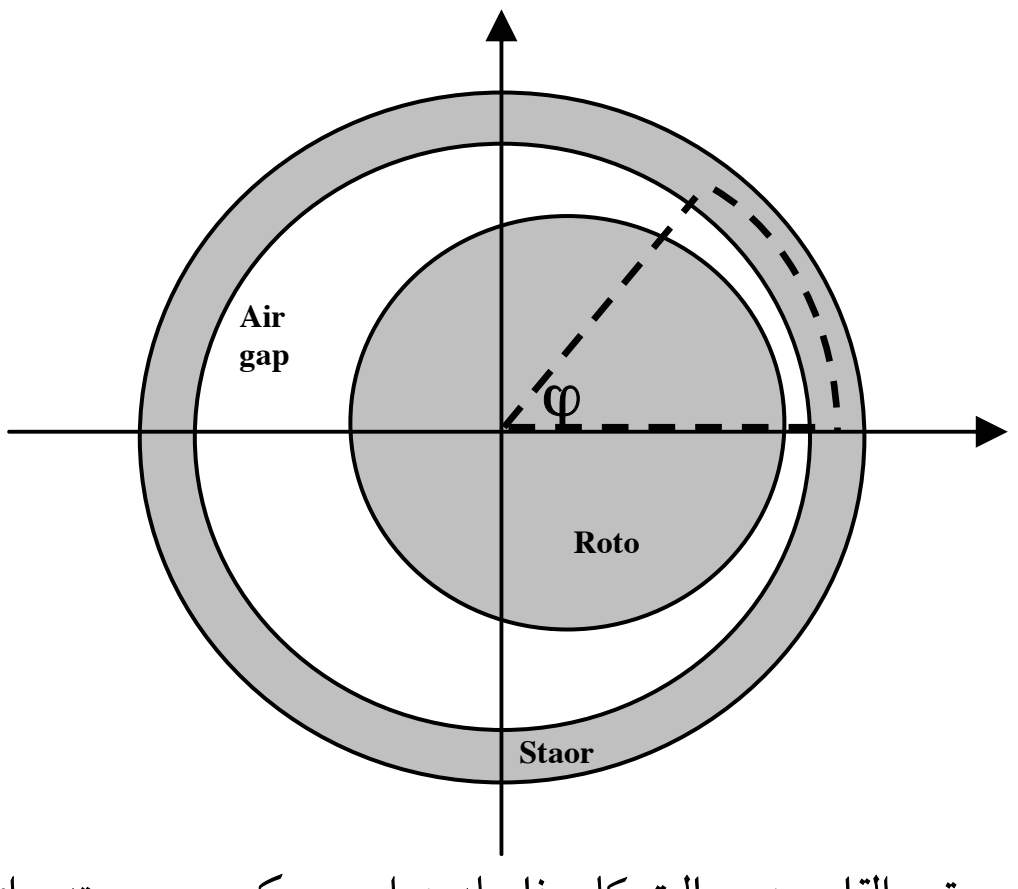

شکل 3-1 : سطح مقطع موتور القایی در حالت کلی فاصله هوایی و یک مسیر بسته برای اعمال قانون مداری آمپر

سطح بسته $S$ ناحیه محصور بین دو خط $\theta = 0$ و $\theta = \varphi$ در دستگاه مرجع استاتور است. تابع دور سیم پیچی $x$ بصورت زیر تعریف می شود.

$$n_x(\varphi) \equiv \frac{\int_S \vec{J}_x(r,\varphi)d\vec{s}}{i_x} \tag{2-3}$$

این رابطه بایستی در هر لحظه از زمان صادق باشد. اگر توزیع جریان سیم بندیهای موتور بصورت گسسته و نقطه‌ای در نظر گرفته شود، می توان رابطه فوق را بصورت زیر در نظر گرفت.

$$n_x(\varphi) = \frac{\sum_s (-1)^{x(\varphi)} i_x}{i_x} \tag{3-3}$$

بنابراین:

$$n_x(\varphi) = \sum_s (-1)^{x(\varphi)} \tag{4-3}$$

در این رابطه $x(\varphi)$ بصورت زیر تعریف می شود:

$x(\varphi) = 0$ : اگر جهت جریان در هادی قرار گرفته در زاویه $\varphi$ در جهت مثبت باشد.

$x(\varphi) = 1$: اگر جهت جریان در هادی قرار گرفته در زاویه $\varphi$ در جهت منفی باشد.

برای مثال تابع دور سیم پیچی شکل 3-2 در شکل 3-3 نشان داده شده است. در تابع دور رسم شده در شکل 3-3 فرض شده است که تابع دور در بالای هادیهای استاتور بصورت پله ای افزایش می یابد. اما از آنجا که جریان روی سطح مقطع سیم نیز توزیع می شود بهتر است افزایش تابع دور در بالای هادیها بصورت خطی در نظر

گرفته شود. این مساله باعث می شود تا هارمونیکهای فضایی بیشتری در تابع دور لحاظ شود. در مطالعات دقیقتر می توان تغییرات تابع دور در بالای هادیها را بر اساس شکل هندسی هادیها و شیارهای استاتور محاسبه کرد.

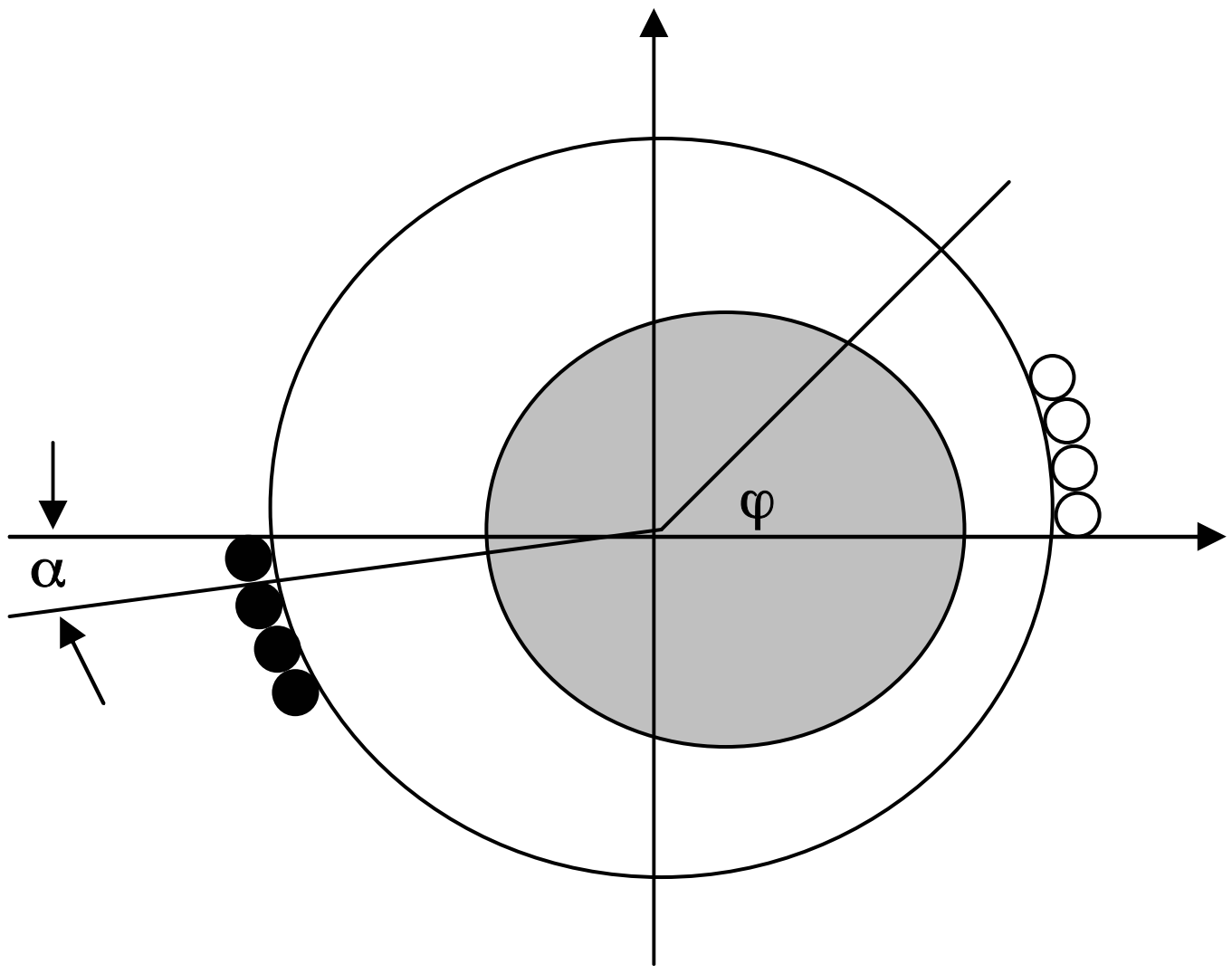

شکل **3-2** : توزیع یک سیم بندی دلخواه در موتور

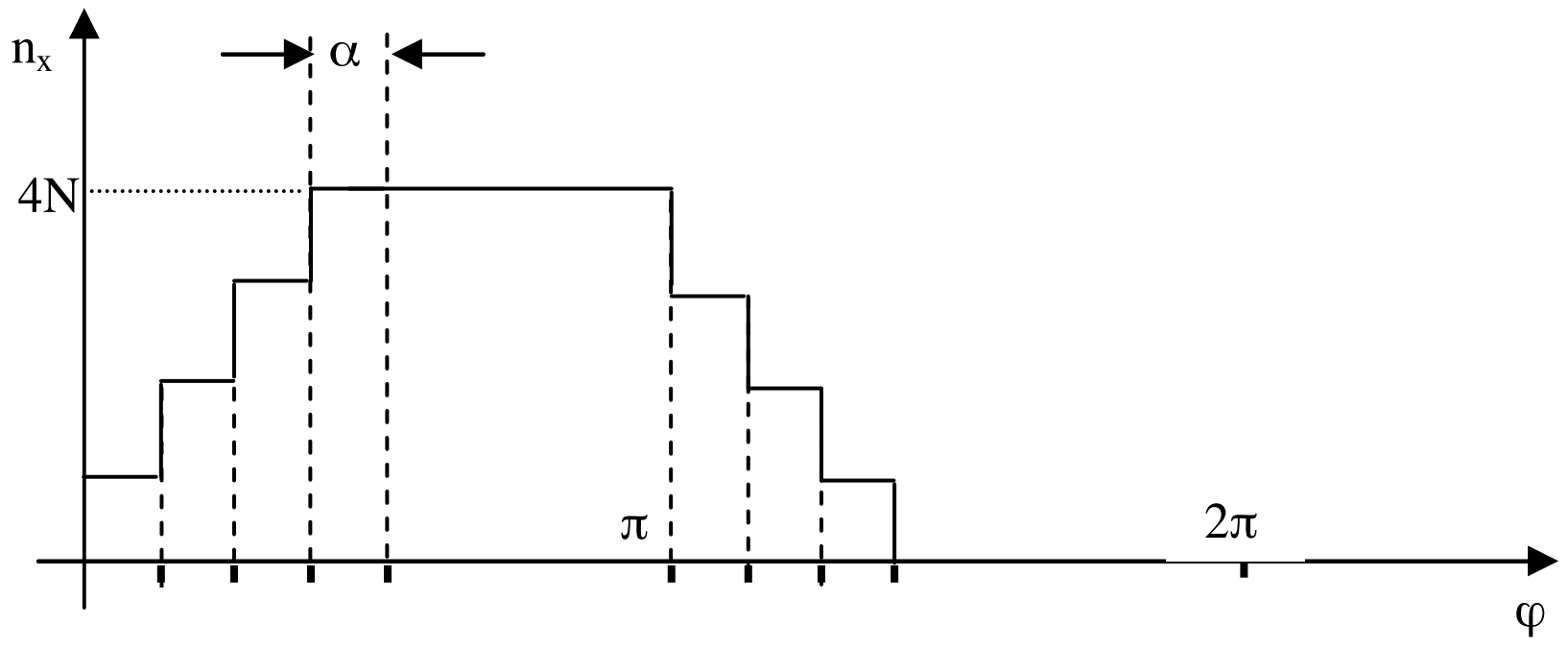

سحں د-3 : تابع دور سیم پیچی نشان داده شده در شکل **3-2** با فرض افزایش پله‌ای چگالی جریان

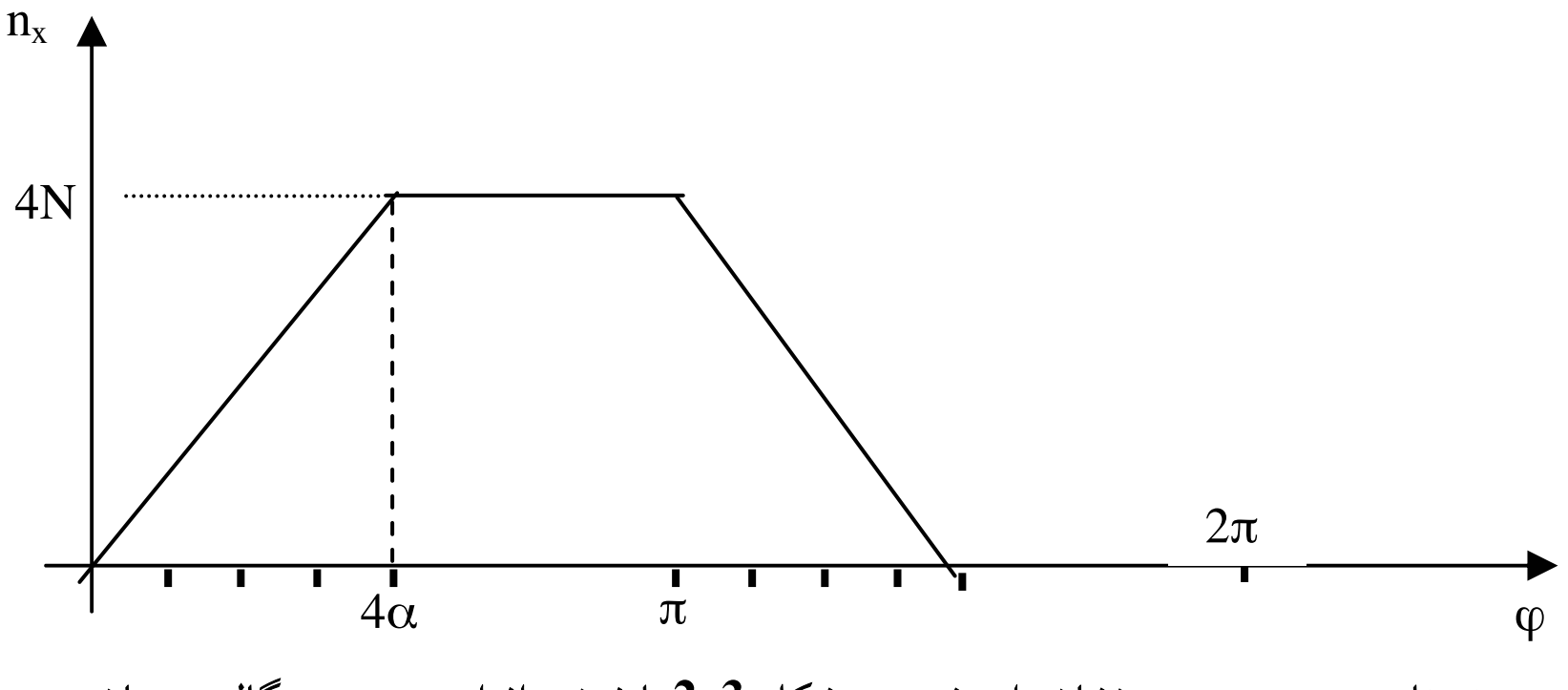

سحں د-4 : تابع دور سیم پیچی نشان داده شده در شکل **3-2** با فرض افزایش خطی چگالی جریان

## 3-3- تعمیم نظریه تابع سیم پیچی

با استفاده از مسیر نشان داده شده در شکل 3-1 و قانون مداری آمپر برای سیم پیچی فرضی $x$ خواهیم داشت:

$$\oint \vec{H}_x(r,\varphi)d\vec{\ell} = \int_S \vec{J}_x(r,\alpha)\varphi\, d\vec{s} \tag{5-3}$$

طبق رابطه (3-2) و تعریف تابع دور داریم :

$$\int_S \vec{J}_x(r,\varphi)d\vec{s} = n_x(\varphi)i_x \tag{6-3}$$

از آنجا که پرمابیلیته نسبی آهن به مراتب بیشتر از پرمابیلیته نسبی فاصله هوایی(یک) است، می توان از افتMMF در داخل آهن نسبت به افت آن در فاصله هوایی صرفنظر کرد. همچنین به علت کوچکی طول فاصله هوایی فرض می‌کنیم در زاویه دلخواه $\varphi$ شدت میدان مغناطیسی مستقل از شعاع و برابر مقدار آن در وسط فاصله هوایی است. بنابراین خواهیم داشت:

$$\oint \vec{H}_x(r,\varphi)d\vec{\ell} = \int_{airgap} \vec{H}_x(\varphi)d\ell = H_x(\varphi)g(\varphi) - H_x(0)g(0) \tag{7-3}$$

از ترکیب روابط اخیر خواهیم داشت:

$$H_x(\varphi)g(\varphi) - H_x(0)g(0) = n_x(\varphi)i_x \tag{8-3}$$

بنابراین :

$$H_x(\varphi) = \frac{n_x(\varphi)i_x + H_x(0)g(0)}{g(\varphi)} \qquad (9\text{-}3)$$

جهت محاسبه $H_x(0)$ سطح گوسی استوانه‌ای‌شکلی را در عمق فاصله هوایی در نظر می‌گیریم در اینصورت داریم:

$$\int_S \mu_0 H_x(\varphi) ds = 0 \qquad (10\text{-}3)$$

با ترکیب دو رابطه اخیر داریم :

$$\int_S \frac{n_x(\varphi)i_x + H_x(0)g(0)}{g(\varphi)} ds = 0 \qquad (11\text{-}3)$$

با توجه به اینکه در مختصات استوانه ای $ds = r\ell d\varphi$ می باشد، داریم:

$$\int_0^{2\pi} \frac{n_x(\varphi)}{g(\varphi)} r_{av}(\varphi) i_x d\varphi + g(0)H_x(0)\int_0^{2\pi} \frac{r_{av}(\varphi)}{g(\varphi)} d\varphi = 0 \qquad (12\text{-}3)$$

با توجه به رابطه اخیر $H_x(0)$ بصورت زیر بدست می آید:

$$H_x(0) = \frac{-\int_0^{2\pi} \frac{r_{av}(\varphi) n_x(\varphi)}{g(\varphi)} i_x d\varphi}{g(0)\int_0^{2\pi} \frac{r_{av}(\varphi)}{g(\varphi)} d\varphi} \qquad (13\text{-}3)$$

با ترکیب روابط (3-9) و (3-13) داریم:

$$H_x(\varphi) = \left( n_x(\varphi) - \frac{\langle \rho n_x \rangle}{\langle \rho \rangle} \right) \rho(\varphi) i_x \qquad (14\text{-}3)$$

در این رابطه $\rho(\varphi)$ توزیع حاصلضرب شعاع متوسط فاصله هوائی در عکس طول فاصله هوایی است و عملگر $\langle f \rangle$ بصورت میانگین تابع $f$ در بازه $[0,2\pi]$ به صورت زیر تعریف شده است:

$$\langle f \rangle = \frac{1}{2\pi} \int_0^{2\pi} f(\alpha) d\alpha \qquad (15\text{-}3)$$

با توجه به رابطه (3-14) تابع سیم پیچی سیم پیچ $x$ یا $N_x(\alpha)$ بصورت زیر تعریف می‌شود:

$$N_x(\varphi) = n_x(\varphi) - \frac{\langle \rho n_x \rangle}{\langle \rho \rangle} \qquad (16\text{-}3)$$

در حالتی که توزیع فاصله هوایی یکنواخت باشد رابطه اخیر بصورت زیر ساده خواهد شد:

$$N_x(\varphi) = n_x(\varphi) - \langle n_x \rangle \qquad (17\text{-}3)$$

اما این رابطه در حالت فاصله هوایی غیر یکنواخت صادق نیست. این در حالی است که در [24] و [25] جهت تحلیل موتور القایی روتور قفسی در شرایط غیر یکنواختی فاصله هوایی از این رابطه استفاده شده است. با توجه به رابطه (3-14) چگالی شار مغناطیسی ناشی از سیم پیچی $x$ عبارت است از:

$$B_x(\varphi) = \mu_0 N_x(\varphi)\rho(\varphi) i_x \tag{18-3}$$

بنابراین شار پیوندی گذرنده از سیم پیچ دیگری مانند $y$ که از سیم پیچی $x$ ناشی شده، بصورت رابطه زیر قابل محاسبه خواهد بود:

$$\lambda_{yx}(\varphi) = l i_x \int_0^{2\pi} P(\varphi) N_x(\varphi) n_y(\varphi) d\varphi \tag{19-3}$$

$$P(\varphi) = \frac{\mu_0 r_{av}(\varphi)}{g(\varphi)} \tag{20-3}$$

در این رابطه $P(\varphi)$ توزیع هدایت مغناطیسی واحد طول فاصله هوایی ، $l$ طول روتور و $r_{av}(\varphi)$ شعاع متوسط فاصله هوایی است. از آنجا که نسبت طول فاصله هوایی به شعاع متوسط آن بسیار کوچک است با تقریب بسیاز خوبی می‌توان از شعاع متوسط فاصله هوایی در حالت فاصله هوایی یکنواخت استفاده کرد.

با توجه به رابطه (3-16) اندوکتانس متقابل دو سیم پیچی $x$ و $y$ بصورت زیر خواهد بود:

$$L_{yx} = l \int_0^{2\pi} P(\varphi) N_x(\varphi) n_y(\varphi) d\varphi \tag{21-3}$$

با ترکیب رابطه (3-16) با رابطه اخیر داریم:

$$L_{yx} = 2\pi l \langle P n_x n_y \rangle - 2\pi l \frac{\langle P n_x \rangle \langle P n_y \rangle}{\langle P \rangle} \tag{22-3}$$

رابطه (3-22) نسبت به $x$ و $y$ خاصیت جابجایی دارد بنابراین همواره داریم:

$$L_{yx} = L_{xy} \qquad (23\text{-}3)$$

این رابطه در کلیه مدارهای مغناطیسی خطی وجود دارد. این در حالی است که اگر مانند نظریه تابع سیم‌پیچی بکار رفته در [24] و [25] عمل کنیم این تساوی نتیجه نمی‌شود.

## 3-4- وضعیت فاصله هوایی در حالت های مختلف ناهم محوری

در حالت خطای ناهم محوری روتور و استاتور به دلایلی که عمدتا مکانیکی هستند محورهای تقارن استاتور، تقارن روتور و چرخش روتور نسبت به یکدیگر جابجا می شوند.

با توجه به چگونگی عدم تطابق محورهای تقارن استاتور، تقارن روتور و چرخش روتور خطای ناهم محوری روتور و استاتور به سه دسته کلی ناهم محوری ایستا، ناهم محوری پویا و ناهم محوری مرکب تقسیم می‌شود.

### 3-4-1- ناهم محوری ایستا(S.E)

در این حالت محور چرخش روتور منطبق بر محور تقارن آن است ولی نسبت به محور تقارن استاتور جابجا شده است. در این حالت توزیع فاصله هوایی در اطراف روتور یکنواختی خود را از دست می‌دهد اما متغیر با زمان نیست[20].

درجه ناهم محوری ایستا یا $\delta_s$ بصورت رابطه زیر تعریف می شود.

$$\delta_s = \frac{\left|\overrightarrow{O_S O_\omega}\right|}{g_0} \qquad (24\text{-}3)$$

در این رابطه $O_S$ مرکز تقارن استاتور ، $O_r$ مرکز تقارن روتور، $O_\omega$ مرکز چرخش روتور و $g_0$ طول فاصله هوایی یکنواخت می باشد. $\alpha_S$ زاویه اولیه ناهم محوری ایستا و بردار $\overrightarrow{O_S O_\omega}$ بردار انتقال ایستا نامیده می‌شود این بردار به ازاء تمامی موقعیتهای زاویه ای روتور ثابت است. شکل 3-5 وضعیت سطح مقطع روتور و استاتور را در این حالت نشان می دهد.

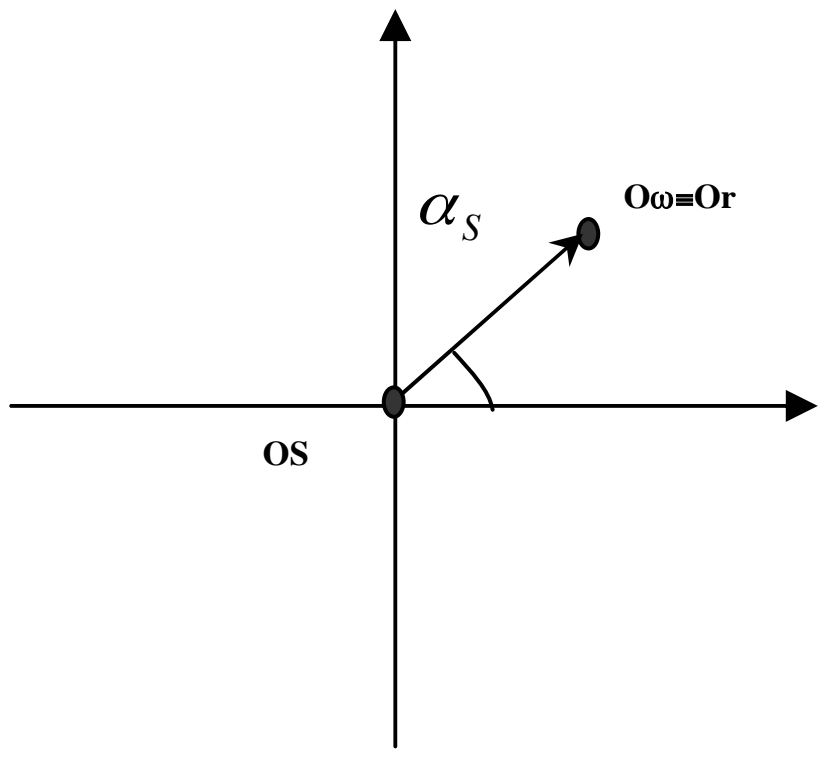

شکل **3-5**: وضعیت سطح مقطع روتور و استاتور در حالت ناهم محوری ایستا در دستگاه مرجع استاتور

### 3-4-2- ناهم محوری پویا(D.E)

در این حالت محور تقارن استاتور وچرخش روتور بر هم منطبق است اما محور تقارن روتور نسبت به آنها جابجا شده است. در این حالت توزیع فاصله هوایی در اطراف روتور بصورت غیریکنواخت و متغییر با زمان خواهد بود

درجه ناهم محوری پویا یا $\delta_d$ بصورت رابطه زیر تعریف می شود.

$$\delta_d = \frac{\left|\overrightarrow{O_\omega O_r}\right|}{g_0} \qquad (25\text{-}3)$$

شکل 3-6 وضعیت سطح مقطع روتور و استاتور را در این حالت نشان می‌دهد.

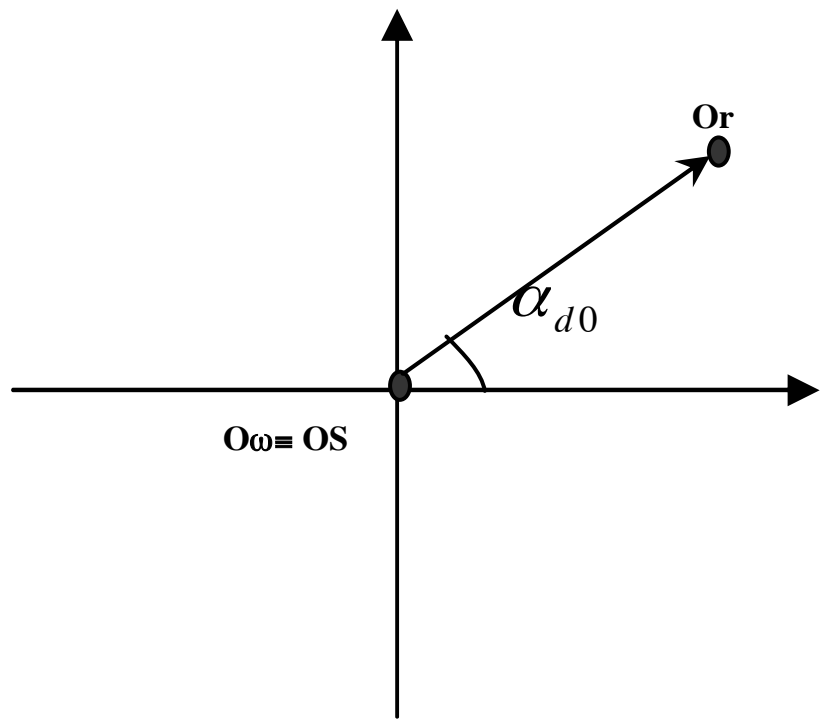

شکل 3-6: وضعیت سطح مقطع روتور و استاتور در حالت ناهم محوری پویا در دستگاه مرجع استاتور

$\alpha_{d0}$ را زاویه اولیه ناهم محوری پویا و بردار $\overrightarrow{O_\omega O_r}$ بردار انتقال پویا نامیده می شود. اندازه این بردار به ازاء همه موقعیتهای زاویه‌ای روتور ثابت است اما زاویه آن تغییر می‌کند.

### 3-4-3- ناهم محوری مرکب(M.E)

در این حالت هر سه محور تقارن روتور، تقارن استاتور و چرخش روتور نسبت به یکدیگر جابجا می‌گردند. این نوع ناهم محوری نتیجه اعمال برآیند یک بردار انتقال ایستا و یک بردار انتقال پویا می‌باشد. درجه ناهم محوری مرکب یا $\delta(\theta)$ تابعی از زاویه چرخش مکانیکی روتور است که بصورت رابطه زیر تعریف می‌شود.

$$\delta(\theta)=\frac{\left|\overrightarrow{O_S O_r}\right|}{g_0}=\frac{\left|\overrightarrow{O_S O_\omega}+\overrightarrow{O_\omega O_r}\right|}{g_0} \qquad (26-3)$$

شکل 3-7 وضعیت سطح مقطع روتور و استاتور را درحالت ناهم محوری مرکب نشان می دهد.

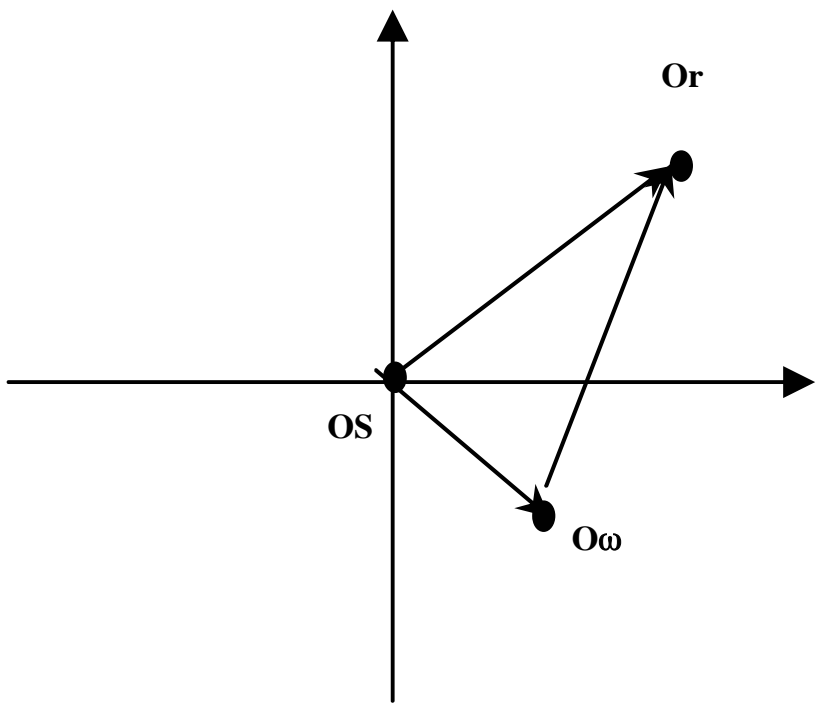

شکل 3-7: وضعیت سطح مقطع روتور و استاتور در حالت ناهم محوری مرکب در دستگاه مرجع استاتور

بردار $\overrightarrow{O_S O_r}$ را بردار انتقال ناهم محوری مرکب و زاویه آنرا زاویه انتقال ناهم محوری مرکب می‌نامند. اندازه و زاویه این بردار تابعی از زاویه مکانیکی روتور است.

## 5-3- محاسبه هدایت مغناطیسی فاصله هوایی در حالت ME

بنا به تعریف، هدایت مغناطیسی فاصله هوایی عکس رلوکتانس مغناطیسی فاصله هوایی است که با توجه به ابعاد هندسی فاصله هوایی بصورت زیر قابل محاسبه است.

$$P(\varphi) = \mu_0 \frac{r_{av}(\varphi)}{g(\varphi)} \tag{27-3}$$

جهت محاسبه این کمیت ابتدا بایستی توزیع فاصله هوایی پیرامون روتور و تغییرات آن نسبت به زاویه مکانیکی روتور را محاسبه کرد. می‌توان نشان داد که زاویه ناهم محوری ایستا همواره ثابت و برابر مقدار اولیه آن است ولی زاویه ناهم محوری پویا با تغییر زاویه مکانیکی روتور تغییر می‌کند.شکل 3-8 وضعیت محورهای تقارن و چرخش روتور را در دستگاه مرجع استاتور را در حالتی نشان می‌دهد که هر دو ناهم‌محوری ایستا و پویا بین روتور و استاتور وجود دارد. در این شکل بردار $\overrightarrow{O_s O_\omega}$ با اندازه $\delta_s$ و زاویه $\alpha_s$ بردار ناهم محوری ایستا و بردار $\overrightarrow{O_\omega O_r}$ با اندازه $\delta_d$ و زاویه $\alpha_d$ بردار ناهم محوری ایستا است. با توجه به شکل 3-8 بردار ترکیبی ناهم محوری یا $\overrightarrow{O_s O_r}$ بصورت زیر تعریف می‌شود:

$$\overrightarrow{O_s O_r} = \overrightarrow{O_s O_\omega} + \overrightarrow{O_\omega O_r} \tag{28-3}$$

در حالتی که زاویه ناهم محوری ایستا و زاویه اولیه ناهم محوری پویا صفر باشد داریم:

$$\delta = \left|\overrightarrow{O_s O_r}\right| = \sqrt{\delta_s^2 + \delta_r^2 + 2\delta_s \delta_r \cos\theta} \quad (29\text{-}3)$$

$$\alpha = \measuredangle\, \overrightarrow{O_s O_r} = \tan^{-1} \frac{\delta_d \sin\theta}{\delta_s + \delta_d \cos\theta} \quad (30\text{-}3)$$

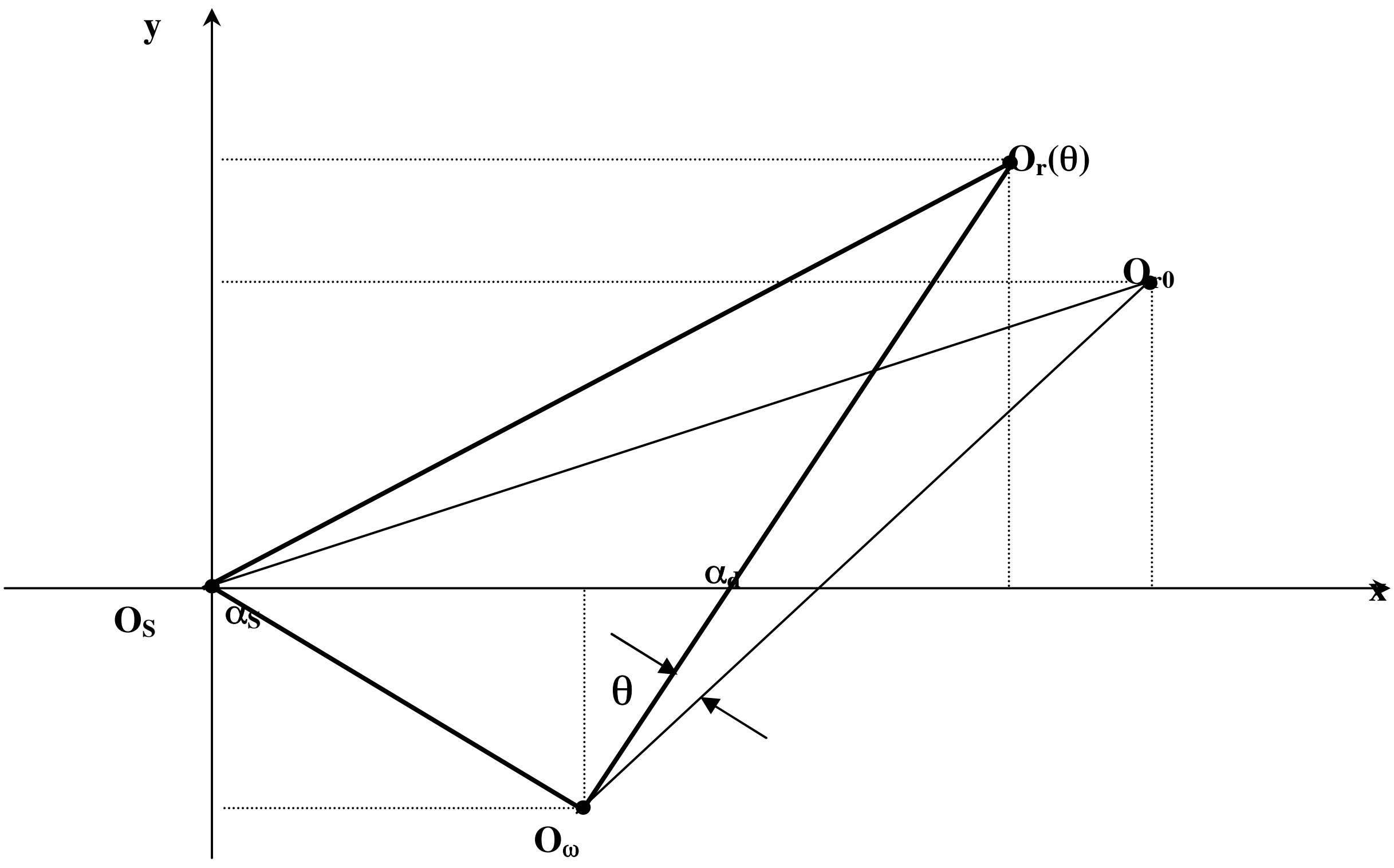

شکل 3-8: وضعیت محورهای تقارن و چرخش روتور در دستگاه مرجع استاتور در حالت ناهم محوری ترکیبی

با توجه به روابط بالا می توان در هر لحظه از زمان مختصات محورهای تقارن روتور و استاتور را به دست آورد و با کمک آنها معادله دایره‌های سطح مقطع روتور و استاتور را در دستگاه مختصات استاتور تشکیل داد. با کمک این معادله ها طول دقیق فاصله هوایی بصورت زیر قابل محاسبه است:

$$g(\varphi) = R_s - \delta g_0 \cos(\varphi - \alpha) - \sqrt{R_r^{\,2} - \delta^2 g_0^2 \sin(\varphi - \alpha)} \qquad (31\text{-}3)$$

با توجه به رابطه اخیر شعاع متوسط فاصله هوایی برابر است با:

$$r_{av}(\varphi) = R_s - \tfrac{1}{2} g(\varphi) \qquad (32\text{-}3)$$

بنابراین هدایت مغناطیسی فاصله هوایی در حالت ناهم محوری ترکیبی بصورت زیر می باشد:

$$P(\varphi) = \mu_0 \frac{0.5R_s + 0.5\delta g_0 \cos(\varphi - \alpha) + 0.5\sqrt{R_r^{\,2} - \delta^2 g_0^2 \sin(\varphi - \alpha)}}{R_s - \delta g_0 \cos(\varphi - \alpha) - \sqrt{R_r^{\,2} - \delta^2 g_0^2 \sin(\varphi - \alpha)}} \qquad (33\text{-}3)$$

از آنجا که طول فاصله هوایی نسبت به شعاع روتور خیلی کوچکتر است با تقریب بسیار خوبی می توان نوشت :

$$\sqrt{R_r^{\,2} - \delta^2 g_0^2 \sin(\varphi - \alpha)} \cong R_r \qquad (34\text{-}3)$$

بنابراین هدایت مغناطیسی فاصله هوایی در حالت ناهم محوری ترکیبی بصورت زیر ساده می شود.

$$P(\varphi) = \mu_0 \frac{0.5(R_s + R_r)}{(R_s - R_r) - \delta g_0 \cos(\varphi - \alpha)} \qquad (35\text{-}3)$$

از آنجا که در حالت متقارن برای طول و شعاع متوسط فاصله هوایی داریم :

(36-3)

$$P_0 = \mu_0 \frac{r_0}{g_0}$$

(37-3)

$$0.5(R_s + R_r) = r_0$$

(38-3)

$$R_s - R_r = g_0$$

می توان نوشت:

(39-3)

$$P(\varphi) = P_0 \frac{1}{1 - \delta \cos(\varphi - \alpha)}$$

رابطه اخیر بصورت سری نامتناهی زیر بسط داده می شود:

(40-3)

$$P(\varphi) = \frac{P_0}{\sqrt{1-\delta^2}} + \frac{2P_0}{\sqrt{1-\delta^2}} \sum_{k=0}^{\infty} \left( \frac{1-\sqrt{1-\delta^2}}{\delta} \right)^k \cos k(\phi - \alpha)$$

تاکنون در مطالعات انجام گرفته بسط نامتناهی فوق با دو جمله تقریب زده شده است [24] و [25]. در بررسی حاضر برای نخستنی بار جهت افزایش دقت محاسبات این بسط تا سه جمله تقریب زده می شود از این رو هدایت مغناطیسی فاصله هوایی بصورت رابطه زیر در می‌آید..

(41-3)

$$P(\varphi) = AP_0 + BP_0 \cos(\varphi - \alpha) + CP_0 \cos 2(\varphi - \alpha)$$

$$A = \frac{1}{\sqrt{1-\delta^2}}$$

$$B = \frac{2}{\sqrt{1-\delta^2}} \left( \frac{1-\sqrt{1-\delta^2}}{\delta} \right)$$

$$C = \frac{2}{\sqrt{1-\delta^2}}\left(\frac{1-\sqrt{1-\delta^2}}{\delta}\right)^2$$

فصل چهارم

# محاسبه اندوکتانسهای موتور القایی سه فاز قفس سنجابی در حالتهای مختلف ناهم‌محوری

نظریه تابع سیم پیچی قادر به محاسبه اندوکتانسهای مغناطیسی موتور است. در این فصل با توجه به تعمیم ارائه شده در فصل دوم به محاسبه اندوکتانسهای موتور القایی سه فاز روتور قفسی تحت شرایط مختلف ناهم محوری روتور و استاتور می پردازیم☜. از ویژگیهای بارز روش ارائه شده در این فصل میتوان به موارد زیر اشاره نمود:

1- روش ارائه شده کاملا تحلیلی است و مقدار اندوکتانسها را با عبارت جبری بسته ای ارائه می دهد .

بنابراین می توان کلیه اندوکتانسهای مغناطیسی موتور را با سرعت بالایی محاسبه نمود. با توجه به اینکه در حالت کلی ناهم محوری کلیه اندوکتانسهای مغناطیسی موتور یعنی اندوکتانسهای خودی و متقابل فازهای استاتور، اندوکتانسهای خودی و متقابل حلقه های روتور و اندوکتانسهای متقابل میان فازهای استاتور و حلقه های روتور تابعی از موقعیت زاویه ای روتور هستند، سایر روشها برای محاسبه مقدار هر اندوکتانس در موقعیت زاویه ای مورد نظر بایستی انتگرال رابطه (2-21) فصل دوم را بصورت عددی در چندین موقعیت زاویه ای روتورحل کرده و درحافظه کامپیوتر ذخیره کنند. در اینصورت روال نرم افزاری تحلیل موتور بایستی برای تعیین اندوکتانسهای مغناطیسی موتور در هر لحظه بر روی اطلاعات ذخیره شده عملیات درونیابی را اجرا نماید[24].

به عنوان مثال یک موتور القایی 3 فاز روتور قفسی با 40 میله روتور 1729 اندوکتانس مغناطیسی دارد بنابراین روال نرم افزاری تحلیل موتور باید برای تحلیل عملکرد موتور در هرگام زمانی 1729 بار عملیات درونیابی را انجام دهد. این مساله باعث وقت گیر بودن روشهای ذکر شده می گردد. روش ارائه شده در این فصل نیازی به درونیابی ندارد و مقدار اندوکتانس در هر لحظه از زمان را با انجام عملیات جبری ساده به دست می آورد.

2- روش ارائه شده هدایت مغناطیسی فاصله هوایی را دقیق تر از سایر روشها محاسبه می کند.

این مساله باعث افزایش دقت در محاسبه اندوکتانسهای موتور می شود. سایر روشها تنها از دو جمله اول سری نامتناهی هدایت مغناطیسی فاصله هوایی استفاده می کنند. روش ارائه شده سری نامتناهی مذکور را با سه جمله اول آن تقریب می زند. در فصل قبل به این موضوع پرداخته شد.

3- روش ارائه شده اثر افزایش خطی mmf فاصله هوایی در بالای شیارها را لحاظ می کند.

این مساله باعث می شود تا در محاسبه اندوکتانسهای مغناطیسی هارمونیکهای فضایی بیشتری لحاظ شوند.

☜ در اين گزارش كليه محاسبات بر اساس موتور القايي معرفي شده در ضميمه مي باشد.

4- روش ارائه شده از نظریه تابع سیم پیچی تعمیم یافته استفاده می کند.

تاکنون کلیه مقالاتی که با استفاده از این نظریه به تحلیل موتور القایی در حالتهای غیریکنواختی فاصله هوایی پرداخته‌اند از نظریه تابع سیم پیچی استفاده کرده‌اند که فقط در حالت فاصله هوایی یکنواخت معتبر است. این مساله باعث شده تا در محاسبات آنها به طور مشهودی تساوی $M_{12} = M_{21} = M$ وجود نداشته باشد.

با استفاده از تعمیم ارائه شده درفصل دوم خواهیم دید همانطور که از همه مدارهای مغناطیسی خطی انتظار میرود همواره با هر توزیع فاصله هوایی تساوی $M_{12} = M_{21} = M$ بر قرار است. در اینصورت نسبت به آنچه در سایر مقالات ارائه شده تردید اساسی وجود دارد[25],[24].

بنابراین روش ارائه شده در این گزارش به مراتب سریعتر و دقیقتر از سایر روشهایی که تاکنون ارائه شدهاند، میباشد. برای اینکه بتوان روابط اندوکتانسهای موتور را بصورت فرمولهای بسته ارائه کرد ابتدا چند تابع اساسی برای فازهای استاتور و حلقه های روتور تعریف می کنیم.

## 4-1- تابع اساسی دور استاتور (AS)

یک دوره تناوب از تابع دور فاز اول استاتور را تابع دور اساسی استاتور می نامیم. تابع دور بقیه فازهای استاتور از تناوب و جابجایی این تابع به دست می آید. برای یک موتور القایی سه فاز چهارقطبی در صورتیکه سیم‌پیچی هر قطب فازهای استاتور $2N$ دور داشته باشد، تابع اساسی دور استاتور بصورت شکل 4-1 می‌باشد.

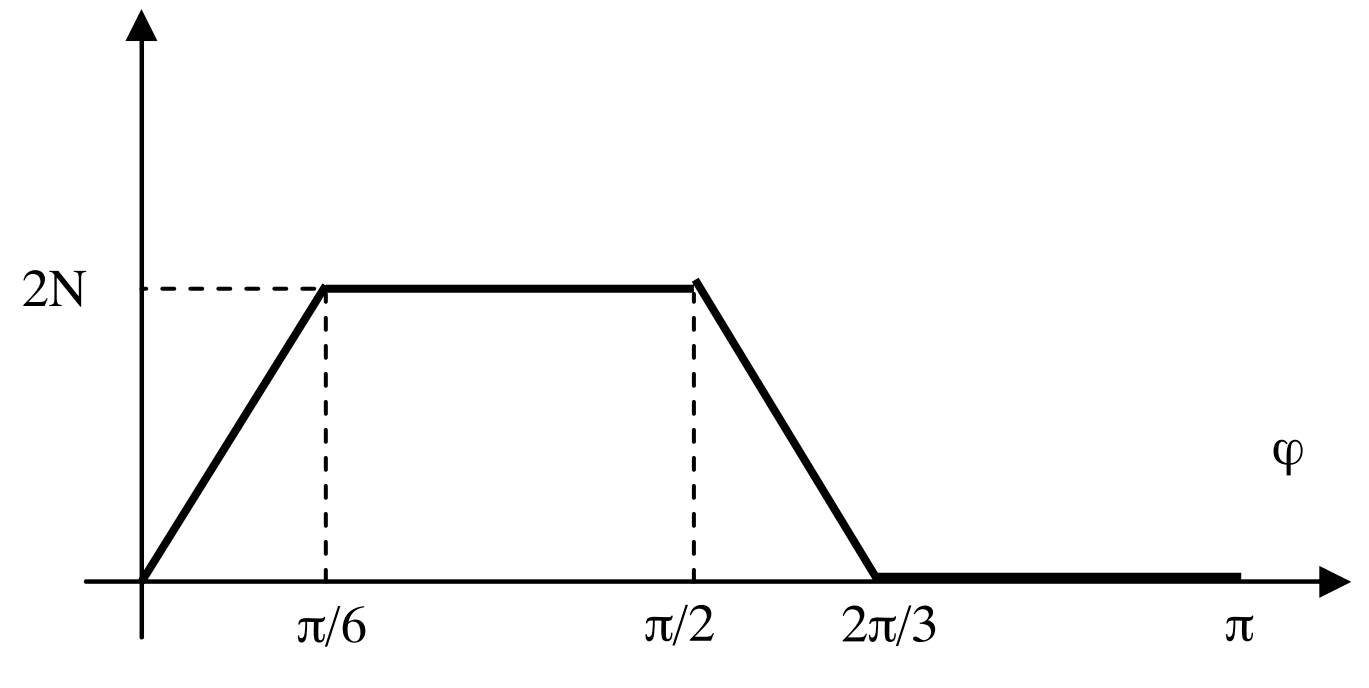

**شکل 4-1: تابع اساسی دور استاتور یک موتور القایی سه فاز چهار قطب**

می توان این تابع را با عبارت جبری زیر نشان داد:

$$AS(\varphi)=\begin{cases}\frac{12N}{\pi}\varphi & when\, 0\le\varphi<\frac{\pi}{6}\\ 2N & when\,\frac{\pi}{6}\le\varphi<\frac{\pi}{2}\\ 8N-\frac{12N}{\pi}\varphi & when\,\frac{\pi}{2}\le\varphi<\frac{2\pi}{3}\\ 0 & when\,\frac{2\pi}{3}\le\varphi<\pi\end{cases} \tag{1-4}$$

تابع دور فازهای استاتور با توجه به تابع $AS(\varphi)$ و باکمک رابطه زیر به دست می آید.

$$AS_i(\varphi)=AS\left(\varphi-\frac{2\pi(i-1)}{3}\right)\qquad i=1,2,3 \tag{2-4}$$

تابع اساسی دور استاتور را میتوان با سری فوریه زیر نشان داد:

$$AS(\varphi)=a_0^{AS}+\sum_{k=1}^{\infty}a_k^{AS}\cos(2k\varphi)+b_k^{AS}\sin(2k\varphi) \tag{3-4}$$

در اینصورت ضرایب سری فوق بصورت زیر می باشند:

$$a_0^{AS}=N \tag{4-4}$$

$$a_k^{AS}=-\frac{6N}{(\pi k)^2}\left[1-(-1)^k\right]\left(1-\cos\frac{k\pi}{3}\right) \tag{5-4}$$

$$b_k^{AS} = -\frac{6N}{(\pi k)^2}\left[1-(-1)^k\right]\sin\frac{k\pi}{3} \qquad (6\text{-}4)$$

محاسبه این ضرایب نشان می دهد که سری فوریه فوق فقط دارای ضرایب $6k \pm 1$ می باشد.

## 2-4- تابع اساسی مربع دور استاتور (BS)

این تابع با توجه به رابطه زیر تعریف می شود:

$$BS(\varphi) = AS^2(\varphi) \qquad (7\text{-}4)$$

بنابراین :

$$BS(\varphi) = \begin{cases} \frac{144N^2}{\pi^2}\varphi^2 & when\ 0 \le \varphi < \frac{\pi}{6} \\ 4N^2 & when\ \frac{\pi}{6} \le \varphi < \frac{\pi}{2} \\ \left(8N - \frac{12N}{\pi}\varphi\right)^2 & when\ \frac{\pi}{2} \le \varphi < \frac{2\pi}{3} \\ 0 & when\ \frac{2\pi}{3} \le \varphi < \pi \end{cases} \qquad (8\text{-}4)$$

شکل زیر نمودار این تابع را نشان می دهد.

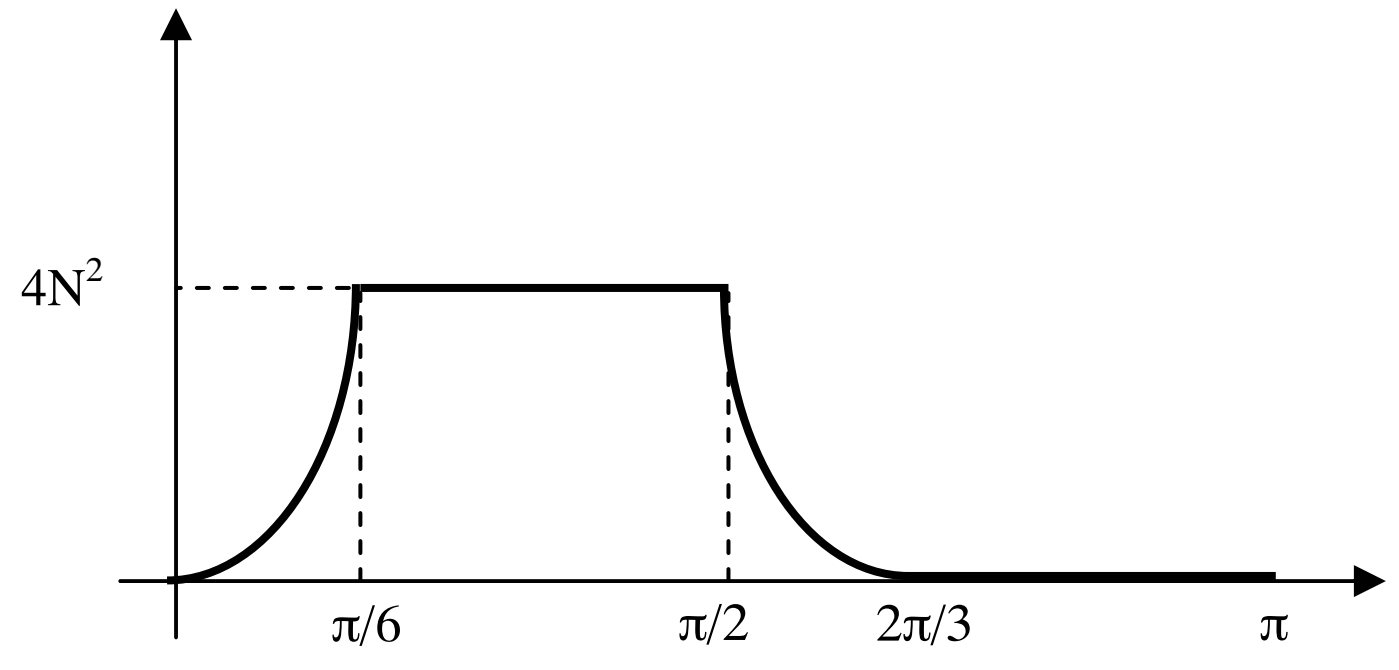

شکل4-2: تابع اساسی مربع دور استاتور یک موتور القایی سه فاز چهار قطب

تابع مربع دور فازهای استاتور با توجه به تابع $BS(\varphi)$ و باکمک رابطه زیر به دست می آید.

$$BS_i(\varphi)=BS\left(\varphi-\frac{2\pi(i-1)}{3}\right)\quad i=1,2,3 \tag{9-4}$$

تابع اساسی مربع دور استاتور را میتوان با سری فوریه  زیر نشان داد:

$$BS(\varphi)=a_0^{BS}+\sum_{k=1}^{\infty}a_k^{BS}\cos(2k\varphi)+b_k^{BS}\sin(2k\varphi) \tag{10-4}$$

در اینصورت ضرایب سری فوق بصورت زیر می باشند:

$$a_0^{BS}=\frac{16}{9}N^2 \tag{11-4}$$

$$a_k^{BS}=\frac{24N^2}{(\pi k)^2}\left(\cos\frac{k\pi}{3}+\cos k\pi\right)-\frac{72N^2}{(\pi k)^3}\left[1+(-1)^k\right]\sin\frac{k\pi}{3} \tag{12-4}$$

$$b_k^{BS}=\frac{24N^2}{(\pi k)^2}\sin\frac{k\pi}{3}-\frac{72N^2}{(\pi k)^3}\left[1+(-1)^k\right]\left(1-\cos\frac{k\pi}{3}\right) \tag{13-4}$$

## 3-4- تابع اساسی حاصلضرب دور استاتور (CS)

این تابع با توجه به رابطه زیر تعریف می شود:

$$CS(\varphi)=AS(\varphi)\times AS\left(\varphi-\frac{2\pi}{3}\right) \tag{14-4}$$

بنابراین :

$$CS(\varphi)=\begin{cases}\frac{24N^2}{\pi}\varphi & when\,0\le\varphi<\frac{\pi}{6}\\ 8N^2-\frac{24N^2}{\pi}\varphi & when\,\frac{\pi}{6}\le\varphi<\frac{\pi}{3}\\ 0 & when\,\frac{\pi}{3}\le\varphi<\pi\end{cases} \quad (15\text{-}4)$$

شکل 4-3 نمودار این تابع را نشان می دهد.

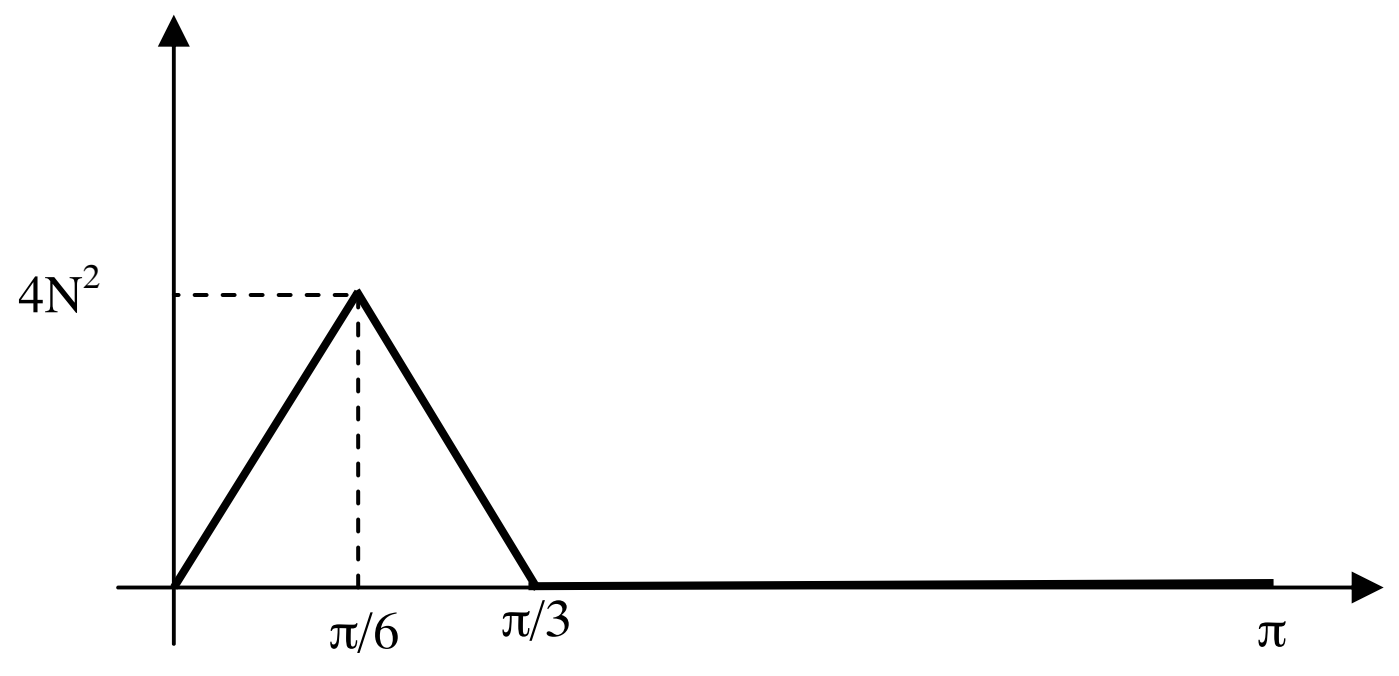

شکل **4-3**: تابع اساسی حاصلضرب دور استاتور یک موتور القایی سه فاز چهار قطب

تابع حاصلضرب دور فازهای $i$ و $j$ استاتور با توجه به تابع $CS(\varphi)$ و باکمک رابطه زیر به دست می آید.

$$CS_{i,j}(\varphi)=CS\left(\varphi-\frac{(i+j)\pi}{3}+\pi\right)\qquad i,j=1,2,3\quad i\neq j \quad (16\text{-}4)$$

تابع اساسی حاصلضرب دور استاتور را میتوان با سری فوریه زیر نشان داد:

$$CS(\varphi) = a_0^{CS} + \sum_{k=1}^{\infty} a_k^{CS} \cos(2k\varphi) + b_k^{CS} \sin(2k\varphi) \tag{17-4}$$

در اینصورت ضرایب سری فوق بصورت زیر می باشند:

$$a_0^{CS} = \frac{2}{3}N^2 \tag{18-4}$$

$$a_k^{CS} = \frac{-12N^2}{(\pi k)^2}\left(1 - 2\cos\frac{k\pi}{3} + \cos\frac{2k\pi}{3}\right) \tag{19-4}$$

$$b_k^{CS} = \frac{-12N^2}{(\pi k)^2}\left(-2\sin\frac{k\pi}{3} + \sin\frac{2k\pi}{3}\right) \tag{20-4}$$

## 4-4- تابع اساسی دور روتور (AR)

یک دوره تناوب از تابع دور حلقه اول روتور را تابع دور اساسی روتور می نامیم. تابع دور بقیه حلقه های روتور از تناوب و جابجایی این تابع به دست می آید. برای یک موتور القایی با $n$ میله در قفس روتور تابع اساسی دور روتور بصورت شکل 4-4 میباشد.

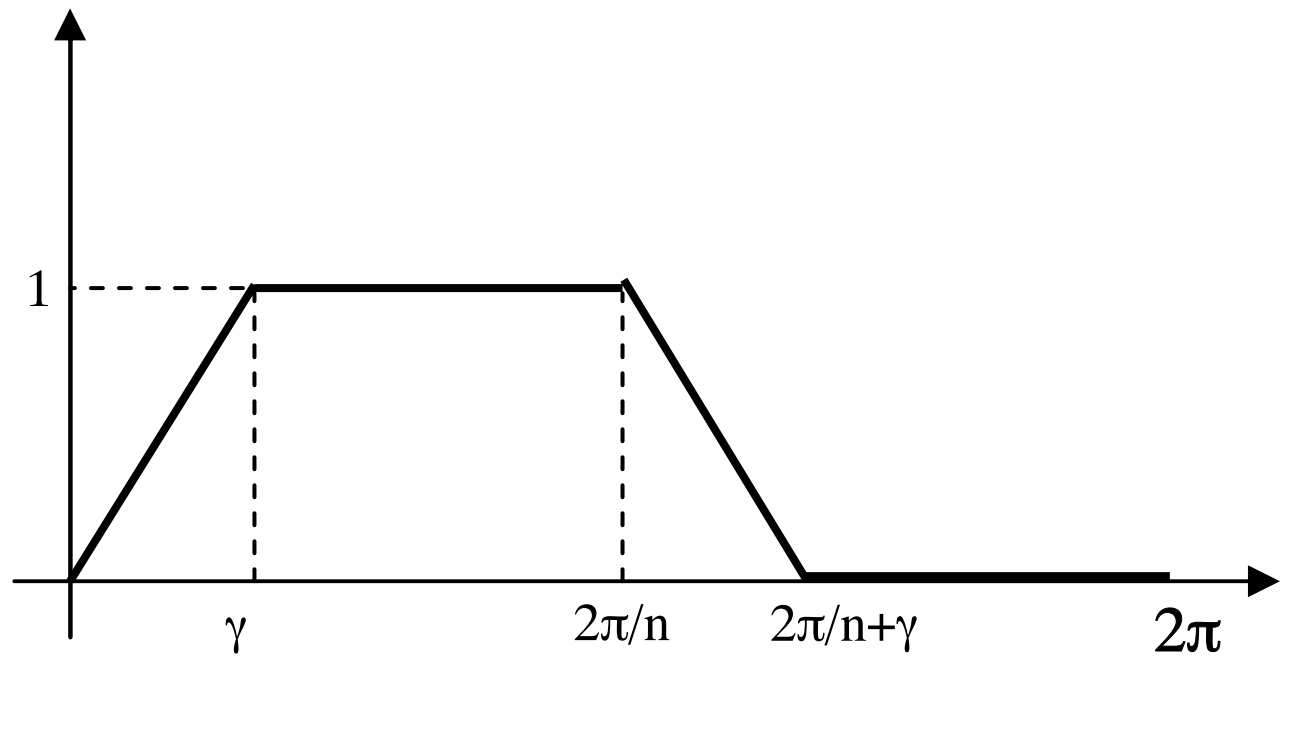

**شکل 4-4**: تابع اساسی دور روتور یک موتور القایی

می‌توان این تابع را با عبارت جبری زیر نشان داد:

$$AR(\varphi)=\begin{cases}\frac{1}{\gamma}\varphi & when\,0\le\varphi<\gamma\\ 1 & when\,\frac{1}{\gamma}\le\varphi<\frac{2\pi}{n}\\ 1-\frac{1}{\gamma}\left(\varphi-\frac{2\pi}{n}\right) & when\,\frac{2\pi}{n}\le\varphi<\frac{2\pi}{n}+\gamma\\ 0 & when\,\frac{2\pi}{n}+\gamma\le\varphi<2\pi\end{cases} \qquad (21\text{-}4)$$

تابع دور حلقه های روتور با توجه به تابع $AR(\varphi)$ و باکمک رابطه زیر به دست می آید.

$$AR_i(\varphi)=AR\left(\varphi-\frac{2\pi(i-1)}{n}\right)\qquad i=1,2,3,...,n \qquad (22\text{-}4)$$

تابع اساسی دور روتور را میتوان با سری فوریه  زیر نشان داد:

$$AR(\varphi)=a_0^{AR}+\sum_{k=1}^{\infty}a_k^{AR}\cos(k\varphi)+b_k^{AR}\sin(k\varphi) \qquad (23\text{-}4)$$

در اینصورت ضرایب سری فوق بصورت زیر می باشند:

$$a_0^{AR}=\frac{1}{n} \qquad (24\text{-}4)$$

$$a_k^{AR}=\frac{4}{\pi k^2\gamma}\sin\frac{k\gamma}{2}\sin\frac{k\pi}{n}\cos\left(\frac{k\pi}{n}+\frac{k\gamma}{2}\right) \qquad (25\text{-}4)$$

$$b_k^{AR}=\frac{4}{\pi k^2\gamma}\sin\frac{k\gamma}{2}\sin\frac{k\pi}{n}\sin\left(\frac{k\pi}{n}+\frac{k\gamma}{2}\right) \qquad (26\text{-}4)$$

## 5-4- تابع اساسی مربع دور استاتور (BR)

این تابع با توجه به رابطه زیر تعریف می شود:

$$BR(\varphi)=AR^2(\varphi) \qquad (27\text{-}4)$$

بنابراین :

$$BR(\varphi)=\begin{cases}\left(\frac{1}{\gamma}\varphi\right)^2 & \text{when } 0\le\varphi<\gamma\\ 1 & \text{when } \frac{1}{\gamma}\le\varphi<\frac{2\pi}{n}\\ \left(1-\frac{1}{\gamma}\left(\varphi-\frac{2\pi}{n}\right)\right)^2 & \text{when } \frac{2\pi}{n}\le\varphi<\frac{2\pi}{n}+\gamma\\ 0 & \text{when } \frac{2\pi}{n}+\gamma\le\varphi<2\pi\end{cases} \quad (28\text{-}4)$$

شکل 4-5 نمودار این تابع را نشان می دهد.

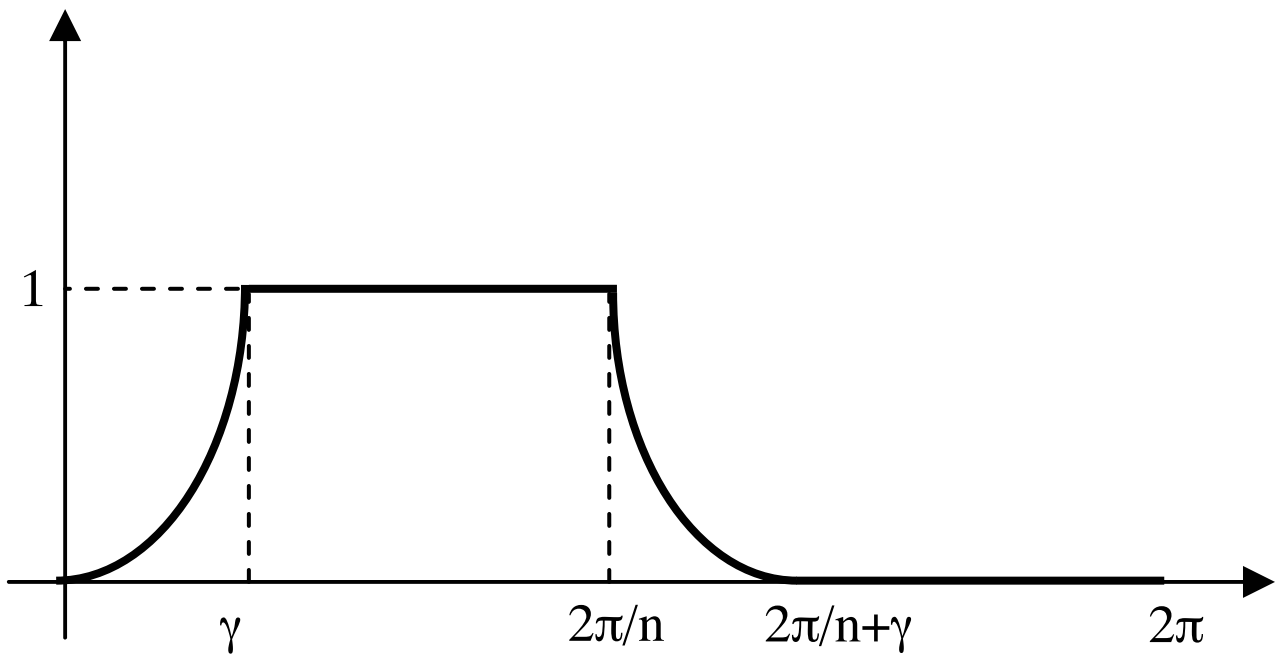

**شکل 4-5**: تابع اساسی مربع دور روتور یک موتور القایی

تابع مربع دور حلقه های روتور با توجه به تابع $BR(\varphi)$ و باکمک رابطه زیر به دست می آید.

$$BR_i(\varphi)=BR\left(\varphi-\frac{2\pi(i-1)}{3}\right) \quad i=1,2,3,\ldots,n \quad (29\text{-}4)$$

تابع اساسی مربع دور روتور را میتوان با سری فوریه زیر نشان داد:

$$BR(\varphi)=a_0^{BR}+\sum_{k=1}^{\infty}a_k^{BR}\cos(k\varphi)+b_k^{BR}\sin(k\varphi) \quad (30\text{-}4)$$

در اینصورت ضرایب سری فوق بصورت زیر می باشند:

$$a_0^{BR} = \frac{1}{n} - \frac{\gamma}{6\pi} \quad (31\text{-}4)$$

$$a_k^{BR} = \frac{4}{\pi k^2 \gamma}\left[cos\left(\frac{k\pi}{n} - \frac{k\gamma}{2}\right) - \frac{2}{k\gamma} cos\frac{k\pi}{n} sin\frac{k\gamma}{2}\right] cos\left(\frac{k\pi}{n} + \frac{k\gamma}{2}\right) \quad (32\text{-}4)$$

$$b_k^{BR} = \frac{4}{\pi k^2 \gamma}\left[cos\left(\frac{k\pi}{n} - \frac{k\gamma}{2}\right) - \frac{2}{k\gamma} cos\frac{k\pi}{n} sin\frac{k\gamma}{2}\right] sin\left(\frac{k\pi}{n} + \frac{k\gamma}{2}\right) \quad (33\text{-}4)$$

## 6-4- تابع اساسی حاصلضرب دور مجاوری روتور (CR)

این تابع با توجه به رابطه زیر تعریف می شود:

$$CR(\varphi) = AR(\varphi) \times AR\left(\varphi + \frac{2\pi}{n}\right) \quad (34\text{-}4)$$

بنابراین :

$$CR(\varphi) = \begin{cases} \frac{1}{\gamma}\varphi - \frac{1}{\gamma^2}\varphi^2 & when\ 0 \le \varphi < \gamma \\ 0 & when\ \gamma \le \varphi < 2\pi \end{cases} \quad (35\text{-}4)$$

شکل 4-6 نمودار این تابع را نشان می دهد.

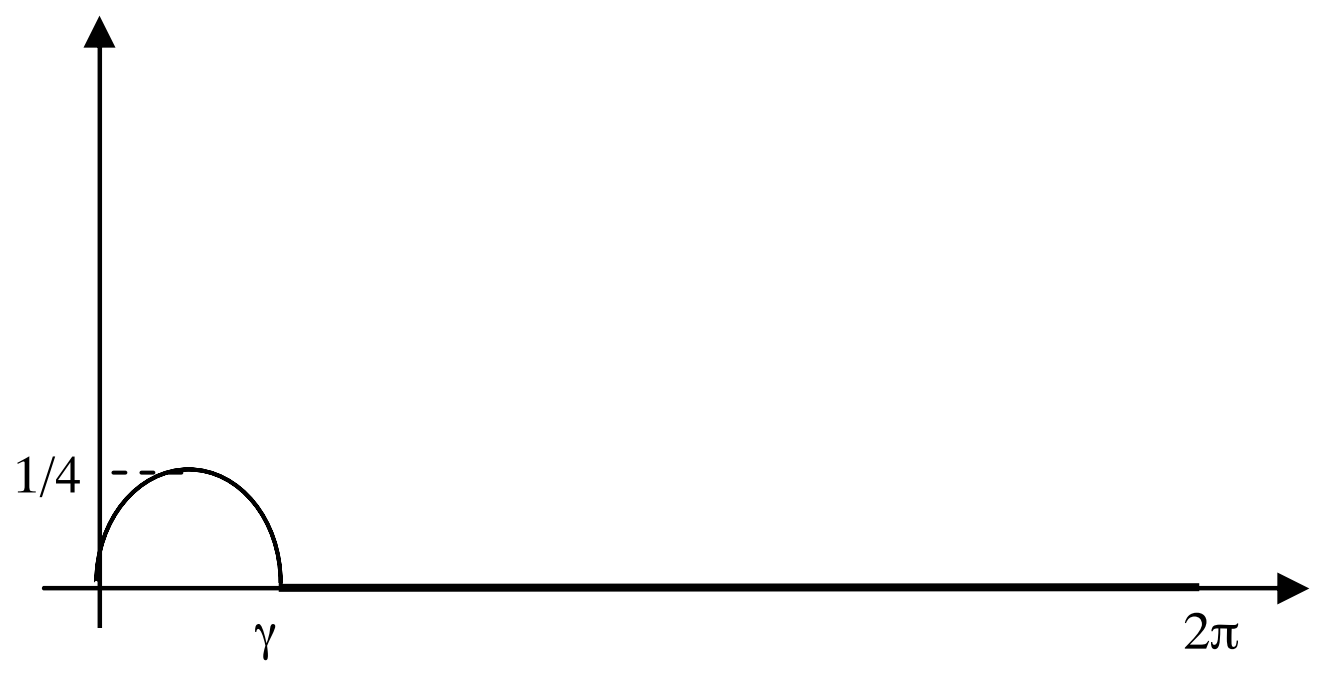

شکل **4-6**: تابع اساسی حاصلضرب دور روتور یک موتور القایی سه فاز چهار قطب

تابع حاصلضرب دور حلقه های های $i$ و $i-1$ روتور با توجه به تابع $CR(\varphi)$و باکمک رابطه زیر به دست میآید.

$$CR_i(\varphi) = CR\left(\varphi - \frac{(i-1)2\pi}{n}\right) \qquad i = 1,2,3,\dots,n \tag{36-4}$$

تابع اساسی حاصلضرب دور مجاوری روتور را میتوان با سری فوریه  زیر نشان داد:

$$CR(\varphi) = a_0^{CR} + \sum_{k=1}^{\infty} a_k^{CR} \cos(k\varphi) + b_k^{CR} \sin(k\varphi) \tag{37-4}$$

در اینصورت ضرایب سری فوق بصورت زیر می باشند:

$$a_0^{CR} = \frac{\gamma}{12\pi} \tag{38-4}$$

$$a_k^{CR} = \frac{2}{\pi\gamma k^2}\left(\frac{2}{k\gamma}\sin\frac{k\gamma}{2} - \cos\frac{k\gamma}{2}\right)\cos\frac{k\gamma}{2} \tag{39-4}$$

$$b_k^{CR} = \frac{2}{\pi\gamma k^2}\left(\frac{2}{k\gamma}\sin\frac{k\gamma}{2} - \cos\frac{k\gamma}{2}\right)\sin\frac{k\gamma}{2} \tag{40-4}$$

تابع اساسی حاصلضرب دور غیر مجاوری روتور همواره متحد با صفر است.

## 4-7- محاسبه اجزاء اصلی اندوکتانسهای موتور در حالتهای ناهم محوری

اندوکتانس متقابل دو سیم پیچی دلخواه $x$ و $y$ از طریق رابطه زیر محاسبه می شود.

$$L_{yx} = 2\pi l\langle Pn_x n_y\rangle - 2\pi l\frac{\langle Pn_x\rangle\langle Pn_y\rangle}{\langle P\rangle} \qquad (41\text{-}4)$$

همچنین در فصل سوم هدایت مغناطیسی فاصله هوایی بصورت زیر در نظر گرفته شد:

$$P(\varphi) = AP_0 + BP_0\cos(\varphi-\alpha) + CP_0\cos 2(\varphi-\alpha) \qquad (42\text{-}4)$$

صرایب $A$، $B$ و $C$ و همچنین زاویه $\alpha$ نیز در فصل دوم تعریف شدند، همچنین نشان داده شده که این پارامترها در حالت کلی تابعی از $\theta$ می باشند.

برای محاسبه تحلیلی رابطه (4-41) بایستی اجزاء اصلی این رابطه را یعنی $\langle Pn_x\rangle$ و $\langle Pn_x n_y\rangle$ را بصورت تحلیلی برای تمامی فازهای استاتور و حلقه های روتور محاسبه کنیم. بنابراین به عنوان قدم بعدی بایستی عبارت ذکر شده در جدول زیر را بصورت تحلیلی محاسبه کرد.

جدول 4-1 : اجزاء اصلی اندوکتانسهای مغناطیسی موتور القایی سه فاز قفس سنجابی

| | |
|---|---|
| میانگین وزندار تابع دور فازهای استاتور با وزن هدایت مغناطیسی فاصله هوایی | $\langle P.AS_i \rangle$ |
| میانگین وزندار تابع مربع دور فازهای استاتور با وزن هدایت مغناطیسی فاصله هوایی | $\langle P.BS_i \rangle$ |
| میانگین وزندار تابع حاصلضرب دور فازهای استاتور با وزن هدایت مغناطیسی فاصله هوایی | $\langle P.CS_{i,j} \rangle$ |
| میانگین وزندار تابع دور حلقه های روتور با وزن هدایت مغناطیسی فاصله هوایی | $\langle P.AR_i \rangle$ |
| میانگین وزندار تابع مربع دور حلقه های روتور با وزن هدایت مغناطیسی فاصله هوایی | $\langle P.BR_i \rangle$ |
| میانگین وزندار تابع حاصلضرب دور مجاوری حلقه های روتور با وزن هدایت مغناطیسی فاصله هوایی | $\langle P.CR_i \rangle$ |
| میانگین وزندار حاصلضرب تابع دور حلقه های روتور در تابع دور فازهای استاتور با وزن هدایت مغناطیسی فاصله هوایی | $\langle P.AS_i.AR_j \rangle$ |

در ادامه ضمن معرفی یک قضیه به محاسبه توابع ذکر شده پرداخته می شود.

*قضيه اساسي:*

*فرض :*

*تابع* $f(x)$ *بصورت سري فوريه* $f(x) = a_0 + \sum_{k=1}^{\infty} a_k \cos pkx + b_k \sin pkx$ *بسط داده مي شود.*

*تابع* $P(x)$ *بصورت كلي* $P(x) = P_0 + \sum_{m=1}^{M} P_i \cos m(x-\alpha)$ *است.*

*حكم :*

*ميانگين وزندار تابع* $f(x-x_0)$ *با وزن* $P(x)$ *بصورت زير قابل محاسبه است:*

$$\langle P(x).f(x-x_0) \rangle = P_0 a_0 + \sum_{i=1}^{M} P_{i\times p} a_i \cos(ip\alpha - ipx_0) + P_{i\times p} b_i \sin(ip\alpha - ipx_0)$$

الف : محاسبه $\langle P.AS_i \rangle$

$$AS_i(\varphi) = AS(\varphi - \frac{2(i-1)\pi}{3}) \Rightarrow \begin{cases} p = 2 \\ \varphi_0 = \frac{2(i-1)\pi}{3} \end{cases} \tag{43-4}$$

$$P(\varphi) = AP_0 + BP_0 \cos(\varphi - \alpha) + CP_0 \cos(2\varphi - 2\alpha) \Rightarrow M = 2$$

بنابراین با توجه به قضیه اصلی و رابطه (4-3) خواهیم داشت:

$$\langle P.AS_i \rangle = P_0 A a_0^{AS} + \frac{P_0 C}{2}\left(a_1^{AS} \cos 2(\alpha - \varphi_0) + b_1^{AS} \sin 2(\alpha - \varphi_0)\right) \tag{44-4}$$

با توجه به روابط (4-4) ،(4-5) ، (4-6) ، (4-43) و اعمال ساده سازیهای جبری و مثلثاتی خواهیم داشت :

$$\langle P.AS_i \rangle = P_0 AN + \frac{6NP_0 C}{\pi^2} \cos 2\left(\alpha - \frac{(2i-1)\pi}{3}\right) \tag{45-4}$$

ب : محاسبه $\langle P.BS_i \rangle$

$$BS_i(\varphi) = BS(\varphi - \frac{2(i-1)\pi}{3}) \Rightarrow \begin{cases} p = 2 \\ \varphi_0 = \frac{2(i-1)\pi}{3} \end{cases} \tag{46-4}$$

$$P(\varphi) = AP_0 + BP_0 \cos(\varphi - \alpha) + CP_0 \cos(2\varphi - 2\alpha) \Rightarrow M = 2$$

بنابراین با توجه به قضیه اصلی و رابطه (4-10) خواهیم داشت:

$$\langle P.BS_i \rangle = P_0 A a_0^{BS} + \frac{P_0 C}{2}\left(a_1^{BS} \cos 2(\alpha - \varphi_0) + b_1^{BS} \sin 2(\alpha - \varphi_0)\right) \tag{47-4}$$

با توجه به روابط (4-11) ، (4-12) ، (4-13) ، (4-47) و اعمال ساده سازیهای جبری و مثلثاتی خواهیم داشت :

$$\left\langle P.BS_i \right\rangle = \frac{16P_0AN^2}{9} + \frac{12N^2P_0C}{\pi^2}\cos 2\left(\alpha - \frac{(2i-1)\pi}{3}\right) \tag{48-4}$$

ج : محاسبه $\left\langle P.CS_{i,j} \right\rangle$

$$CS_{i,j}(\varphi) = CS(\varphi - \frac{(i+j)\pi}{3} + \pi) \Rightarrow \begin{cases} p = 2 \\ \varphi_0 = \frac{(i+j)\pi}{3} - \pi \end{cases} \tag{49-4}$$

$$P(\varphi) = AP_0 + BP_0\cos(\varphi - \alpha) + CP_0\cos(2\varphi - 2\alpha) \Rightarrow M = 2$$

بنابراین با توجه به قضیه اصلی و رابطه (4-17) خواهیم داشت:

$$\left\langle P.CS_{i,j} \right\rangle = P_0Aa_0^{CS} + \frac{P_0C}{2}\left(a_1^{CS}\cos 2(\alpha - \varphi_0) + b_1^{CS}\sin 2(\alpha - \varphi_0)\right) \tag{50-4}$$

با توجه به روابط (4-18) ، (4-19) ، (4-20) ، (4-50) و اعمال ساده سازیهای جبری و مثلثاتی خواهیم داشت :

$$\left\langle P.CS_{i,j} \right\rangle = \frac{2P_0AN^2}{3} - \frac{6N^2P_0C}{\pi^2}\cos 2\left(\alpha - \frac{(i+j-1)\pi}{3}\right) \tag{51-4}$$

د : محاسبه $\langle P.AR_i \rangle$

(52-4)

$$AR_i(\varphi) = AR(\varphi - \frac{(i-1)2\pi}{n}) \Rightarrow \begin{cases} p = 1 \\ \varphi_0 = \dfrac{(i-1)2\pi}{n} \end{cases}$$

$$P(\varphi) = AP_0 + BP_0 \cos(\varphi - \alpha) + CP_0 \cos(2\varphi - 2\alpha) \Rightarrow M = 2$$

بنابراین با توجه به قضیه اصلی و رابطه (4-23) خواهیم داشت:

(53-4)

$$\langle P.AR_i \rangle = P_0 A a_0^{AR} + \frac{P_0 B}{2}\left(a_1^{AR} \cos(\alpha - \varphi_0) + b_1^{AR} \sin(\alpha - \varphi_0)\right) + \frac{P_0 C}{2}\left(a_2^{AR} \cos 2(\alpha - \varphi_0) + b_2^{AR} \sin 2(\alpha - \varphi_0)\right)$$

با توجه به روابط (4-24) ، (4-25) ، (4-26) ، (4-53) و اعمال ساده سازیهای جبری و مثلثاتی خواهیم داشت :

(54-4)

$$\langle P.AR_i \rangle = \frac{P_0 A}{n} + \frac{2P_0 B}{\pi\gamma} \sin\frac{\gamma}{2} \sin\frac{\pi}{n} \cos\left(\alpha - \frac{2(i-1)\pi}{n} - \frac{\pi}{n} - \frac{\gamma}{2}\right) + \frac{P_0 C}{2\pi\gamma} \sin\gamma \sin\frac{2\pi}{n} \cos 2\left(\alpha - \frac{2(i-1)\pi}{n} - \frac{\pi}{n} - \frac{\gamma}{2}\right)$$

ه : محاسبه $\langle P.BR_i \rangle$

(55-4)

$$BR_i(\varphi) = BR(\varphi - \frac{(i-1)2\pi}{n}) \Rightarrow \begin{cases} p = 1 \\ \varphi_0 = \dfrac{(i-1)2\pi}{n} \end{cases}$$

$$P(\varphi) = AP_0 + BP_0 \cos(\varphi - \alpha) + CP_0 \cos(2\varphi - 2\alpha) \Rightarrow M = 2$$

بنابراین با توجه به قضیه اصلی و رابطه (4-40) خواهیم داشت:

$$\langle P.BR_i\rangle = P_0 A a_0^{BR} + \frac{P_0 B}{2}\left(a_1^{BR}\cos(\alpha-\varphi_0)+b_1^{BR}\sin(\alpha-\varphi_0)\right) + \frac{P_0 C}{2}\left(a_2^{BR}\cos 2(\alpha-\varphi_0)+b_2^{BR}\sin 2(\alpha-\varphi_0)\right) \tag{56-4}$$

با توجه به روابط (4-31) ، (4-32) ، (4-33) ، (4-56) و اعمال ساده سازیهای جبری و مثلثاتی خواهیم داشت

:

$$\langle P.BR_i\rangle = \left(\frac{1}{n}-\frac{\gamma}{6\pi}\right)P_0 A + \frac{2P_0 B}{\pi\gamma}\left[\cos\left(\frac{\pi}{n}-\frac{\gamma}{2}\right)-\frac{2}{\gamma}\cos\frac{\pi}{n}\sin\frac{\gamma}{2}\right]\cos\left(\alpha-\frac{2(i-1)\pi}{n}-\frac{\pi}{n}-\frac{\gamma}{2}\right) + \frac{P_0 C}{2\pi\gamma}\left[\cos\left(\frac{2\pi}{n}-\gamma\right)-\frac{1}{\gamma}\cos\frac{2\pi}{n}\sin\gamma\right]\cos 2\left(\alpha-\frac{2(i-1)\pi}{n}-\frac{\pi}{n}-\frac{\gamma}{2}\right) \tag{57-4}$$

و: محاسبه $\langle P.CR_i\rangle$

$$CR_i(\varphi) = CR\left(\varphi-\frac{(i-1)2\pi}{n}\right) \Rightarrow \begin{cases} p=1 \\ \varphi_0 = \dfrac{(i-1)2\pi}{n} \end{cases} \tag{58-4}$$

$$P(\varphi) = AP_0 + BP_0\cos(\varphi-\alpha)+CP_0\cos(2\varphi-2\alpha) \Rightarrow M=2$$

بنابراین با توجه به قضیه اصلی و رابطه (4-37) خواهیم داشت:

$$\langle P.CR_i\rangle = P_0 A a_0^{CR} + \frac{P_0 B}{2}\left(a_1^{CR}\cos(\alpha-\varphi_0)+b_1^{CR}\sin(\alpha-\varphi_0)\right) + \frac{P_0 C}{2}\left(a_2^{CR}\cos 2(\alpha-\varphi_0)+b_2^{CR}\sin 2(\alpha-\varphi_0)\right) \tag{59-4}$$

با توجه به روابط (4-38) ، (4-39) ، (4-40) ، (4-59) و اعمال ساده سازیهای جبری و مثلثاتی خواهیم داشت:

$$\left\langle P.CR_i \right\rangle = \frac{P_0 A\gamma}{12\pi} + \frac{P_0 B}{\pi\gamma}\left[\frac{2}{\gamma}\sin\frac{\gamma}{2} - \cos\frac{\gamma}{2}\right]\cos\left(\alpha - \frac{2(i-1)\pi}{n} - \frac{\gamma}{2}\right) + \frac{P_0 C}{4\pi\gamma}\left[\frac{1}{\gamma}\sin\gamma - \cos\gamma\right]\cos 2\left(\alpha - \frac{2(i-1)\pi}{n} - \frac{\gamma}{2}\right) \tag{60-4}$$

ز: محاسبه $\left\langle P.AS_i.AR_j \right\rangle$

$$\left\langle P.AS_i.AR_j \right\rangle = \frac{1}{2\pi}\int_0^{2\pi} P(\varphi)AS\left(\varphi - \frac{2(i-1)\pi}{3}\right)AR\left(\varphi - \frac{2(j-1)\pi}{n} - \theta\right)d\varphi = \frac{1}{2\pi}\int_e^{f} P(\varphi)AS\left(\varphi - \frac{2(i-1)\pi}{3}\right)AR\left(\varphi - \frac{2(j-1)\pi}{n} - \theta\right)d\varphi \tag{61-4}$$

که در این رابطه :

$$e = \frac{2(j-1)\pi}{n} + \theta \tag{62-4}$$

$$f = \frac{2j\pi}{n} + \theta + \gamma \tag{63-4}$$

از آنجا که $f - e = \frac{2\pi}{n} + \gamma$ مقدار کوچکی است می توان در انتگرال رابطه (4-61) از تغییرات $P(\varphi)$ در بازه $[e,f]$ صرفنظر کرد و به جای آن از مقدار این تابع در نقطه $\frac{f+e}{2}$ استفاده نمود. بنابراین خواهیم داشت:

$$\left\langle P.AS_i.AR_j \right\rangle = \frac{1}{2\pi} P\left(\theta + \frac{(2j-1)\pi}{n} + \frac{\gamma}{2}\right) \int_0^{2\pi} AS\left(\varphi - \frac{2(i-1)\pi}{3}\right) AR\left(\varphi - \frac{2(j-1)\pi}{n} - \theta\right) d\varphi \quad (64\text{-}4)$$

حال باید انتگرال موجود در رابطه (4-64) محاسبه گردد.

$$X_{i.j} = \int_0^{2\pi} AS\left(\varphi - \frac{2(i-1)\pi}{3}\right) AR\left(\varphi - \frac{2(j-1)\pi}{n} - \theta\right) d\varphi \quad (65\text{-}4)$$

با توجه به رابطه (4-1) و (4-22) و قانون انتگرال گیری جزء به جزء خواهیم داشت:

$$\begin{aligned} X_{i,j} = \frac{12P_0N}{\pi} \Bigg[ & - m_r\left(\frac{\pi}{6} - \theta + \frac{2(i-1)\pi}{3} - + \frac{2(j-1)\pi}{n}\right) + m_r\left(-\theta + \frac{2(i-1)\pi}{3} - + \frac{2(j-1)\pi}{n}\right) \\ & + m_r\left(\frac{2\pi}{3} - \theta + \frac{2(i-1)\pi}{3} - + \frac{2(j-1)\pi}{n}\right) - m_r\left(\frac{\pi}{2} - \theta + \frac{2(i-1)\pi}{3} - + \frac{2(j-1)\pi}{n}\right) \\ & - m_r\left(\frac{7\pi}{6} - \theta + \frac{2(i-1)\pi}{3} - + \frac{2(j-1)\pi}{n}\right) + m_r\left(\pi - \theta + \frac{2(i-1)\pi}{3} - + \frac{2(j-1)\pi}{n}\right) \\ & - m_r\left(\frac{5\pi}{3} - \theta + \frac{2(i-1)\pi}{3} - + \frac{2(j-1)\pi}{n}\right) + m_r\left(\frac{3\pi}{2} - \theta + \frac{2(i-1)\pi}{3} - + \frac{2(j-1)\pi}{n}\right) \Bigg] \\ & - P_0Nk_r(2\pi - \theta) + \frac{4\pi N}{n} \end{aligned} \quad (66\text{-}4)$$

که در این رابطه :

$$k_r(\varphi) = \int_0^{\varphi} AR(\varphi')d\varphi' \quad (67\text{-}4)$$

$$m_r(\varphi) = \int_0^{\varphi} k_r(\varphi')d\varphi' \quad (68\text{-}4)$$

## 4-8- محاسبه اندوکتانسهای موتور در حالتهای ناهم محوری

با توجه به آنچه در بخشهای قبل گفته شد می توان کلیه اندوکتانسهای موتور را ازطریق روابط زیرمحاسبه نمود :

### 4-8-1- اندوکتانس خودی فاز $i$ استاتور

$$L_i^S = 2\pi l\langle P.BS_i\rangle - 2\pi l\frac{\langle P.AS_i\rangle^2}{\langle P\rangle} \tag{69-4}$$

همانطور که در شکلهای4-7 تا 4-10 دیده می شود اندوکتانس خودی فازهای استاتور در حالتهای ناهم محوری پویا و مرکب تابعی از موقعیت زاویه ای روتور است ولی در حالت ناهم محوری ایستا مانند حالت سالم این اندوکتانسها ثابت هستند. دلیل این امر این است که در حالت ناهم محوری ایستا روتور در جای خود می چرخد لذا مدار معادل مغناطیسی دیده شده برای فازهای استاتور همواره ثابت است. اما در حالتهای ناهم محوری پویا و مرکب که روتور در داخل استاتور علاوه بر حرکت چرخشی یک حرکت انتقالی نیز دارد این مدار معادل بسته به موقعیت روتور تغییر می کند. نکته دیگری که در این شکلها مشهود است افزایش وابستگی این اندوکتانس به موقعیت زاویه ای روتور در حالت ناهم محوری مرکب نسبت به سایر حالتهاست. این مساله به دلیل وابستگی ضرایب $A,B$ و $C$ به موقعیت زاویه ای روتور در این حالت است. همانطور که در فصل دوم نشان داده شد این ضرایب تنها در حالت ناهم محوری مرکب تابعی از موقعیت زاویه‌ای روتور هستند و در حالتهای ناهم محوری ایستا و پویا مستقل از این متغیرند.

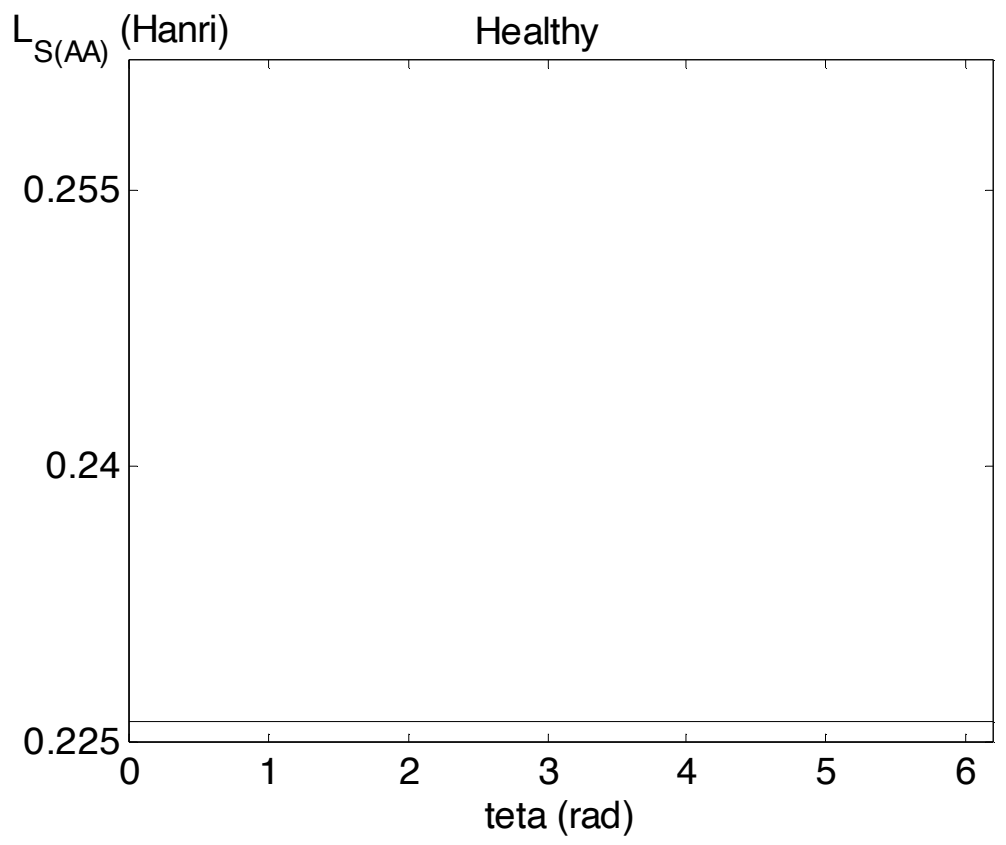

شکل**4-7**: اندوکتانس خودی فاز اول استاتور در حالت سالم

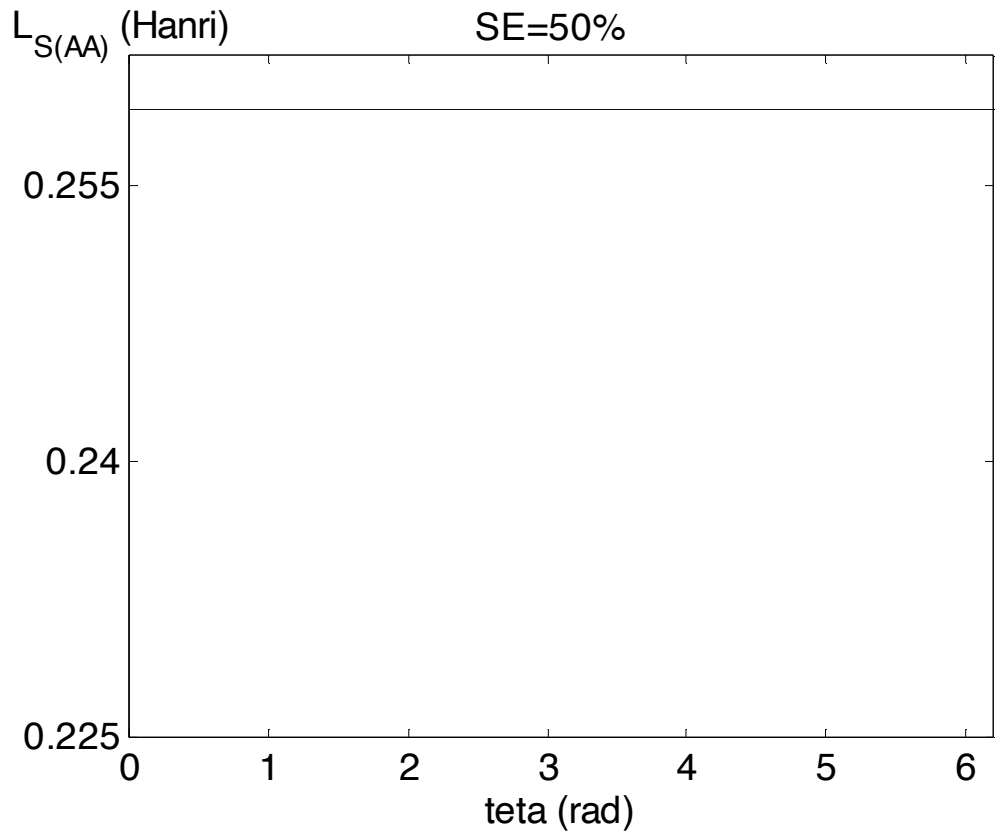

شکل**4-8**: اندوکتانس خودی فاز اول استاتور در حالت: **50**٪ ناهم محوری ایستا

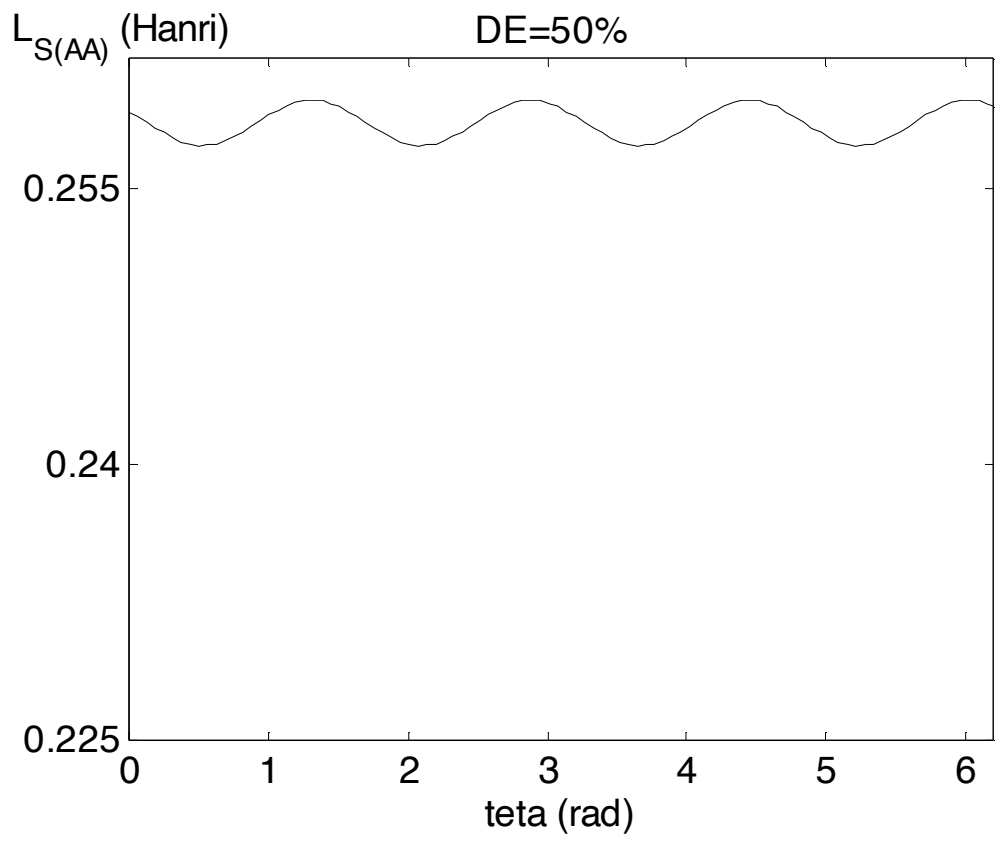

شکل**4-9**: اندوکتانس خودی فاز اول استاتور در حالت **50**٪ ناهم محوری پویا

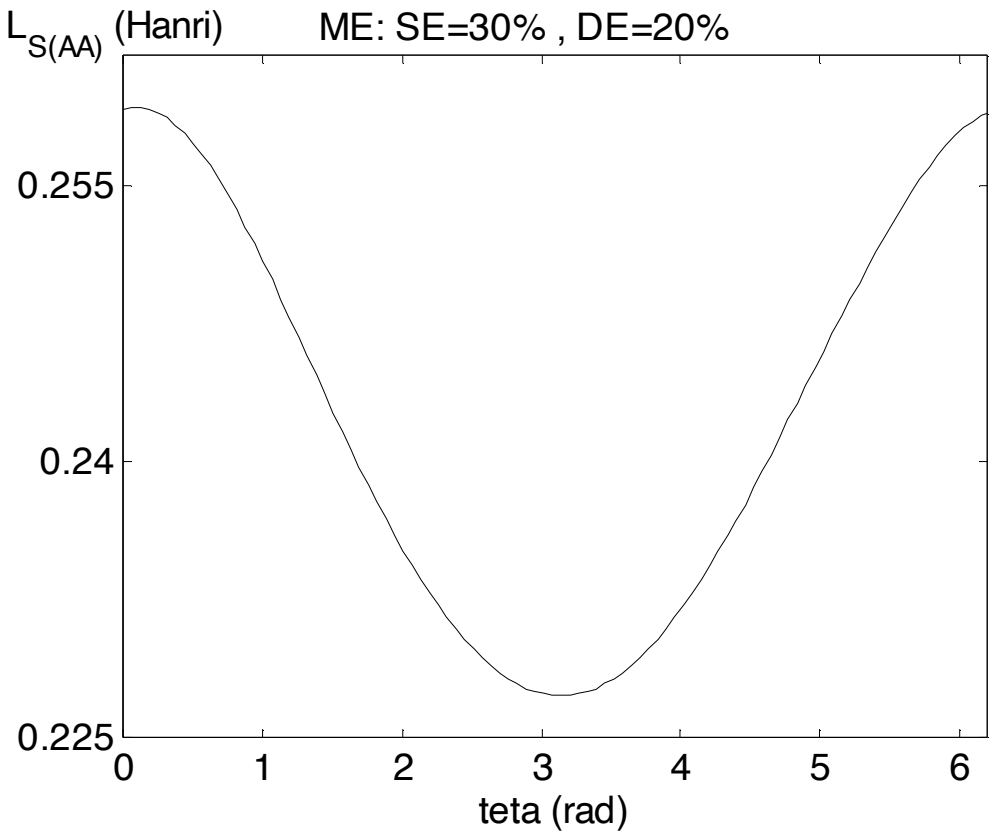

شکل**4-10**: اندوکتانس خودی فاز اول استاتور در حالت **30**٪ ناهم محوری پویا و **20**٪ ناهم محوری ایستا

## **4-8-2**- اندوکتانس متقابل فاز $i$ و $j$ استاتور

$$L_{i,j}^{S} = 2\pi l\left\langle P.CS_{i,j}\right\rangle - 2\pi l\frac{\left\langle P.AS_i\right\rangle\left\langle P.AS_j\right\rangle}{\left\langle P\right\rangle} \qquad (70\text{-}4)$$

شکلهای 4-11 تا 4-14 اندوکتانس متقابل فاز اول و دوم استاتور را برای حالت سالم، ناهم محوری ایستا، ناهم محوری پویا و ناهم محوری مرکب نشان می دهد.

همانطور که شکلهای 4-11 تا 4-14 نشان می دهند رفتار این اندوکتانس و وابستگی آن به موقعیت زاویه ای روتور نیز در حالتهای مختلف ناهم محوری مانند اندوکتانس خودی فازهای استاتور است.

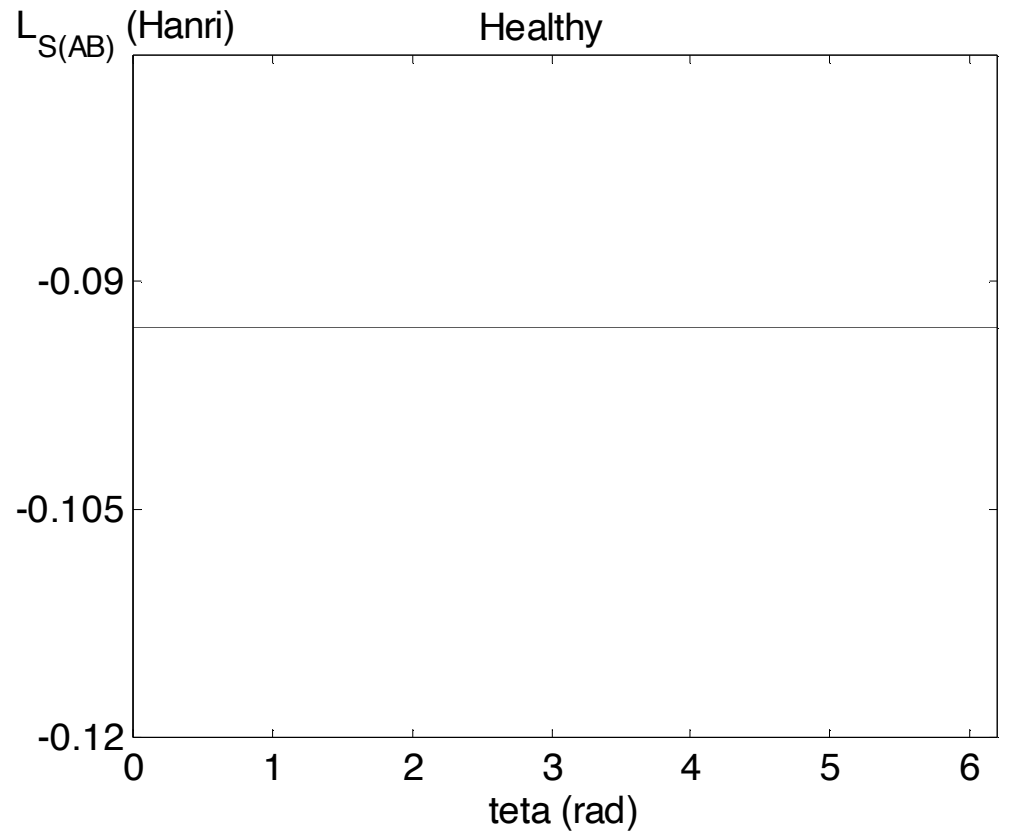

شکل **4-11**: اندوکتانس متقابل فاز اول و دوم استاتور در حالت سالم

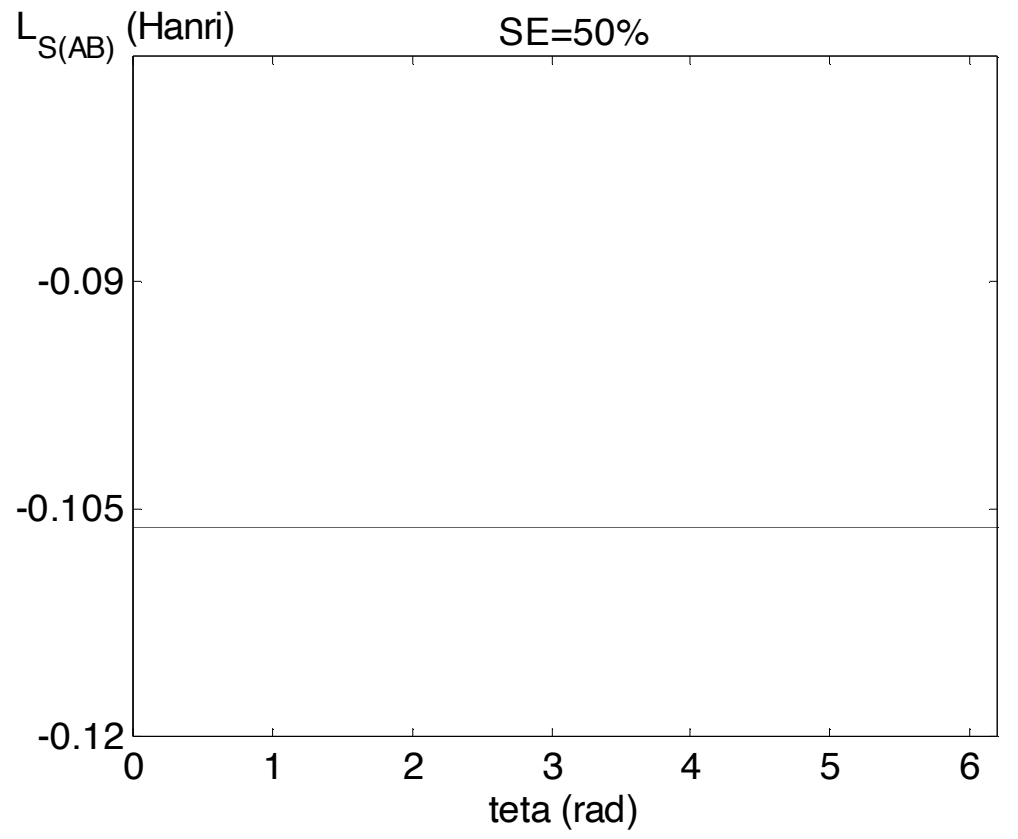

شکل **4-12**: اندوکتانس متقابل فاز اول و دوم استاتور در حالت: **50**٪ ناهم محوری ایستا

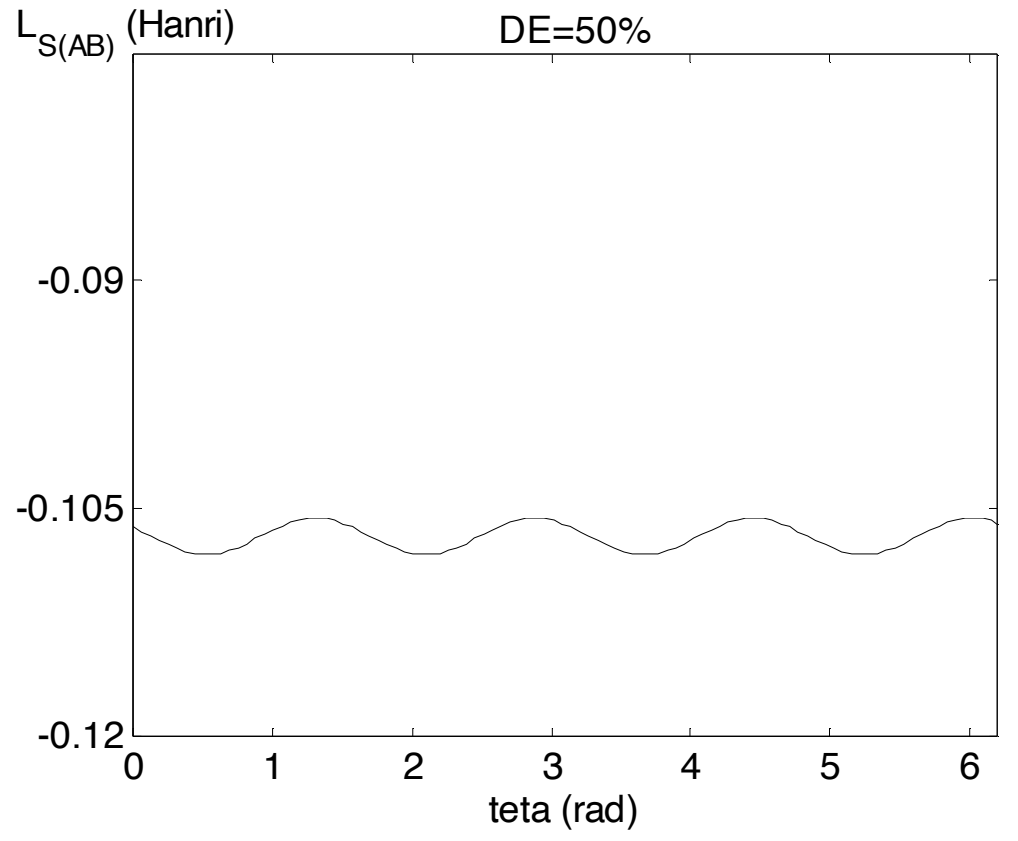

شکل**4-13**: اندوکتانس متقابل فاز اول و دوم استاتور در حالت **50**٪ ناهم محوری پویا

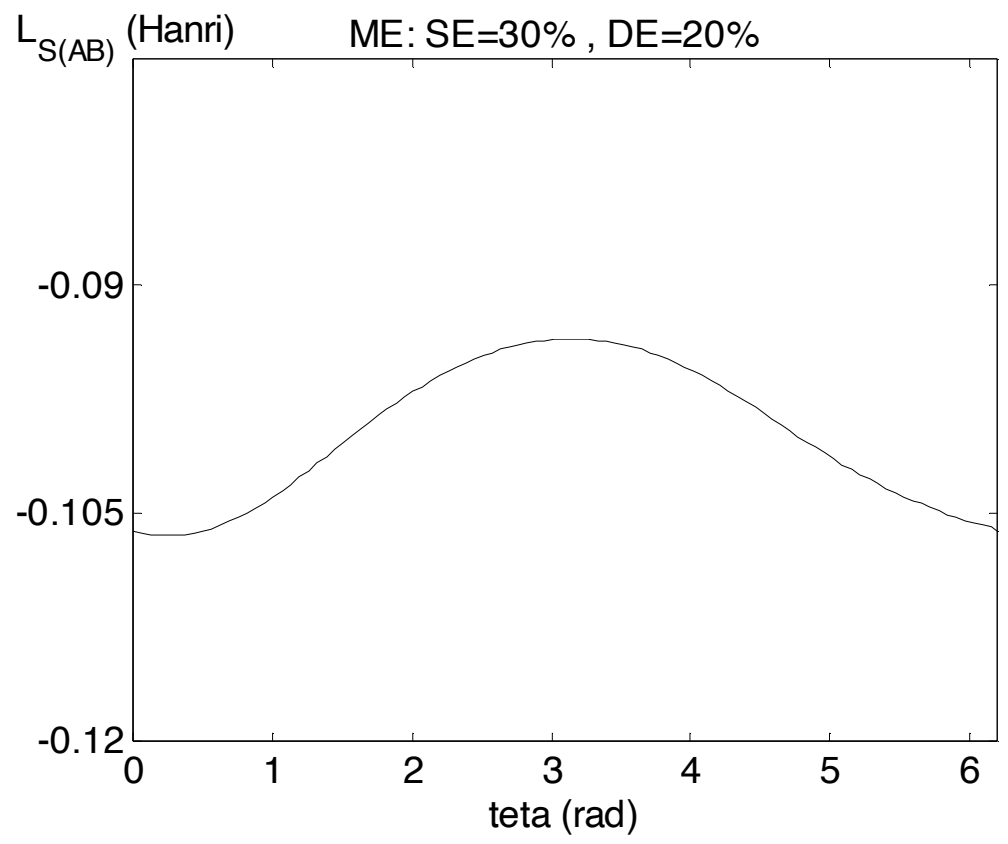

شکل **4-14**: اندوکتانس متقابل فاز اول و دوم استاتور در حالت **30**٪ ناهم محوری پویا و **20**٪ ناهم محوری ایستا

### **4-8-3**- اندوکتانس خودی حلقه $i$ روتور

$$L_i^R = 2\pi l\langle P.BR_i\rangle - 2\pi l\frac{\langle P.AR_i\rangle^2}{\langle P\rangle} \qquad (71\text{-}4)$$

شکلهای 4-15 تا 4-18 اندوکتانس خودی حلقه اول روتور را برای حالتهای مختلف نشان می دهد. همانطور که در این شکلها ملاحظه می شود. اندوکتانسهای روتور در حالتهای ناهم محوری ایستا و مرکب تابعی از موقعیت زاویه‌ای روتور بوده و در حالت ناهم محوری پویا مانند حالت سالم ثابت هستند. دلیل این امر این است که در حالت ناهم محوری پویا همراه با تغییرات مدارمعادل مغناطیسی با چرخش روتور موقعیت روتور نیز عوض می شود بنابراین مدار معادل مغناطیسی دیده شده توسط یک حلقه روتور همواره ثابت است.

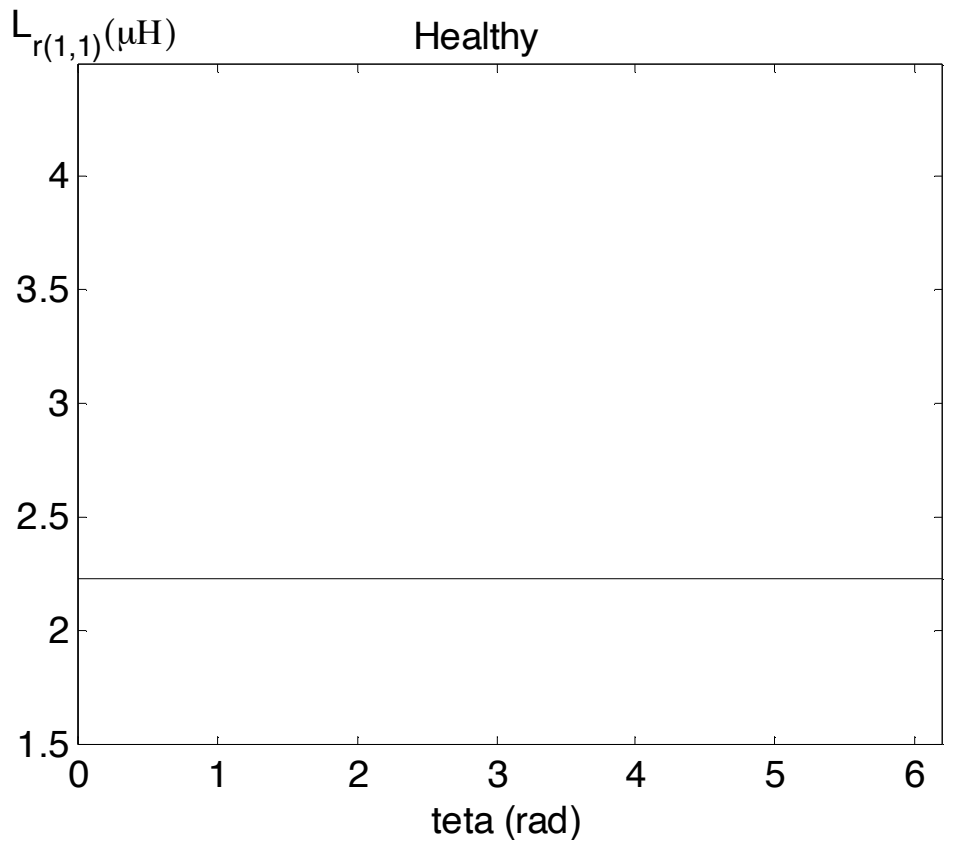

شکل 4-15: اندوکتانس خودی حلقه اول روتور در حالت سالم

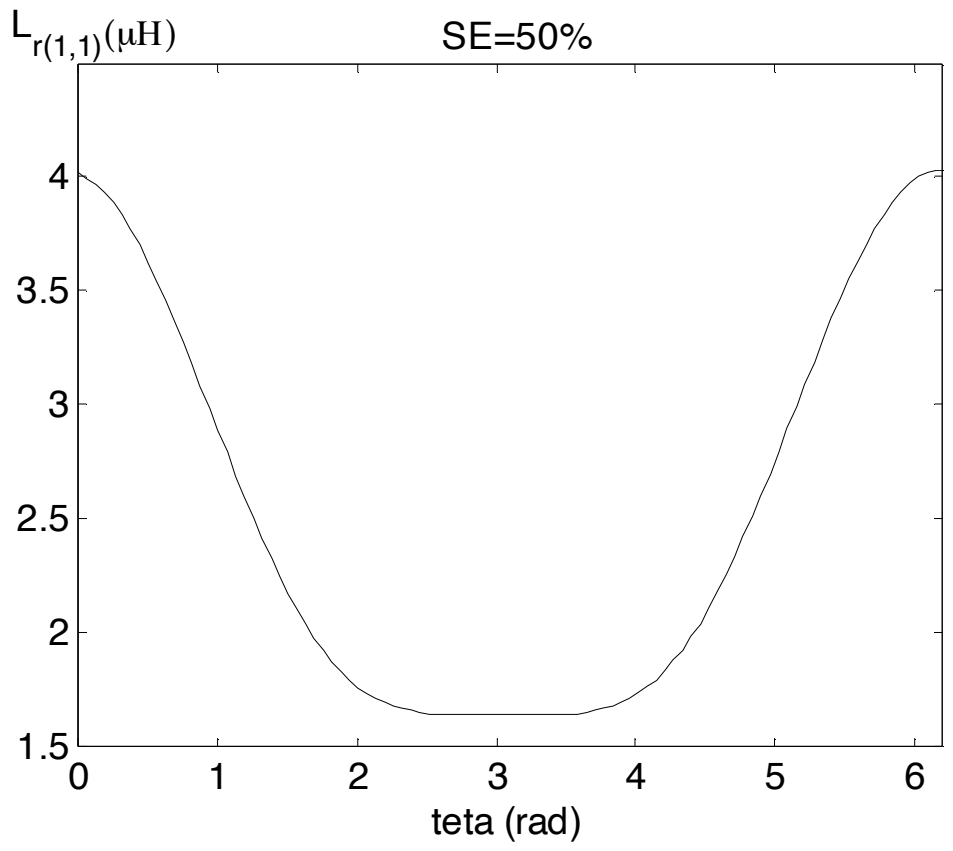

شکل 4-16: اندوکتانس خودی حلقه اول روتور در حالت: 50٪ ناهم محوری ایستا

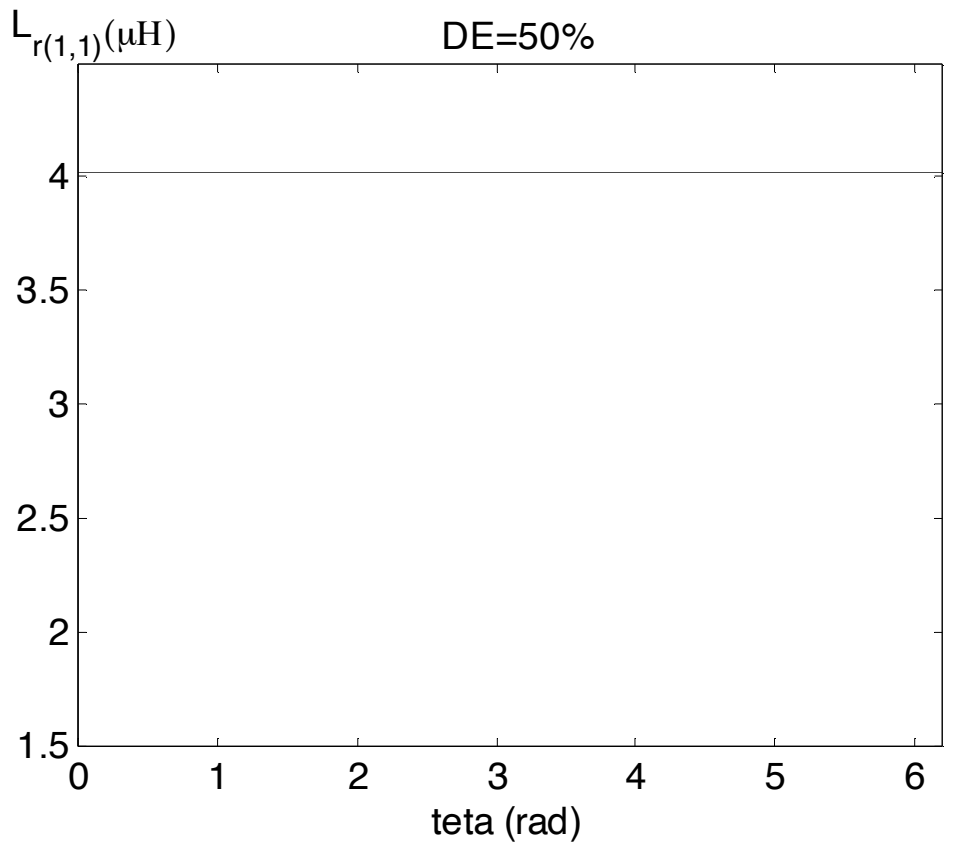

شکل4-17: اندوکتانس خودی حلقه اول روتور در حالت 50٪ ناهم محوری پویا

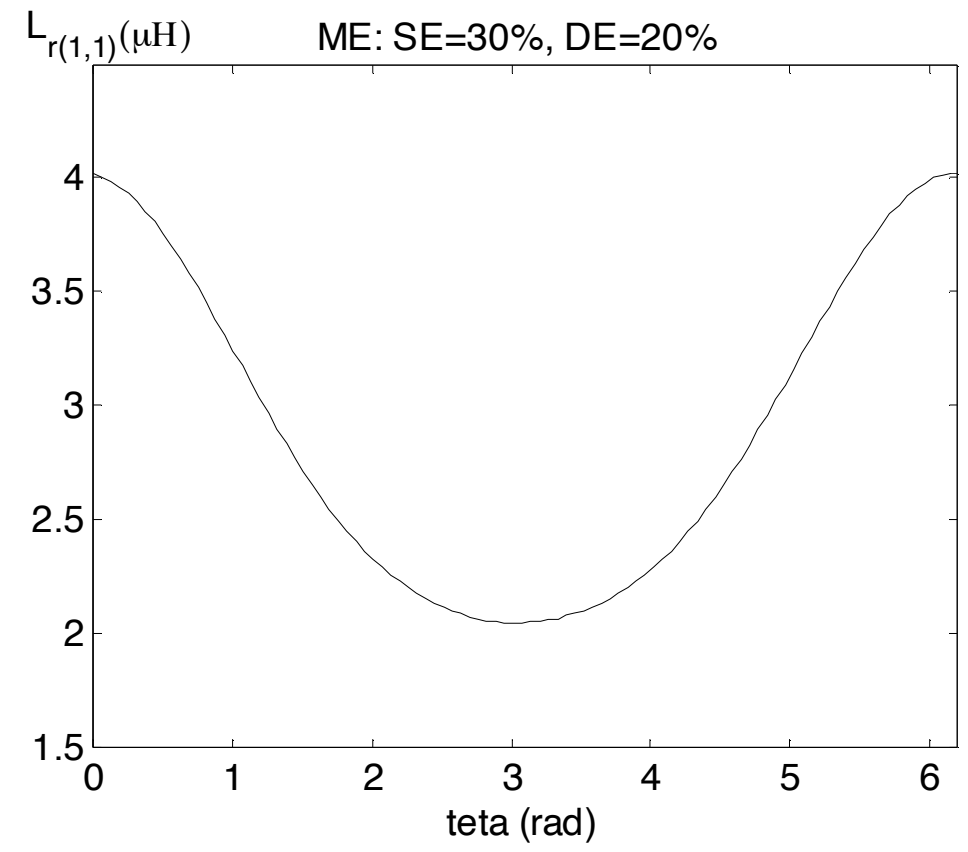

شکل **4-18**: اندوکتانس خودی حلقه اول روتور در حالت **30**٪ ناهم محوری پویا و **20**٪ ناهم محوری ایستا

**4-8-4**- : اندوکتانس متقابل حلقه $i$ و $i-1$ (مجاور) روتور

$$L^{R}_{i,i-1} = 2\pi l \langle P.CR_{i,j} \rangle - 2\pi l \frac{\langle P.AR_i \rangle \langle P.AR_j \rangle}{\langle P \rangle} \qquad (71\text{-}4)$$

شکلهای 4-19 تا 4-22 اندوکتانس متقابل حلقه اول و دوم روتور را برای حالتهای مختلف نشان می دهد. وابستگی و عدم وابستگی این اندوکتانسها به موقعیت زاویه ای روتور مانند اندوکتانسهای خودی روتور است.

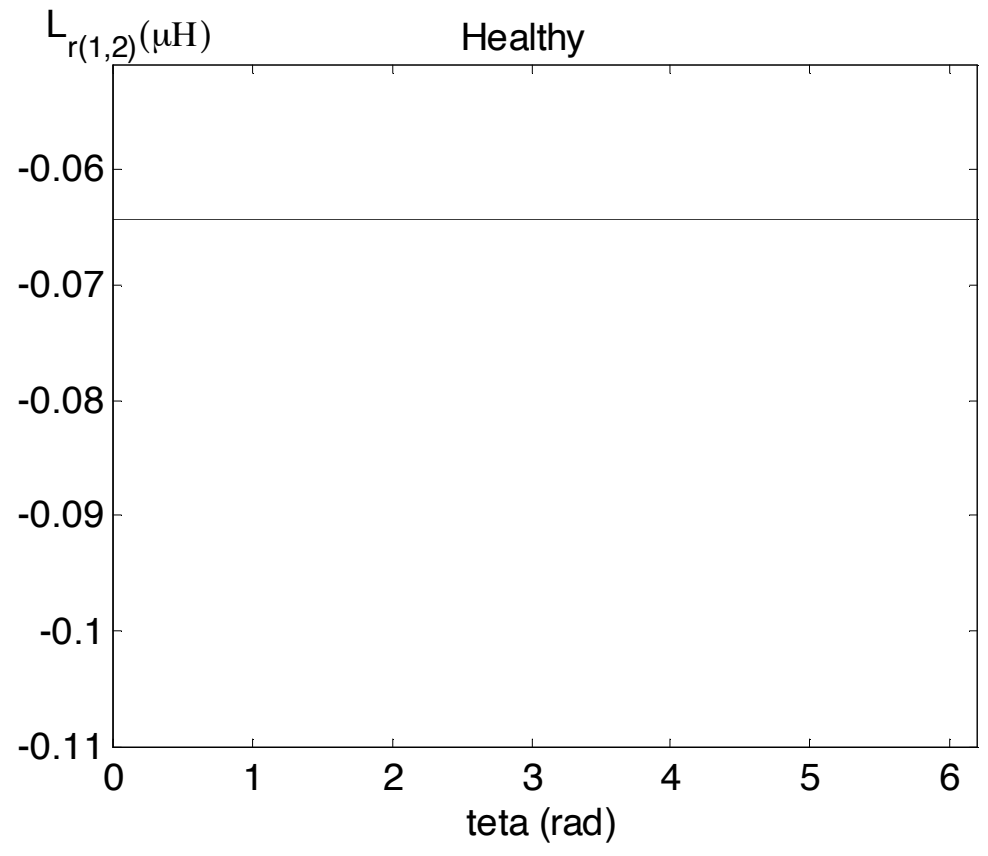

شکل **4-19**: اندوکتانس متقابل حلقه اول و دوم روتور در حالت سالم

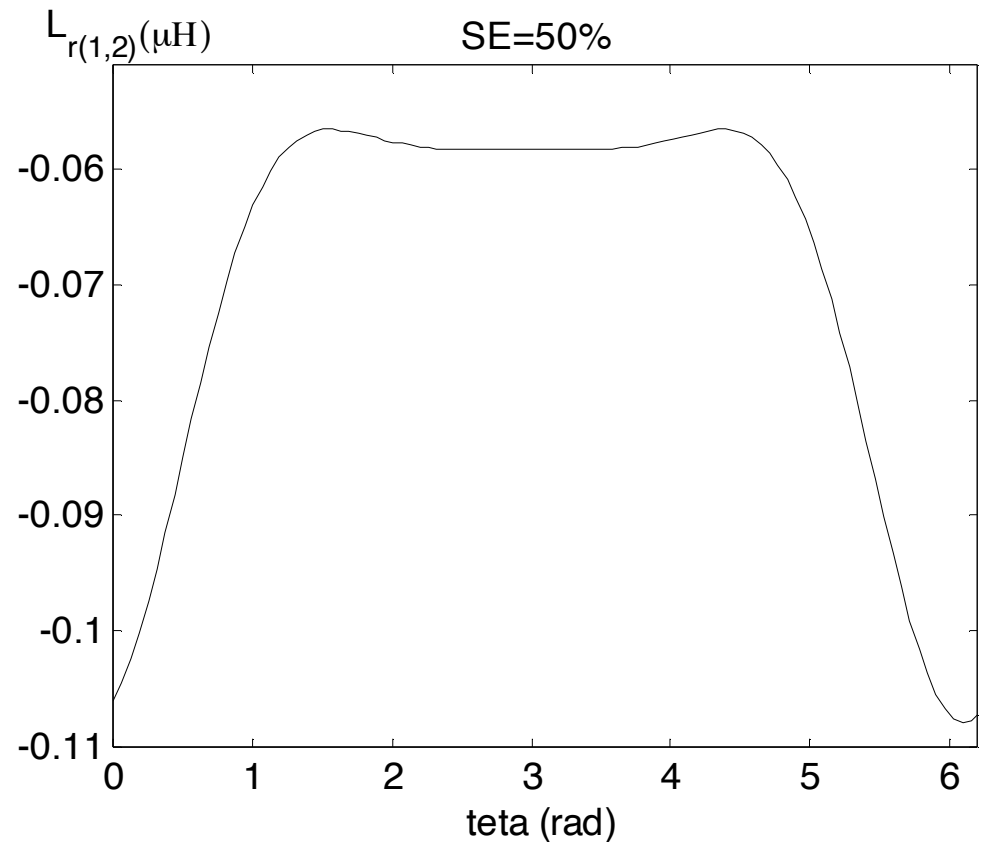

شکل 4-20: اندوکتانس متقابل حلقه اول و دوم روتور در حالت: 50٪ ناهم محوری ایستا

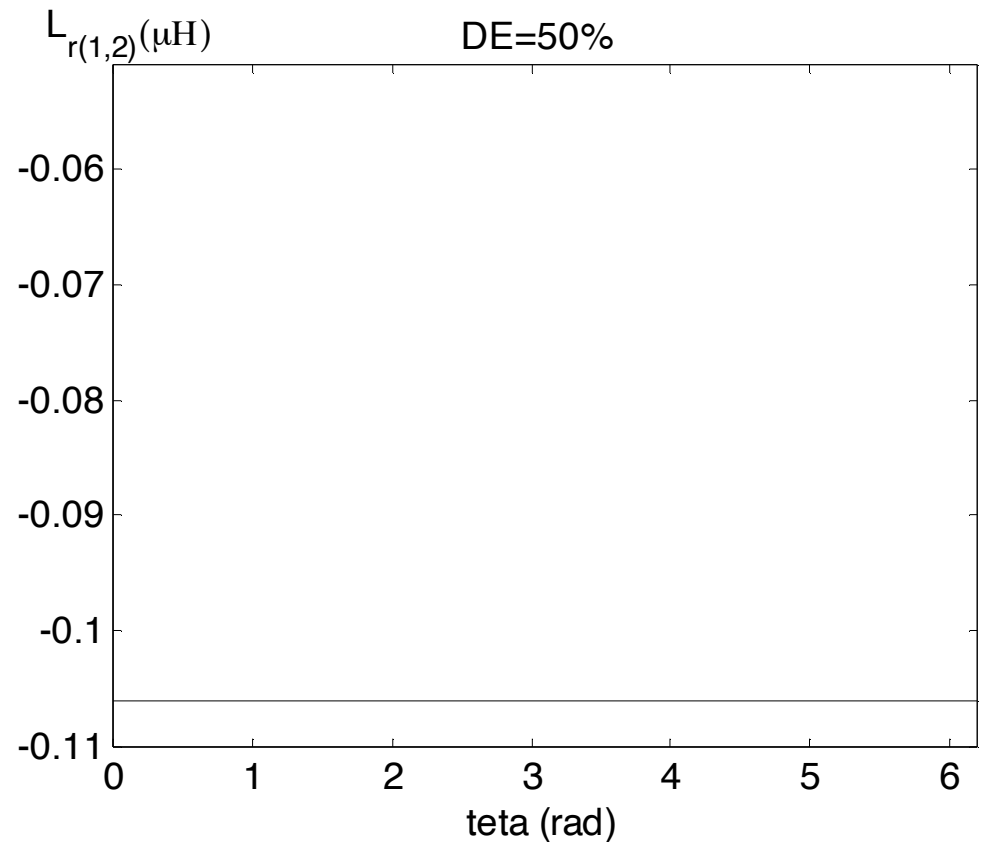

شکل4-21: اندوکتانس متقابل حلقه اول و دوم روتور در حالت 50٪ ناهم محوری پویا

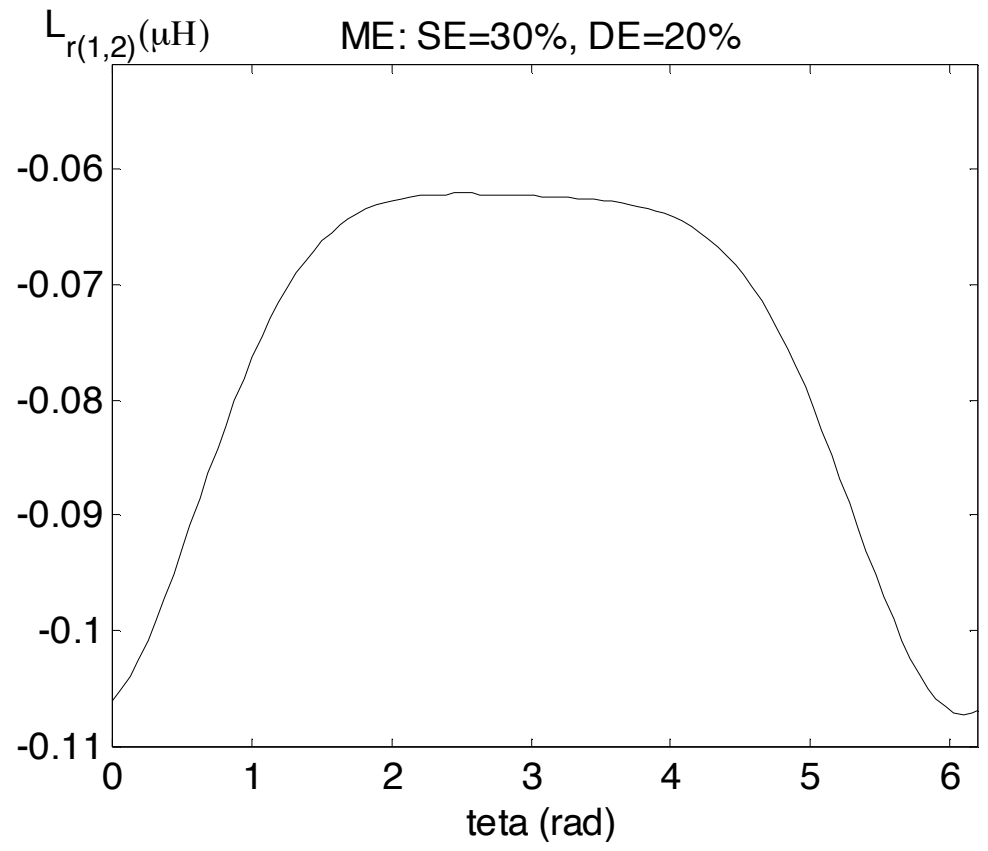

شکل 4-22: اندوکتانس متقابل حلقه اول و دوم روتور در حالت 30٪ ناهم محوری پویا و 20٪ ناهم محوری ایستا

**4-8-5**- اندوکتانس متقابل  دو حلقه غیرمجاور روتور

$$L^{R}_{i,j} = -2\pi d \frac{\langle P.AR_i \rangle \langle P.AR_j \rangle}{\langle P \rangle} \qquad (72\text{-}4)$$

شکلهای 4-23 تا 4-26 اندوکتانس متقابل حلقه اول و سوم روتور را برای حالتهای مختلف نشان می دهد. وابستگی و عدم وابستگی این اندوکتانسها به موقعیت راویه ای روتور مانند اندوکتانسهای خودی روتور است.

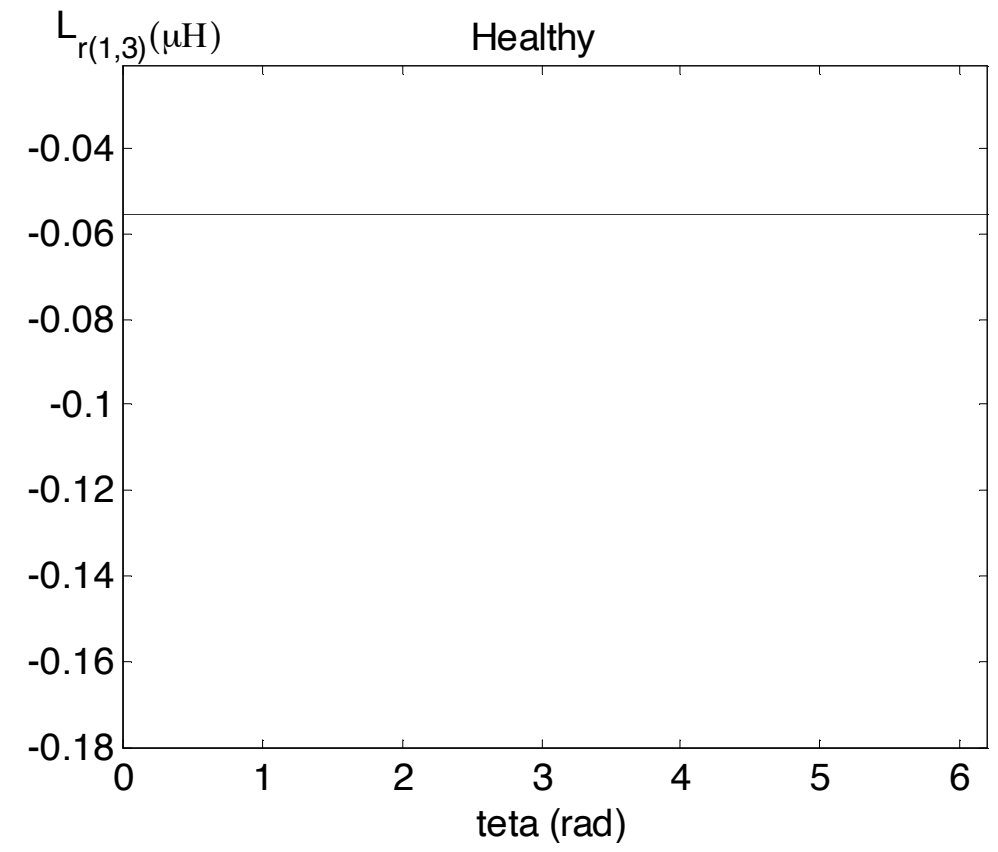

شکل **4-23**: اندوکتانس متقابل حلقه اول و سوم روتور در حالت سالم

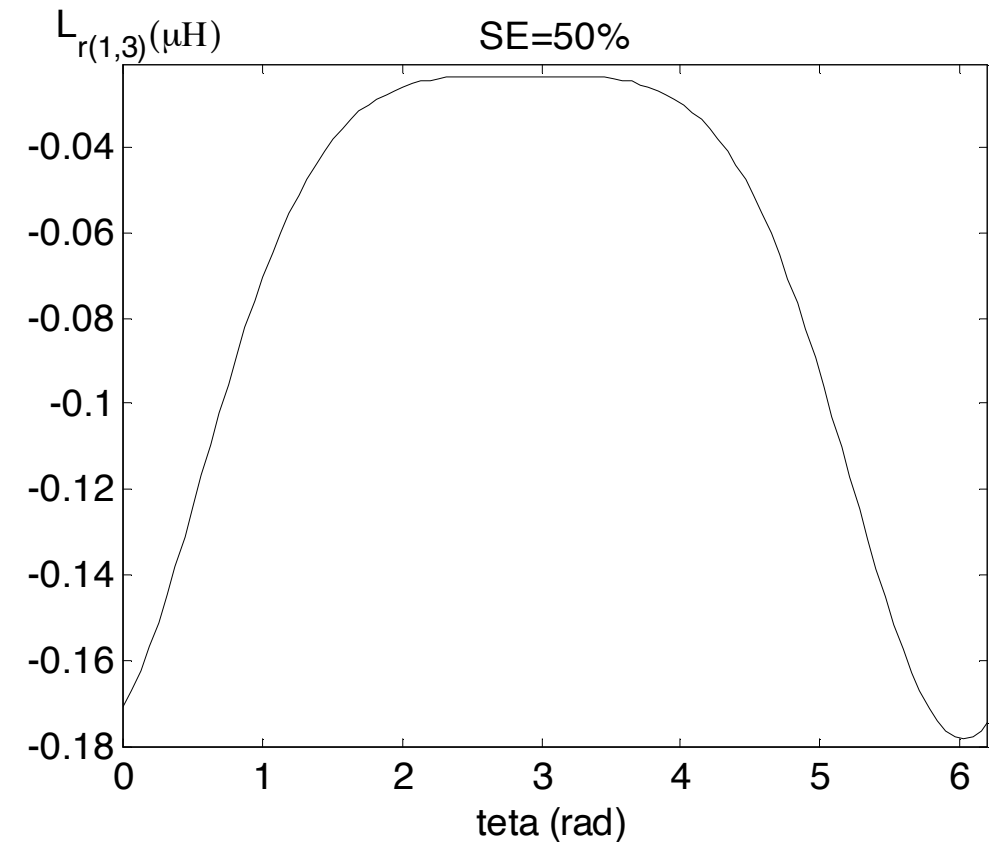

شکل **4-24**: اندوکتانس متقابل حلقه اول و سوم روتور در حالت: **50**٪ ناهم محوری ایستا

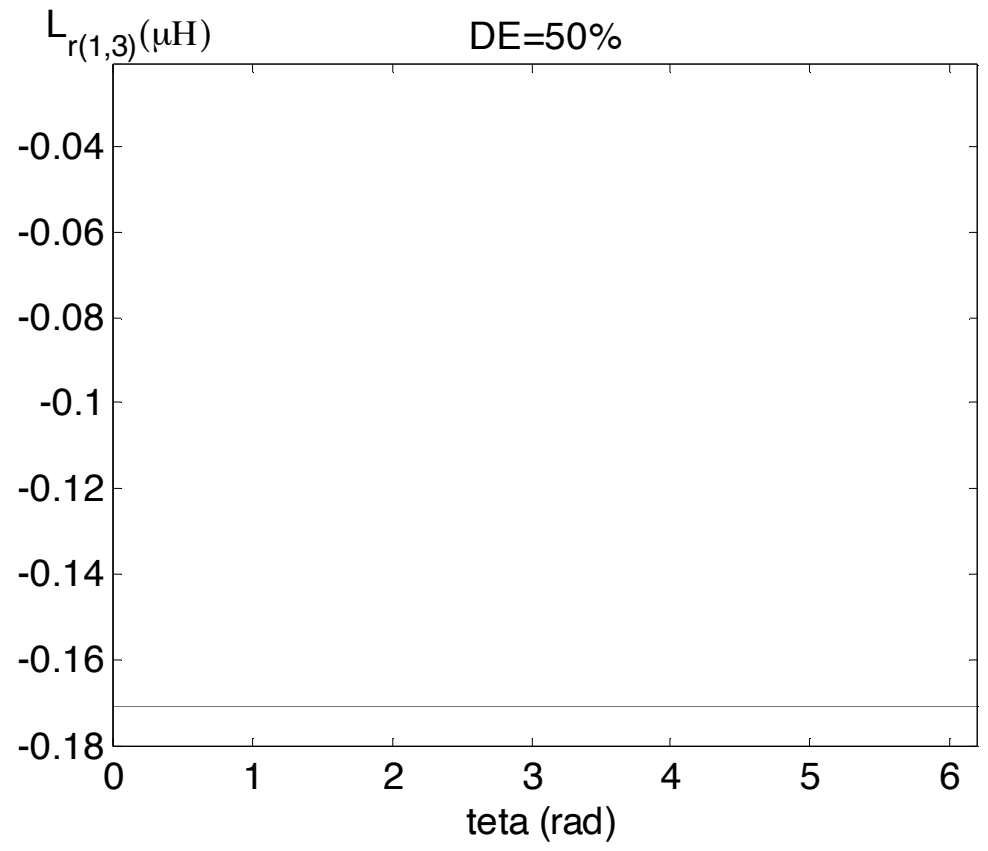

شکل **4-25**: اندوکتانس متقابل حلقه اول و سوم روتور در حالت **50**٪ ناهم محوری پویا

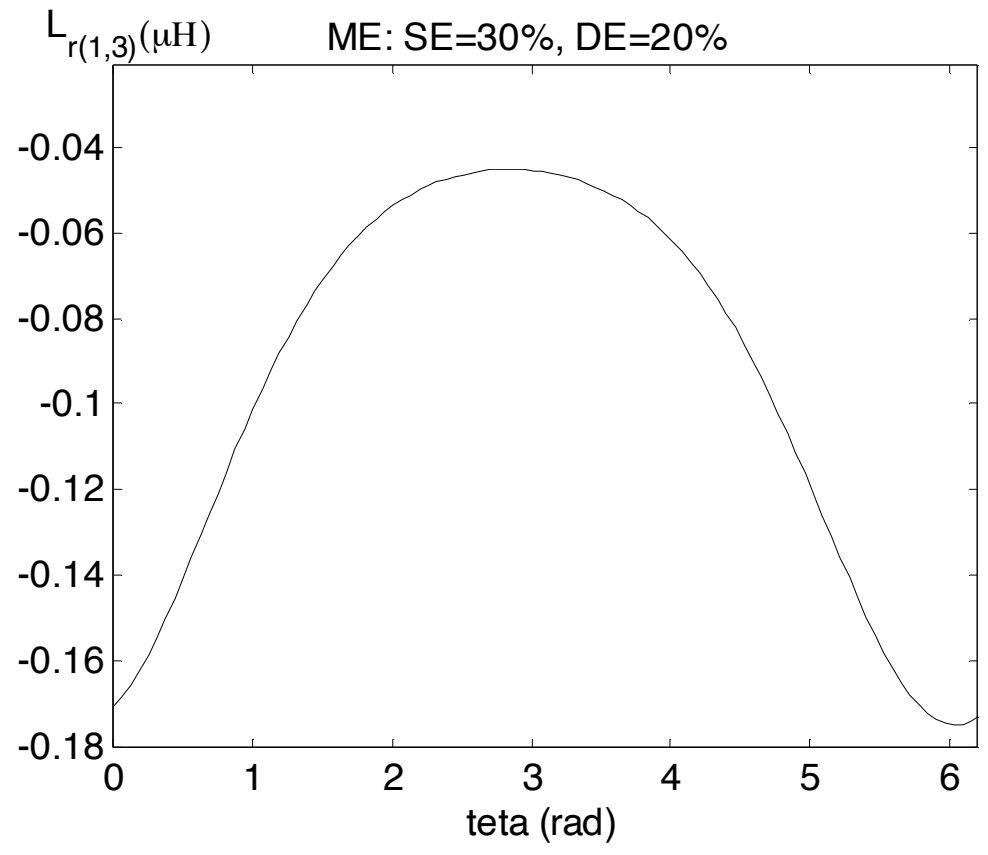

شکل **4-26**: اندوکتانس متقابل حلقه اول و دوم روتور در حالت **30**٪ ناهم محوری پویا و **20**٪ ناهم محوری ایستا

### **4-8-6**- اندوکتانس متقابل فاز $i$ استاتور و حلقه $j$ روتور

$$L_{i,j}^{SR} = 2\pi l\left\langle P.AS_i.AR_j\right\rangle - 2\pi l\frac{\left\langle P.AS_i\right\rangle\left\langle P.AR_j\right\rangle}{\left\langle P\right\rangle} \quad (73\text{-}4)$$

شکلهای4-27 تا 4-30 اندوکتانس متقابل فاز اول استاتور و حلقه اول روتور را برای حالتهای مختلف نشان می دهد.

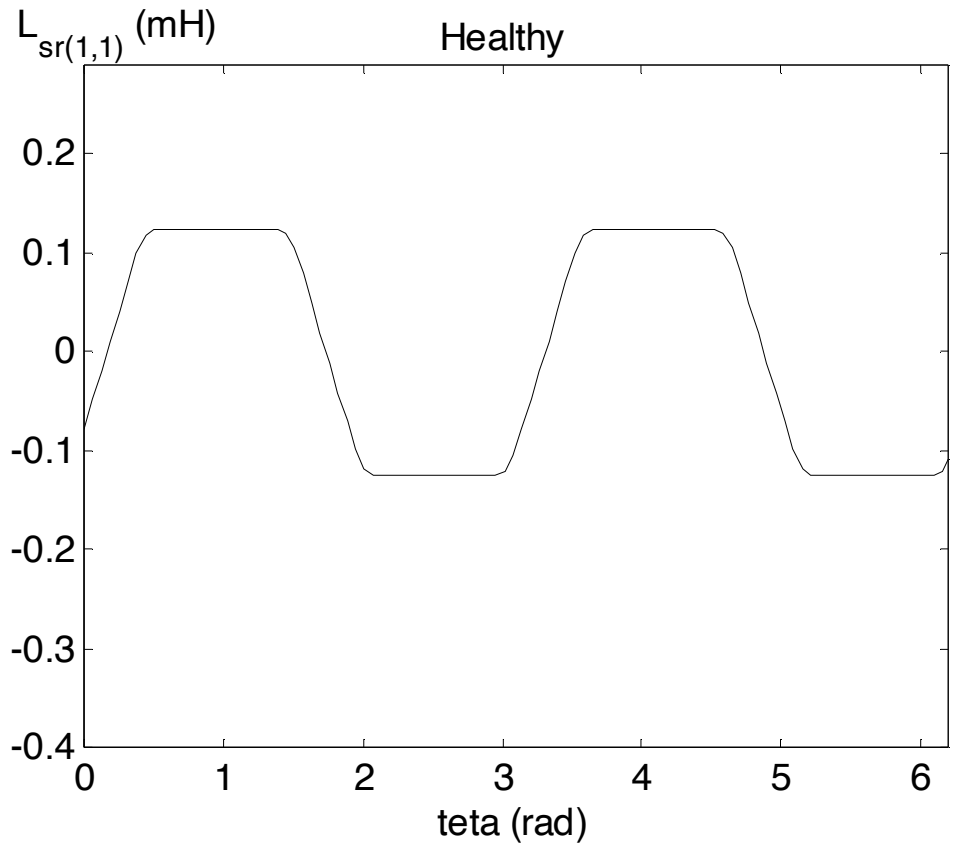

شکل **4-27**: اندوکتانس متقابل فاز اول استاتور و حلقه اول روتور در حالت سالم

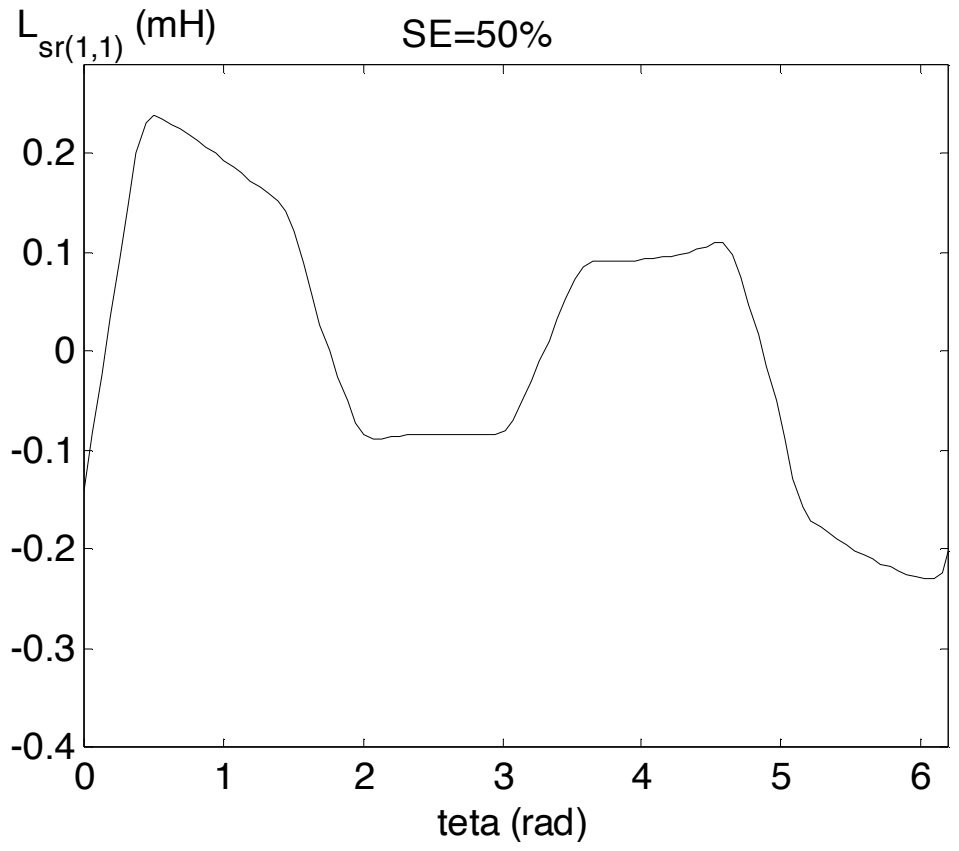

شکل **4-28**: اندوکتانس متقابل فاز اول استاتور و حلقه اول روتور  در حالت: **50**٪ ناهم محوری ایستا

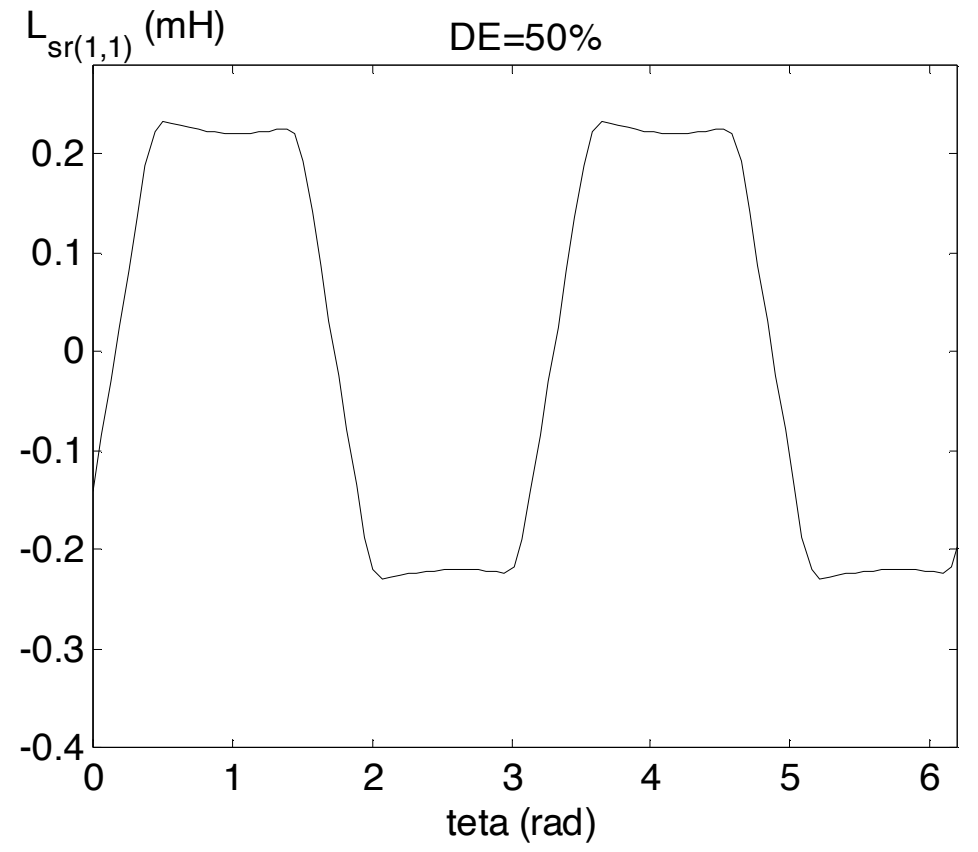

شکل 4-29: اندوکتانس متقابل فاز اول استاتور و حلقه اول روتور در حالت: 50٪ ناهم محوری پویا

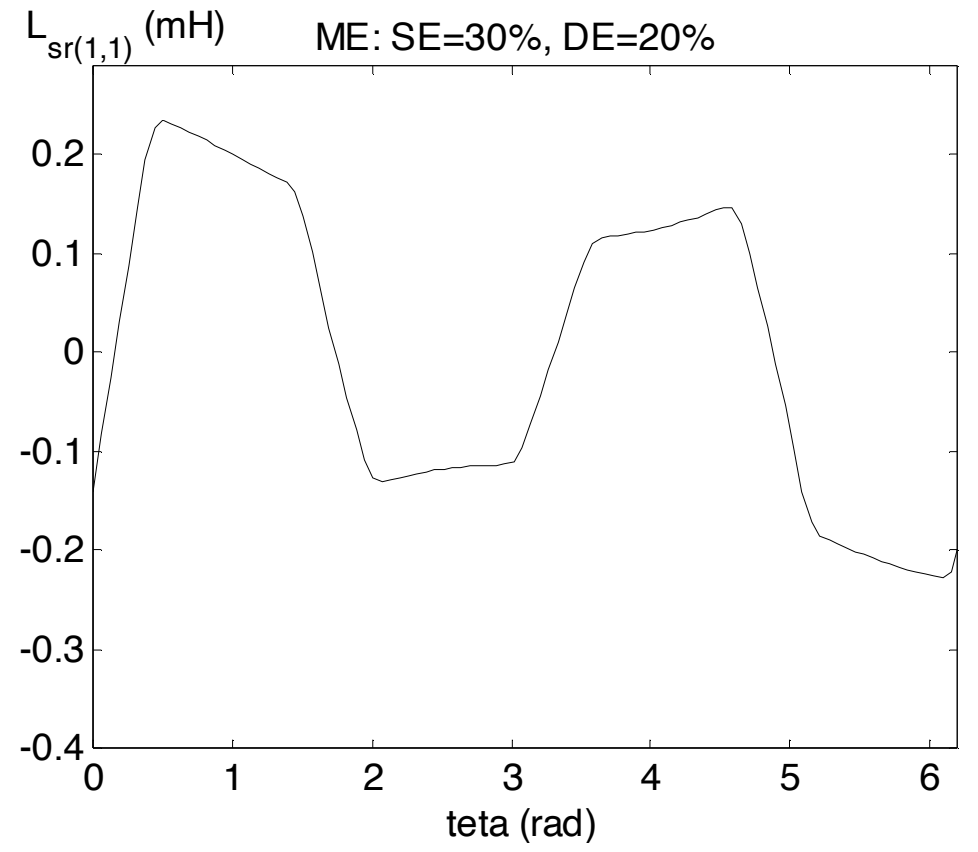

شکل 4-30: اندوکتانس متقابل فاز اول استاتور و حلقه اول روتور در حالت 30٪ ناهم محوری پویا و 20٪ ناهم محوری ایستا

## 4-9- بررسی اثر تورب میله های روتور بر روی اندوکتانسهای موتور

در صورتیکه که $\gamma$ زاویه تورب میله های روتور، $\ell$ طول موتور و $L'(\theta)$ یکی از اندوکتانسهای مغناطیسی موتور در واحد طول باشد. با در نظر گرفتن اثر تورب میله های روتور اندوکتانس مغناطیسی در کل طول موتور بصورت زیر محاسبه می شود:

(74-4)

$$L_{skew}(\theta)=\int_{-\frac{\ell}{2}}^{\frac{\ell}{2}} L'\left(\theta-\xi\frac{\gamma}{\ell}\right)d\xi$$

با استفاده از تغییر متغیر در انتگرالگیری، انتگرال فوق بصورت زیر در خواهد آمد:

(75-4)

$$L_{skew}(\theta)=\int_{-\frac{\gamma}{2}}^{\frac{\gamma}{2}} \frac{\ell}{\gamma}L'(\theta-x)dx$$

$\ell L'(\theta)$ اندوکتانس مغناطیسی مذکور در کل طول موتور بدون در نظر گرفتن اثر تورب میله های موتور است که با $L(\theta)$ نشان داده می شود. در اینصورت داریم:

(76-4)

$$L_{skew}(\theta)=\frac{1}{\gamma}\int_{-\frac{\gamma}{2}}^{\frac{\gamma}{2}} L(\theta-x)dx$$

در صورتیکه که $\gamma$ نسبتا کوچک باشد با تقریب بسیار خوبی از طریق بکارگیری قانون انتگرالگیری ذوزنقه ای حاصل انتگرال بالا بصورت زیر خواهد بود:

(77-4)

$$L_{skew}(\theta)=\frac{L\left(\theta-\frac{\gamma}{2}\right)+2L(\theta)+L\left(\theta+\frac{\gamma}{2}\right)}{4}$$

بنابراین پس از استفاده از روشی که در بخش قبل برای محاسبه اندوکتانسهای مغناطیسی موتور بدون در نظرگرفتن تورب میله ها ارائه شد می توان رابطه (4-77) را بکار برد تا این اثر لحاظ گردد. اما از آنجا که در بررسی حاضر

اثر خطی افزایش تابع دور سیم پیچهای در نظر گرفته شده است ملاحظه می شود که تورب میله های روتور اثر قابل توجهی بر روی اندکتانسهای مغناطیسی نمی گذارند. شکل 4-31 اندوکتانس متقابل فاز اول استاتور و حلقه اول روتور را با در نظر گرفتن و بدون در نظر گرفتن تورب میله های روتور در حالت ناهم‌محوری مرکب نشان می دهد. همانطور که ملاحظه می کنید این دو اندوکتانس تقریبا بر هم منطبقند. شکل 4-32 با بزرگنمایی بیشتری قسمتهای شکستگی این اندوکتانس را نشان می دهد. در این نواحی اثر تورب میله های روتور اندکی به چشم می خورد، اما قابل ملاحظه نیست.

شکلهای 4-33 و 4-34 دو شکل فوق را برای حالتی که موتور سالم است نشان می دهد. در این شکلها نیز اثر تورب میله های روتور ناچیز است.

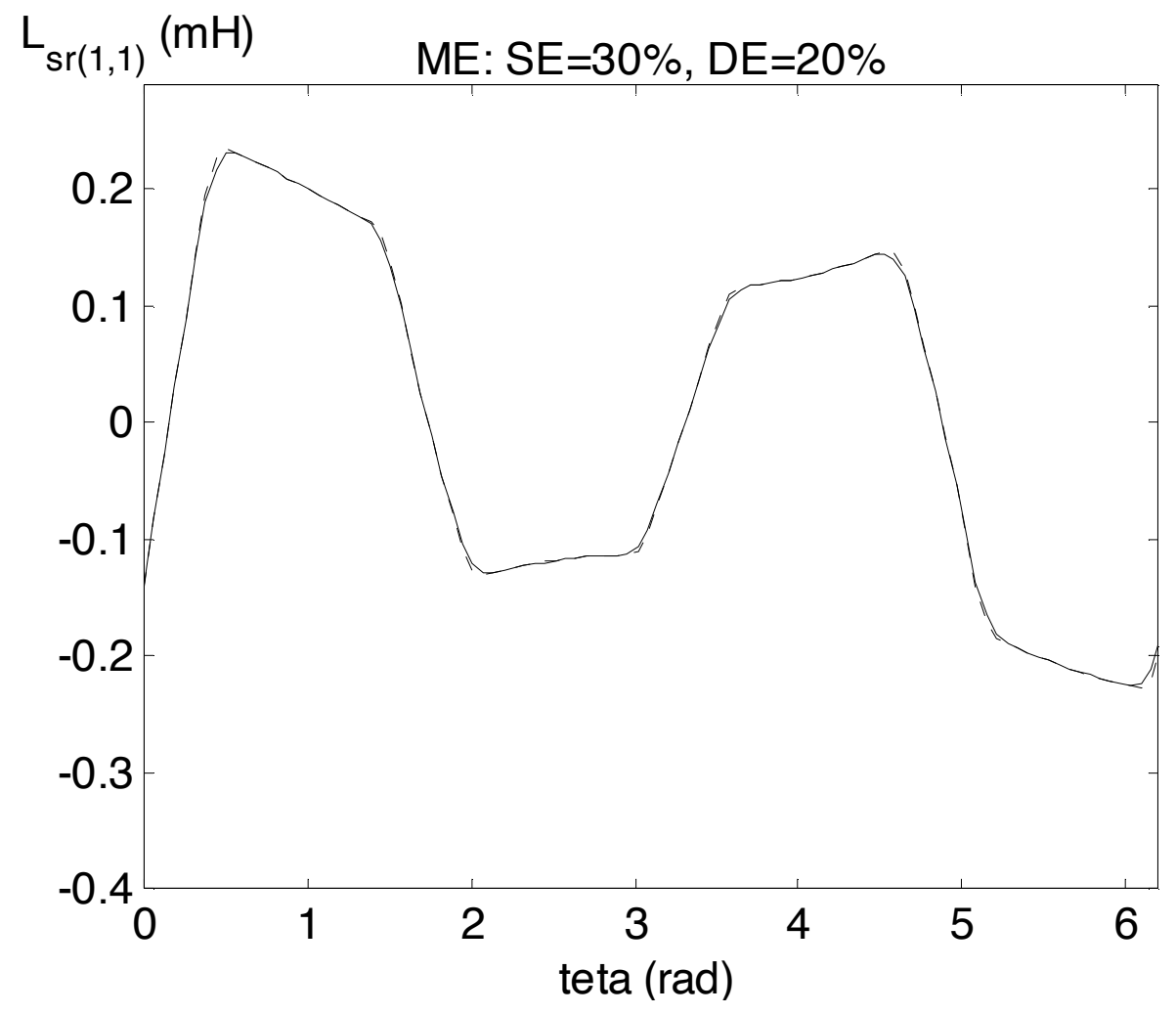

شکل **4-31**: اندوکتانس متقابل فاز اول استاتور و حلقه اول روتور در حالت ناهم محوری مرکب

خط چین: بدون در نظر گرفتن تورب میله های روتور. توپر : با در نظر گرفتن تورب میله های روتور

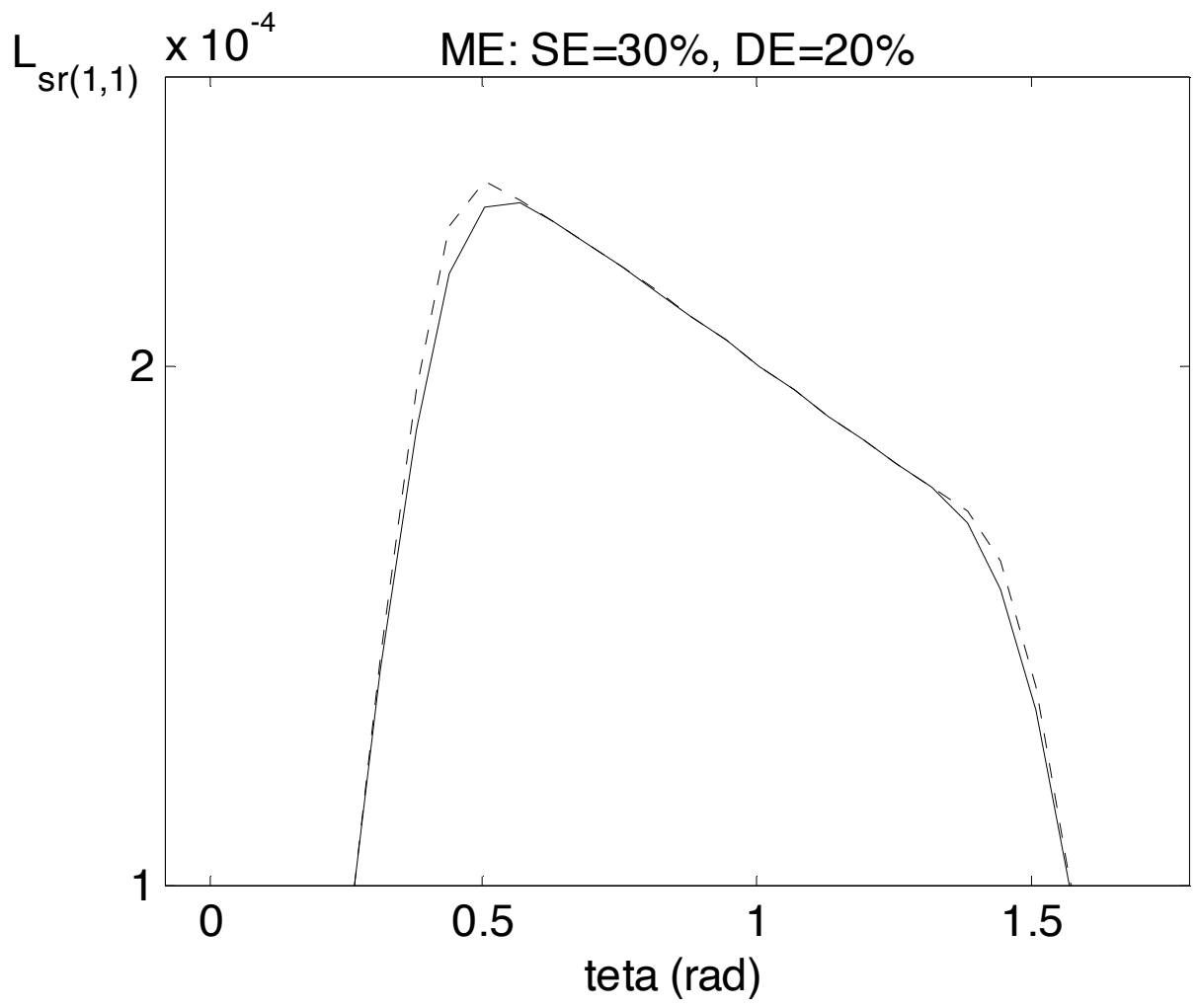

شکل **4-32**: اندوکتانس متقابل فاز اول استاتور و حلقه اول روتور در حالت ناهم محوری مرکب با بزرگنمایی بیشتر

خط چین: بدون در نظر گرفتن تورب میله های روتور. توپر : با در نظر گرفتن تورب میله های روتور

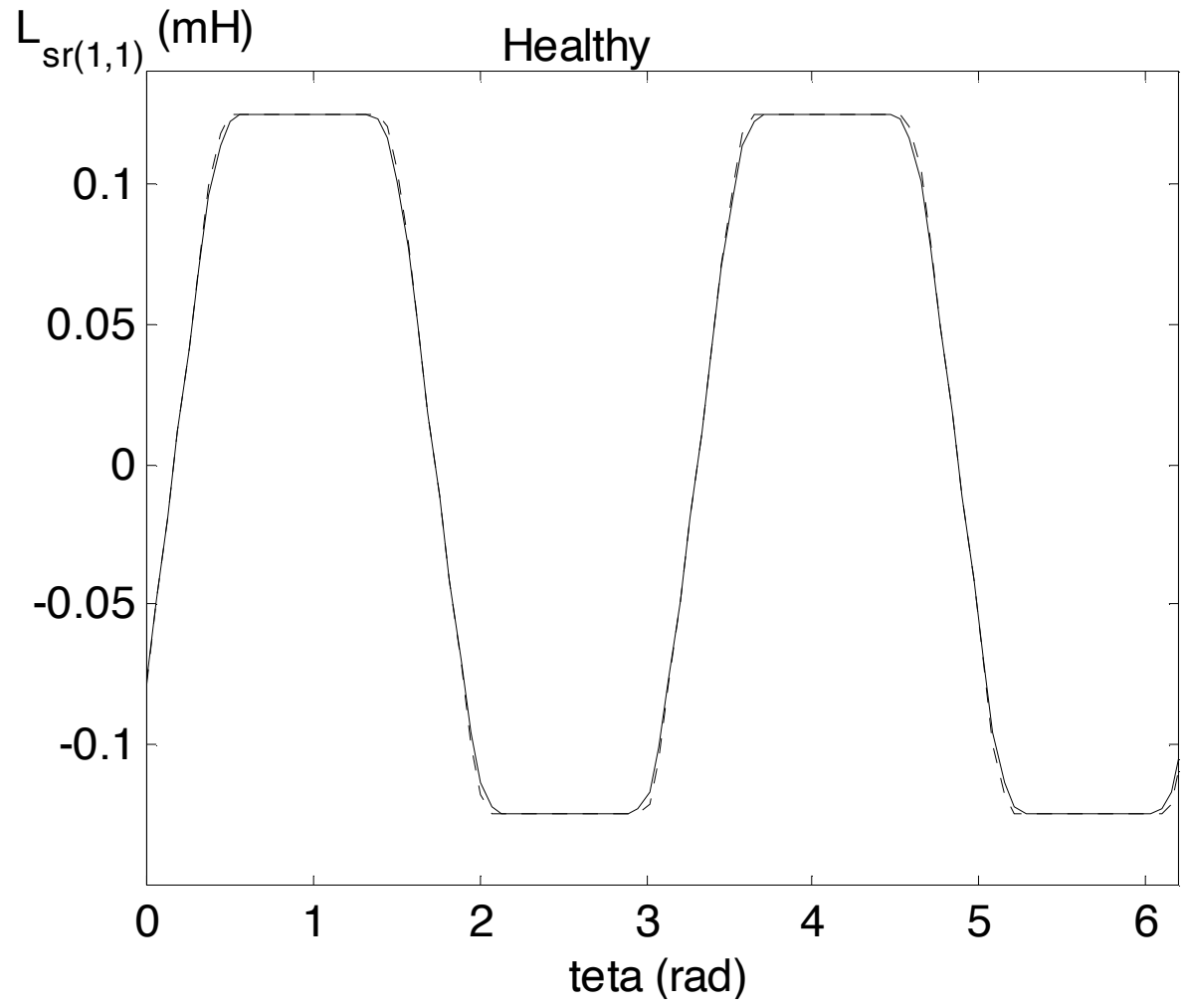

شکل **4-33**: اندوکتانس متقابل فاز اول استاتور و حلقه اول روتور در حالت سالم

خط چین: بدون در نظر گرفتن تورب میله های روتور. توپر : با در نظر گرفتن تورب میله های روتور

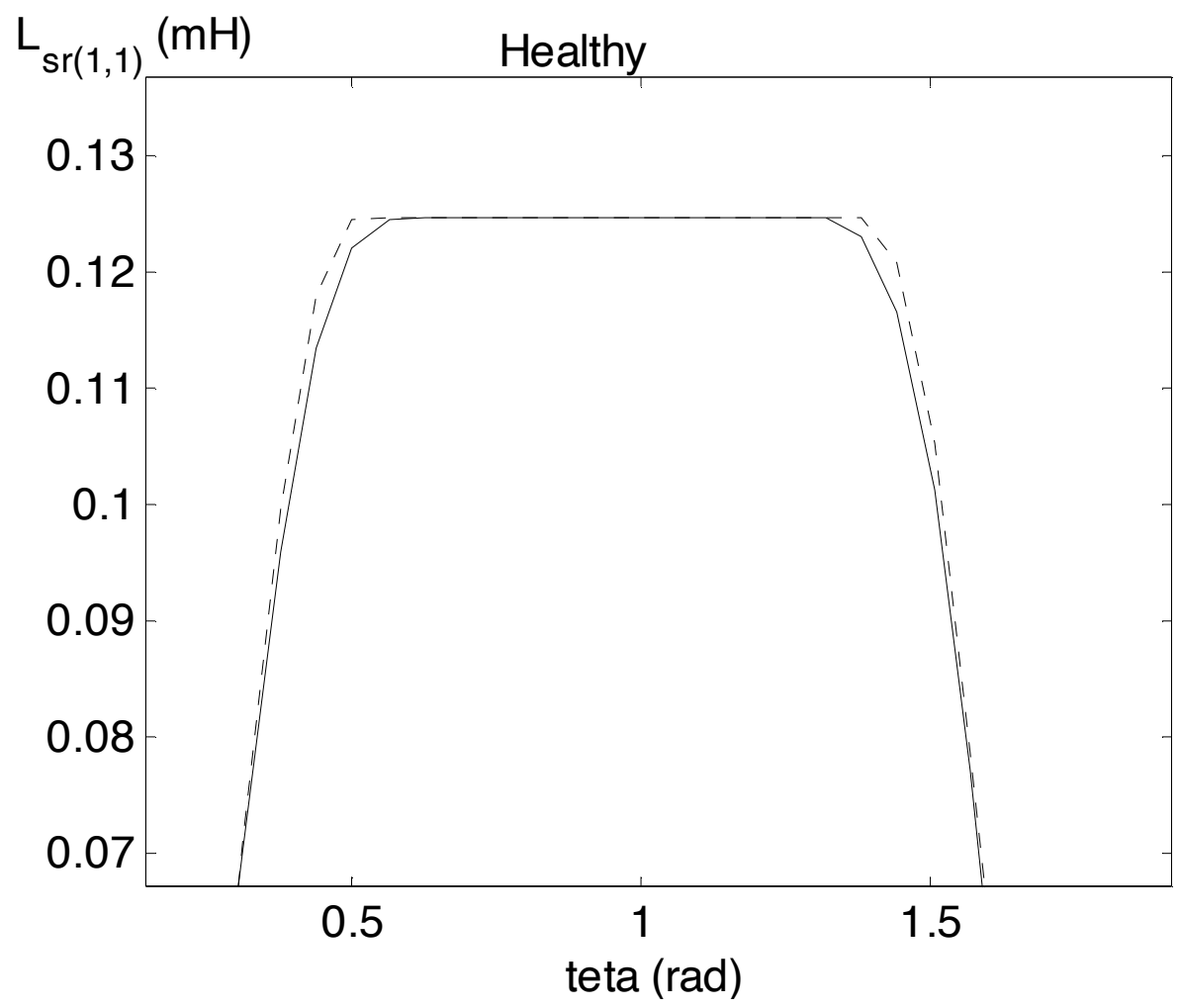

شکل **4-34**: اندوکتانس متقابل فاز اول استاتور و حلقه اول روتور در حالت سالم با بزرگنمایی بیشتر

خط چین: بدون در نظر گرفتن تورب میله های روتور. توپر : با در نظر گرفتن تورب میله های روتور

# فصل پنجم

# تحلیل کامپیوتری موتور القایی سه فاز قفس سنجابی در حالت ناهم محوری بین روتور و استاتور

بکارگیری نظریه تابع سیم پیچی مستلزم محاسبه کلیه اندوکتانسهای خودی و متقابل سیم‌پیچهای استاتور و مدارهای روتور می باشد. در شرایط متقارن این اندوکتانسها به سادگی قابل محاسبه هستند، به عنوان مثال می‌توان اندوکتانس متقابل فاز اول استاتور و حلقه اول روتور را محاسبه و سپس با انجام عملیات جبری ساده‌ای اندوکتانس متقابل بین بقیه فازهای استاتور و حلقه های روتور را بدست آورد. اما با وجود ناهم محوری بین روتور و استاتور وضعیت پیچیده تر است.

ناهم محوری روتور و استاتور که منجر به فاصله هوایی غیر یکنواخت در پیرامون روتور می‌شود، تغییرات هدایت مغناطیسی فاصله هوایی نسبت به موقعیت زاویه‌ای روتور را ایجاد می‌کند. تغییرات هدایت مغناطیسی فاصله هوایی بر روی اندوکتانسهای مختلف موتور نیز تاثیر می گذارد. این تغییرات پیچیدگی معادلات توصیف موتور را به دنبال خواهد داشت، به طوری که در این شرایط تحلیل مدل موتور منجر به صرف زمان طولانی در شبیه سازی کامپیوتری می‌شود. اغلب این شبیه سازیها برای بدست آوردن بررسی تغییرات طیف فرکانسی جریان استاتور است. ناهم محوری روتور و استاتور موجب تولید هارمونیکهای اضافی با فرکانسهای معین در طیف فرکانسی جریان فازهای استاتور می شود.

در اغلب موارد کاربرد نظریه تابع سیم پیچی، به منظور کاهش زمان تحلیل، ماتریس اندوکتانس موتور پیشاپیش محاسبه و مقدار آن در موقعیتهای زاویه ای مختلف روتور ذخیره می شود و تحلیل بر اساس درونیابی از این اطلاعات انجام می گیرد. این مساله در حالتیکه تحلیل موتور در حالت خطای ناهم محوری روتور و استاتور مدنظر است غیر قابل اجتناب است. تا کنون هیچ رابطه تحلیلی میان اندوکتانسهای مختلف فازهای استاتور و حلقه های روتور در حالت ناهم محوری میان روتور و استاتور ارائه نشده و جهت محاسبه هر اندوکتانس در هر موقعیت راویه ای روتور گریزی از حل عددی یک انتگرال نسبتا پیچیده نبوده است.

در فصل چهارم این گزارش یک روش تحلیلی جدید جهت محاسبه اثر حالتهای مختلف ناهم محوری روتور و استاتور بر روی اندوکتانسهای موتور القایی قفس سنجابی ارائه شد. استفاده از این روش تحلیلی باعث تسریع تحلیل مدل کامپیوتری موتور القایی در شرایط ناهم محوری می‌شود. در روشهای دیگر بایستی ماتریس اندوکتانس در موقعیت های مختلف زاویه‌ای روتور ذخیره و از طریق درونیابی مقدار این ماتریس در لحظه مورد

نظر محاسبه شود. این عمل ضمن کاهش دقت در محاسبات مستلزم عملیات حجیم درونیابی دو بعدی است. روش ارائه شده در این گزارش تنها با ذخیره کردن اندوکتانس یکی از فازهای استاتور و حلقه های روتور در حالت سالم و انجام محاسبات تحلیلی ساده قادر به محاسبه کلیه اندوکتانسهای بین فازهای استاتور و حلقه های روتور در هر موقعیت دلخواه زاویه‌ای روتور و هر وضعیت ناهم محوری میان روتور و استاتور می‌باشد. با استفاده از این اندوکتانسها و معادلات تزویج الکترومغناطیسی و الکترومکانیکی موتور می توان عملکرد موتور را تحلیل و بررسی نمود.

## 5-1- معادلات توصیف موتور در حالت ناهم محوری روتور و استاتور

معادلات دیفرانسیل توصیف موتور القایی در شرایط ناهم محوری روتور و استاتور شامل معادلات تزویج الکترومغناطیسی و معادلات الکترومکانیکی بصورت زیرند:

$$[\psi_s] = [L_s][I_s] + [L_{sr}][I_r] \qquad (1\text{-}5)$$

$$[\psi_r] = [L_{rs}][I_s] + [L_s][I_s] \qquad (2\text{-}5)$$

$$[V_s] = [R_s][I_s] + \tfrac{d}{dt}[\psi_s] \qquad (3\text{-}5)$$

$$[V_r] = [R_r][I_r] + \tfrac{d}{dt}[\psi_r] \qquad (5\text{-}5)$$

$$\mathrm{T_m} = [\mathrm{I_s}]'\left(\tfrac{\mathrm{d}}{\mathrm{dt}}[\mathrm{L_{sr}}]\right)[\mathrm{I_r}] + 0.5[\mathrm{I_s}]'\left(\tfrac{\mathrm{d}}{\mathrm{dt}}[\mathrm{L_s}]\right)[\mathrm{I_s}] + 0.5[\mathrm{I_r}]'\left(\tfrac{\mathrm{d}}{\mathrm{dt}}[\mathrm{L_r}]\right)[\mathrm{I_r}] \qquad (5\text{-}5)$$

$$J\,\tfrac{d}{dt}\omega_r = T_m - T_l \qquad (6\text{-}5)$$

$$\tfrac{d}{dt}\theta_r = \omega_r \qquad (7\text{-}5)$$

در این معادلات :

$$\left[F_s\right]=\left[f_1^s \quad f_2^s \quad f_3^s\right]^T \tag{8-5}$$

$$\left[F_r\right]=\left[f_1^r \quad f_2^r \quad \dots \quad f_n^r\right]^T \tag{9-5}$$

$n$ تعداد میله های روتور است و $F$ می تواند $I$ ، $V$ یا $\psi$ باشد.

ماتریسهای مقاومت روتور و استاتور بصورت زیر محاسبه می‌شوند.

$$\left[R_S\right]=\begin{bmatrix} r_S & 0 & 0 \\ 0 & r_S & 0 \\ 0 & 0 & r_S \end{bmatrix} \tag{10-5}$$

$$\left[R_r\right]=\left[r_{ij}\right]_{n\times n}:\begin{cases} r_{ij}=2(r_b+r_e) & i=j \\ r_{ij}=-r_b & (i=j\pm 1)\,or\,(i=1,j=n)\,or\,(i=n,j=1) \\ r_{ij}=0 & else \end{cases} \tag{11-5}$$

ماتریس $R_r$ با نوشتن معادلات $KVL$ در حلقه های روتور بدست می‌آید. انواع مختلف ناهم محوری بر روی اندوکتانسهای استفاده شده در روابط (5-1) تا (5-7) تاثیر می‌گذارد. روش محاسبه این اندوکتانسها در فصل سوم ارائه شد.

## 5-2- تحلیل کامپیوتری موتور درحالتهای مختلف ناهم محوری

به منظور حل معادلات (5-1) تا (5-7) از روش رانگ کوتای مرتبه 4 و 5 نرم افزار MATLAB استفاده شده است. همچنین موارد زیر به منظور ساده سازی مدنظر قرار گرفته‌اند:

1. ماشین در ناحیه خطی منحنی مغناطیسی کار می‌کند و به اشباع نمی‌رود.
2. توزیع جریان هر فاز بصورت یکنواخت است.
3. موتور توسط یک منبع سینوسی سه فاز متقارن تغذیه می شود.

4. بار مکانیکی ثابت فرض شده است.

موارد زیر نقاط قوت روش مذکور هستند.

1. امکان تحلیل عملکرد موتور در حالت ناهم محوری مرکب میان روتور و استاتور
2. در نظر گرفتن اثر تورب میله های روتور
3. عدم استفاده از عملیات درونیابی عددی
4. عدم نیاز به حافظه کامپیوتر جهت ذخیره اندوکتانسهای موتور

در ادامه نتایج حاصل از شبیه سازی کامپیوتری موتور در حالت سالم و ناهم محوریهای مختلف میان روتور و استاتور بررسی، و با یکدیگر مقایسه می‌شوند. کلیه این شکل موجها برای شرایط تغذیه و بارمکانیکی زیر می‌باشند :

**جدول 5-1: شرایط تحلیل کامپیوتری موتور**

| فرکانس تغذیه | دامنه ولتاژ تغذیه | ممان اینرسی بار مکانیکی | گشتاور بار مکانیکی |
|---|---|---|---|
| $50Hz$ | $380Volt$ | $0.05\, {N.m}/{s^2}$ | $20N.m$ |

شکلهای 5-1 تا 5-5 به ترتیب منحنیهای سرعت مکانیکی روتور، جریان سه فاز استاتور و گشتاور تولیدی موتور را تا رسیدن به حالت دائمی نشان می دهند. شکلهای 5-6 تا 5-10، 5-11 تا 5-15 و 5-16 تا 5-20 به ترتیب این منحنیها را برای حالت 20٪ ناهم محوری ایستا، 15٪ ناهم محوری پویا و ناهم محوری مرکب شامل 20٪ ناهم محوری ایستا و 15٪ ناهم محوری پویا نشان می‌دهند.

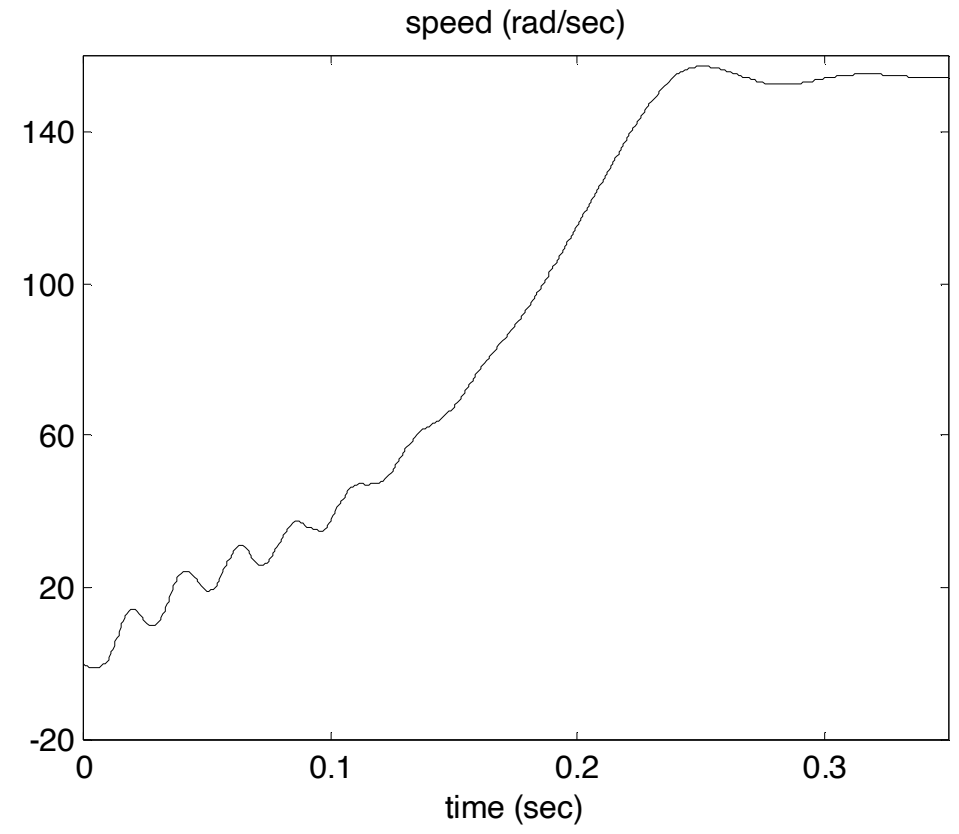

شکل **5-1** : منحنی تغییرات سرعت موتور در حالت سالم

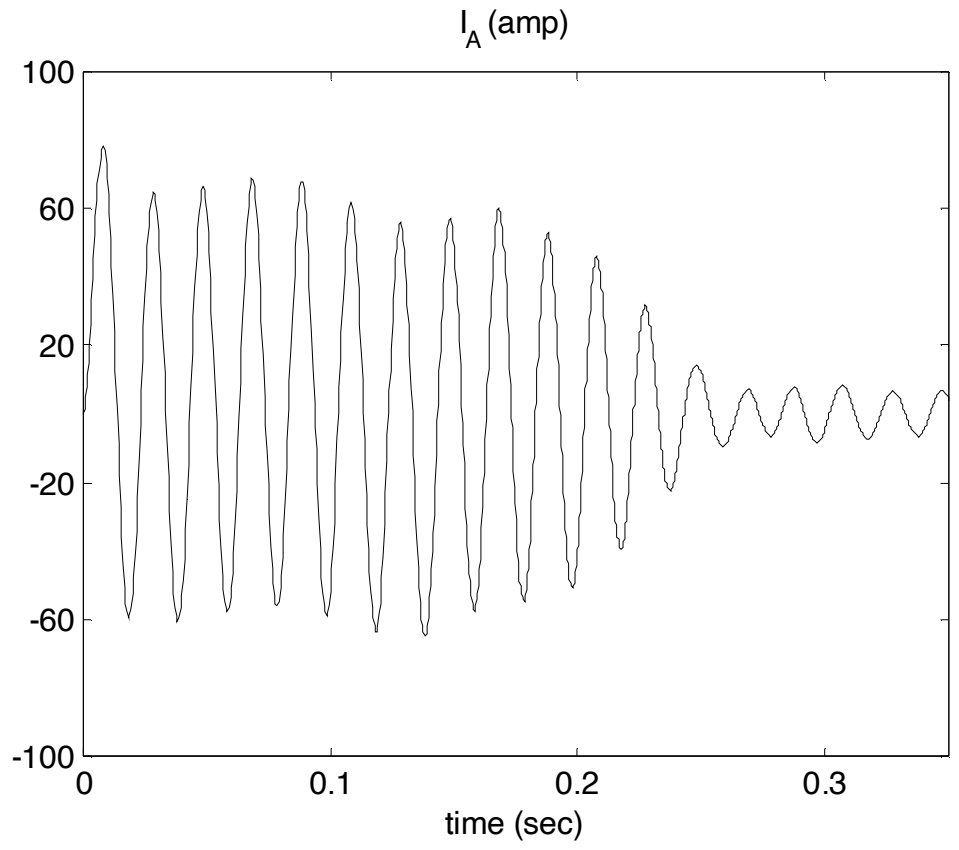

شکل **5-2** : منحنی تغییرات جریان فاز اول موتور در حالت سالم

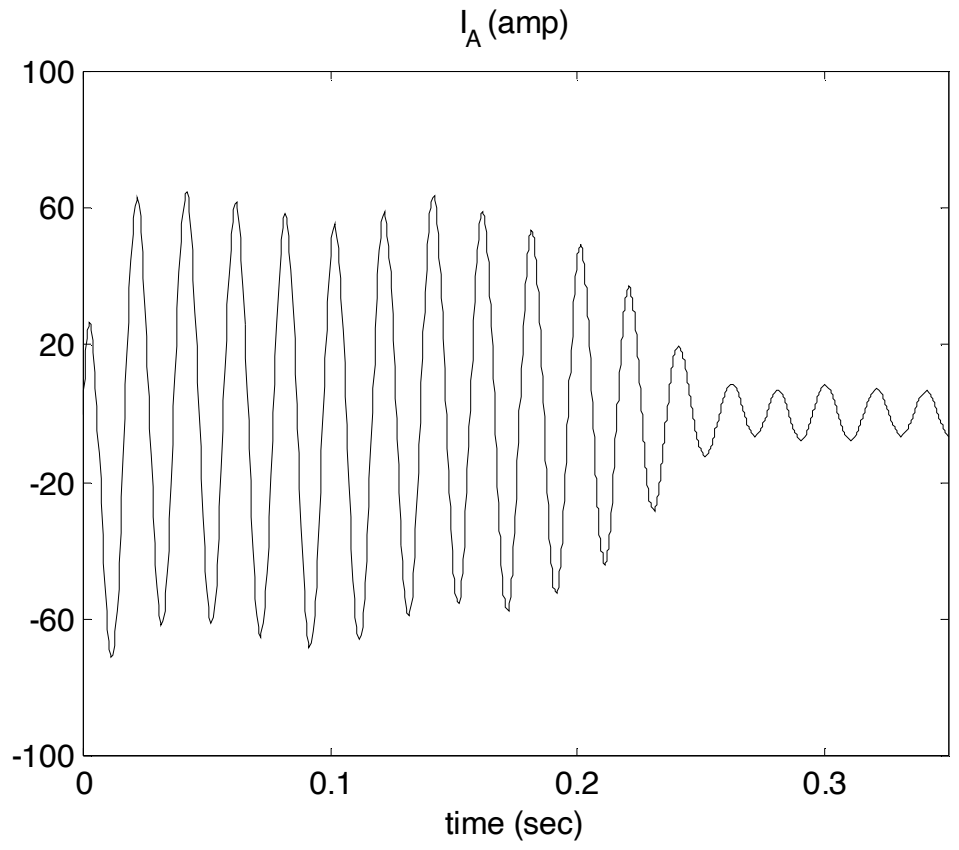

شکل **5-3** : منحنی تغییرات جریان فاز دوم موتور در حالت سالم

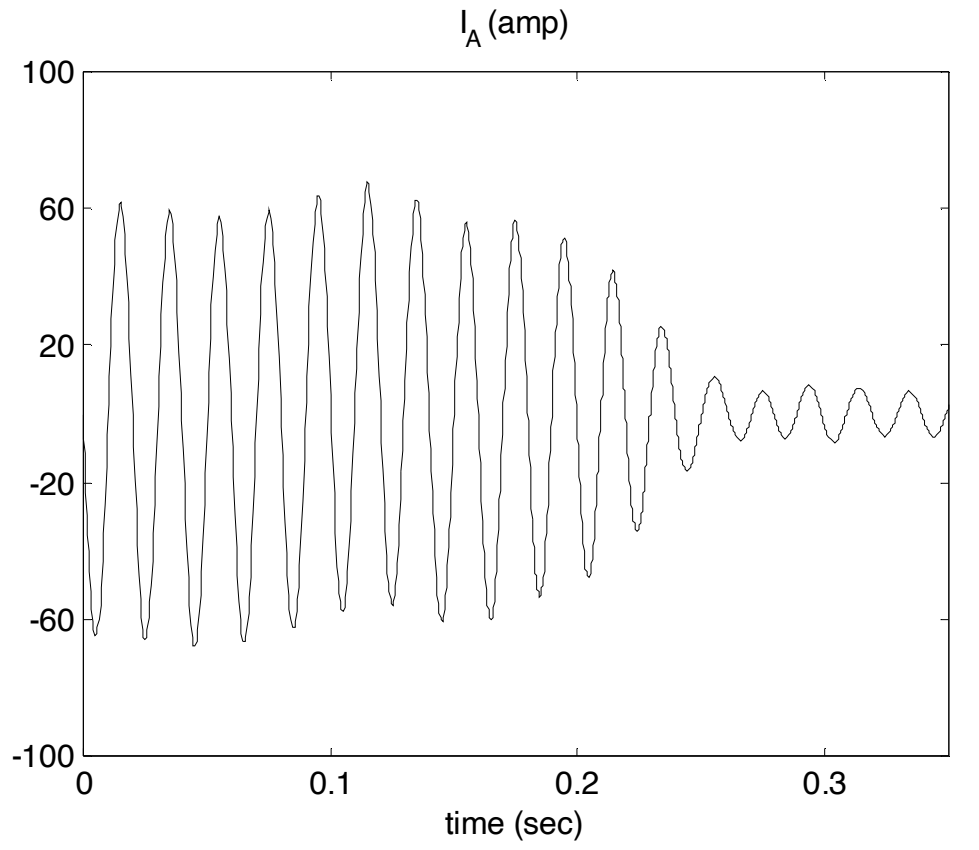

شکل 5-5 : منحنی تغییرات جریان فاز سوم موتور در حالت سالم

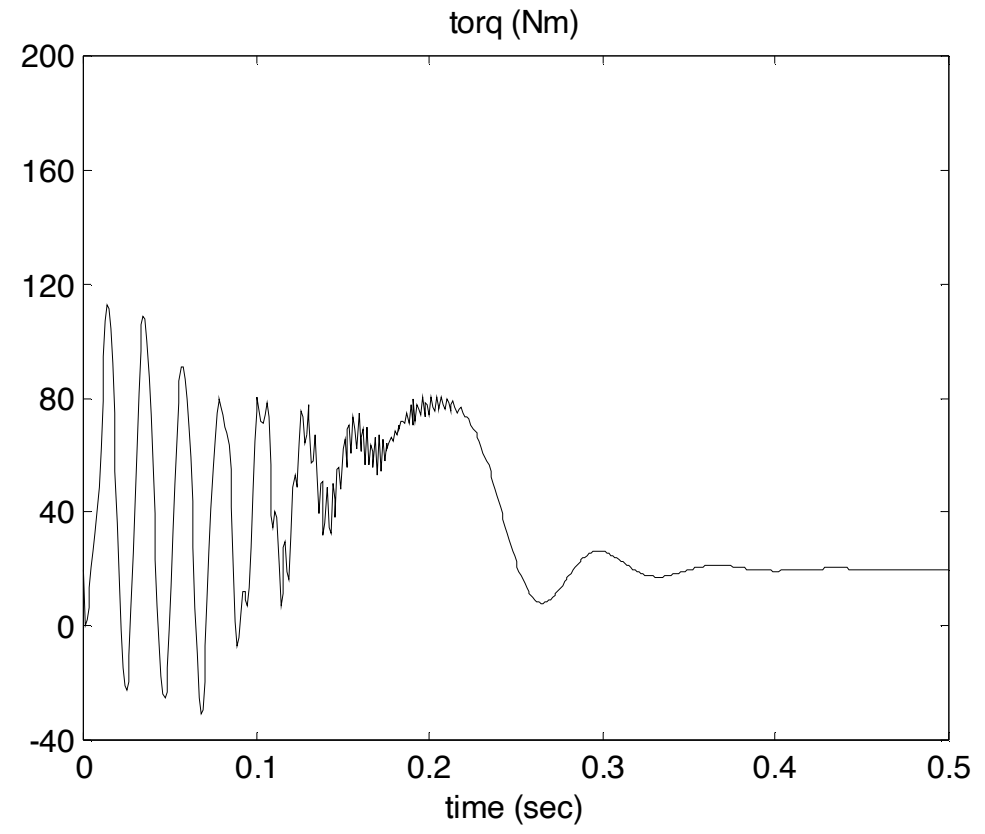

شکل 5-5 : منحنی تغییرات گشتاور موتور در حالت سالم

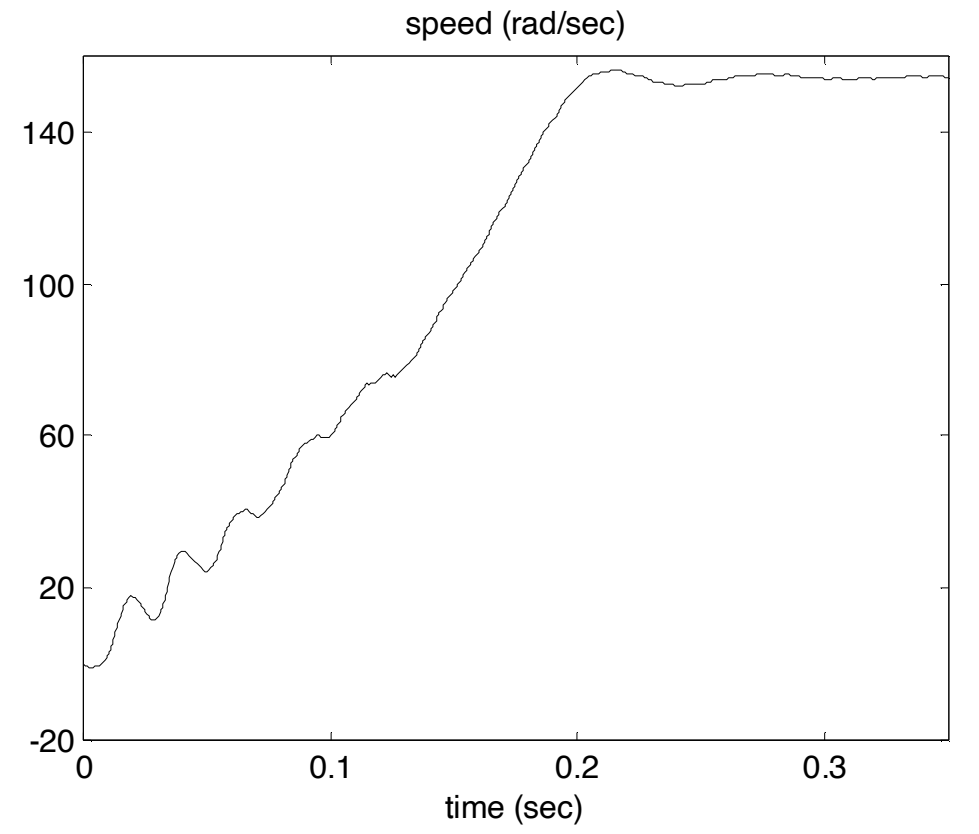

شکل 5-6 : منحنی تغییرات سرعت موتور در حالت 20٪ ناهم محوری ایستا

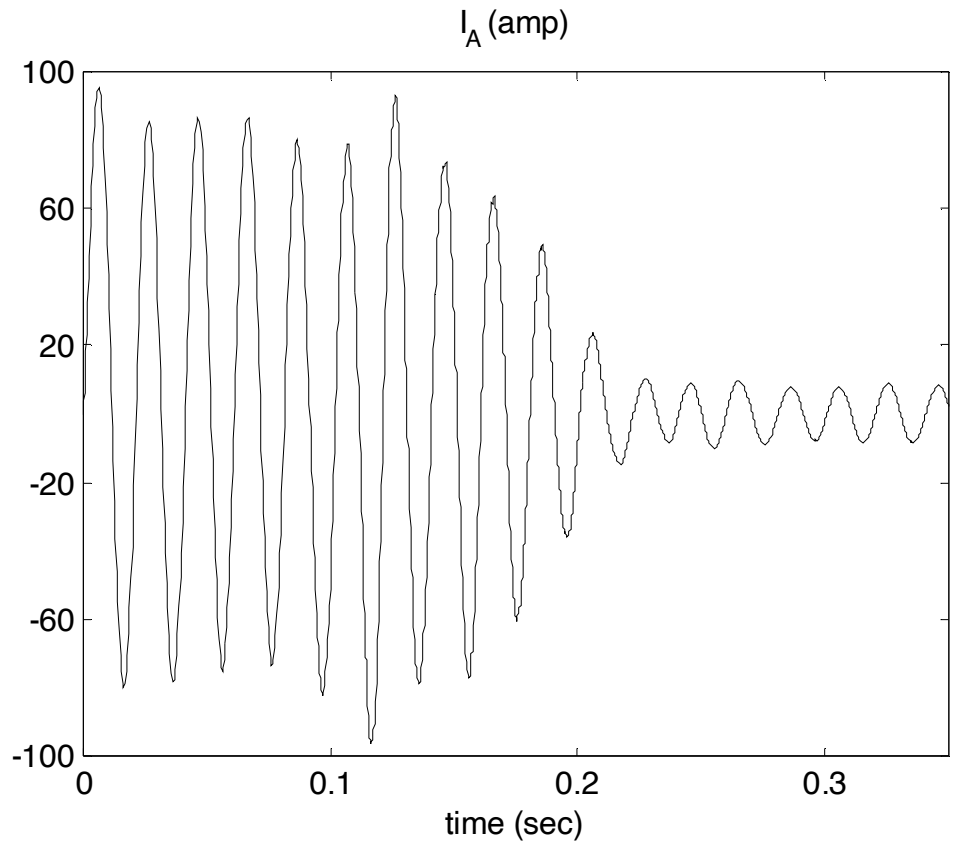

شکل 5-7 : منحنی تغییرات جریان فاز اول موتور در حالت 20٪ ناهم محوری ایستا

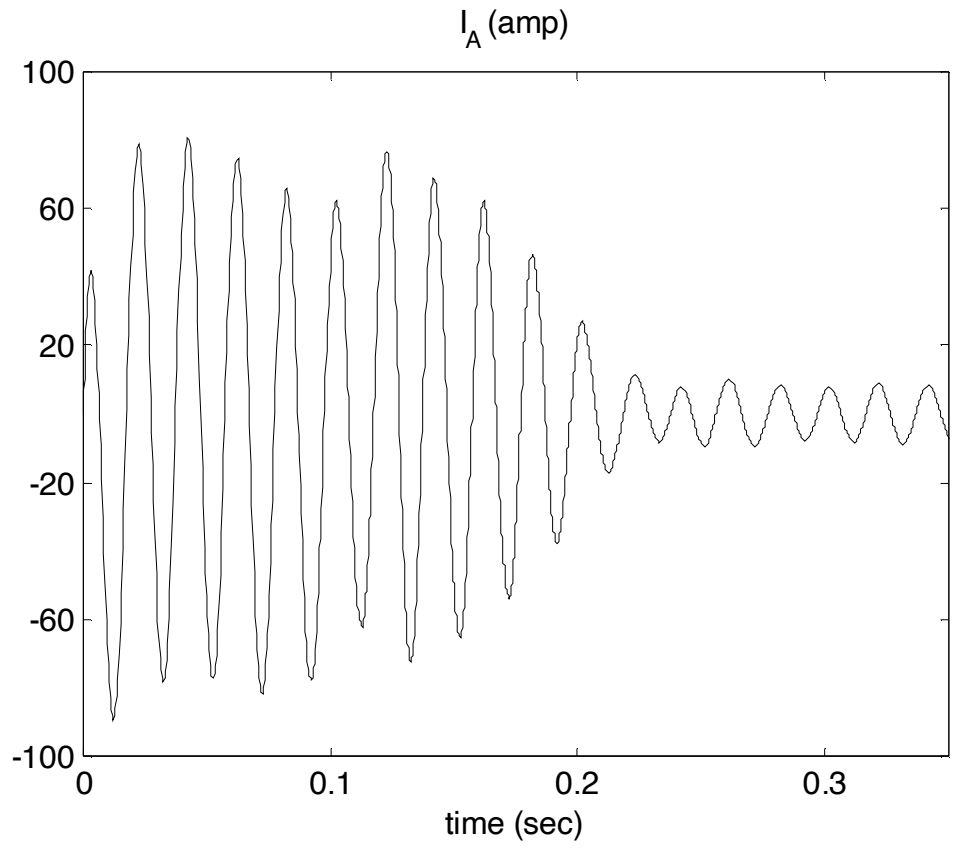

شکل 5-8 : منحنی تغییرات جریان فاز دوم موتور در حالت 20٪ ناهم محوری ایستا

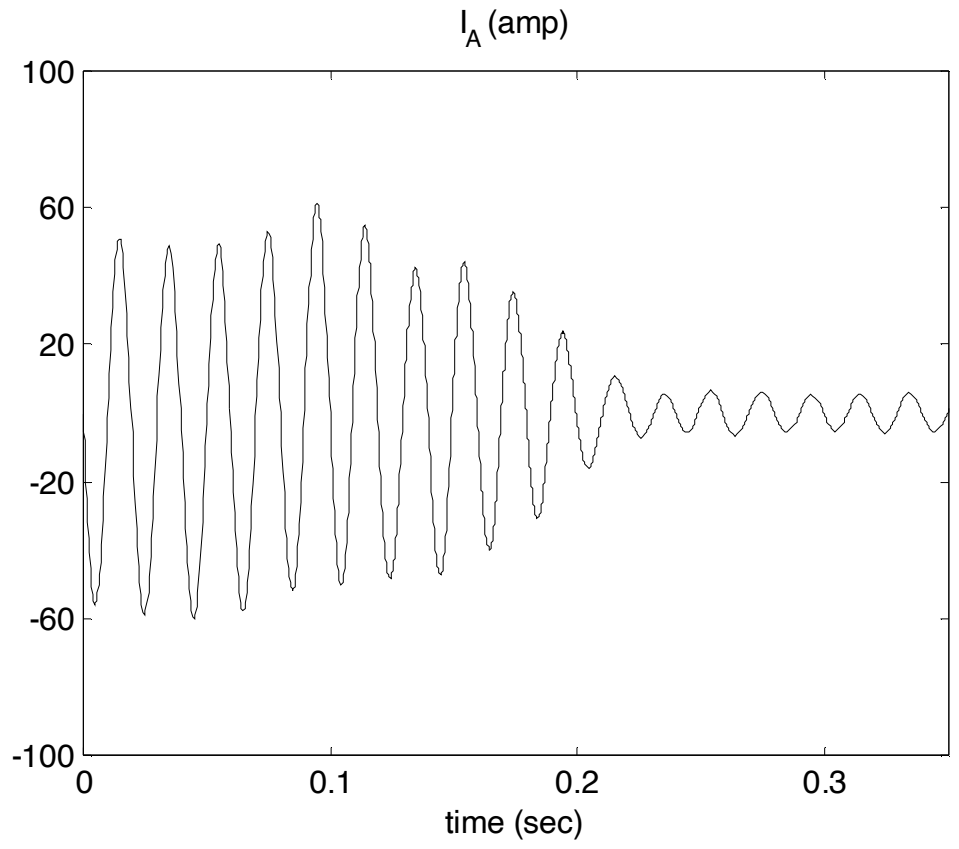

شکل 5-9 : منحنی تغییرات جریان فاز سوم موتور در حالت 20٪ ناهم محوری ایستا

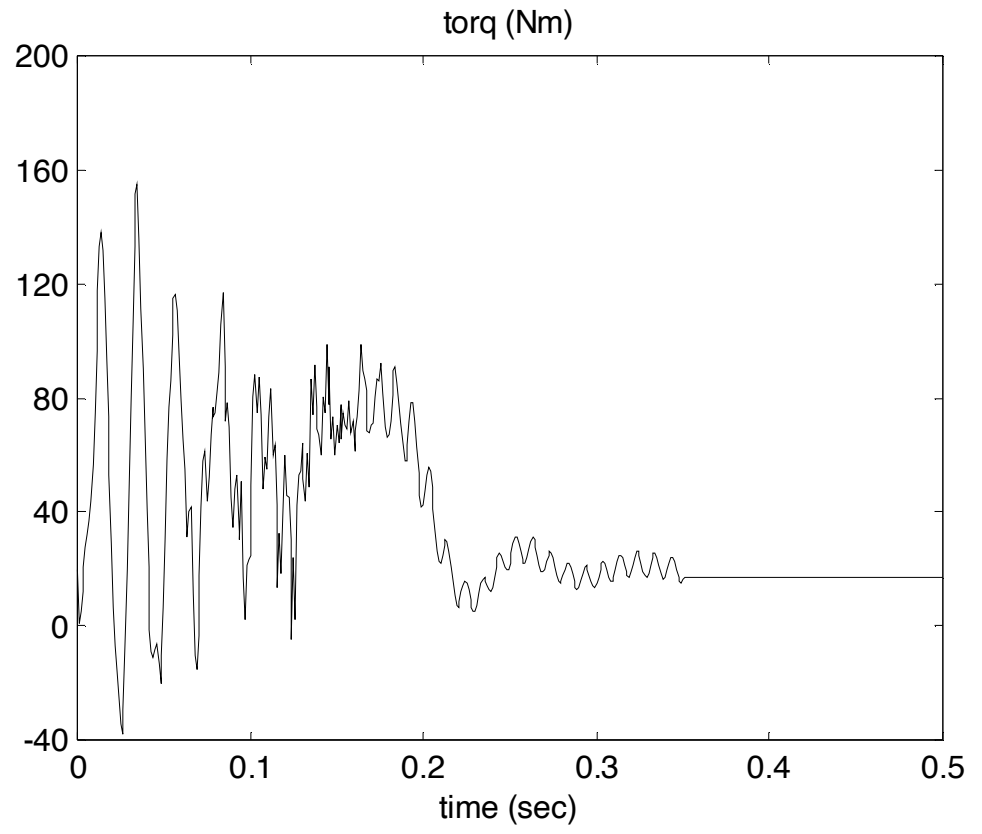

شکل **5-10** : منحنی تغییرات گشتاور موتور در حالت **20**٪ ناهم محوری ایستا

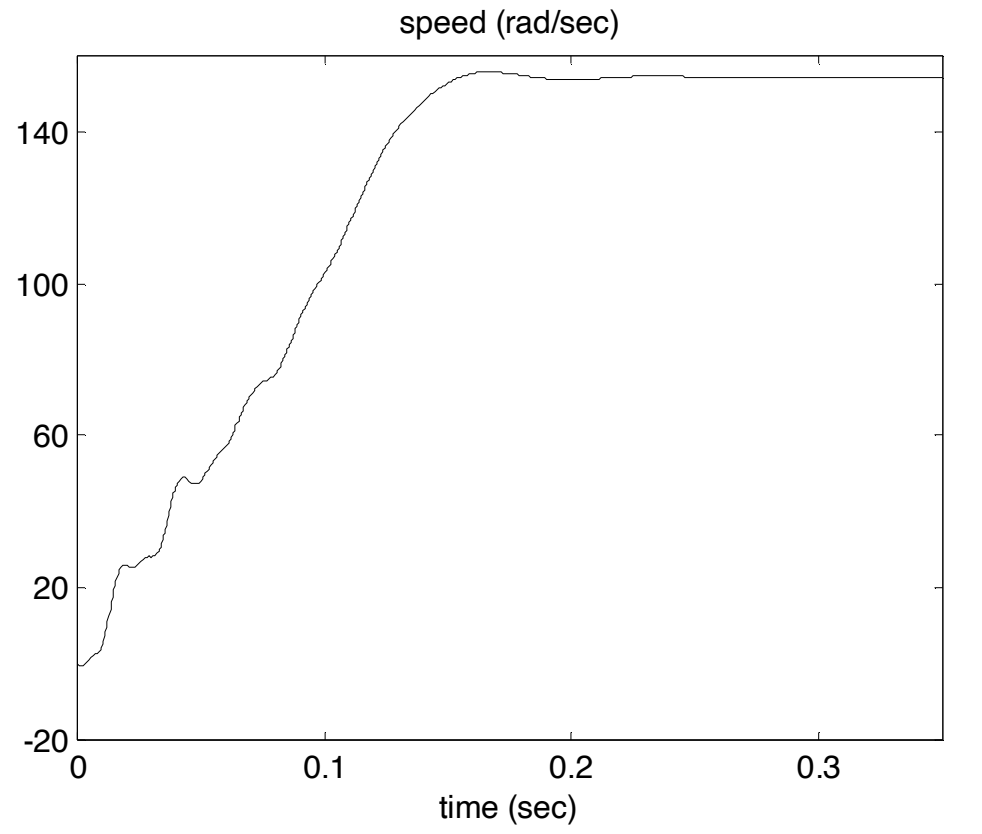

شکل **5-11** : منحنی تغییرات سرعت موتور در حالت **15**٪ ناهم محوری پویا

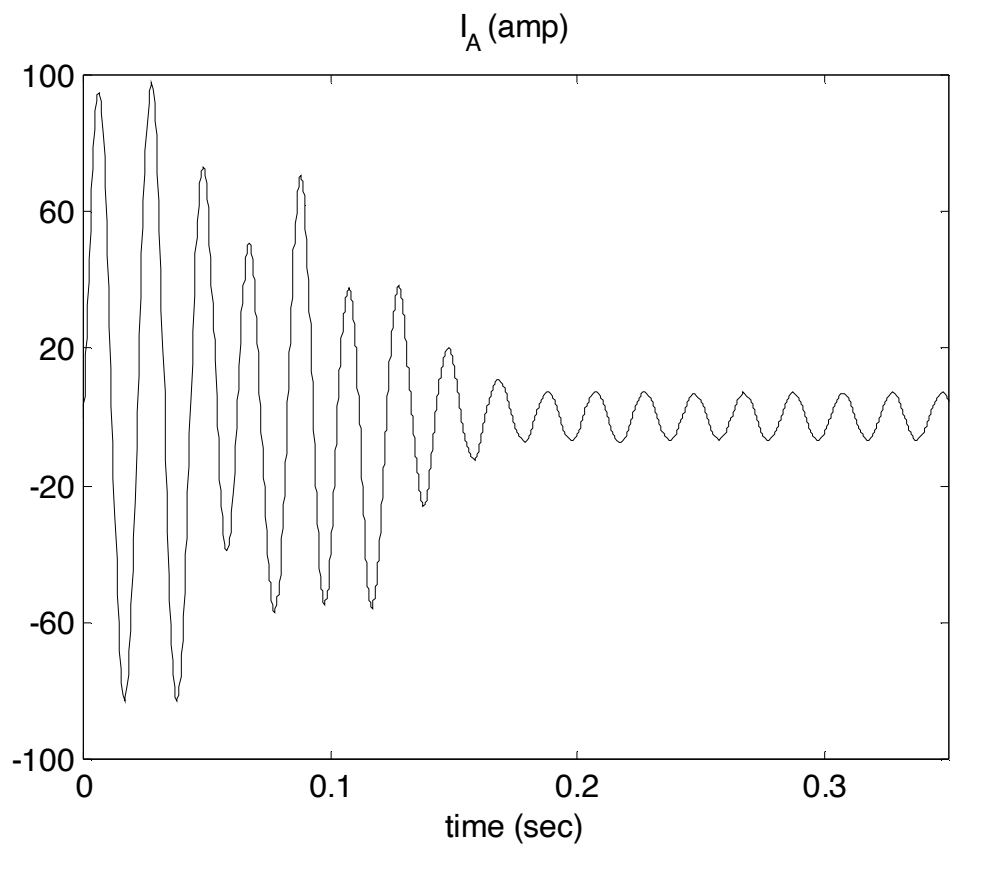

شکل **5-12** : منحنی تغییرات جریان فاز اول موتور در حالت **15**٪ ناهم محوری پویا

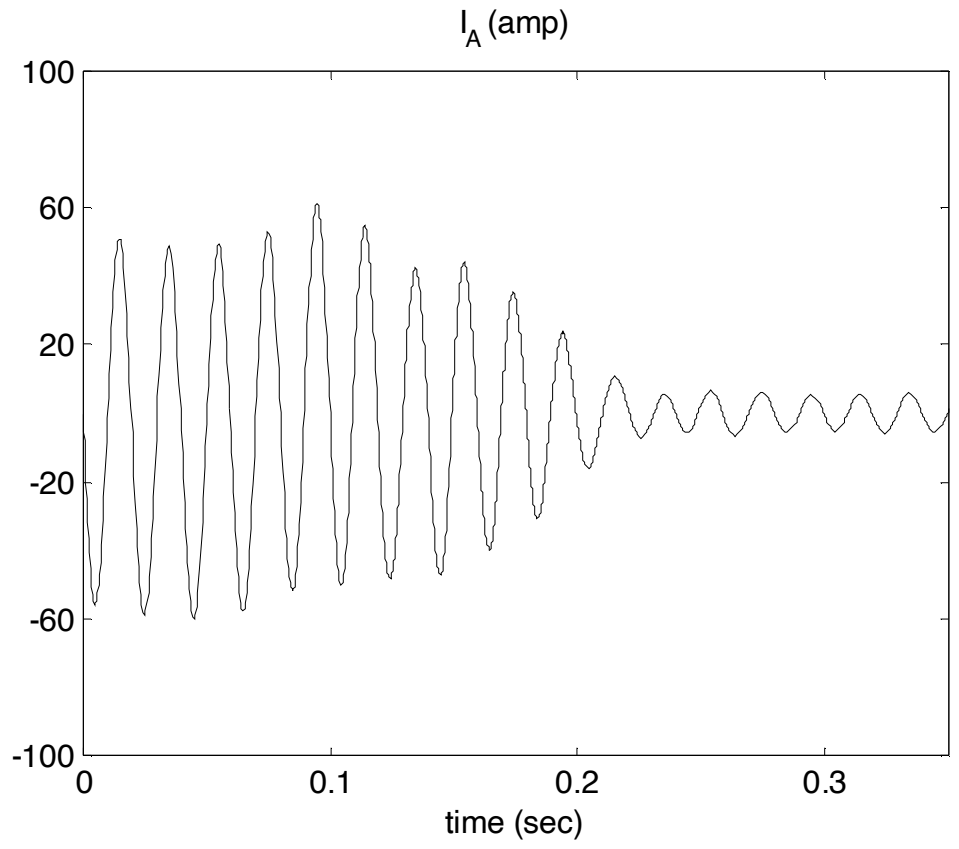

شکل **5-13** : منحنی تغییرات جریان فاز دوم موتور در حالت **15**٪ ناهم محوری پویا

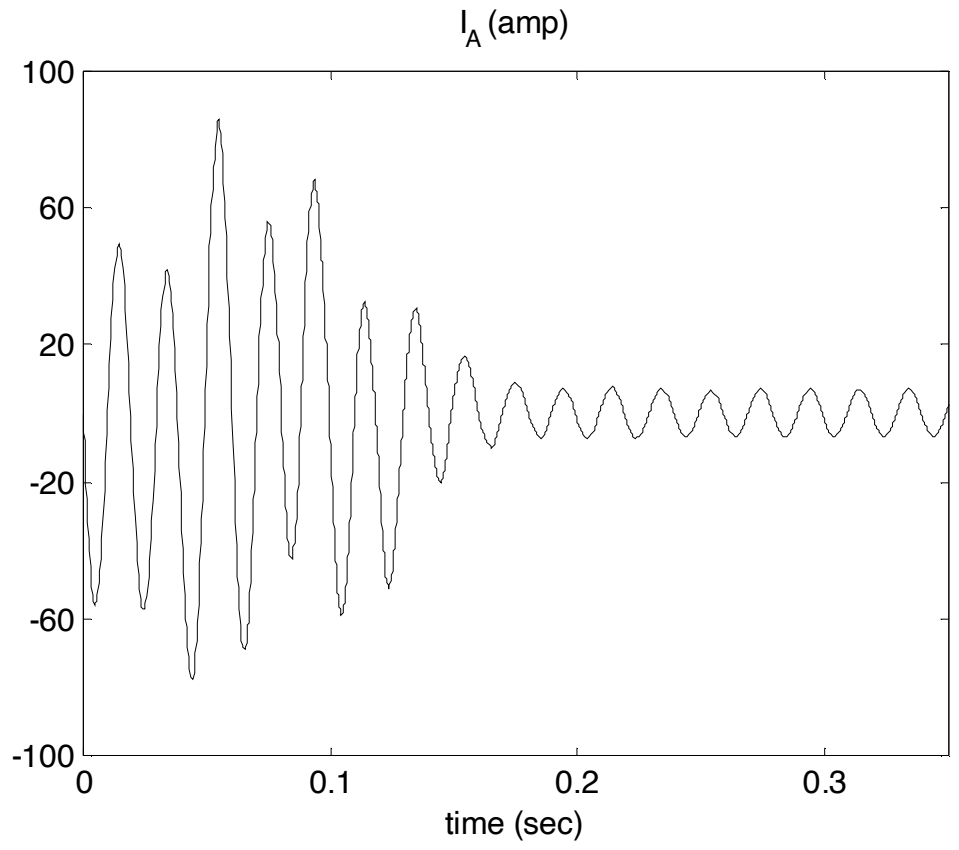

شکل 5-14 : منحنی تغییرات جریان فاز سوم موتور در حالت 15٪ ناهم محوری پویا

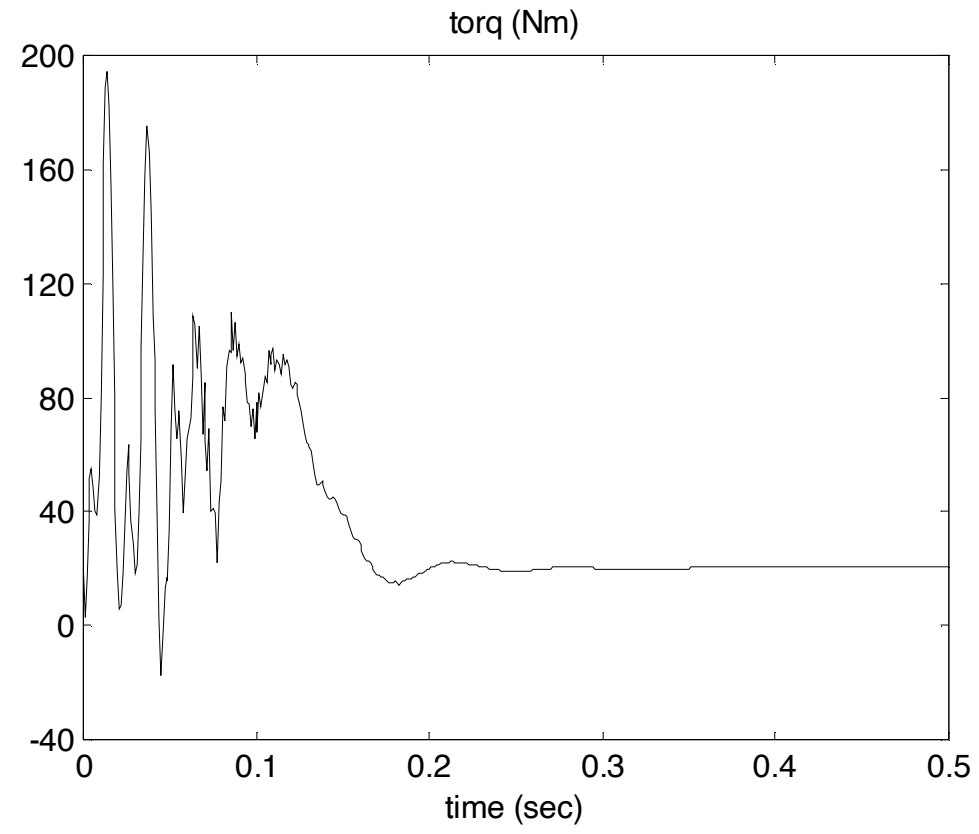

شکل 5-15 : منحنی تغییرات گشتاور موتور در حالت 15٪ ناهم محوری پویا

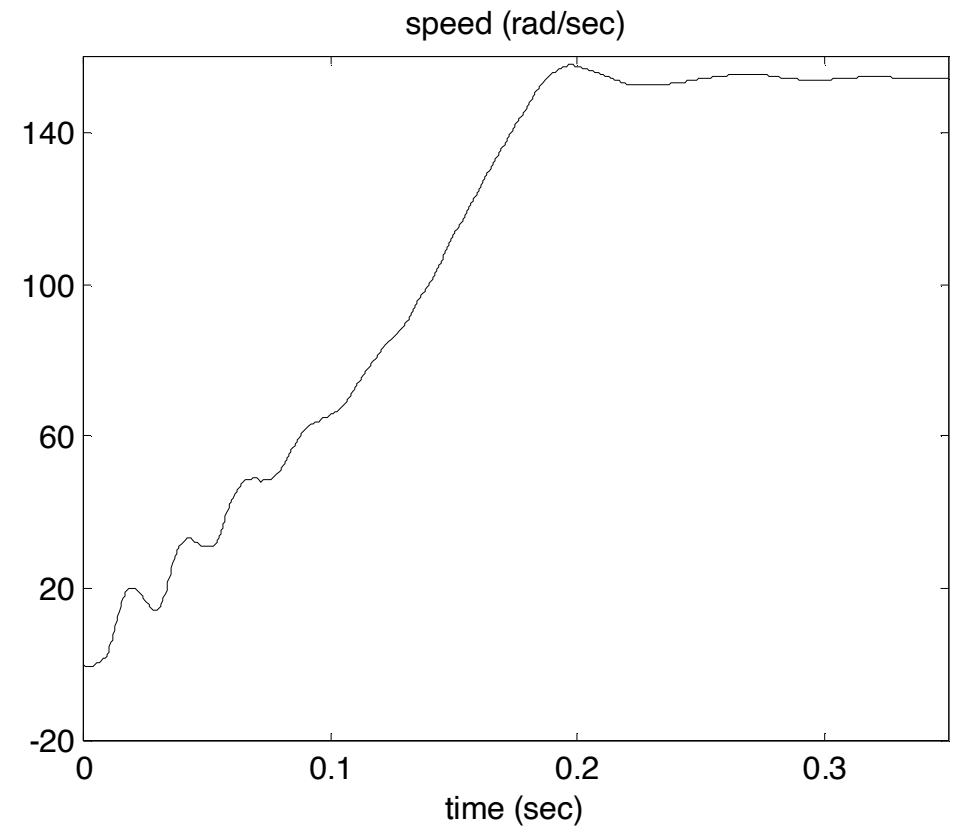

شکل 5-16 : منحنی تغییرات سرعت موتور در حالت 20٪ ناهم محوری ایستا و 15٪ ناهم محوری پویا

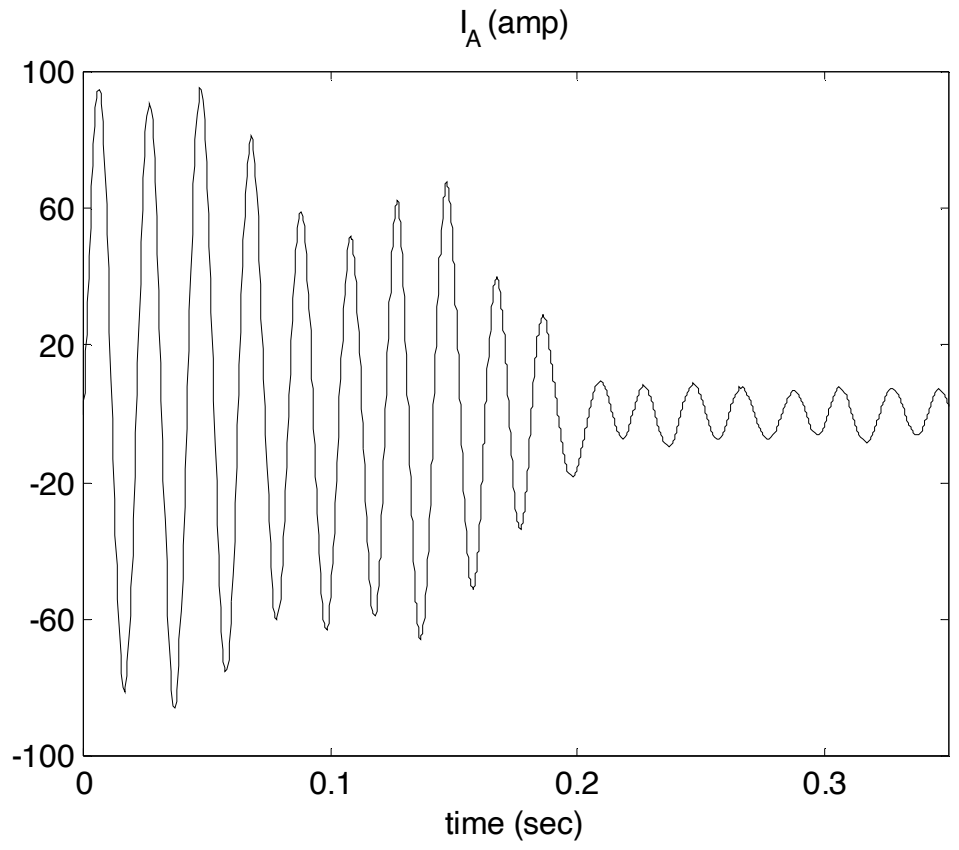

شكل 5-17 : منحنی تغییرات جریان فاز اول موتور در حالت 20٪ ناهم محوری ایستا و 15٪ ناهم محوری پویا

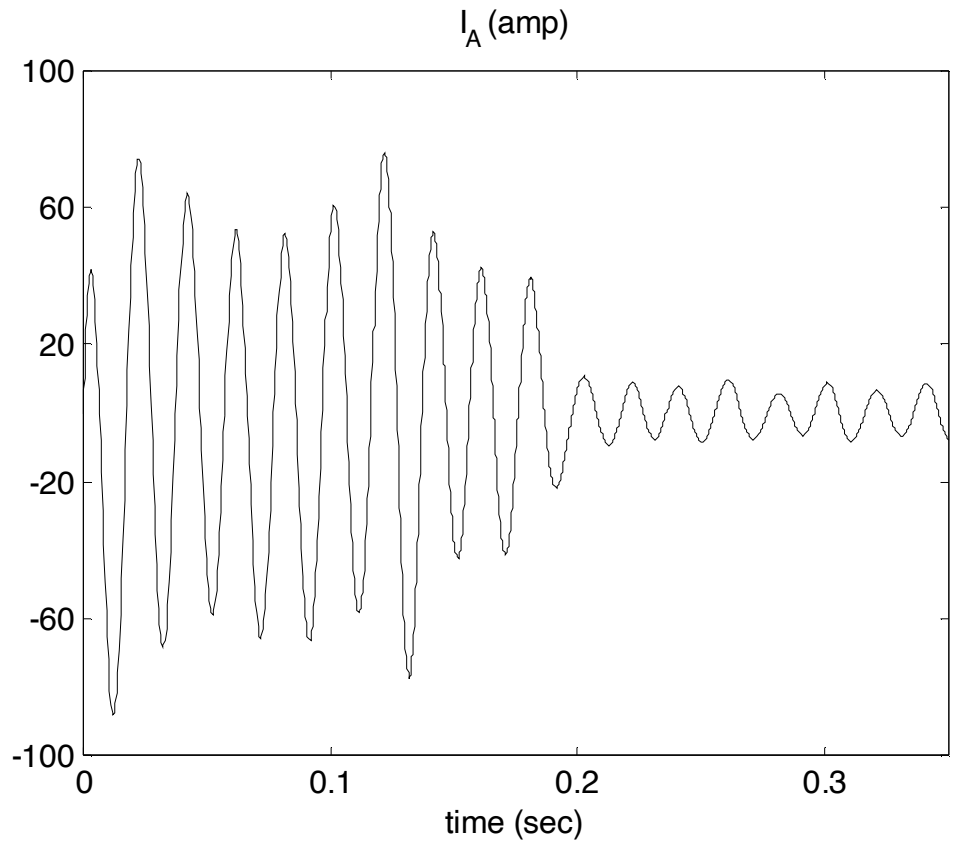

شكل 5-18 : منحنی تغییرات جریان فاز دوم موتور در حالت 20٪ ناهم محوری ایستا و 15٪ ناهم محوری پویا

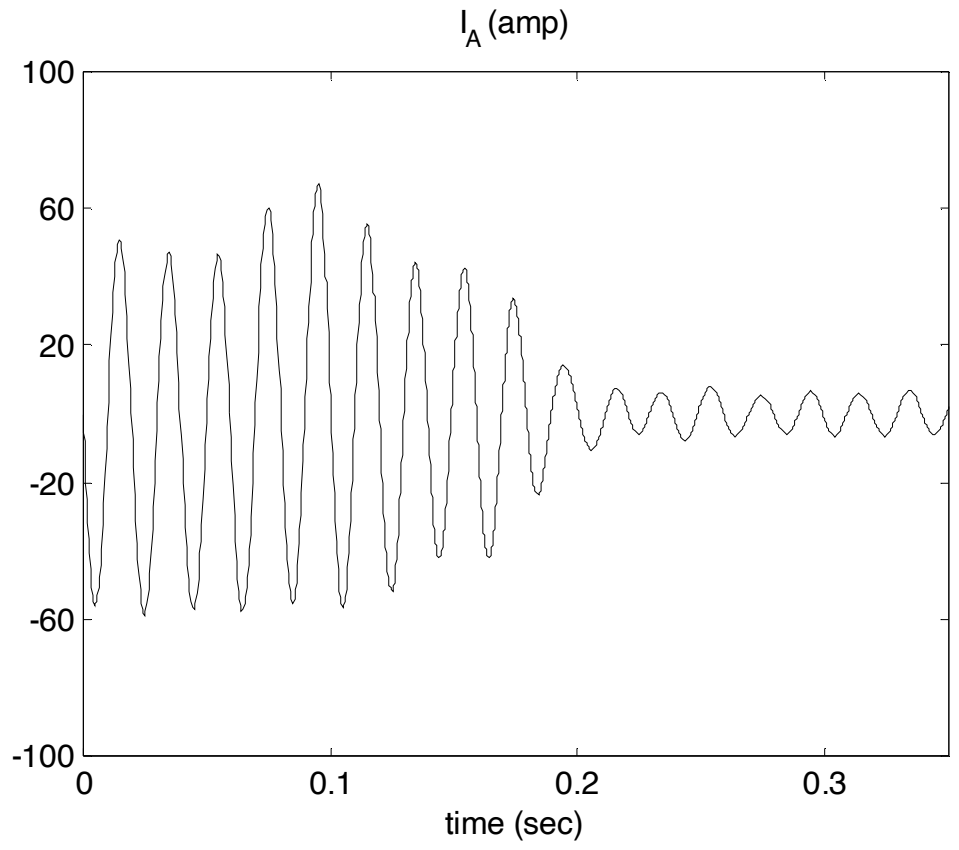

شكل 5-19 : منحنی تغییرات جریان فاز سوم موتور در حالت 20٪ ناهم محوری ایستا و 15٪ ناهم محوری پویا

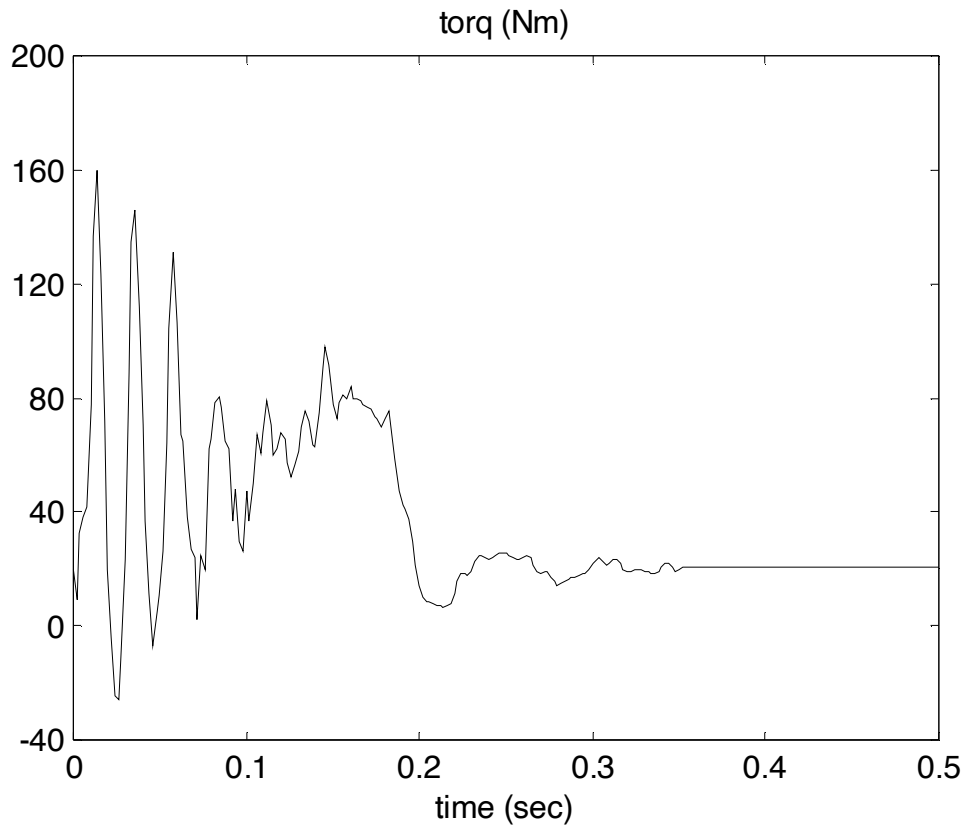

**شکل 5-20** : منحنی تغییرات گشتاور موتور در حالت 20٪ ناهم محوری ایستا و 15٪ ناهم محوری پویا

نتایج حاصل از شبیه سازی نشان می‌دهند که در حالتهای مختلف ناهم محوری، موتور زودتر از حالت سالم به وضعیت ماندگار می‌رسد. این امر باعث می‌شود که جریان و گشتاور راه اندازی موتور در این حالتها نسبت به حالت سالم افزایش پیدا کند. منحنیهای (5-1) تا (5-12) به خوبی گویای این فرضیه می باشند.

در حالت ناهم محوری مرکب موتور در زمانی بین زمان رسیدن به حالت ماندگار حالت ناهم محوری ایستا و پویا به حالت ماندگار می‌رسد.  با اجرای برنامه شبیه سازی در حالتهای مختلف ناهم محوری ایستا و پویا می‌توان پیش بینی کرد که تحت شرایط مختلف ناهم محوری ایستا و یا پویا موتور پس از چه زمانی به وضعیت ماندگار خود می‌رسد و یا اینکه با اندازه‌گیری این زمان میزان ناهم محوری ایستا و پویا را محاسبه کرد. منحنیهای شکل (5-21) به ترتیب، زمان اولین فراجهش موتور را بر حسب درجه ناهم محوری ایستا و پویا برای موتور شرایط جدول 5-1 نشان می‌دهند. اطلاعات موجود در این نمودارها پس از تحلیل کامپیوتری موتور در این حالتها و بررسی منحنی سرعت آن محاسبه شده اند.

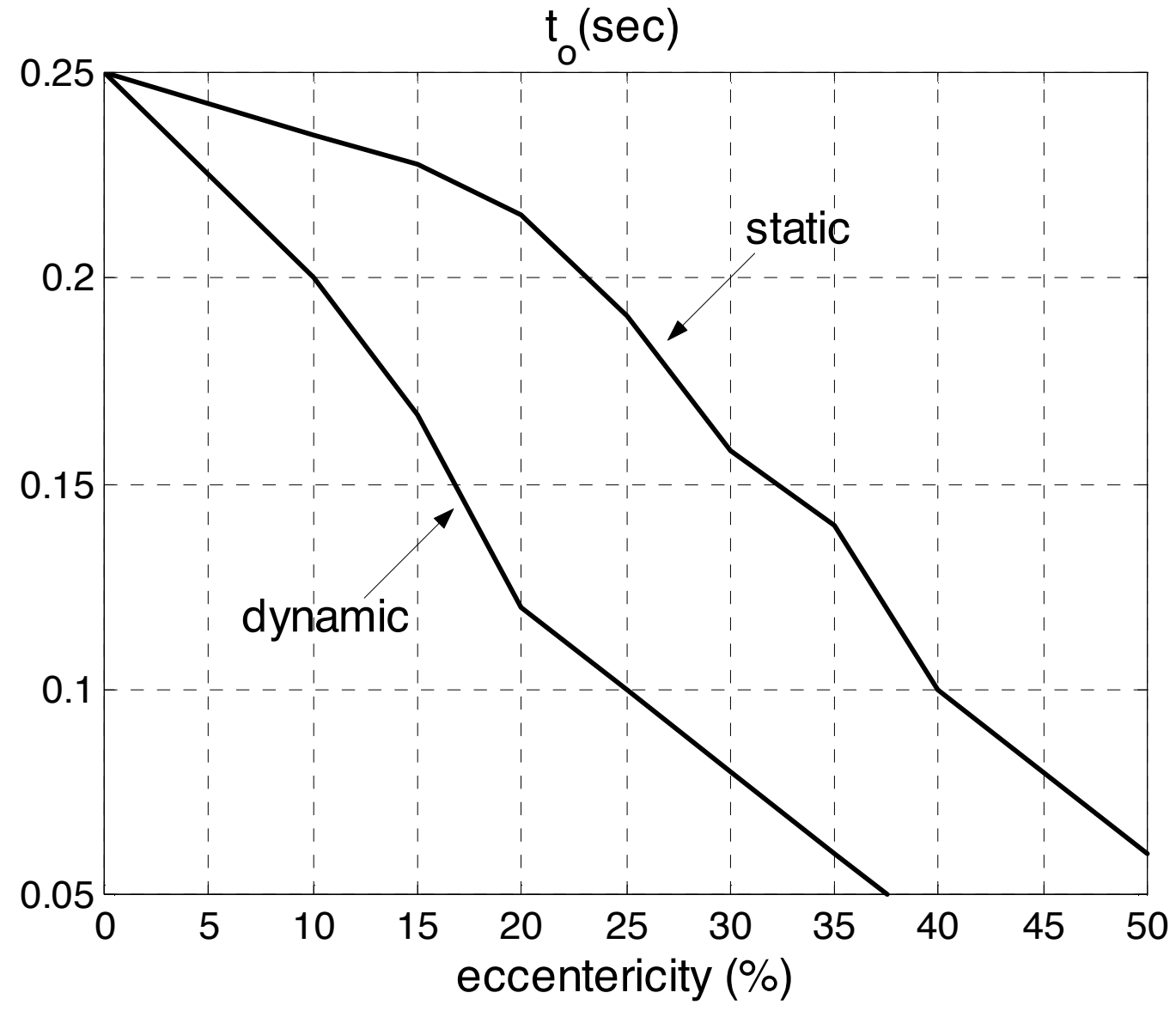

شکل **5-21** : زمان اولین فراجهش در حالتهای ناهم محوری ایستا و پویا بر حسب درجه ناهم محوری

همچنین این نتایج نشان می دهند که موتوری که دارای ناهم محوری پویای x است زودتر از موتوری که دارای درجه ناهم محوری ایستای x است به حالت دائمی می‌رسد. این نکته در منحنیهای (5-5) و (5-7) نیز با اینکه در جه ناهم‌محوری ایستا بیشتر از ناهم محوری پویا است مشاهده می‌شود.

معیار نشان داده شده در شکل (5-21) و یا هر معیار دیگری که مربوط به زمان رسیدن پارامترهای موتور به نقطه و یا ناحیه معینی است- در عمل نمی‌توانند معیار خوبی برای تشخیص ناهم محوری باشند. اولین و مهمترین دلیل این امر وابستگی این زمانها به پارامترهای مکانیکی بار موتور است. منحنی نشان داده شده در شکل (5-21) زمانهای فراجهش را برای یک بار مکانیکی با ممان اینرسی و گشتاور بار ثابت نشان می‌دهد.

## 5- 3- بررسی طیف فرکانسی استاتور در حالتهای مختلف ناهم محوری

هارمونیکهای ناشی از شیارگذاری قفس روتور و ناهم‌محوری میان روتور و استاتور از رابطه زیر محاسبه می‌شوند[46].

$$f_{har} = f_h(k, n_d, s, \eta) = \left[ (kR \pm n_d)\frac{1-s}{p} \pm \eta \right] f_0 \tag{12-5}$$

در این رابطه:

$f_0$ : فرکانس اصلی تغذیه

$s$ : لغزش روتور

$p$ : تعداد زوج قطب استاتور

$R$ : تعداد شیارهای روتور

$n_d$ : مرتبه ناهم محوری پویا

ناهم محوری ایستا : $n_d = 0$

ناهم محوری پویا : $n_d = 1,2,3,...$

$k$ : یک عدد صحیح مثبت $k = 0,1,2,...$

$\eta$ : مرتبه هارمونیک استاتور $\eta = 1,3,5,...$

در صورتیکه $n_d = 0$، $k = 1$ و $\eta = 1$ باشد، رابطه (5-12) بصورت زیر ساده می‌شود. این هارمونیک، مشخصه شیارگذاری موتور یا PSH(Principle Slot Harmonics) نام دارد.

$$f_{psh} \equiv f_h(1,0,s,1) = \left( \frac{1-s}{p} R \pm 1 \right) f_0 \Rightarrow \begin{cases} f_{PSH+} = \left( \frac{1-s}{p} R + 1 \right) f_0 \\ f_{PSH-} = \left( \frac{1-s}{p} R - 1 \right) f_0 \end{cases} \tag{13-5}$$

و در صورتیکه $n_d = 1$، $k = 0$ و $\eta = 1$ باشد رابطه (5-12) بصورت رابطه زیر ساده می‌شود.

$$f_{ecc0} \equiv f_h(0,1,s,1) = \left( 1 \pm \frac{1-s}{p} \right) f_0 \tag{14-5}$$

این هارمونیکها که در اطراف فرکانس اصلی تغذیه قرار دارند ناشی از ناهم‌محوری بین روتور و استاتور می‌باشند. هارمونیکهای دیگر ناشی از ناهم‌محوری بین روتور و استاتور از رابطه زیر و با اعمال $n_d = 1$، $k = 1$ و $\eta = 1$ قابل محاسبه می‌شوند. این هارمونیکها در اطراف فرکانس PSH هستند.

$$f_{ecc1} \equiv f_h(1,1,s,1) = \left(\frac{R \pm 1}{p}(1-s) \pm 1\right) f_0 \qquad (15\text{-}5)$$

لازم به ذکر است که طیف فرکانسی جریان استاتور در حالت ناهم محوری بین روتور و استاتور علاوه بر هارمونیکهای مشخص شده در روابط اخیر هارمونیکهای اضافی دیگری نیز دارد، اما این هارمونیکها نسبتا ضعیف و آشکارسازی آنها بسیار دشوار است.

به منظور محاسبه طیف فرکانسی جریان استاتور پس از شبیه سازی عملکرد موتور در حالت دائم، شکل موج جریان را در یک فایل MAT ذخیره کرده و با استفاده از جعبه ابزار SIMULINK نرم افزار MATLAB از یک ZERO ORDER HOLD عبور می‌دهیم. . شکل 5-22 بلوک دیاگرام این عمل را نشان می‌دهد.

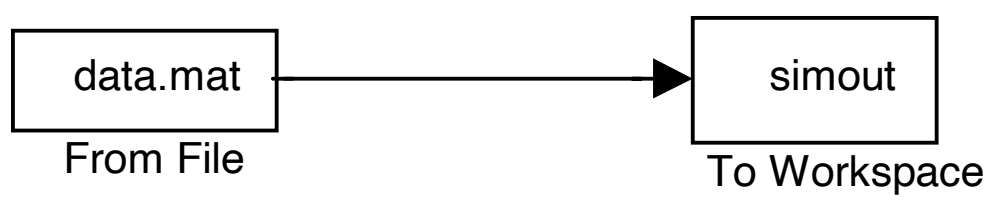

**شکل5-22: عبور سیگنال جریان از ZERO ORDER HOLD با استفاده از بلوکهای جعبه ابزار SIMULINK**

سپس با اجرای دستورات زیر در محیط MATLAB طیف فرکانسی جریان بدست می‌آید:

```
%fast fourier transformation using MATLAB commands
```

```
Tf1=fft([simout]);
Tf2=abs(Tf1);
Tf3=2*Tf2/length(Tf2);
Tf4=20*log(Tf3)/log(10);
plot((0:length(Tf4)-1)/1.5,Tf4,'black');
```

شکلهای 5-23 تا 5-26 تبدیل فوریه شکل موج جریان فاز اول استاتور را در حالت سالم و چند حالت مختلف ناهم محوری ایستا نشان می‌دهند. این نتایج بر مبنای عملکرد موتور در شرایط ارائه شده در جدول 5-1 است.

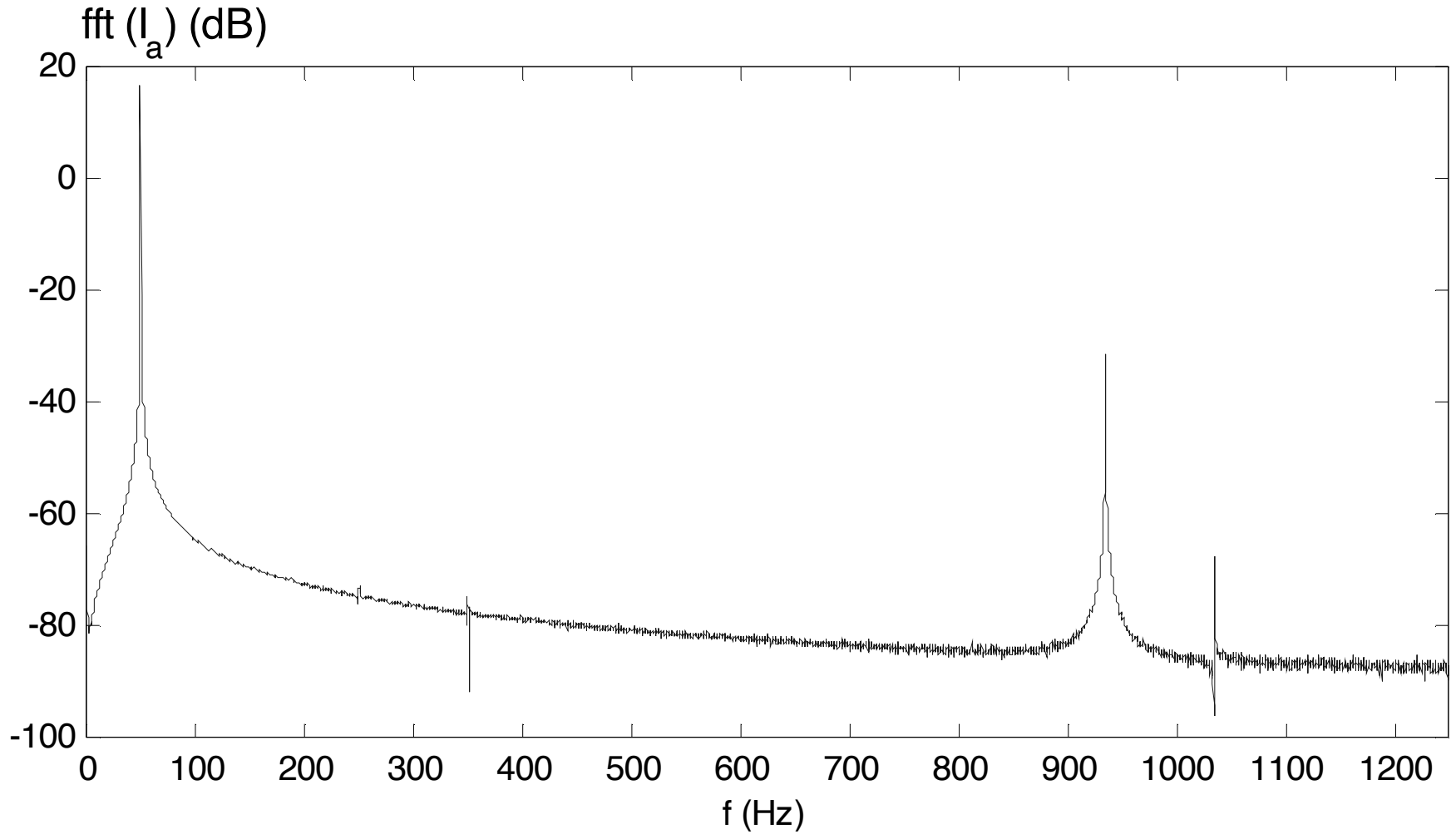

شکل **5-23**: طیف فرکانسی جریان استاتور در حالت سالم

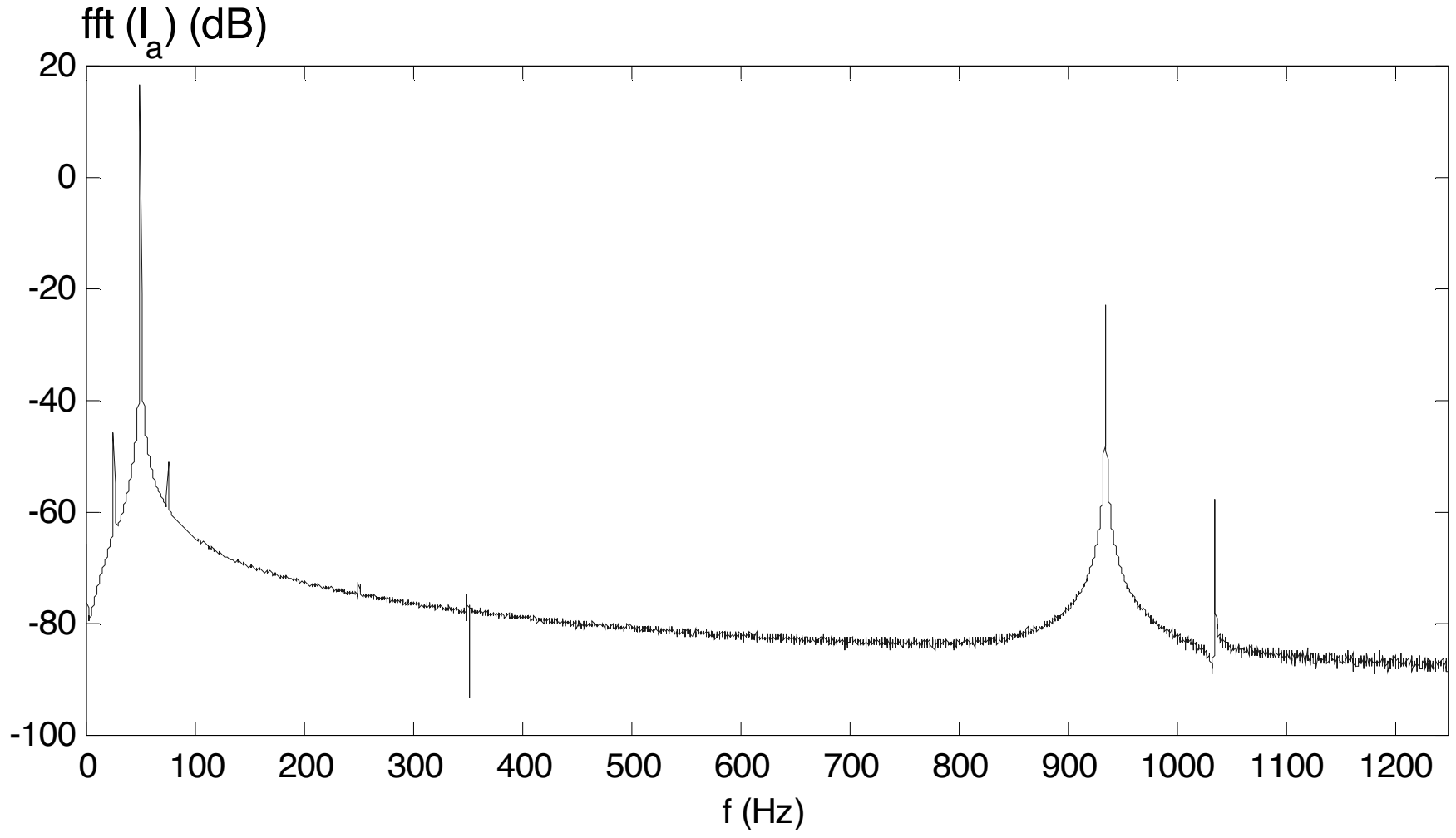

شکل 5-24 : طیف فرکانسی جریان استاتور در حالت 5٪ ناهم محوری ایستا

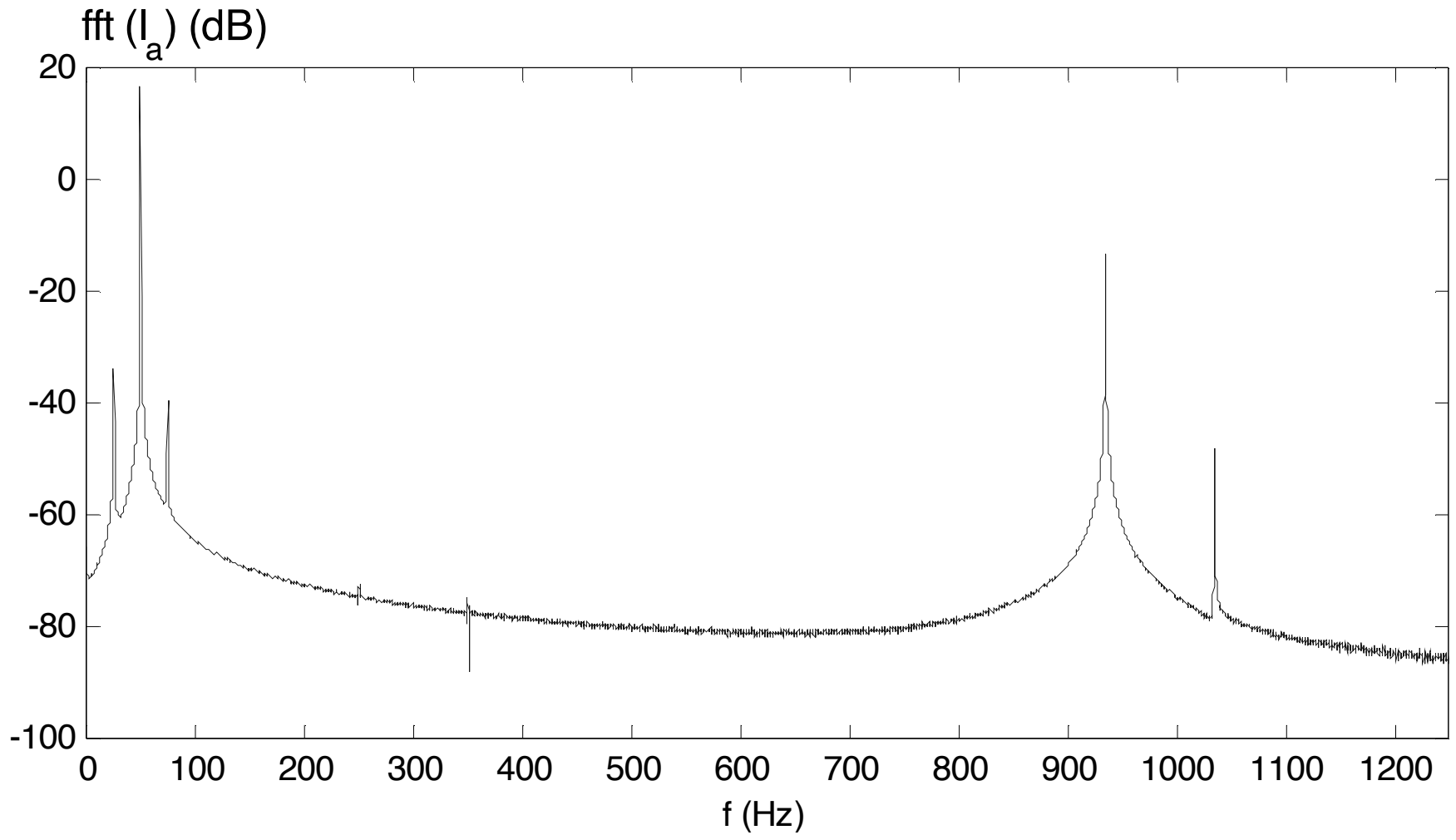

شکل 4-25 : طیف فرکانسی جریان استاتور در حالت 15٪ ناهم محوری ایستا

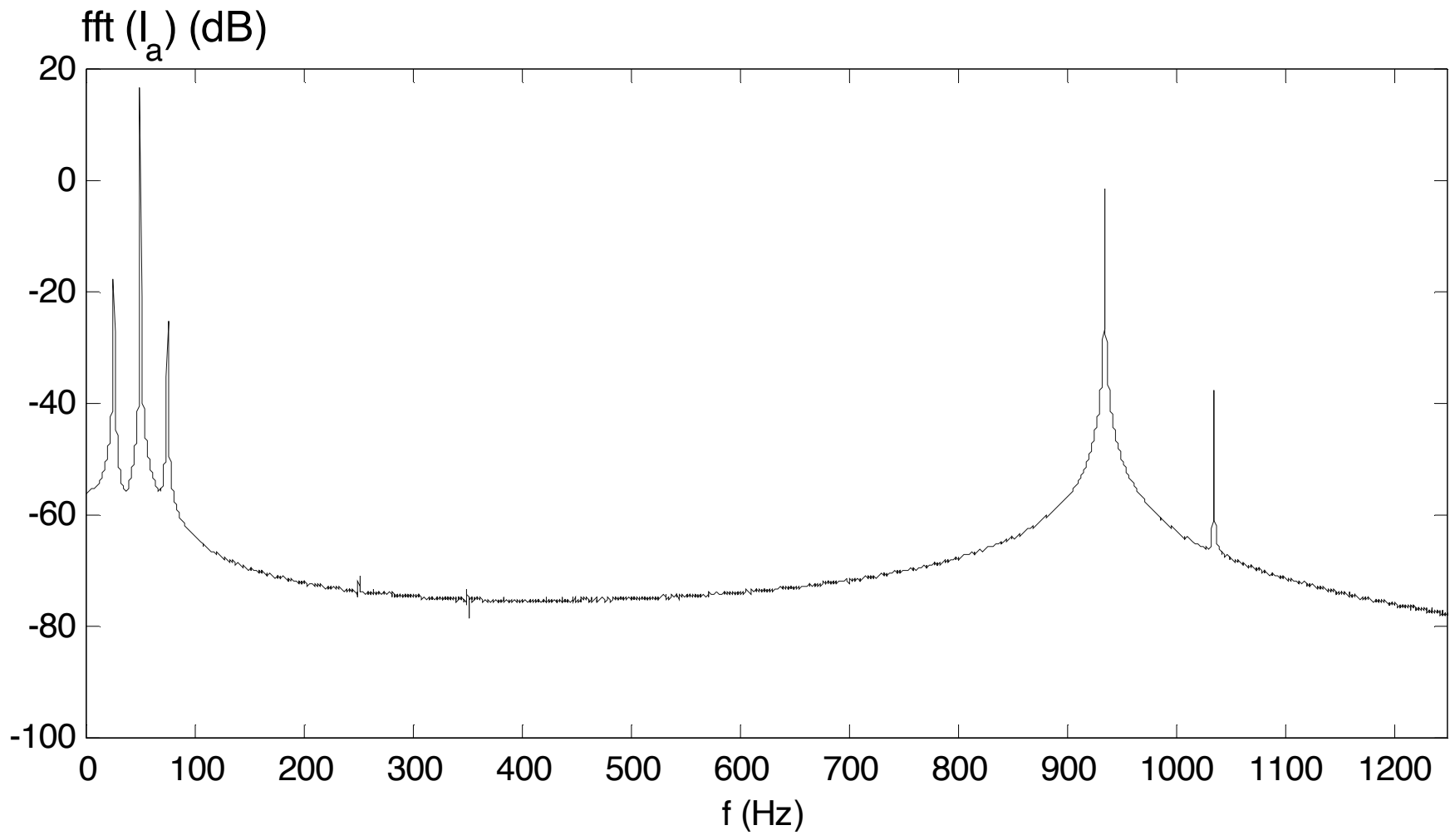

شکل 5-26 : طیف فرکانسی جریان استاتور در حالت 35٪ ناهم محوری ایستا

این نکته قابل ذکر است که نتایج شبیه سازی نشان می دهدکه در حالت دائمی مستقل از وضعیت ناهم محوری میان روتور و استاتور لغزش موتور حدود $0.015$ است. لذا فرکانس هارمونیکهای نشان داده شده در شکلهای فوق با استفاده لز جاگذاری این مقدار در رابطه (5-14) بصورت زیر محاسبه شده است.

$$f_{ecc0-} = 25.375Hz, \quad f_{ecc0+} = 74.625Hz$$

همچنین با توجه به رابطه (5-15) انتظار می رود که در حالتهای مختلف ناهم محوری طیف فرکانسی اطراف هارمونیک شیارگذاری تغییر کند. با توجه به مقدار لغزش محاسبه شده و رابطه (5-13) فرکانسهای هارمونیکهای شیار گذاری بصورت زیر قابل محاسبه هستند:

$$f_{psh-} = 935Hz, \quad f_{psh+} = 1035Hz$$

همانطور که در این شکلها دیده می‌شود هارمونیکهای متناظر با $f_{ecc0}$ در حالت ناهم محوری ایستا تقویت می‌شوند. با استفاده از نتایج شبیه سازی می‌توان نشان داد که در حالت ناهم محوری ایستا همواره هارمونیک با فرکانس کمتر از فرکانس اصلی دامنه کوچکتری نسبت به هارمونیک با فرکانس بیشتر از فرکانس اصلی دارد. با توجه به رابطه (5-12) در حالت ناهم محوری ایستا $(n_d = 0)$ انتظار می‌رود که بدون ایجاد هارمونیکهای اضافی در اطراف هارمونیک شیار گذاری، این هارمونیکها تقویت شوند. شکل 5-24 تا 5-26 صحت این ادعا را نشان می‌دهد. منحنیهای شکلهای 5-27 دامنه هارمونیکهای ناشی از ناهم محوری ایستا را بر حسب درجه های مختلف ناهم‌محوری ایستا نشان می‌دهند.

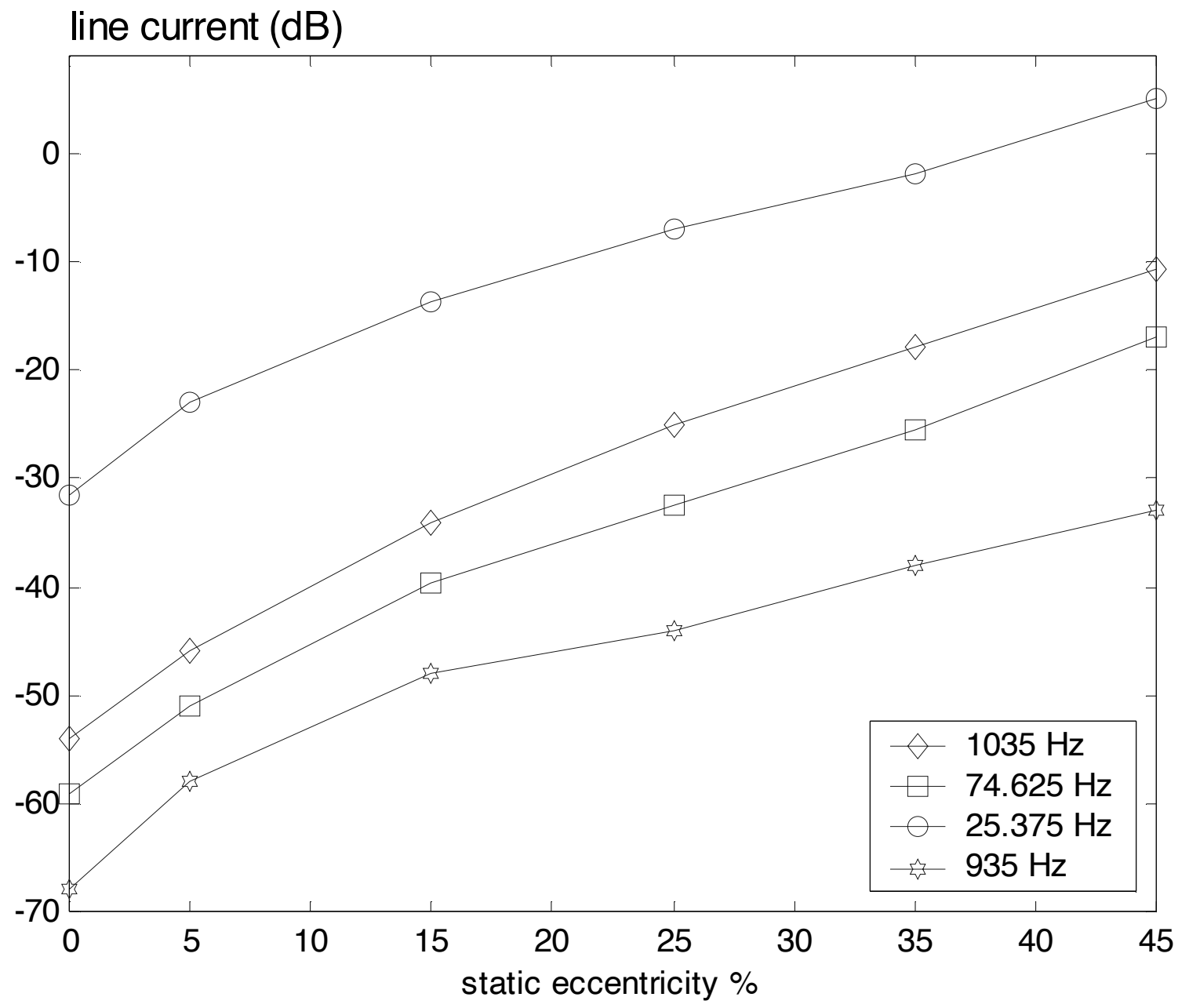

شکل5-27 : تغییرات دامنه هارمونیکهای ناهم محوری اطراف فرکانس اصلی نسبت به تغییرات درجه ناهم محوری ایستا

با توجه به رابطه (5-15) انتظار می رود تا درحالت ناهم‌محوری پویا و مرکب هارمونیکهای اضافی در اطراف هارمونیک شیار گذاری ایجاد شود. با توجه به علامتهای ± این رابطه این هارمونیکها می‌توانند دارای فرکانسهای زیر باشند.

$$f_{1ecc1} = f_h(1,1,s,1) = \left(\frac{R-1}{p}(1-s)-1\right)f_0 = 910.37Hz \quad (16\text{-}5)$$

$$f_{2ecc1} = \left(\frac{R+1}{p}(1-s)-1\right)f_0 = 959.62Hz \quad (17\text{-}5)$$

$$f_{ecc1} = \left(\frac{R-1}{p}(1-s)+1\right)f_0 = 1010.37Hz \quad (18\text{-}5)$$

$$f_{ecc1} = \left(\frac{R+1}{p}(1-s)+1\right)f_0 = 1059.62Hz \quad (19\text{-}5)$$

این هارمونیکها در شکلهای 5-28 و 5-30 که طیف فرکانسی جریان فاز اول استاتور در حالتهای مختلف ناهم محوری پویا ست به چشم‌می‌خورند . در حالتهای مختلف ناهم‌محوری با توجه به شرایط طراحی و عملکرد موتور این هارمونیکها تقویت و یا تضعیف می‌شوند. به هر حال وجود حتی یکی از این هارمونیکها دلیل خوبی برای احتمال وجود انواع مختلف ناهم محوری میان روتور و استاتور می‌باشد.

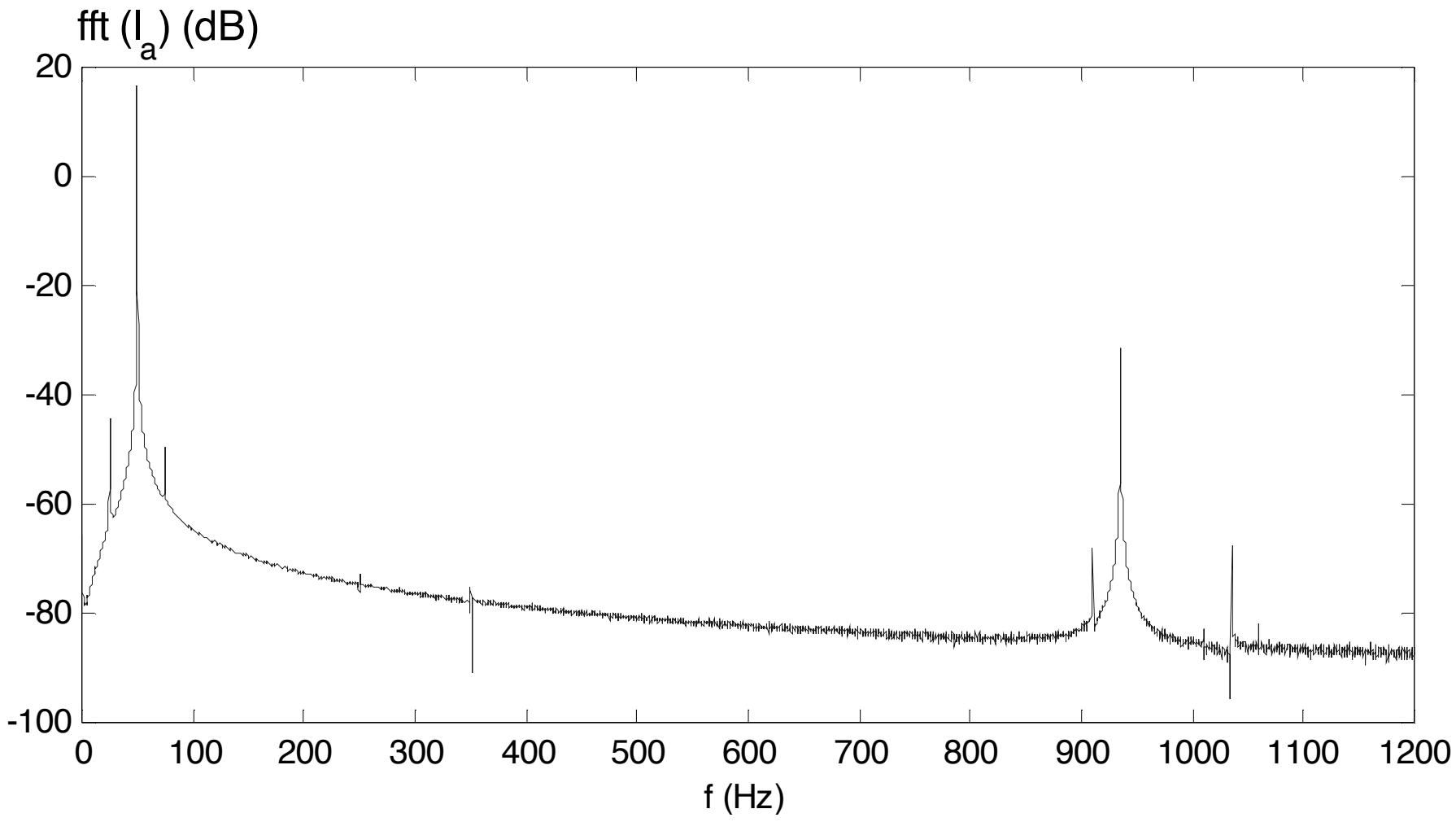

شکل **5-28** : طیف فرکانسی جریان استاتور در حالت **5**٪ ناهم محوری پویا

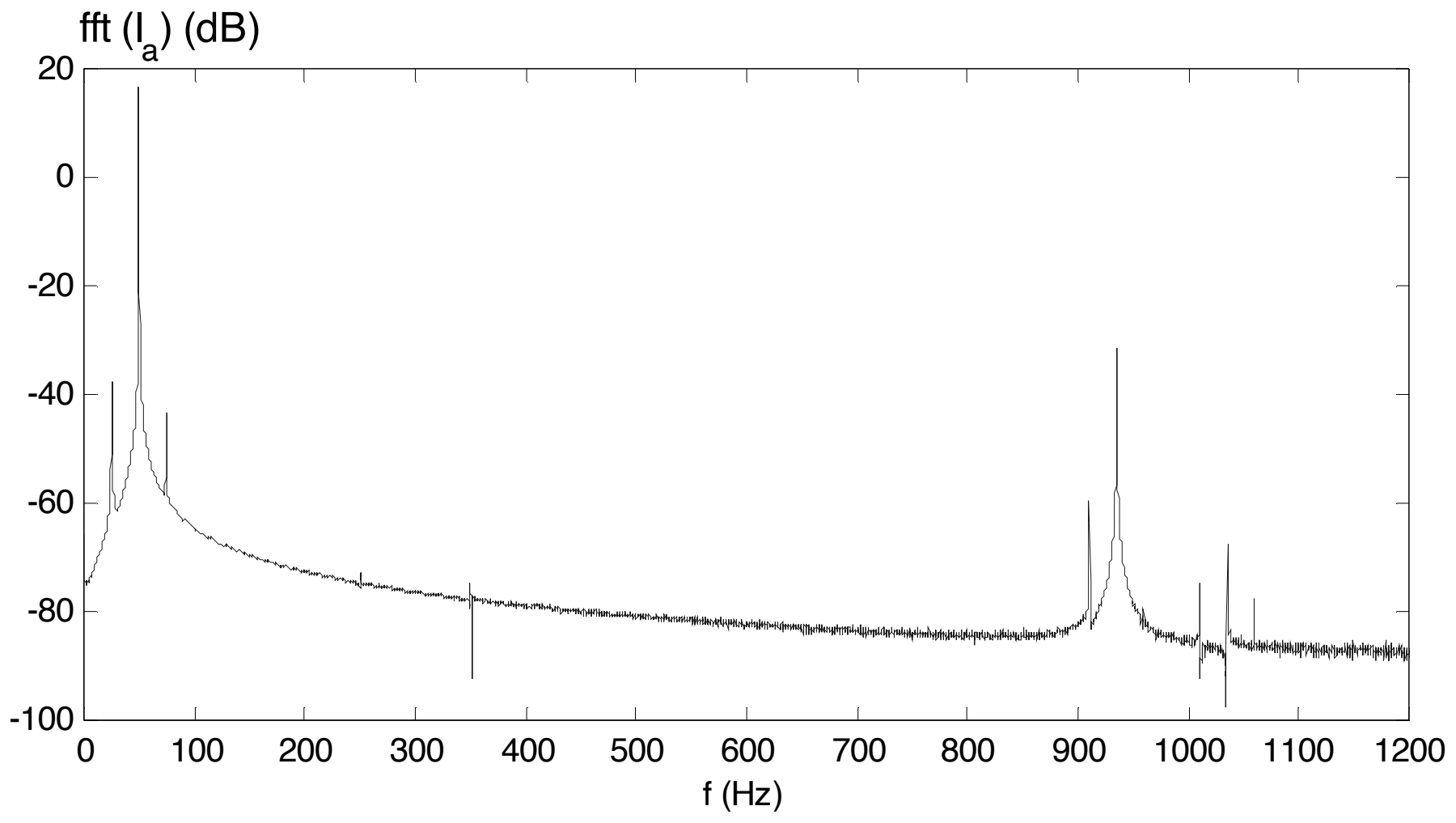

شكل 5-29 : طيف فركانسی جريان استاتور در حالت 15٪ ناهم محوری پويا

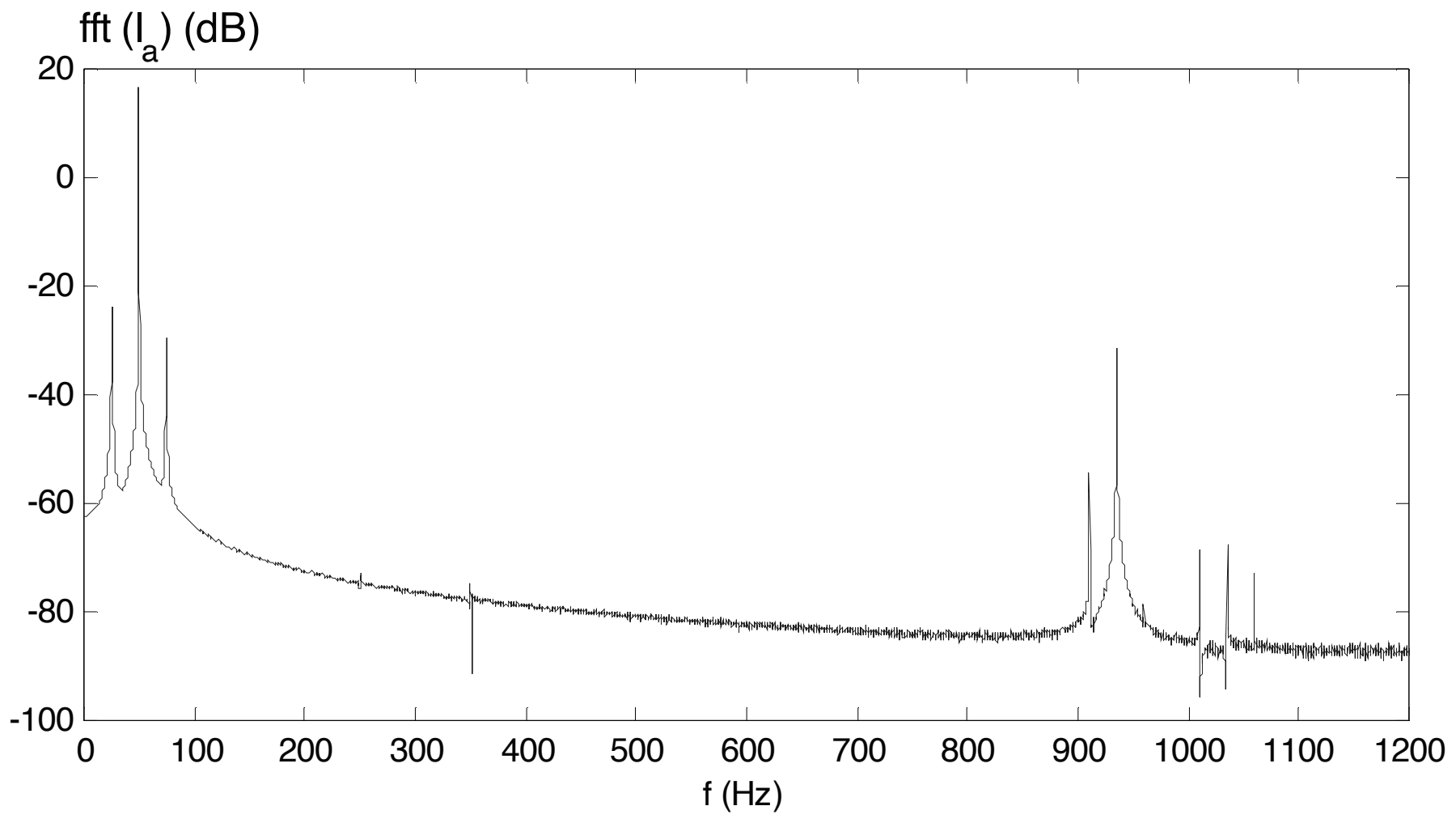

شكل 5-30 : طيف فركانسی جريان استاتور در حالت 25٪ ناهم محوری پويا

منحنیهای شکلهای 5-31 دامنه هارمونیکهای ناشی از ناهم محوری پویا را بر حسب درجه های مختلف ناهم‌محوری پویا نشان می‌دهند.

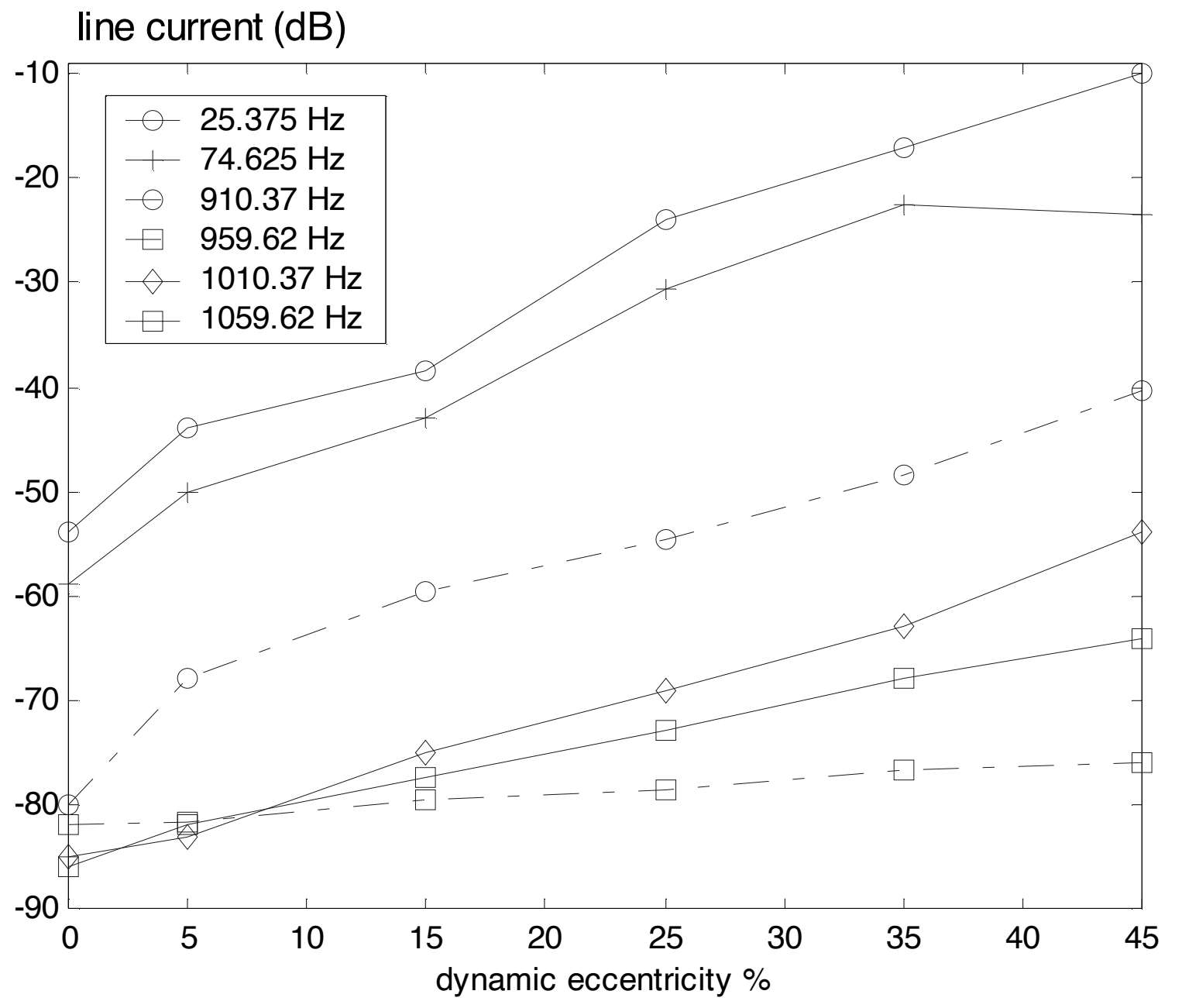

شكل**5-31** : تغييرات دامنه هارمونيكهای ناهم محوری اطراف فركانس اصلی نسبت به تغييرات درجه ناهم محوری ايستا

شكلهای 5-32 و 5-33 تبديل فوريه شكل موج جريان فاز اول استاتور را در چند حالت مختلف ناهم‌محوری مركب نشان می‌دهند. مطابق اين شكلها در حالت ناهم محوری مركب كليه هلرمونكيهای ناشی از ناهم محوری ايستا و پويا وجود دارند.

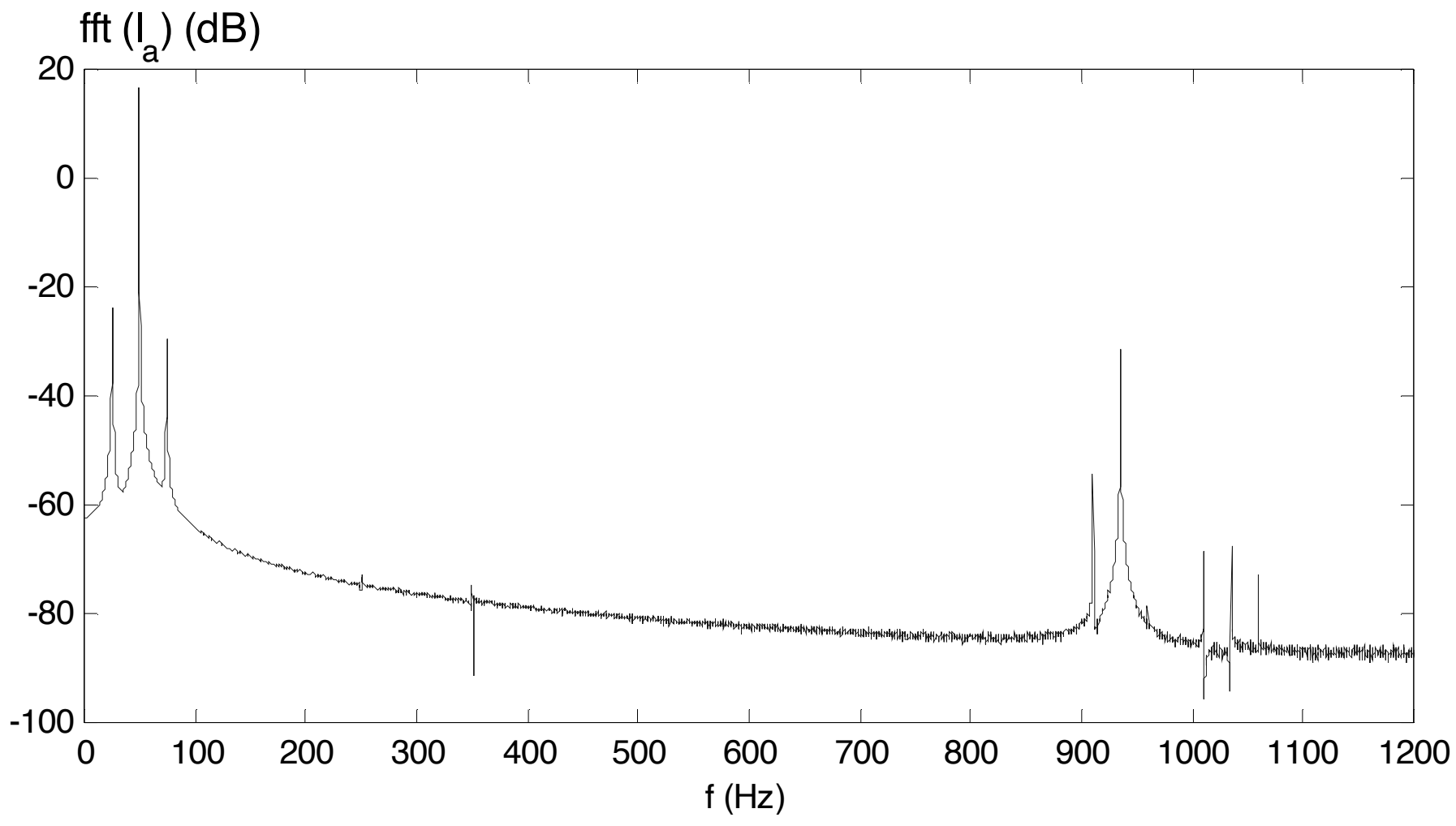

شكل 5-32 : طيف فركانسی جريان استاتور در حالت 5٪ ناهم محوری ايستا و 5٪ ناهم محوری پويا

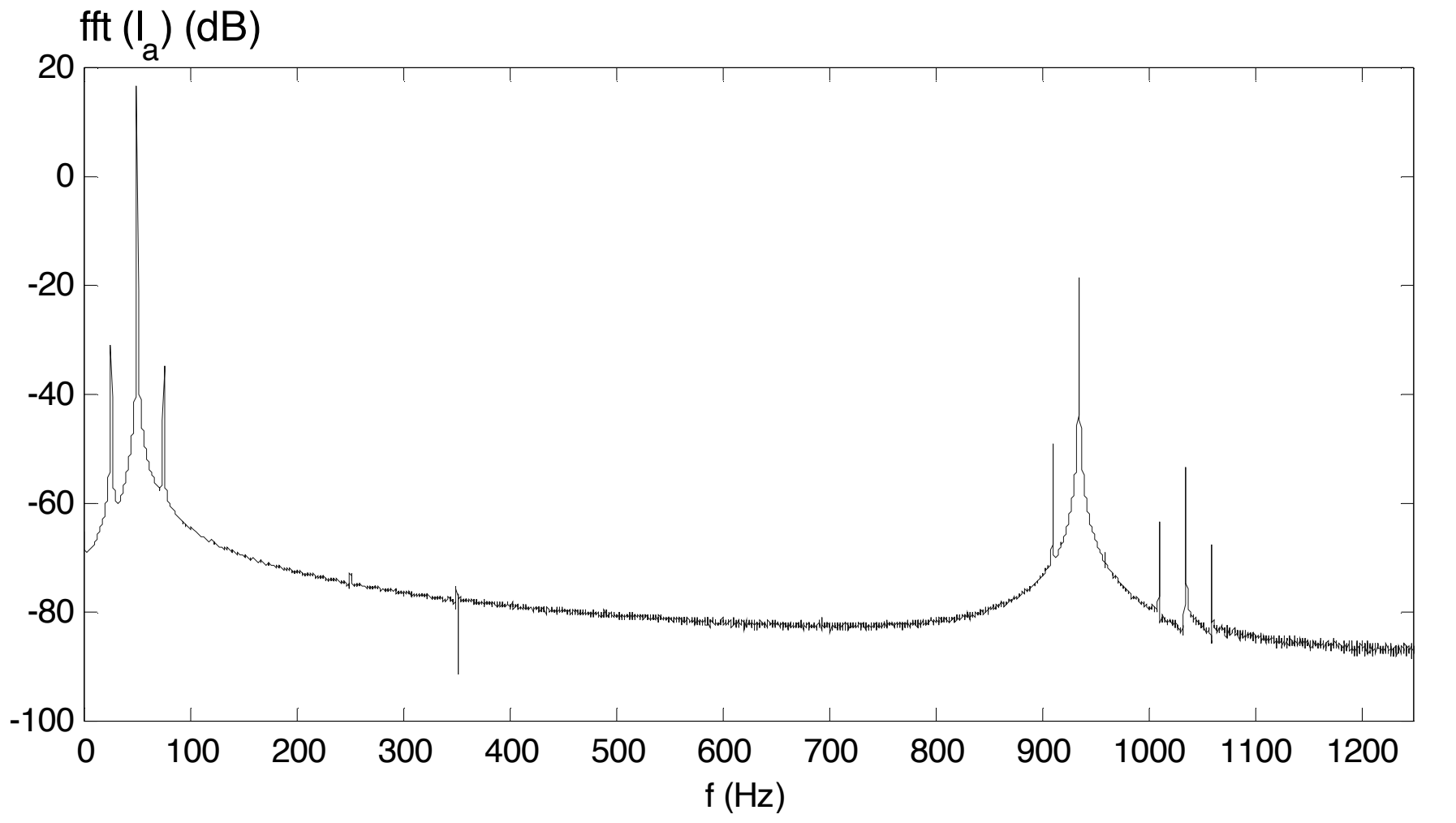

شكل 5-33 : طيف فركانسی جريان استاتور در حالت 5٪ ناهم محوری ايستا و 15٪ ناهم محوری پويا

# فصل ششم

# تحلیل کامپیوتری موتور القایی سه فاز قفس سنجابی در حالت شکستگی میله های روتور

شکستگی میله های روتور در موتورهای القایی قفس سنجابی موجب ایجاد ارتعاشات اضافی، نویز و جرقه در هنگام راه اندازی موتور می شود. این علائم موقعی بروز می کند که خطا گسترش یافته و باعث شکستگی چندین میله شود. بنابراین تشخیص و جلوگیری از توسعه این خطا موجب افزایش کارایی موتور می‌شود.

## 6-1- معادلات توصیف موتور در حالت شکستگی یک یا چند میله قفس روتور

شکل 6-1 مدل گسترده شده قفس روتور را در حالت سالم نشان می دهد. همانطور که در این شکل نشان داده شده هر دو میله مجاور تشکیل یک حلقه می دهد که متناظر با آن می‌توان یک معادله جریان KCL نوشت. بنابراین از آنجا که به ازاء هر میله یک حلقه وجود دارد می توان به تعداد حلقه ها معادله KCL مستقل خطی نوشت.

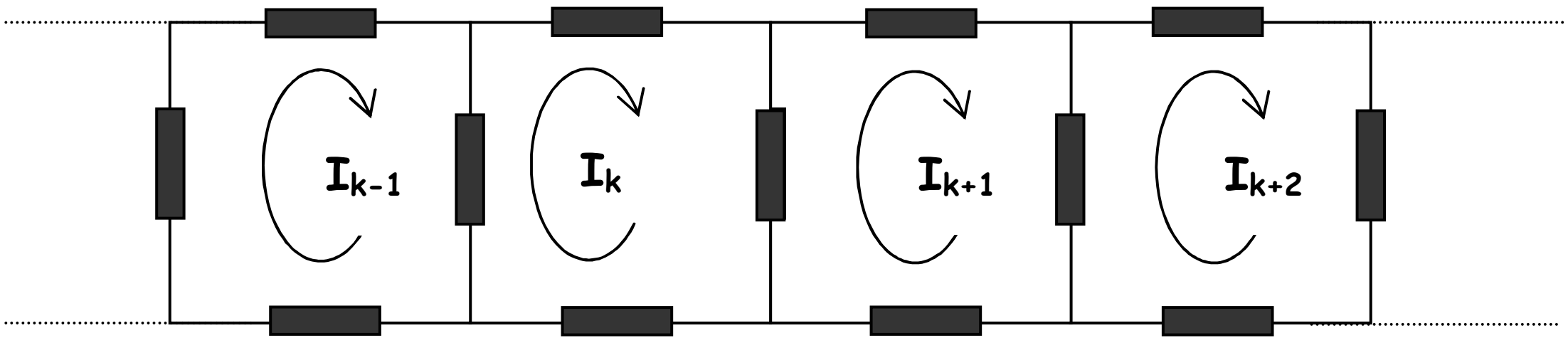

**شکل 6-1 : شکل گسترده قفس روتور در حالت سالم**

شکل 6-2 مدل گسترده شده قفس روتور را در حالت شکستگی یک میله نشان می دهد. همانطور که این شکل نشان می دهد در حالت شکستگی میله جریان دو حلقه أی که میله مورد نظر در آنها وجود دارد برابر است. بنابراین تعداد معادلات مستقل خطی نسبت به حالت سالم یکی کم می شود.

یک راه دیگر برای لحاظ شکستگی میله های روتور در مدل موتور تغییر ماتریس مقاومت روتور است. با این کار می توان حتی ترک خوردگی میله های روتور را در نظر گرفت. با افزایش مقاومت میله شکسته این امر

تحقق می پذیرد. برای مثال اگر میله بین حلقه K ام و K+1 ام روتور شکسته باشد ماتریس مقاومت روتور بصورت زیر تغییر خواهد کرد.

$$R_r(K,K) = r_{brok} \qquad (2\text{-}6)$$

$$R_r(K+1,K+1) = r_{brok} \qquad (3\text{-}6)$$

$$R_r(K+1,K) = R_r(K,K+1) = -r_{brok} \qquad (4\text{-}6)$$

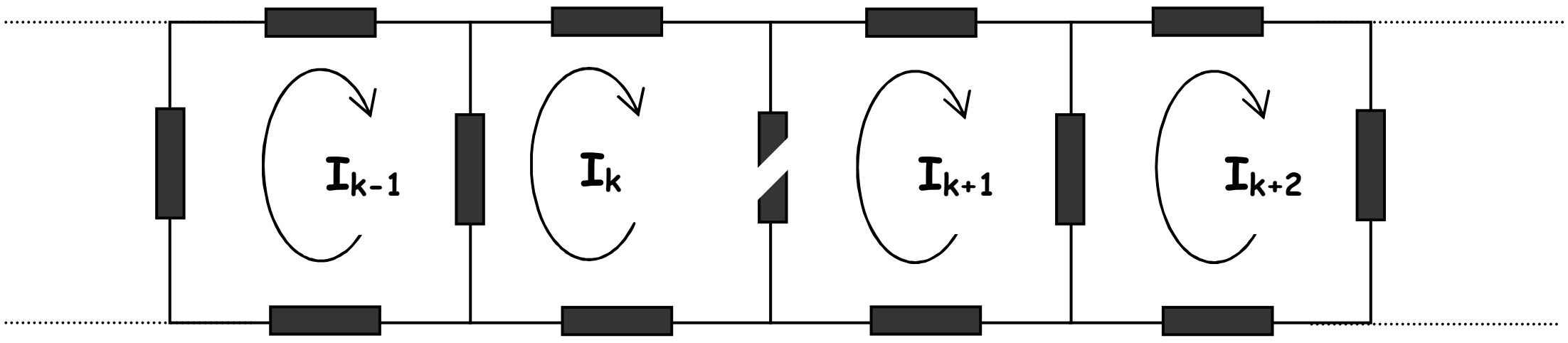

شکل6-2 : شکل گسترده قفس روتور در حالت شکستگی یک میله

در گزارش حاضر مقدار $r_{brok}$ هزار برابر مقاومت میله سالم در نظر گرفته شده است.

## 6-2- تعریف تابع سیم پیچی در حالت شکستگی یک یا چند میله قفس روتور

در این قسمت فرض می کنیم در حالت شکستگی یک یا چند میله قفس روتور هیچکدام از شرایط ناهم محوری بین روتور و استاتور وجود نداشته باشد. به عبارت دیگر فاصله هوایی بصورت یکنواخت در نظر گرفته می شود. بنابراین با توجه به آنچه در فصل دوم گفته شد تعریف تابع سیم پیچی در این حالت بصورت زیر ساده می شود.

$$N_i(\varphi) = n_i(\varphi) - \langle n_i \rangle \qquad (5\text{-}6)$$

$$\langle n_i \rangle = \frac{1}{2\pi}\int_0^{2\pi} n_i(\varphi)d\varphi \qquad (6\text{-}6)$$

در این بخش نیز از توابع دور معرفی شده در فصل سوم برای فازهای استاتور و حلقه های روتور استفاده می‌شود. قدم بعدی محاسبه ماتریس اندوکتانسهای مغناطیسی موتور است.

## 6-3- محاسبه ماتریس اندوکتانس موتور در حالت شکستگی میله های قفس روتور

همانطور که گفته شد برای تحلیل رفتار موتور در حالت شکستگی میله های قفس روتور بایستی دستگاه معادلاتی شامل دستگاه معادلات توصیف موتور در حالت سالم و معادلاتی مانند معادله (6-1) حل گردد. بنابراین ماتریس اندوکتانس موتور در حالت شکستگی میله های قفس روتور مانند حالت سالم است. این مساله به راحتی از روی تعریف تابع دور و آنچه در فصل دوم ذکر شد استنباط می شود.

با توجه به آنچه در فصل دوم گفته شد در حالت کلی اندوکتانس متقابل بین دو سیم پیچی $x$ و $y$ با کمک رابطه زیر قابل محاسبه است.

$$L_{yx} = l\int_0^{2\pi} P(\varphi)N_x(\varphi)n_y(\varphi)d\varphi \qquad (7\text{-}6)$$

از آنجا که در حالت شکستگی میله های قفس روتور هدایت مغناطیسی فاصله هوایی یکنواخت است بنابراین با توجه به رابطه (3-36) رابطه فوق بصورت زیر ساده می شود.

$$L_{yx} = \frac{\mu_0 r_0 l}{g_0}\int_0^{2\pi} N_x(\varphi)n_y(\varphi)d\varphi \qquad (8\text{-}6)$$

برای محاسبه کلیه اندوکتانسهای موتورر از روابط(6-8)و (6-6) و توابع دور معرفی شده در فصل دوم استفاده می‌کنیم. همچنین در این قسمت جهت مدلسازی شکستگی میله های روتور فرض شده است که مقاومت مسله شکسته هزار برابر مقاومت میله سالم است. اندوکتانس نشتی میله شکسته با میله سالم مساوی فرض شده است. یکی از معایب این روش کاهش پارامتر rcond ماتریس اندوکتانس روتور در نرم افزار MATLAB است. هر چقدر این پارامتر به یک نزدیکتر باشد محاسبات دقیقتر و هرچقدر به صفر نزدیک باشد خطای مخاسبات بیشتر است. با اینحال در این فصل نشان داده می شود که دقت این روش کافی است و نتایج شبیه سازی این روش با نتایج شبیه سازی حاصل از روش اول یکسان است. از مزایای این روش نسبت به روش اول سهولت مدلسازی شکستگی میله های مختلف موتور می باشد. در حالیکه در روش دوم این کار به آسانی ممکن نیست.

*الف: ماتریس اندوکتانس استاتور*

اندوکتانس مغناطیسی خودی فازهای A,B,C استاتور همچنین اندوکتانسهای متقابل این فازها بصورت زیر بدست می آید.

$$L_{AA} = L_{BB} = L_{CC} = L_{ss} = \frac{14\pi}{9}\frac{\mu_0 r_0 l_0}{g_0}N^2 + L_{\ell s} \quad (9\text{-}6)$$

$$L_{AB} = L_{BC} = L_{CA} = L_{BA} = L_{CB} = L_{AC} = L_{ms} - \frac{2\pi}{3}\frac{\mu_0 r_0 l_0}{g_0}N^2 \quad (10\text{-}6)$$

در این رابطه $N$ تعداد دور سیم پیچی یک فاز استاتور و $L_{\ell s}$ اندوکتانس پراکندگی یک فاز استاتور می‌باشد. بنابراین ماتریس اندوکتانس استاتور بصورت زیر می باشد:

$$[L_s] = \begin{bmatrix} L_{ss} & L_{ms} & L_{ms} \\ L_{ms} & L_{ss} & L_{ms} \\ L_{ms} & L_{ms} & L_{ss} \end{bmatrix} \quad (11\text{-}6)$$

از آنجا که ماتریس $L_s$ برای محاسبه $[\psi_s]$ از رابطه زیر استفاده می شود.

$$\begin{bmatrix}\psi_A\\ \psi_B\\ \psi_C\end{bmatrix}=\begin{bmatrix}L_{ss} & L_{ms} & L_{ms}\\ L_{ms} & L_{ss} & L_{ms}\\ L_{ms} & L_{ms} & L_{ss}\end{bmatrix}\begin{bmatrix}I_A\\ I_B\\ I_C\end{bmatrix}+\left[L_{sr}\right]_{3\times n}\left[I_r\right]_{n\times 1} \tag{12-6}$$

و اگر جریان فازها سه فاز متقارن باشد :

$$I_A+I_B+I_C=0 \tag{13-6}$$

می توان ماتریس اندوکتانس استاتور را بصورت زیر ساده نمود:

$$\left[L_s\right]=\begin{bmatrix}L_{ss}-L_{ms} & 0 & 0\\ 0 & L_{ss}-L_{ms} & 0\\ 0 & 0 & L_{ss}-L_{ms}\end{bmatrix} \tag{14-6}$$

*ب: ماتریس اندوکتانس روتور*

برای یک موتور القایی با $n$ میله روتور ماتریس اندوکتانس روتور بصورت زیر می‌باشد.

$$\left[L_r\right]=\left[B_{ij}\right]_{n\times n}=\begin{cases}\left(\dfrac{2\pi}{n}-\dfrac{2\pi}{n^2}-\dfrac{\gamma}{3}\right)\dfrac{\mu_0 r_0 l_0}{g_0}+2(l_b+l_e) & i=j\\ \left(-\dfrac{2\pi}{n^2}+\dfrac{\gamma}{6}\right)\dfrac{\mu_0 r_0 l_0}{g_0}-l_b & i=j\pm 1\\ 0 & else\end{cases} \tag{15-6}$$

در این رابطه $l_b$ و $l_e$ به ترتیب اندوکتانس پراکندگی یک میله و یک قسمت از حلقه انتهایی روتور است. با نوشتن معادلات کات‌ست در یکی از حلقه های انتهایی روتور نتیجه می گیریم که مجموع جریانهای میله های روتور صفر است. بنابراین مانند آنچه در مورد ماتریس اندوکتانس استاتور نتیجه گرفته شد ماتریس اندوکتانس روتور نیز بصورت زیر ساده خواهد شد:

$$[L_r] = [B'_{ij}]_{n\times n} = \begin{cases} \left(\frac{2\pi}{n} - \frac{\gamma}{3}\right)\frac{\mu_0 r_0 l_0}{g_0} + 2(l_b + l_e) & i = j \\ \frac{\gamma}{6}\frac{\mu_0 r_0 l_0}{g_0} - l_b & i = j \pm 1 \\ 0 & else \end{cases} \qquad (16\text{-}6)$$

ج : *ماتریس اندوکتانس متقابل استاتور و روتور*

با استفاده از قانون انتگرالگیری جزء به جزء برای محاسبه انتگرال رابطه (6-8) می توان اندوکتانس متقابل میان فاز A استاتور و حلقه اول روتور را برای یک موتور القایی سه فاز چهار قطب بصورت رابطه زیرتنها بر حسب تابع انتگرال اول و دوم تابع دور روتور محاسبه کرد.

$$\begin{aligned} M(\theta) = \frac{\mu_0 r_0 l_0}{g_0}\frac{12N}{\pi}[&-m_r\left(\tfrac{\pi}{6} - \theta\right) + m_r(-\theta) + m_r\left(\tfrac{2\pi}{3} - \theta\right) - m_r\left(\tfrac{\pi}{2} - \theta\right) \\ &- m_r\left(\tfrac{7\pi}{6} - \theta\right) + m_r(\pi - \theta) - m_r\left(\tfrac{5\pi}{3} - \theta\right) + m_r\left(\tfrac{3\pi}{2} - \theta\right)] - P_0 N k_r(2\pi - \theta) \end{aligned} \qquad (17\text{-}6)$$

در این رابطه $K_r$ و $m_r$ به ترتیب انتگرال اول و دوم تابع دور حلقه اول روتور می باشند که از روابط تحلیلی زیر قابل محاسبه هستند.

$$n_r(\varphi)=\begin{cases}\frac{1}{\gamma}\varphi & 0<\varphi<\gamma \\ 1 & \gamma<\varphi<\frac{2\pi}{n} \\ 1-\frac{1}{\gamma}\left(\varphi-\frac{2\pi}{n}\right) & \frac{2\pi}{n}<\varphi<\frac{2\pi}{n}+\gamma \\ 0 & \frac{2\pi}{n}+\gamma<\varphi<2\pi\end{cases} \qquad (18\text{-}6)$$

$$k_r(\varphi)=\int_0^{\varphi} n_r(\varphi')d\varphi' \qquad (19\text{-}6)$$

$$m_r(\varphi)=\int_0^{\varphi} k_r(\varphi')d\varphi' \qquad (20\text{-}6)$$

با استفاده از این رابطه و قانون جابجایی می توان اندوکتانس متقابل بقیه فازهای استاتور و حلقه های روتور را بدست آورد. در اینصورت ماتریس اندوکتانس متقابل استاتور و روتور بصورت زیر می باشد:

$$\left[L_{sr}\right]=\left[C_{ij}\right]_{3\times n}=\left[M\left(\theta-\frac{(i-1)2\pi}{3}+\frac{(j-1)2\pi}{n}\right)\right]_{3\times n} \qquad (35\text{-}6)$$

## 6-4- تحلیل کامپیوتری موتور با روش افزایش مقاومت میله های شکسته

به منظور حل معادلات توصیف موتور در حالت میله شکسته از روش رانگ کوتای مرتبه 4 و 5 نرم افزار MATLAB استفاده شده است. همچنین موارد اساسی زیر به منظور ساده سازی فرض شده اند:

1. ماده مغناطیسی در ناحیه خطی است و موتور در شرایط اشباع قرار نگرفته است.
2. توزیع جریان هر فاز بصورت یکنواخت در نظر گرفته شده است.
3. موتور توسط یک منبع سینوسی سه فاز متقارن تغذیه می شود.
4. بار مکانیکی ثابت در نظر گرفته شده است.

در ادامه نتایج حاصل از شبیه سازی کامپیوتری موتور در حالت سالم و حالتهای مختلف شکستگی میله‌های روتور بررسی، و با یکدیگر مقایسه می‌شوند. کلیه این شکل موجها برای شرایط تغذیه و بارمکانیکی مانند جدول 4-1 می‌باشند.

شکلهای 6-3 تا 6-7 به ترتیب منحنیهای سرعت مکانیکی روتور، جریان فاز اول، دوم و سوم استاتور و گشتاور تولیدی موتور را تا رسیدن به حالت دائمی برای موتور با یک میله شکسته نشان می دهند. شکلهای 6-8 تا 6-12، 6-13 تا 6-17 به ترتیب این منحنیها را برای 2 میله شکسته و 4میله شکسته نشان می‌دهند.

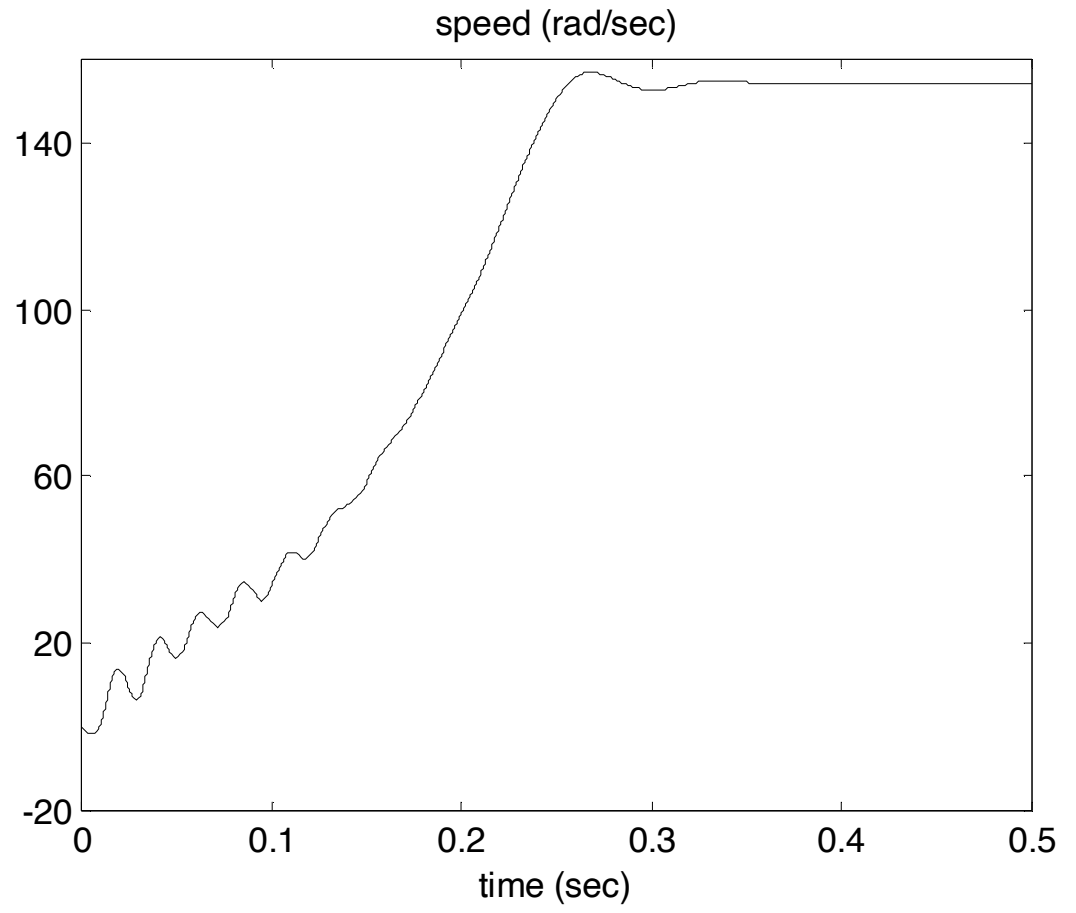

شکل**6-3**: منحنی تغییرات سرعت موتور در حالت **1** میله شکسته

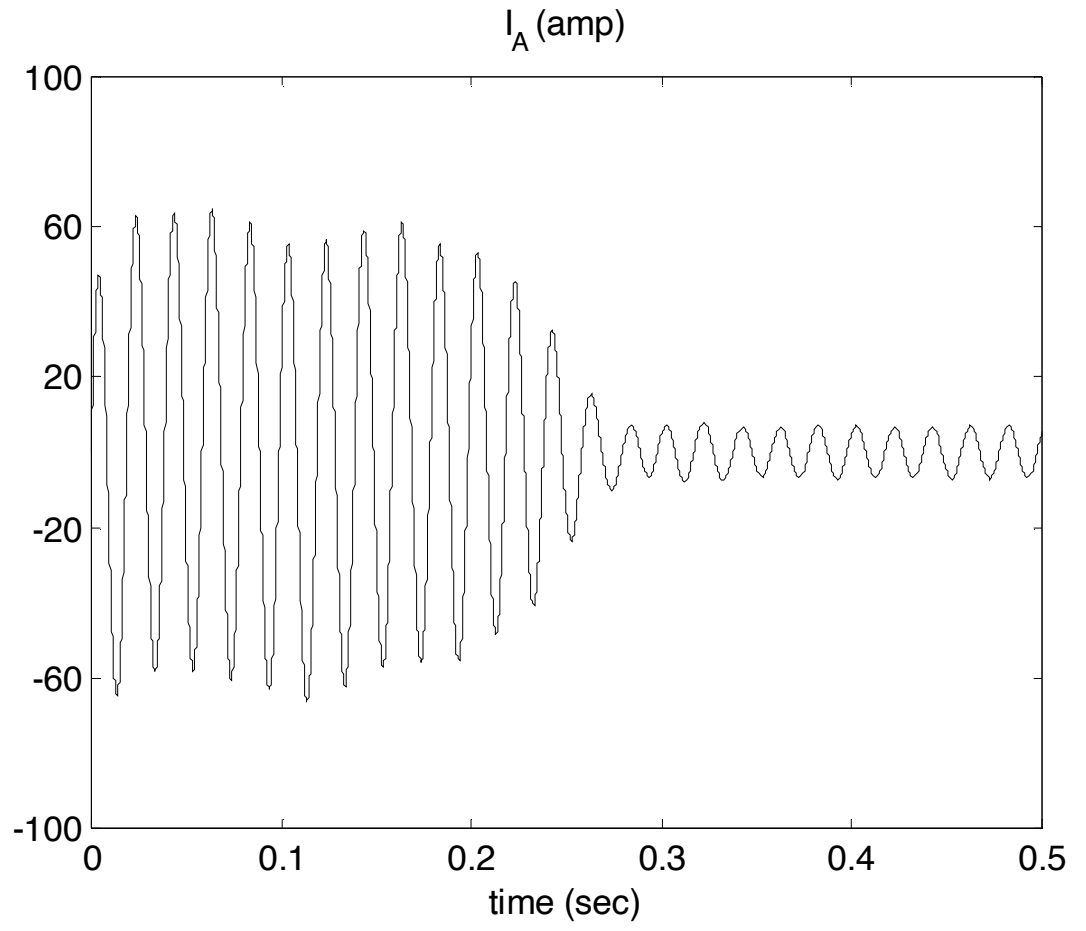

شکل**6-4**: منحنی تغییرات جریان فاز **A** در حالت **1** میله شکسته

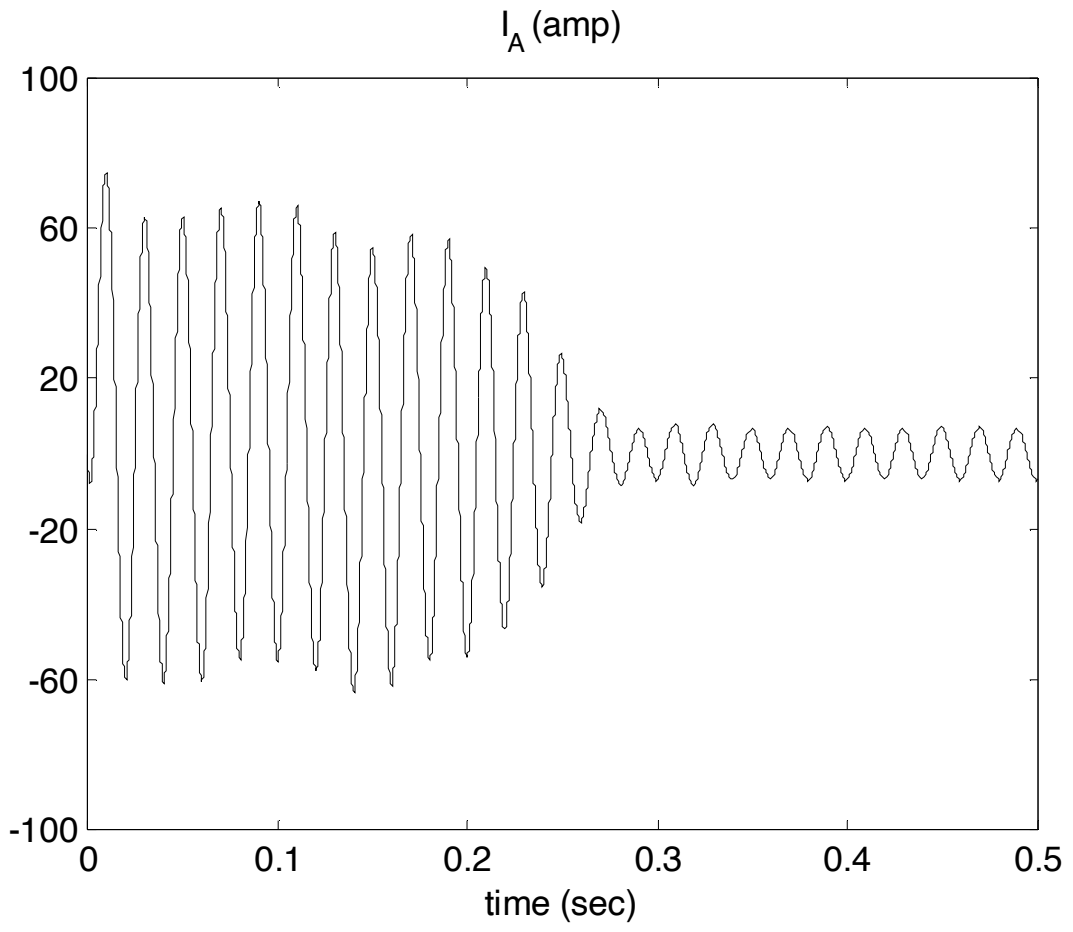

شکل**6-5**: منحنی تغییرات جریان فاز **B** در حالت **1** میله شکسته

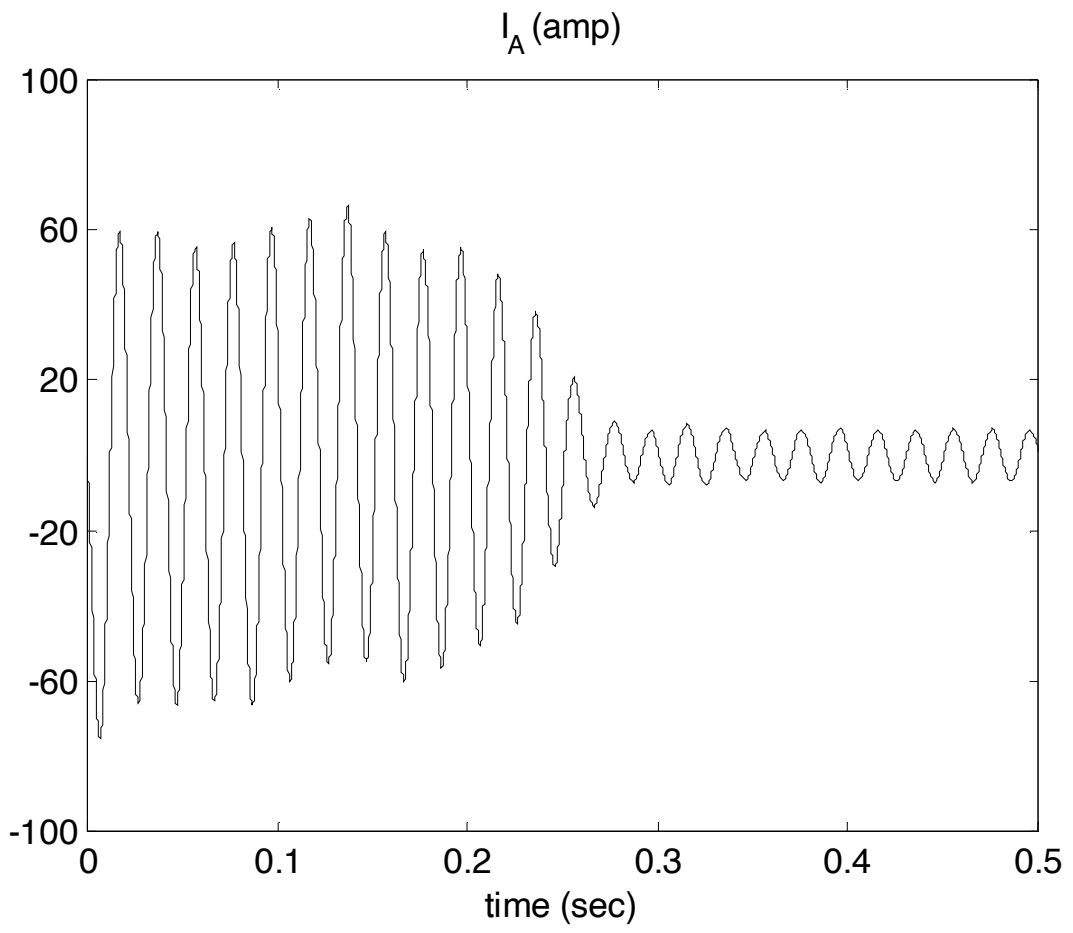

شکل**6-6**: منحنی تغییرات جریان فاز **C** در حالت **1** میله شکسته

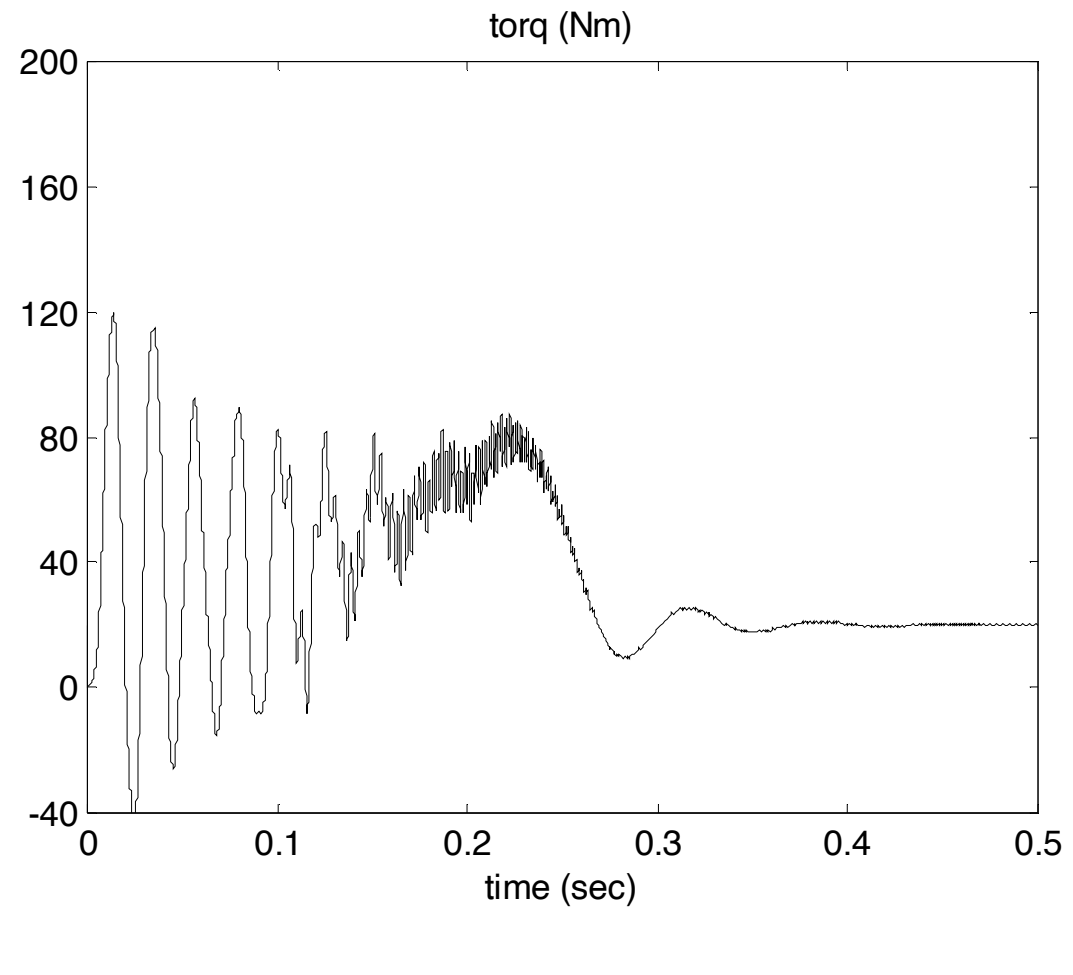

شکل**6-7**:: منحنی تغییرات گشتاور موتور  در حالت **1** میله شکسته

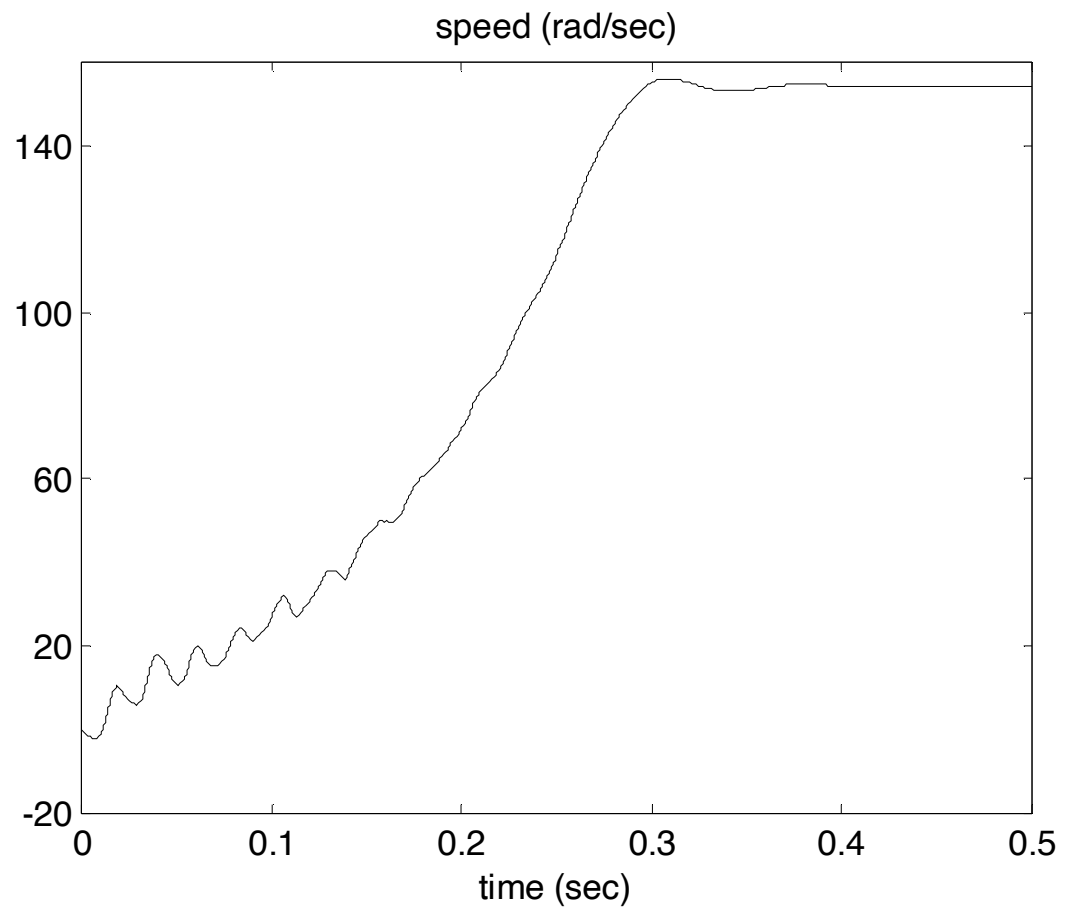

شكل**6-8**: منحنی تغییرات سرعت موتور در حالت **2** میله شکسته

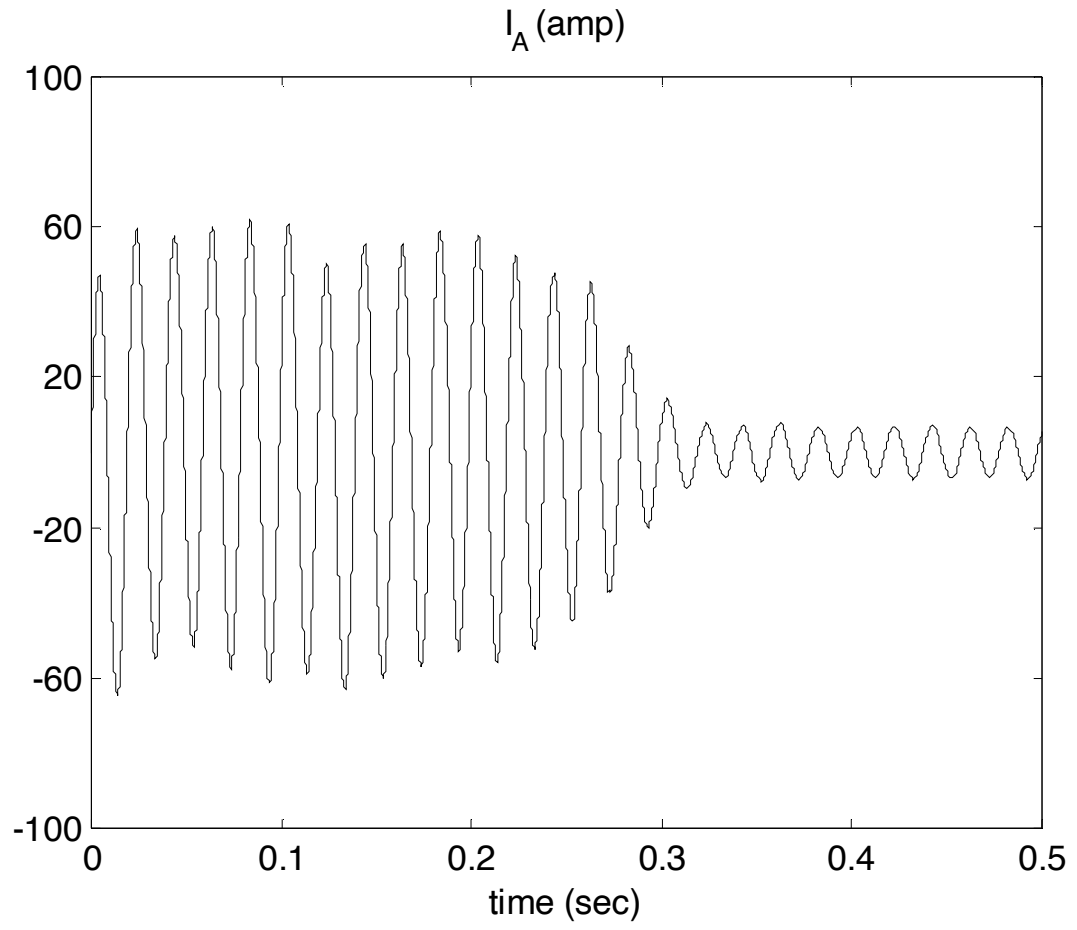

شكل**6-9**: منحنی تغییرات جریان فاز **A** در حالت **2** میله شکسته

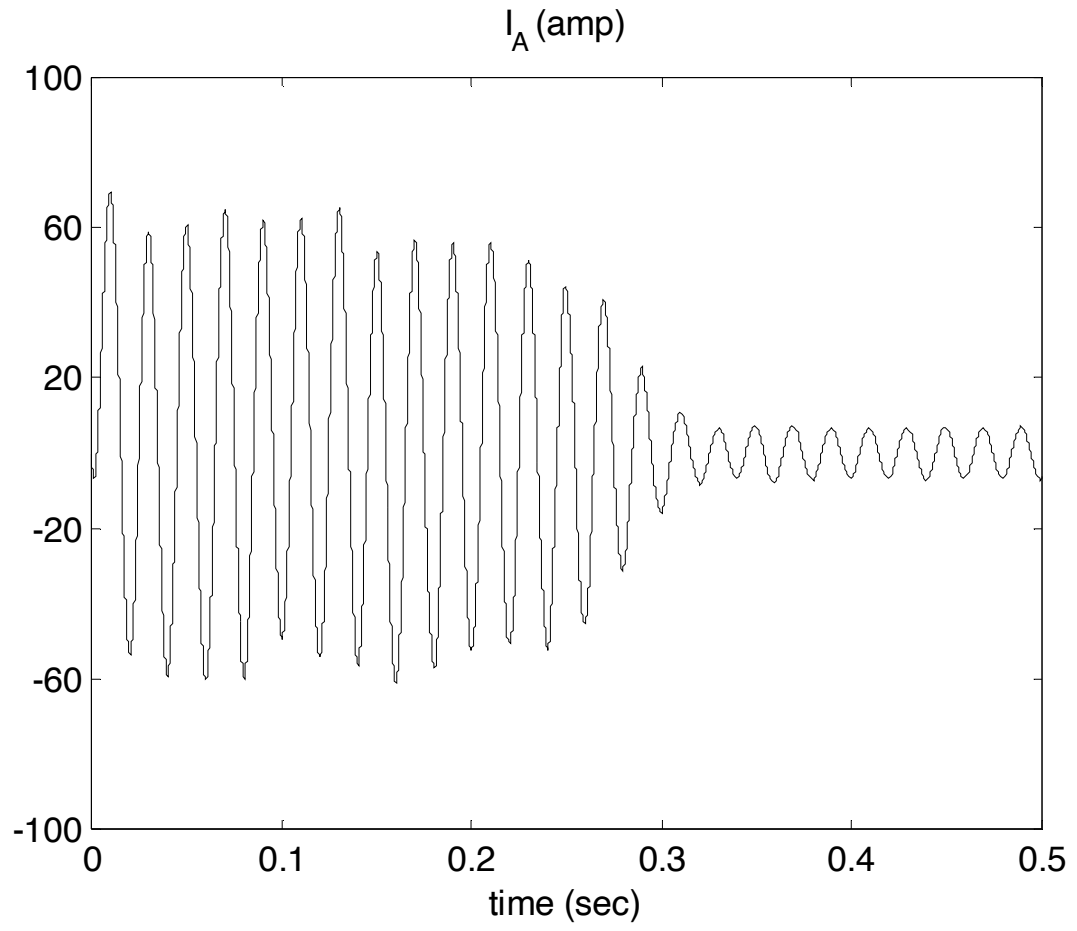

شكل**6-10**: منحنی تغییرات جریان فاز **B** در حالت **2** میله شکسته

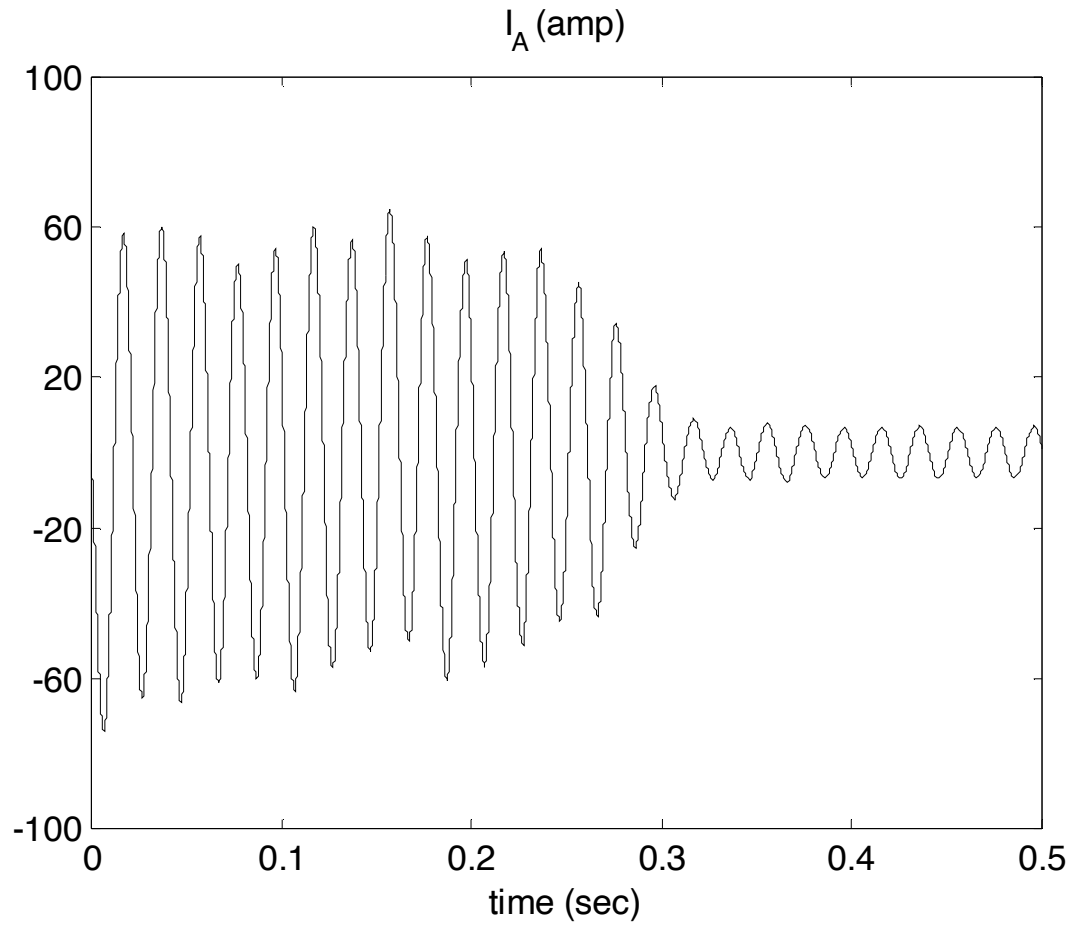

شکل6-11: منحنی تغییرات جریان فاز C در حالت 2 میله شکسته

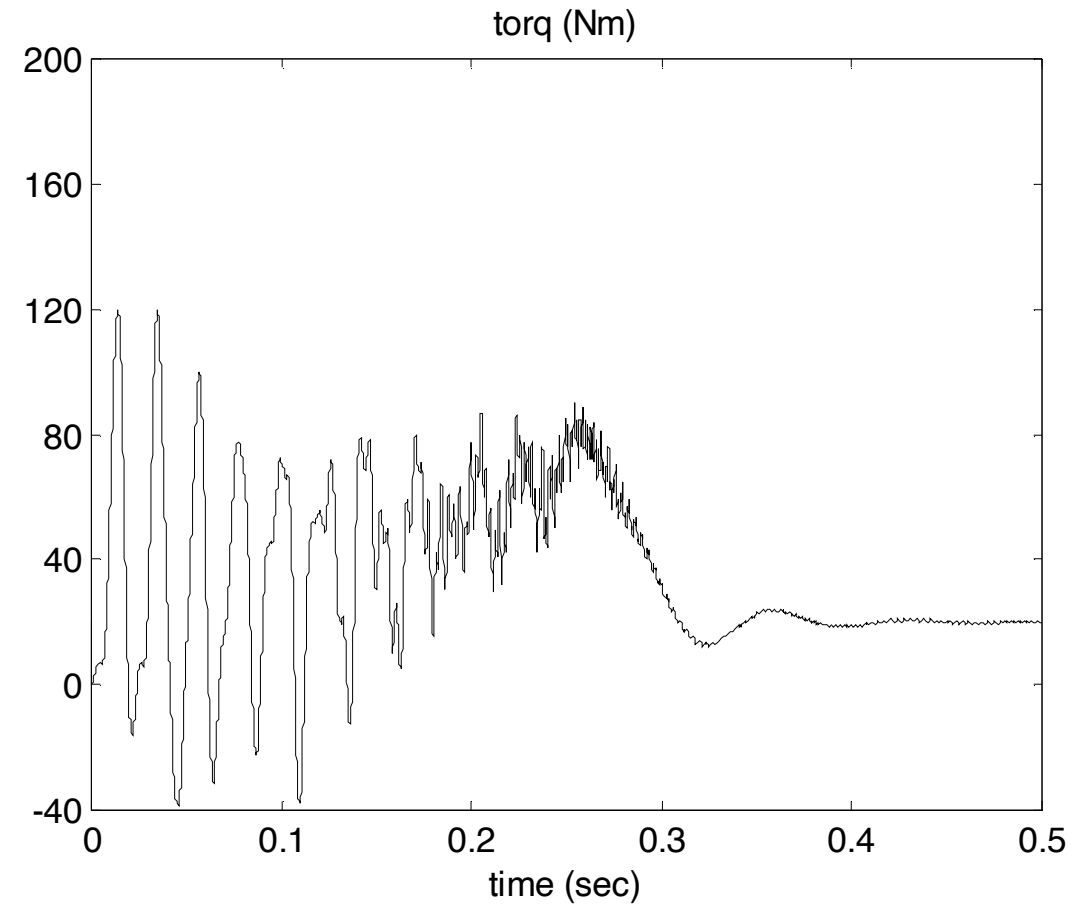

شکل6-12: منحنی تغییرات گشتاور تولیدی موتور در حالت 2 میله شکسته

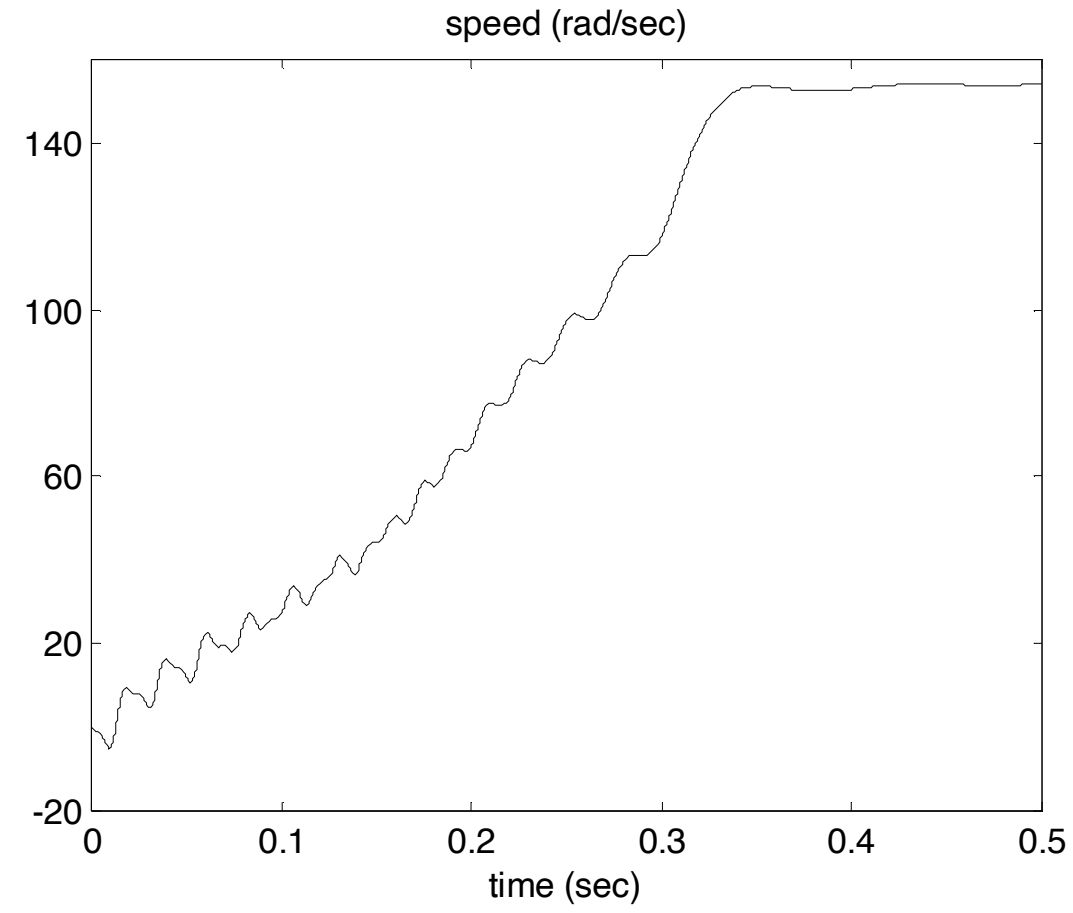

شکل6-13: منحنی تغییرات سرعت موتور در حالت 4 میله شکسته

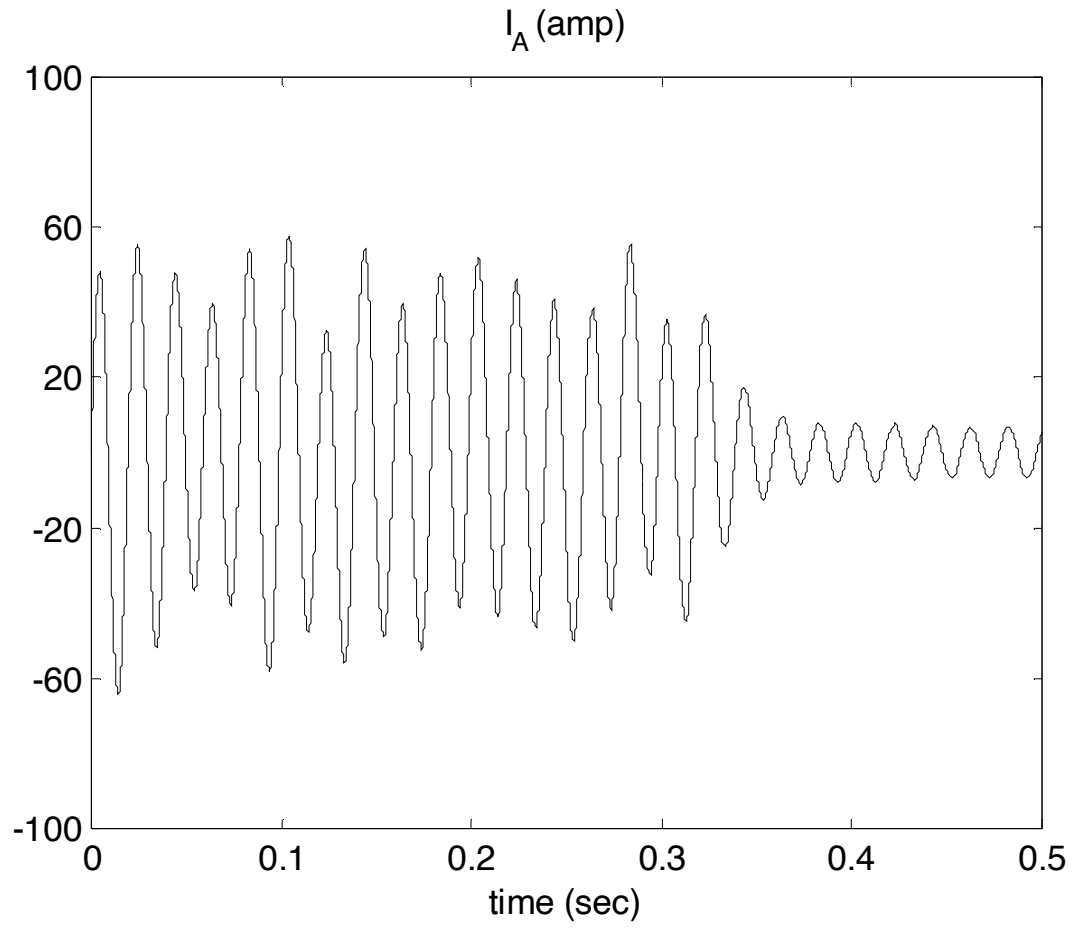

شکل**6-14**: منحنی تغییرات جریان فاز **A** در حالت **4** میله شکسته

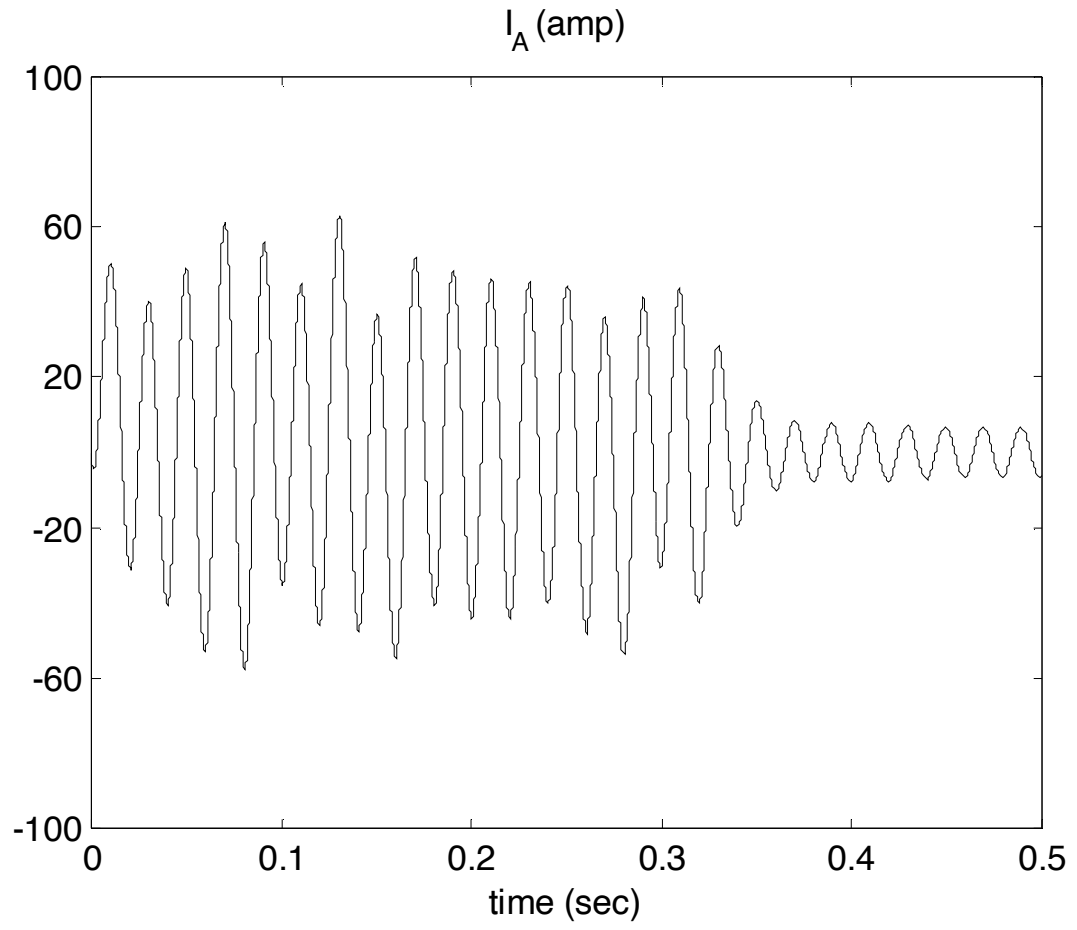

شکل**6-15**: منحنی تغییرات جریان فاز **B** در حالت **4** میله شکسته

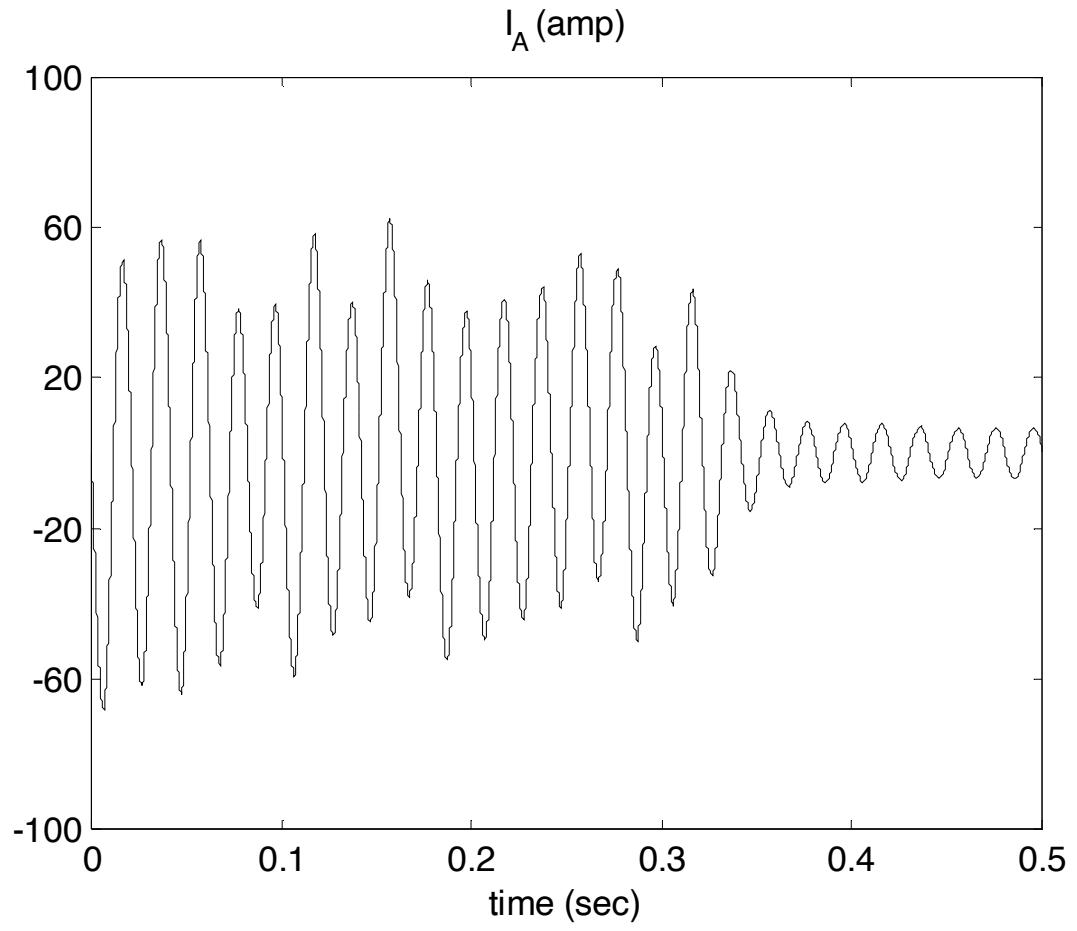

شکل**6-16**: منحنی تغییرات جریان فاز **C** در حالت **4** میله شکسته

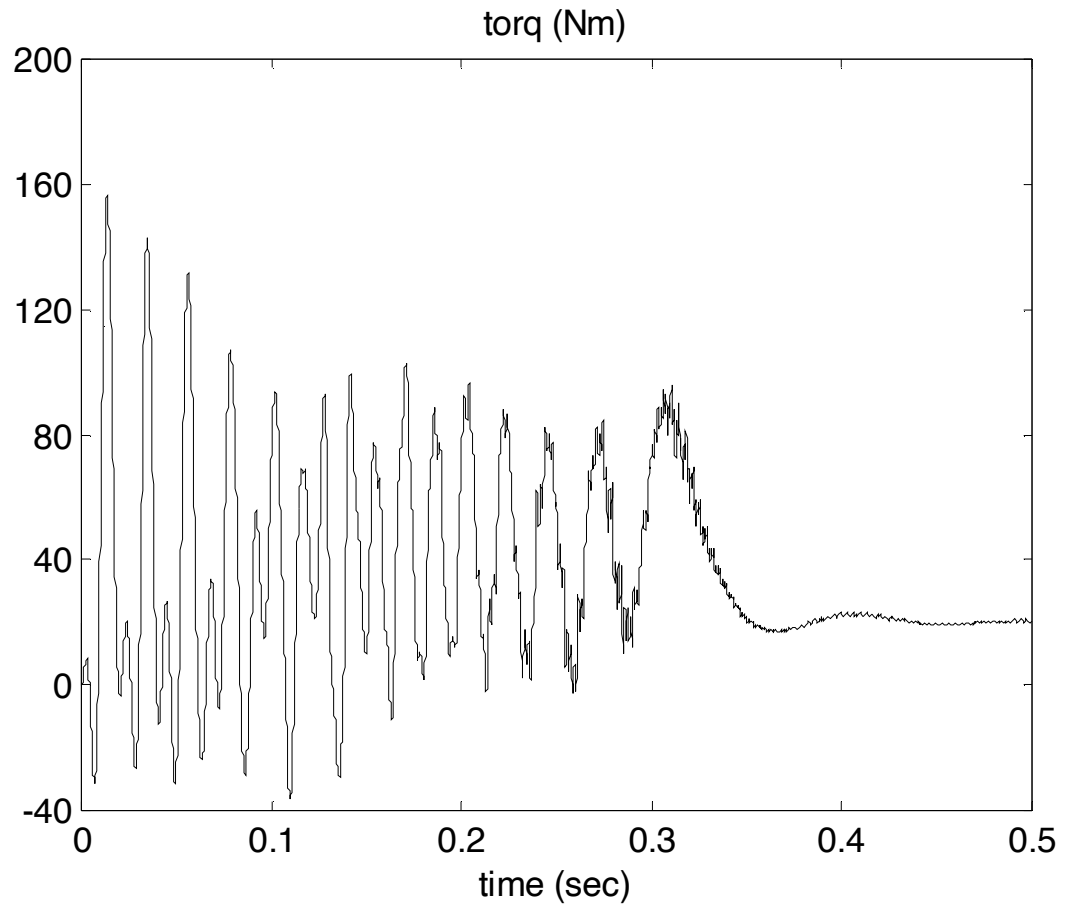

شکل**6-17**: منحنی تغییرات گشتاور موتور در حالت **4** میله شکسته

همانطور که درشکلهای 6-3 تا 6-17 نشان داده شده است هر چقدر تعداد مسله های شکسته روتور بیشتر باشد موتور دیرتر به حالت ماندگار می‌رسد و جریان راه اندازی بیشتری دارد. همچنین بررسی جریان حالت ماندگار موتور در این حالت نشان می دهد که در حالت شکستگی میله های روتور جریان یک منحنی پوش سینوسی با فرکانس کمتر از فرکانس اصلی دارد. شکل 6-18 این وضعیت را برای حالت شکستگی 4 میله روتور نشان می‌دهد. در حالتهای شکستگی یک و یا دو میله این وض عیت وجود دارد اما خیلی محسوس نیست.

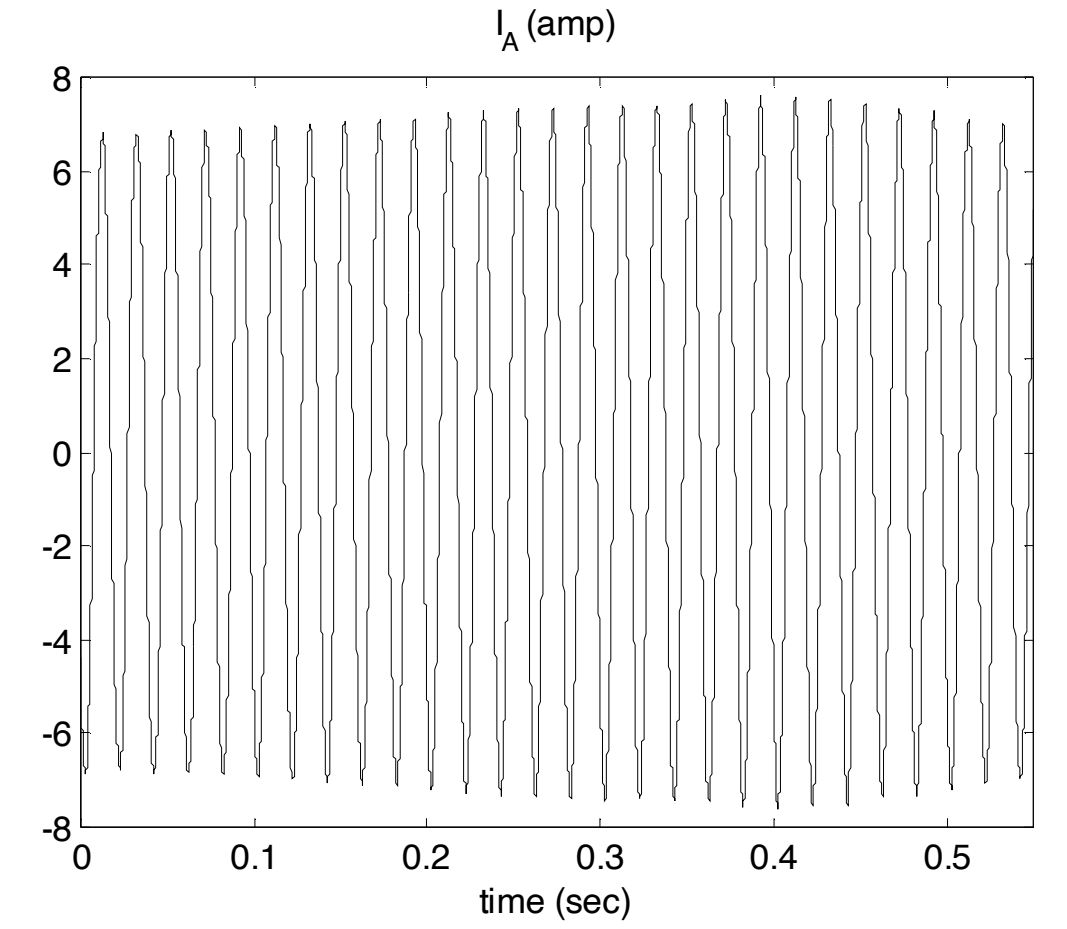

شکل**6-18**: منحنی تغییرات جریان فاز **A** در حالت ماندگار در حالت **4** میله شکسته

## 5-6- تحلیل کامپیوتری موتور با روش حذف حلقه شامل میله های شکسته

فرض کنیم m میله پشت سرهم اول از n میله روتور شکسته باشد در اینصورت تعداد حلقه های روتور به n-m کاهش می یابد. و تغییراتی در ماتریسهای اندوکتانس روتور و اندوکتانس متقابل میان روتور و استاتور پیش می‌آید. در این حالت تابع دور حلقه اول برابر مجموع توابع دور حلقه های حذف شده است. بنابراین اندوکتانس خودی این حلقه تغییر می‌کند. همانطور که در بخش قبل اندوکتانس خودی حلقه ها محاسبه شدند می توان نشان داد که این اندوکتانس بصورت زیر تغییر خواهد کرد:

$$L_r(1,1)=\frac{\mu_0 r_0 l_0}{g_0}\left(2\pi\frac{m+1}{n}\left(1-\frac{m+1}{n}\right)-\frac{\gamma}{3}\right)+2\left(l_b+ml_e+l_e\right) \tag{36-6}$$

همچنین اندوکتانس متقابل کلیه حلقه های روتور با این حلقه بصورت زیر تغییر می کند.

$$L_r(1,2)=L_r(1,n-m)=\frac{\mu_0 r_0 l_0}{g_0}\left(2\pi\frac{m+1}{n^2}+\frac{\gamma}{6}\right)-l_b \tag{37-6}$$

همچنین اندوکتانس متقابل کلیه فازهای استاتور با حلقه اول تغییر می کند. از آنجا که تابع دور حلقه اول در حالت شکسته برابر مجموع توابع دور حلقه های حذف شده است این اندوکتانس در حالت شکسته نیز برابر مجموع اندوکتانسهای متناظر حذف شده است. برای مثال اندوکتانس متقابل بین فاز اول استاتور و حلقه اول روتوری با m میله شکسته برابر مجموع اندوکتانسهای متقابل فاز اول و m حلقه اول روتور در حالت سالم است.

نتایج شبیه سازی حاصل از این روش منطبق بر نتایج شبیه سازی با روش افزایش مقاومت میله های روتور است.

## 6-6- بررسی طیف فرکانسی استاتور در حالتهای مختلف شکستگی میله های روتور

مهمترین هارمونیکهای ناشی از شکستگی میله های روتور دارای فرکانسهای زیر می‌باشند:[46]

$$f_h = (1 \pm 2s) f_0 \tag{38-6}$$

که $s$ لغزش موتور و $f_0$ فرکانس اصلی تغذیه می‌باشد . شکلهای 6-13، 6-14، 6-15 و 6-16 طیف فرکانسی جریان فاز اول استاتور را به ترتیب برای یک میله شکسته، دو میله شکسته و چهار میله شکسته نشان می‌دهند. این طیفها نیز مانند آنچه در فصل قبل گفته شد با نرم افزارMATLAB محاسبه شده اند.

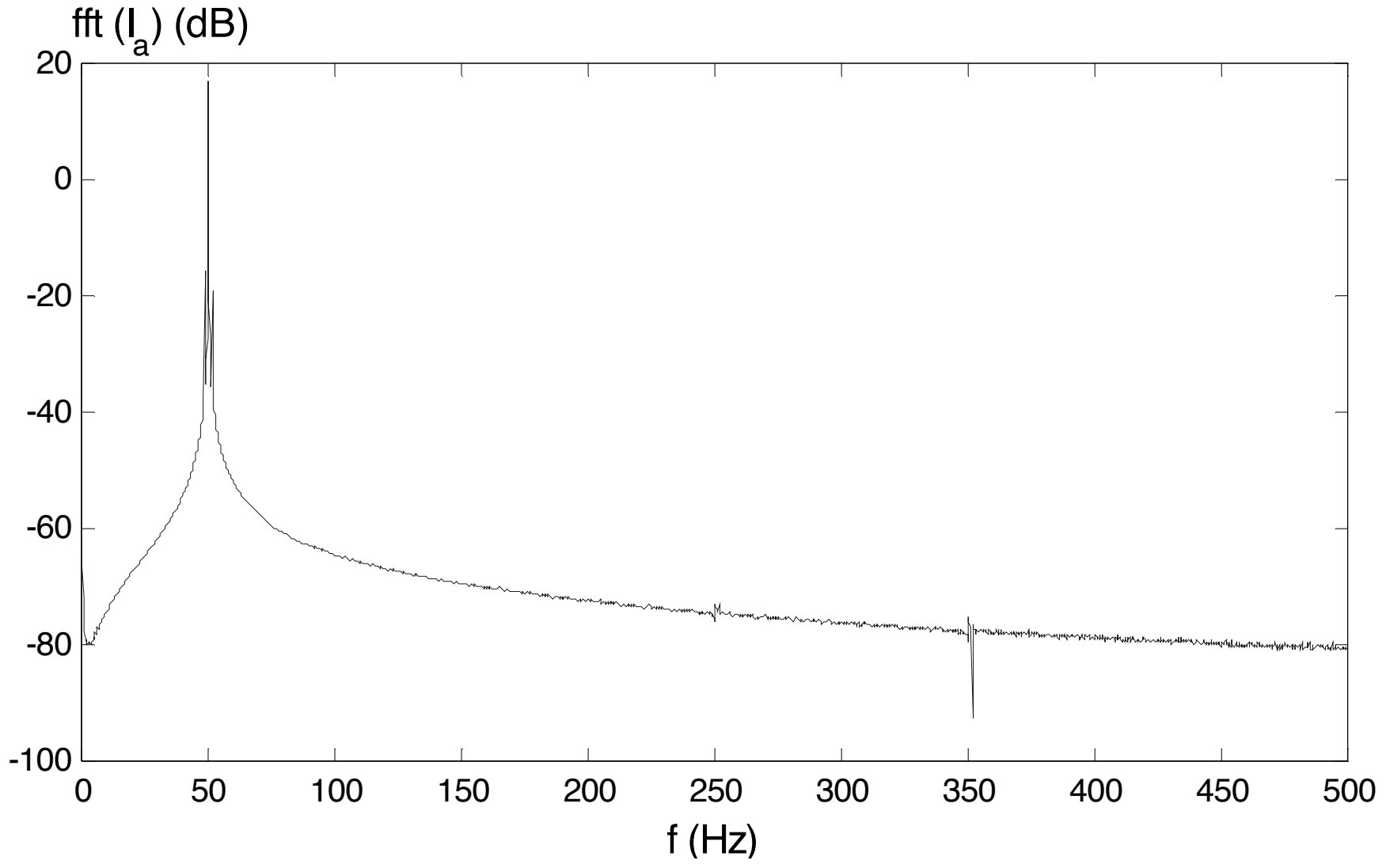

شکل**6-19**: طیف فرکانسی جریان فاز **A** در حالت **1** میله شکسته

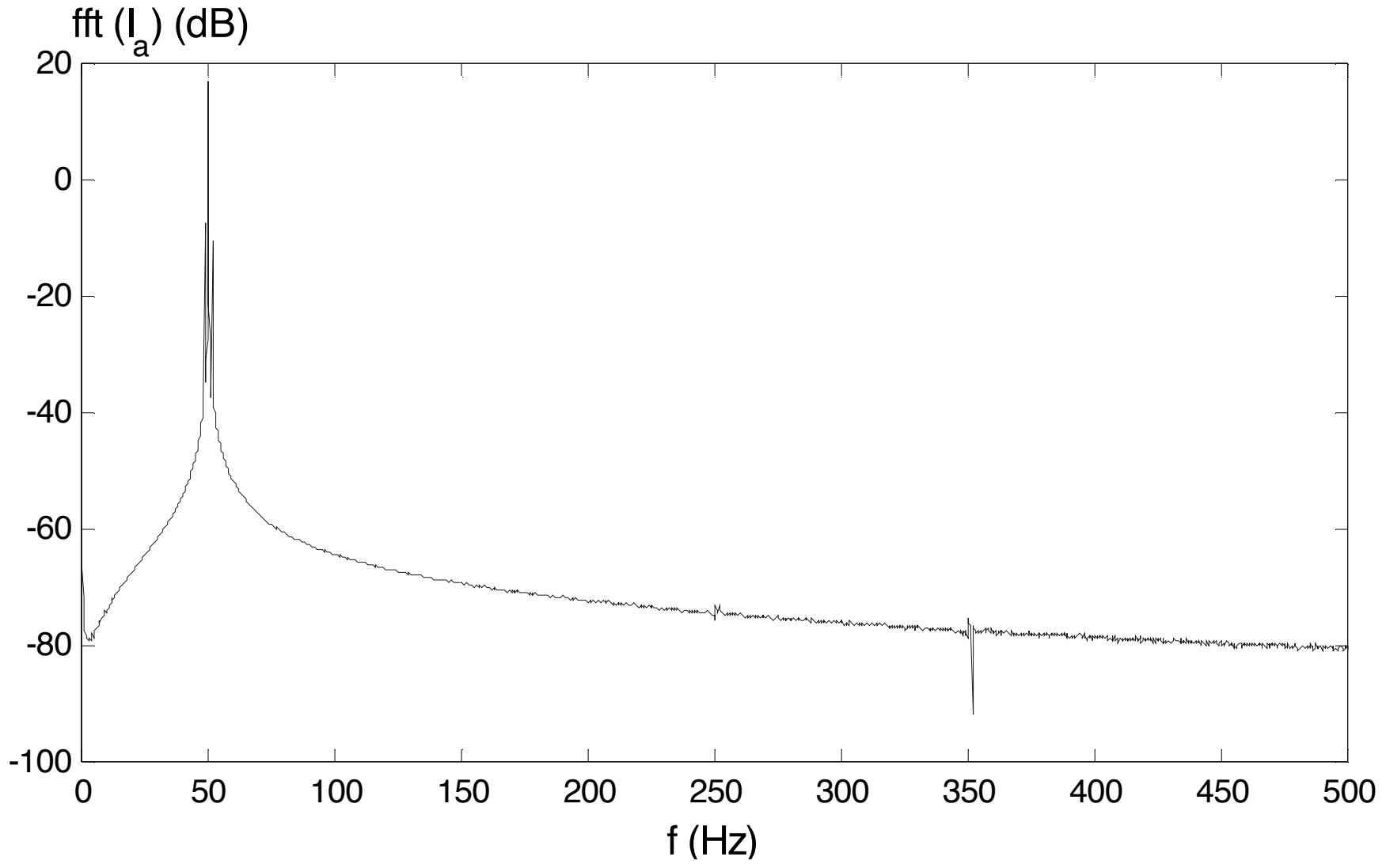

شکل**6-20**: طیف فرکانسی جریان فاز **A** در حالت **2** میله شکسته

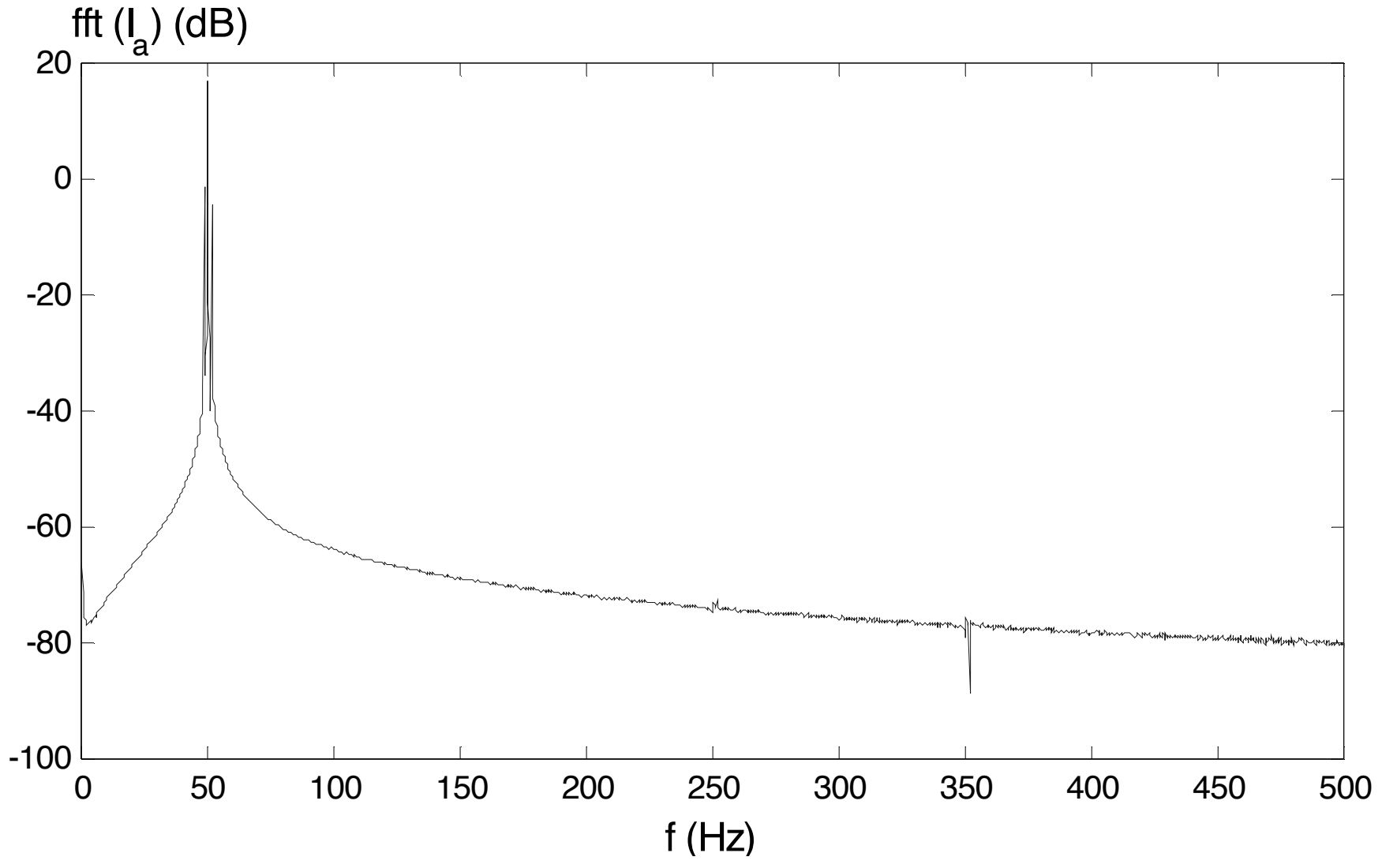

شکل**6-21**: طیف فرکانسی جریان فاز **A** در حالت **4** میله شکسته

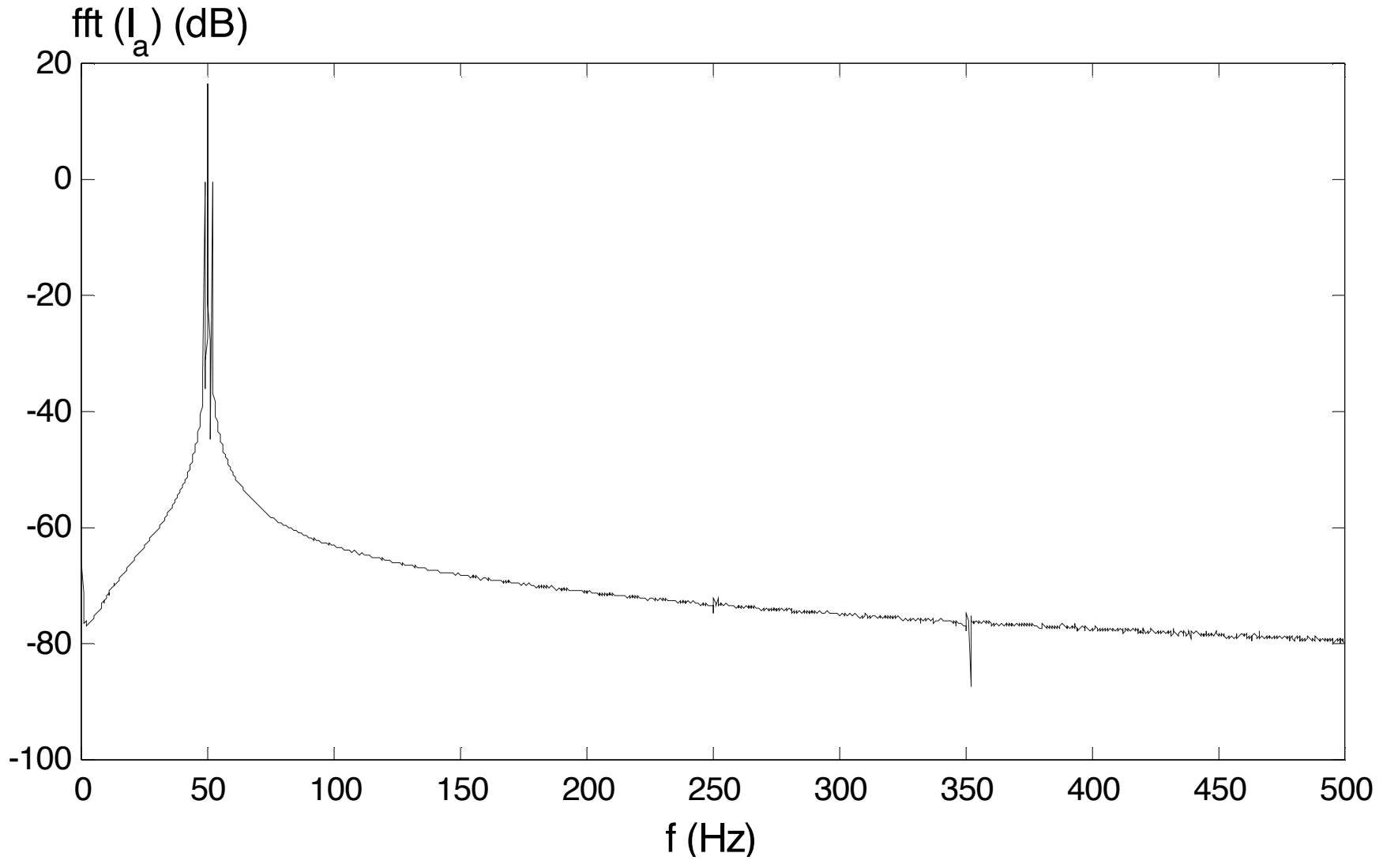

شکل**6-22**: طیف فرکانسی جریان فاز **A** در حالت **6** میله شکسته

شکل 6-23 تا 6-26 این طیف ها را برای محدوده نزدیک فرکانس اصلی نشان می دهد. همانطور که این شکلها نشان می‌دهند با شکستگی میله های بیشتری از روتور هارمونیکهای رابطه (6-38) تقویت می‌شوند.

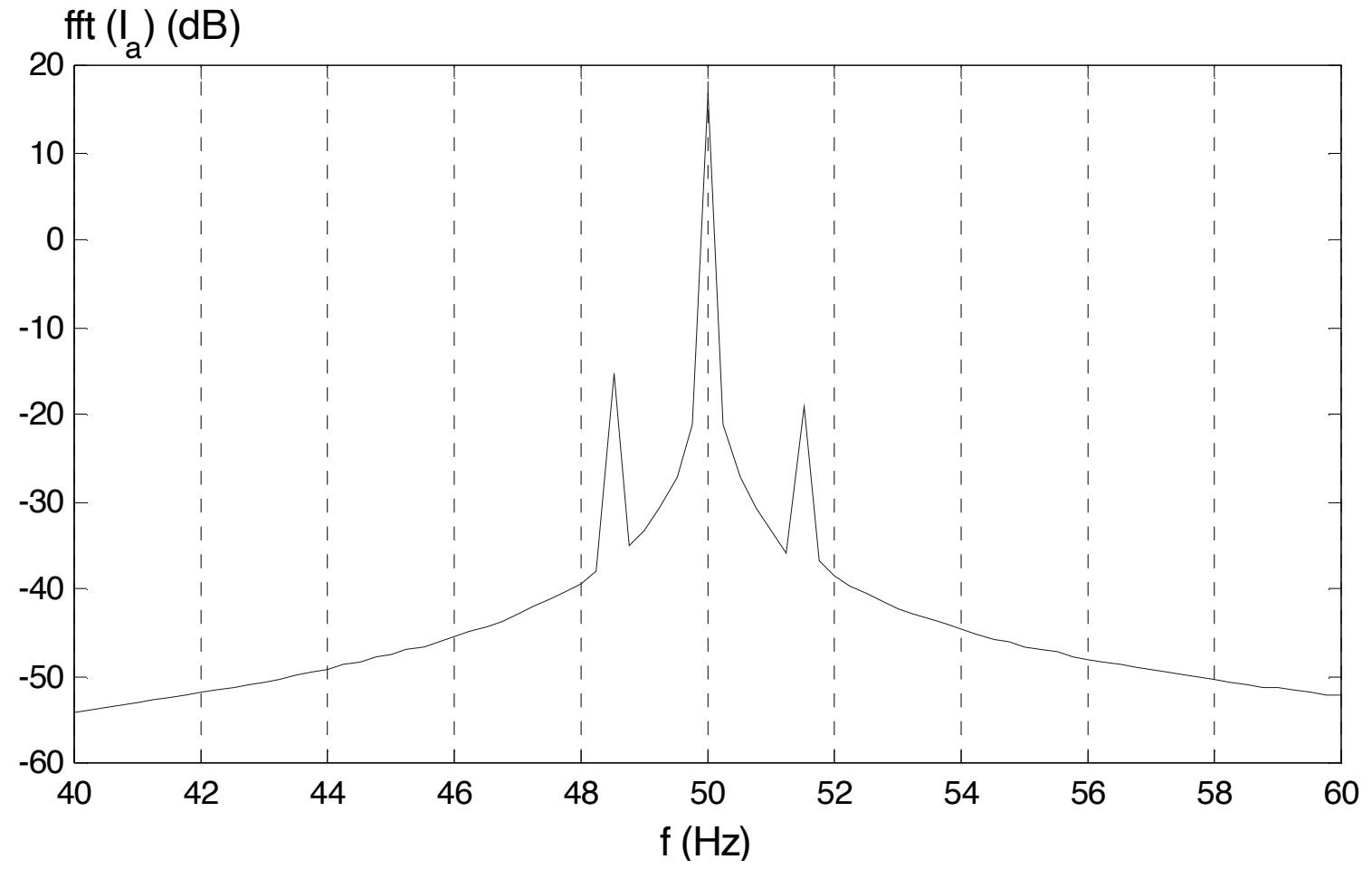

شکل**6-23**: طیف فرکانسی جیان فاز **A** در حالت **1** میله شکسته (در اطراف فرکانس اصلی)

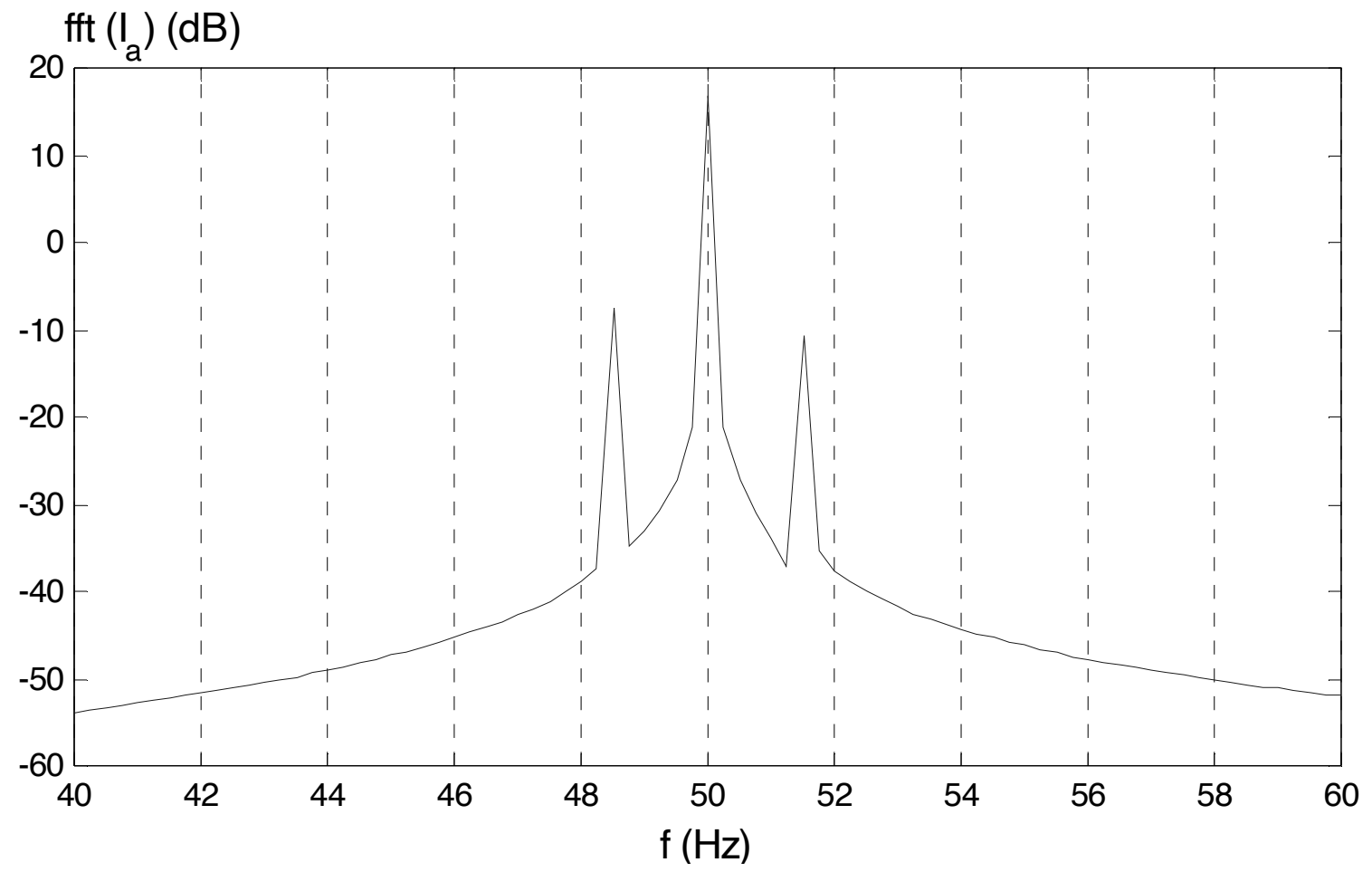

شکل**5-24**: طیف فرکانسی جریان فاز **A** در حالت **2** میله شکسته(در اطراف فرکانس اصلی)

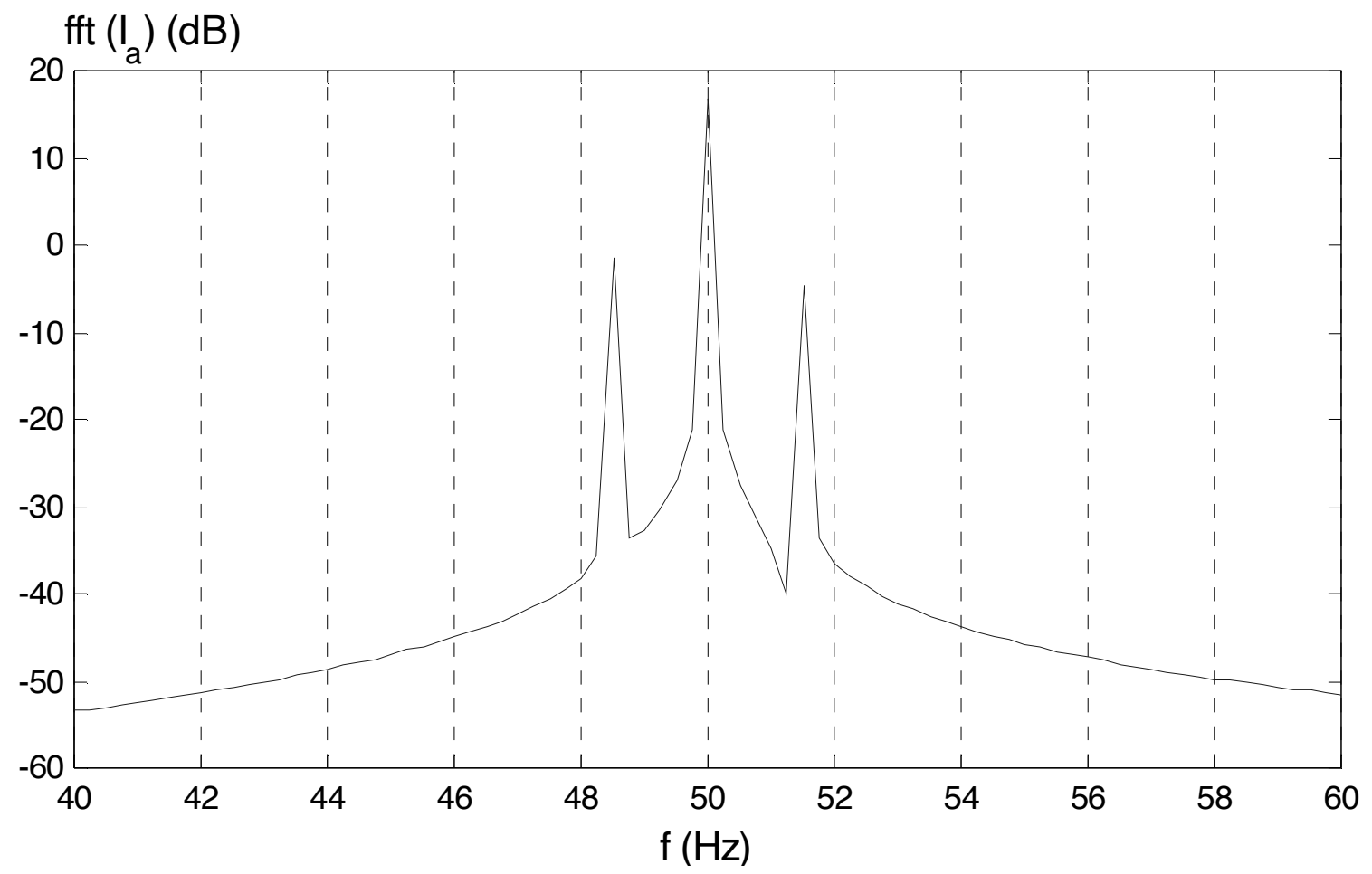

شکل**5-25**: طیف فرکانسی جریان فاز **A** در حالت **4** میله شکسته(در اطراف فرکانس اصلی)

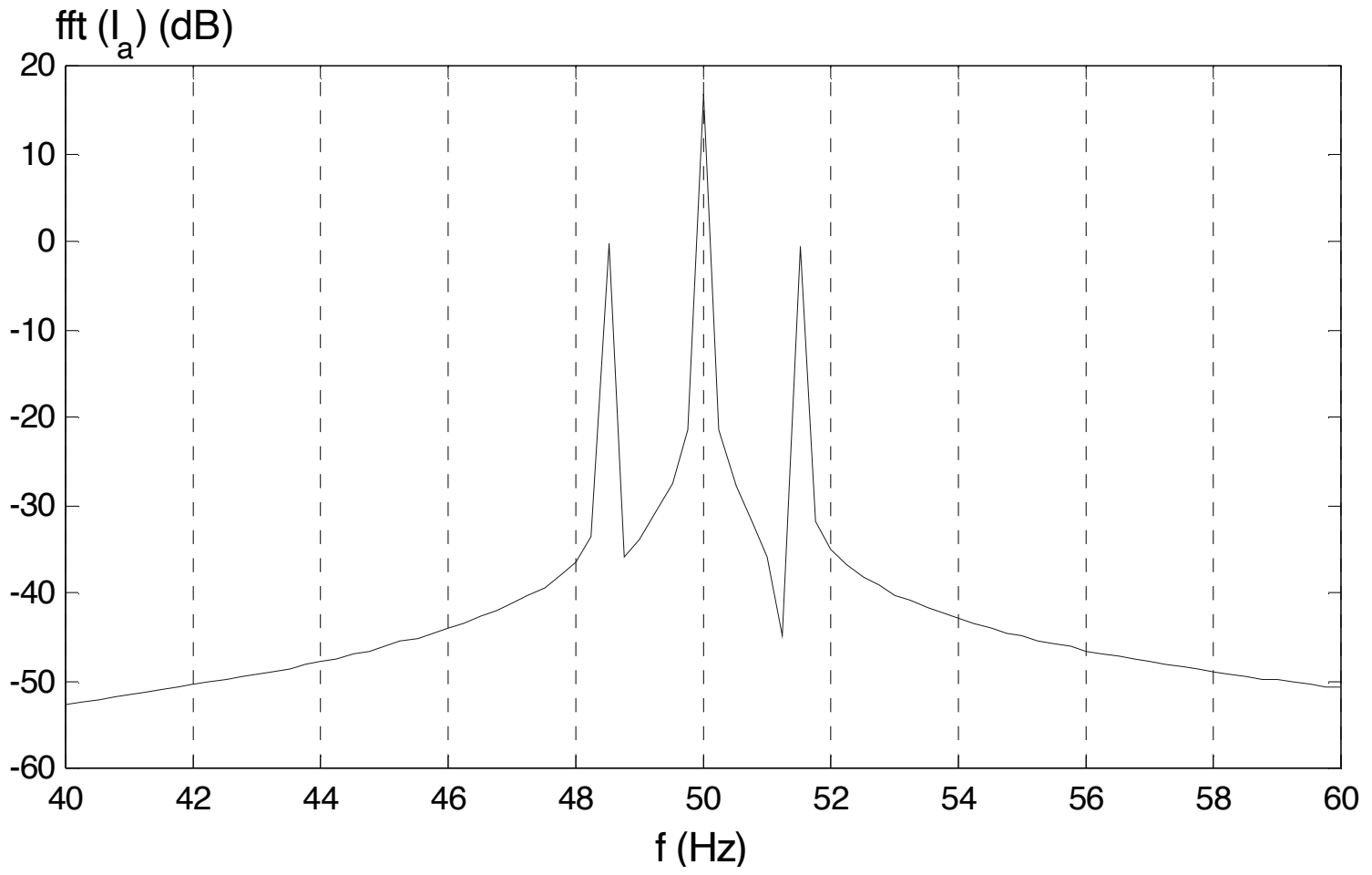

شکل**5-26**: طیف فرکانسی جریان فاز **A** در حالت **6** میله شکسته(در اطراف فرکانس اصلی)

شکل 5-27 نمودار دامنه این هارمونیکها را بر حسب تعداد میله های شکسته نشان می‌دهد.

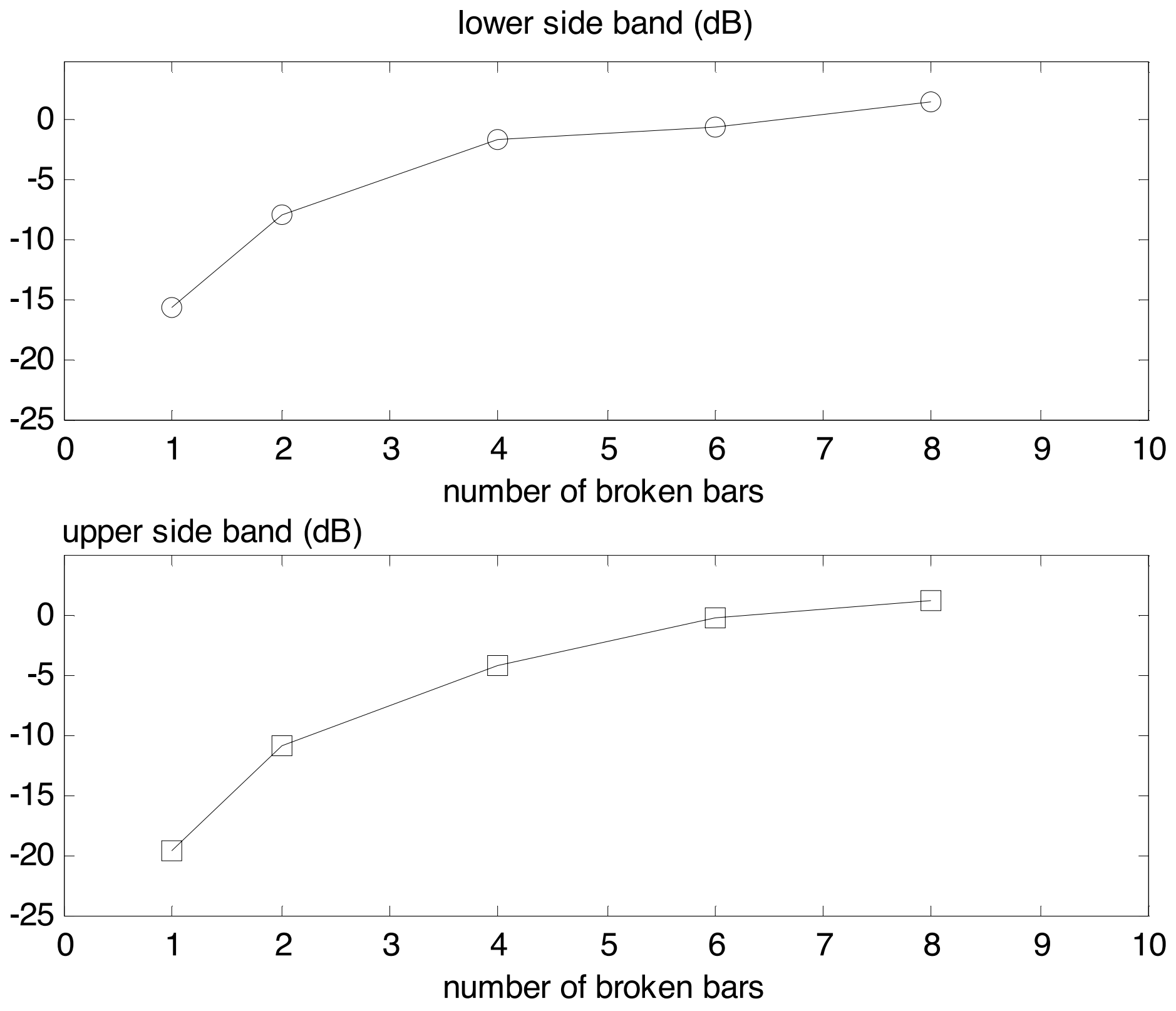

شکل **5-27**: نمودار دامنه این هارمونیکها را بر حسب تعداد میله های شکسته

# نتیجه گیری و پیشنهادات

در این پایاننامه با هدف پیش بینی رفتار دینامیکی موتور القائی سه فاز روتور قفسی تحت شرایط ناهم محوری میان روتور و استاتور و شکستگی میله های روتور به مدلسازی و شبیه سازی کامپیوتری این موتور تحت این شرایط با استفاده از نظریه تابع سیم پیچی پرداخته شد. بدین منظور ضمن تعمیم نظریه تابع سیم پیچی و اصلاح آن در حالت غیر یکنواختی فاصله هوائی مدل دقیقتری برای هدایت مغناطیسی فاصله هوائی ارائه شد. این امر با بازنگری معادلات الکترومغناطیسی نظریه تابع سیم پیچی و محاسبات هندسی موتور صورت گرفت. سپس با استفاده از مدل هندسی موتور برای اولین بار کلیه اندوکتانسهای موتور در حالت کلی ناهم محوری مرکب بصورت تحلیلی با استفاده از انتگرالگیری و محاسبات جبری محاسبه شدند، بگونه‌ای که هر ماتریس اندوکتانس موتور با یک تابع تحلیلی بسته که متغیرهای ورودی آن زاویه مکانیکی روتور و درجه های ناهم محوری ایستا و پویا است محاسبه شود. این نکته قابل ذکر است که این اندوکتانسها با در نظر گرفتن اثر افزایش خطی MMF فاصله هوائی در بالای شیارها و مورب بودن میله های روتور محاسبه شده اند.

محاسبه تحلیلی اندوکتانسهای موتور علاوه بر دقت بیشتر در محاسبات، زمان شبیه سازی کامپیوتری را به نسبت کاهش می‌دهد. پس از شبیه سازی موتور با روش رانگ کوتای مرتبه 4و5 نرم افزار MATLAM عملکرد موتور تحت شرایط مختلف ناهم محوری مورد بررسی قرار گرفت و نشان داده شد که انواع مختلف ناهم‌محوری علاوه بر تاثیر بر روی حالت گذرا و زمان نشست دینامیکی موتور بر روی طیف فرکانسی جریان خط نیز تاثیر گذار هستند، بگونه ای که انواع مختلف ناهم محوری موجب ایجاد هارمونیکهای اضافی معین روی جریان خط می شوند. این هارمونیکها در اطراف فرکانس اصلی و فرکانس شیارگذاری روتور نسبتا مشهود تر ند. فرکانس این هارمونیکها با تغییر بار مکانیکی موتور و یا لغزش تغییر می‌کند.

به منظور مدلسازی موتور تحت شرایط شکستگی میله های روتور دو روش افزایش مقاومت میله های روتور و حذف حلقه میله شکسته شده با یکدیگر مقایسه و نشان داده شد که در روتور با تعداد میله های شکسته اندک این دو روش نتایج کاملا مشابهی دارند. این در حالی است که مدلسازی با استفاده از روش افزایش مقاومت بسیار ساده تر از روش حذف حلقه است اما دقت محاسبات و زمان شبیه سازی روش حذف حلقه بهتر است. در انتها نشان داده شد که

شکستگی میله های روتور نیز علاوه بر ایجاد تغییرات بر روی حالت گذرای موتور و زمان نشست آن، بر روی طیف فرکانسی جریان خط استاتور نیز تاثیر می گذارد. این تاثیر در نواحی اطراف فرکانس اصلی مشهودتر است. فرکانس این هارمونیکها نیز با تغییر بار مکانیکی موتور تغییر می‌کند.

در پایان به منظور تکمیل و ارتقاء مطالعات انجام گرفته موضوعات زیر قابل طرحند:

1. فرکانس هارمونیکهای ناشی از ناهم محوری روتور و استاتور و شکستگی میله های روتور با تغییر بار مکانیکی موتور تغییر می کند. با اینحال پیش بینی رفتار موتور به ازاء تغییرات پارامترهای بار مکانیکی موضوع پیچیده قابل مطالعه است.

2. کلیه شبیه سازیهای بررسی حاضر بر مبنای تغذیه موتور با ولتاژ سه فاز سینوسی خالص انجام گرفته است. این در حالی است که نرم افزار طراحی شده قادر به در نظر گرفتن هر سری زمانی برای ولتاژهای خط است. بررسی اثر هارمونیکهای زمانی روی ناهم محوری روتور و استاتور مقوله دیگری است که در جهت تکمیل کار حاضر قابل بررسی است.

3. هارمونیکهای ناشی از ناهم محوری بین روتور و استاتور دارای فرکانسی متناسب با لغزش موتور هستند. بنابراین بایستی بتوان از روی ناهم محوری ذاتی موجود در موتورهای القائی سرعت موتور را در حین کار تخمین زد.تخمین سرعت با آشکارسازی این هارمونیکها موضوع دیگری است که می تواند در ادامه کار حاضر مطرح گردد.

4. در سالهای اخیر استفاده از تبدیلهای ریاضی دیگری به جز تبدیل فوریه مطرح شده است. یکی از متداولترین این تبدیلها، تبدیل موجک یا (Wavelet) است که نقاط قوت زیادی نسبت به تبدیل فوریه دارد. آشکارسازی اثر ناهم‌محوری روتور و استاتور و یا شکستگی میله های روتور با اعمال تبدیل موجک بر روی سیگنال جریان موضوع دیگری است که در جهت تکمیل بحث حاضر می توان به آن پرداخت.

5. یکی از موضوعات قابل طرح جهت تکمیل نظریه تابع سیم پیچی منظورکردن اثر اشباع است. این مساله به ویژه در حالتهای مختلف ناهم محوری قابل بررسی است.

6. در بررسی حاضر محاسبه تحلیلی اندوکتانسها بصورت دقیق و کامل انجام شده است. از آنجا که اندوکتانسهای روتور و استاتور موتور در حالتهای مختلف ناهم محوری مرکب نسبت به تغییرات زاویه مکانیکی روتور یک رفتار تقریبا سینوسی دارند شاید بتوان این اندوکتانسها را با روابط جبری ساده‌تری تخمین زد تا محاسبه آنها با حجم عملیات جبری کمتر انجام شود.
7. کلیه نرم افزارهای شبیه سازی بااستفاده از نرم افزار MATLAB و روشهای حل معادلات دیفرانسیل این نرم افزار طراحی شده است. پیش بینی می شود که استفاده از نرم افزار C++ و طراحی روالهای حل معادلات دیفرانسیل در این نرم افزار موجب افزایش سرعت محاسبات می‌شود.
8. با استفاده از نرم افزارهای ارائه شده می توان با اندکی تغییر موتور را در حالت وقوع همزمان ناهم محوری بین روتور و استاتور و شکستگی میله های روتور تحلیل کرد. مطالعه تاثیر همزمان این دو عیب بر روی عملکرد موتور و طیف فرکانسی جریان خط یک موضوع جالب و قابل بررسی دیگر است.
9. به منظور افزایش دقت مدلسازی موتور می توان اثر شیارها، افزایش mmf فاصله هوایی در بالای شیارها بر حسب شکل شیار و تغییرات شعاعی شدت میدان مغناطیسی را لحاظ کرد.
10. بررسی عملکرد موتور در حالت انحراف محور موتور یا Missalighment

# مراجع

# ضمیمه یک: مشخصات موتور القایی شبیه سازی شده

| | |
|---|---|
| **number of stator phases** | 3 |
| **number of rotor bars** | 40 |
| **number of pole pairs** | 2 |
| **stator phase resistance** | 1.75Ω |
| **rotor bar resistance** | $31\times10^{-6}\Omega$ |
| **end-ring segment resistance** | $2.2\times10^{-6}\Omega$ |
| **view angle of rotor bars** | π/86**rad** |
| **length of stack** | 0.11**m** |
| **length of air-gap** | 0.0008**m** |
| **leakage inductance of stator phases** | 0.009H |
| **turn number of stator phase** | 56 |
| **rotor bar leakage inductance** | $95\times10^{-9}H$ |
| **end-ring segment leakage inductance** | $18\times10^{-9}H$ |

# ضمیمه دو : راهنمای استفاده از نرم افزارهای نوشته شده

## 1- توضیحات کلی

این نرم افزار به منظور تحلیل موتور القائی سه فاز روتور قفسی تحت شرایط مختلف ناهم‌محوری میان روتور و استاتور اعم از ناهم‌محوری ایستا، ناهم‌محوری پویا و ناهم‌محوری مرکب و شکستگی میله های روتور طراحی شده است. نرم افزار مذکور با استفاده از معادلات تزویج الکترومغناطیسی میان مدارهای الکتریکی موتور به محاسبه جریان فازهای استاتور، جریان میله های روتور و سرعت مکانیکی روتور می‌پردازد. یک قسمت اساسی این نرم افزار مربوط به محاسبه تحلیلی اندوکتانسهای مغناطیسی مدارهای موتور است. فلوچارت کلی این نرم افزار در شکل 1 نشان داده شده است. نرم افزار مذکور در محیط Windows 2000 و تحت نرم افزار MATLAB 5.3 طراحی شده است.

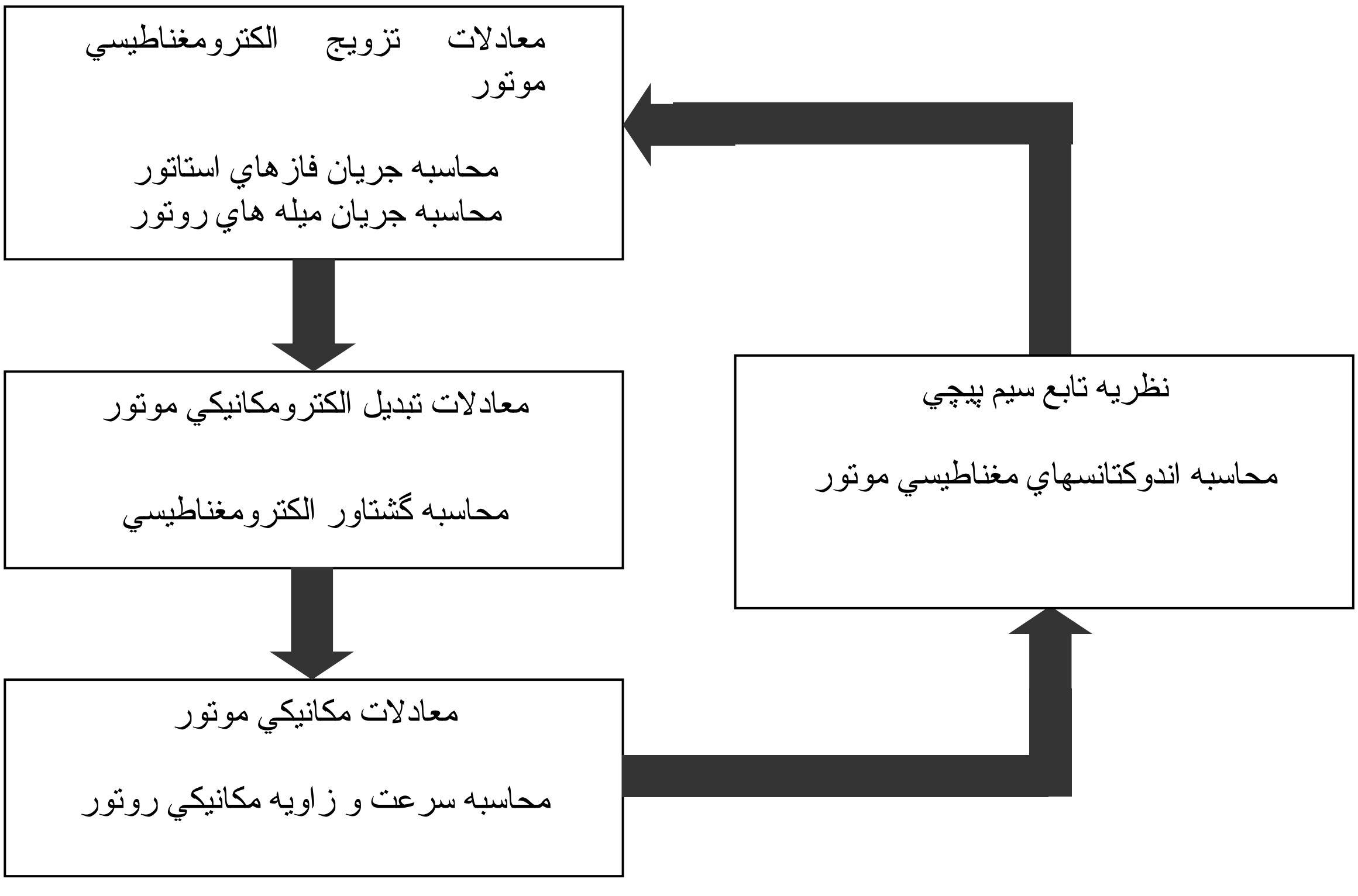

شکل 1 : فلوچارت کلی برنامه

## 2- راهنمای کاربران

✓ نحوه اجرا و نصب نرم افزار :

پس از نصب نرم افزار MATLAB 5.3 در محیط Windows 2000 اجرای برنامه دارای مراحل اصلی زیر است:

الف : قرار دادن فایلهای برنامه در مسیر فعال نرم افزار MATLAB : برای این منظور پس از انتخاب گزینه Set Path از ستون File درMATLAB COMMAND WINDOW ، آدرس محلی که فایلهای برنامه قرار دارد را مشخص کرده و گزینه Save Path را انتخاب می‌کنیم.

ب : تکمیل فایل پارامترهای الکتریکی و مکانیکی موتور: این پارامترها که در قسمتهای بعدی معرفی می شوند در فایل Data.M ذخیره می‌شوند. به منظور مقداردهی به این پارامترها بایستی فایل مذکور را در MATLAB EDITOR تکمیل نمود.

ج : تعیین پارامترهای تحلیل : این پارامترها که شامل زمان تحلیل، مقادیر اولیه متغیرها و تلرانسهای حل معادلات می‌باشند در فایل Run.M مقدار دهی می شوند.

د : اجرای فایل Run.M در MATLAB COMMAND WINDOW

ه : پس از اجرای کامل برنامه کلیه متغیرهای خروجی برنامه در فایلی به نام Result.M ذخیره می‌شوند.

فهرست پارامترهای ورودی برنامه :

پارامترهای ورودی برنامه به دو دسته کلی پارامترهای مدلسازی و پارامترهای تحلیل دسته بندی می شوند. در جدول 1 این پارامترها و واحد آنها معین شده است.

جدول 1 : پارامترهای ورودی برنامه

| نام پارامتر | نماد | واحد |
| --- | --- | --- |
| تعداد فازهای استاتور | m | --- |
| تعداد میله های روتور | n | --- |
| تعداد زوج قطب استاتور | p | --- |
| ممان اینرسی بار مکانیکی | J | N.m/s |
| بار مکانیکی | Tl | N.m |
| ولتاژ تغذیه | Vs | Volt |
| فرکانس تغذیه | Ws | Rad/sec |
| مقاومت یک فاز استاتور | rst | Ohm |
| اندوکتانس پراکندگی یک فاز استاتور | lls | Hanri |
| مقاومت یک میله روتور | r_bar | Ohm |
| مقاومت یک قسمت از طوق انتهایی روتور | r_end | Ohm |
| اندوکتانس پراکندگی یک میله روتور | l_bar | Hanri |
| اندوکتانس پراکندگی یک قسمت از طوق انتهایی روتور | l_end | Hanri |
| تعداد دور سیم بندی فازهای استاتور | Ns | --- |
| زاویه دید میله های روتور از محور استاتور | gama | rad |
| طول فاصله هوائی در حالت متقارن | g | m |
| شعاع روتور | rot_rad | m |
| طول روتور | stack_length | m |
| ضریب نفوذپذیری خلاء | u | --- |
| ضریب ناهم محوری ایستا | es | --- |
| ضریب ناهم محوری پویا | ed | --- |
| زمان شبیه سازی | t_end | sec |
| تلرانس شبیه سازی | RTol | --- |
| زمان شبیه سازی | ATol | --- |

✓ فهرست متغیرهای خروجی برنامه :

متغیرهای خروجی برنامه در فایل Result.M بصورت ماتریسی ذخیره می شود. ماتریس مذکور بصورت [t l y] است. t بردار زمان و y شامل بردارهای متغیرهای خروجی برنامه متناظر با بردار زمان است. شکل کلی y بصورت زیر می باشد:

y= [بردار جریان فاز اول استاتور| بردار جریان فازدوم استاتور | بردار جریان فاز سوم استاتور| بردار جریان حلقه اول روتور| بردار جریان حلقه دوم روتور| ... | بردار جریان حلقه آخر روتور| زاویه مکانیکی روتور | سرعت مکانیکی روتور]

✓ نحوه وارد کردن پارامترهای ورودی برنامه :

همانطور که گفته شد پارامترهای ورودی برنامه در دو فایل Data.M و Run.M وارد می‌شوند. شمای کلی این فایلها در صفحات بعدی آمده است. یک روش دیگربرای وارد کردن پارامترهای ورودی برنامه وارد کردن پارامترهای در هر فایل است، این روش منجر به افزایش سرعت شبیه سازی می‌شود.

در انتهای معرفی نرم افزار پارامترهای لازم برای هر فایل مشخص می‌شود، در اینصورت کاربر می‌تواند با مقدار دهی این پارامترها در اول هر فایل دستوری را که مربوط به خواندن پارامترها از DATA.M است را حذف کند.

Data.M

```
es=0.45;             %static eccentricity coefficient
ed=0;                %dynamic eccentricity coefficient
n=40;                %number of rotor bars
m=3;                 %number of stator phases
P=2;                 %number of rotor pair poles
J=0.05;              %inertia coefficient of load
Tl=1;                %load
rst = 1.75;          %resistance of stator phases
r_bar=31e-6;         %resistance of rotor bars
r_end=2.2e-6;        %resistance of rotor end-ring segments
rot_rad=0.0820;      %rotor radius
stack_length=0.11;   %stack length
g=0.0008;            %length of air-gap
u=4*pi*1e-7;         %air permeability
Vs=380;              %supply voltage
```

```
Ws=314;            %supply frequency
lls=0.009;         %leakage inductance of stator phases
l_bar=95e-9;       %leakage inductance of rotor bars
l_end=18e-9;       %leakage inductance of rotor end-rings
```

**Run.M**

```
clc;clear;
load data
t_end=0.5;
RTol=1e-6;
ATol=1e-6;
option=odeset('RelTol',Rtol,'AbsTol',ATol*ones(m+n+2,1));
[t,y]  = ode23('diffrential',[0 t_end],initial,option);
save result t y
```

✓ نحوه خواندن متغیرهای خروجی برنامه :

با اجرای دستور زیر در محیط متغیرهای خروجی برنامه در حافظه قرار می گیرد:

Load result [t y]

## 3-گزارش فنی نرم افزار

مجموعه فایلهای محاسبه اندوکتانس متقابل فاز اول استاتور و حلقه اول روتور:

همانطور که در فصل سوم توضیح داده شد اندوکتانسهای متقابل میان فازهای اتساتور و حلقه های روتور در حالت ناهم محوری روتور و استاتور از روی اندوکتانس متقابل فاز اول استاتور و حلقه اول روتور حالت سالم قابل محاسبه است. به همین دلیل چند فایل از نرم افزار به منظور محاسبه تحلیلی این اندوکتانس در نظر گرفته شده است. تابع Lsr_a1 با داشتن مقدار teta (زاویه مکانیکی روتور) این اندوکتانس را محاسبه می کند. در فصل سوم در مورد نحوه محاسبه این اندوکتانس از روی تابع دور حلقه اول روتور و توابع انتگرال اول و دوم آن بحث شد. برای پیاده سازی روش بحث شده از فایلهای جدول 3 استفاده شده است.

**جدول 3 : M فایلهای محاسبه اندوکتانس متقابل فاز اول استاتور و حلقه اول روتور حالت سالم**

| ورودی | فایل | خروجی |
|---|---|---|
| زاویه (در دستگاه استاتور) | R_tf.M | مقدار تابع دور روتور |
| زاویه (در دستگاه استاتور) | R_itf.M | مقدار انتگرال تابع دور روتور |
| زاویه (در دستگاه استاتور) | R_iitf.M | مقدار انتگرال دوم تابع دور روتور |
| زاویه مکانیکی روتور | Lsr_a1.M | مقدار اندوکتانس متقابل فاز اول استاتور و حلقه اول روتورد حالت سالم |
| زاویه مکانیکی روتور | DLsr_a1.M | مقدارمشتق‌اندوکتانس متقابل‌فازاول‌استاتور وحلقه اول روتورموتور سالم |

شکل 2 فلوچارت محاسبه این اندوکتانس را با استفاده از فایلهای معرفی شده نشان می دهد.

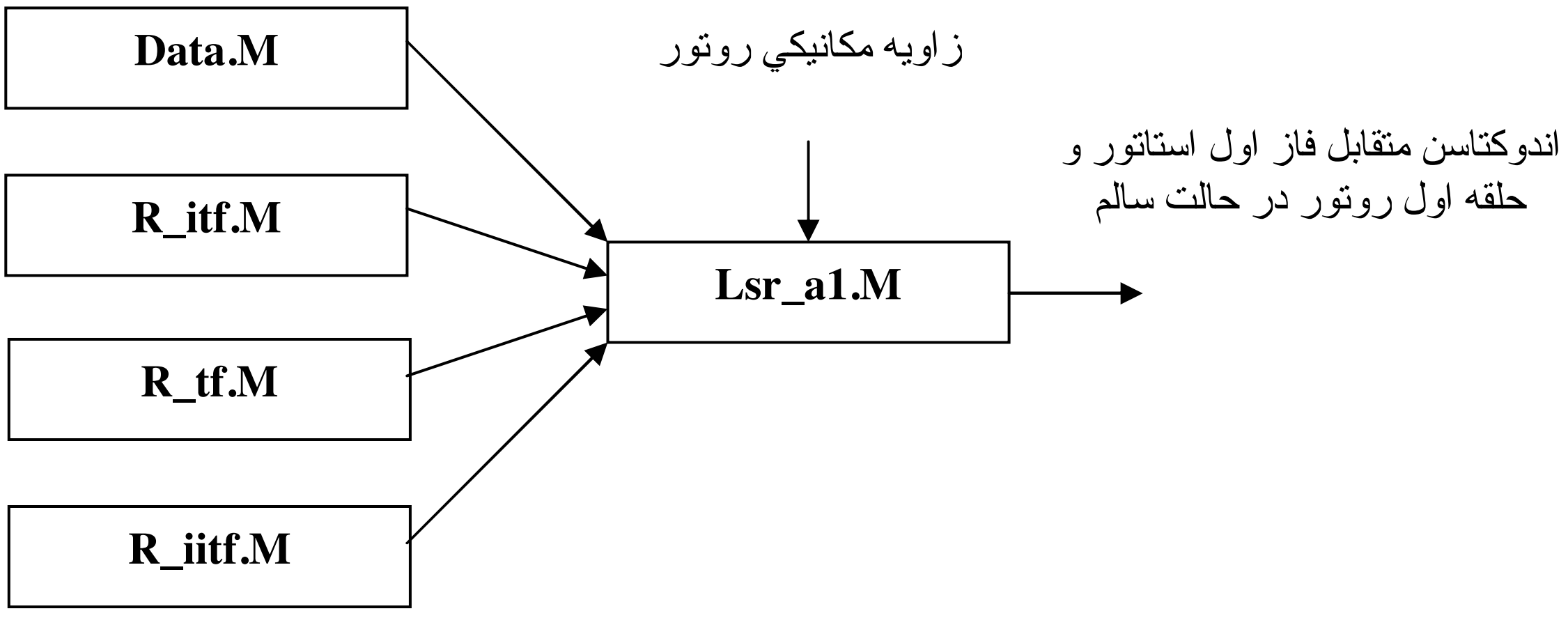

شکل2 : فلوچارت محاسبه اندوکتانس متقابل فاز اول استاتور و حلقه اول روتور در حالت سالم

مجموعه فایلهای محاسبه ماتریسهای اندوکتانس موتور:

جدول 3 این فایلها و ورودی و خروجی آنها را نشان می دهد.

جدول3 :M فایلهای محاسبه ماتریس اندوکتانسها و مشتق اندوکتانسهای موتور در حالتهای مختلف ناهم‌محوری میان روتور و استاتور

| خروجی | فایل | ورودی |
|---|---|---|
| ماتریس اندوکتانس استاتور | Ls.M | زاویه مکانیکی روتور<br>ضریب ناهم محوری ایستا<br>ضریب ناهم محوری پویا |
| ماتریس مشتق اندوکتانس استاتور | dLs.M | |
| ماتریس اندوکتانس روتور | Lr.M | |
| ماتریس مشتق اندوکتانس روتور | dLr.M | |
| ماتریس اندوکتانس میان استاتور و روتور | Lsr.M | |
| ماتریس مشتق اندوکتانس میان استاتور وروتور | dLsr.M | |

فایل محاسبه معادلات فضای حالت موتور:

در فایل diffrential.M این معادلات ذخیره شده است. ورودی این فایل مقدار متغیرهای فضای حالت و خروجی آن مقدار مشتق معادلات فضای حالت است. به این ترتیب معادلات فضای حالت موجود در این فایل از طریق دستور ode قابل حل است. این فایل بااستفاده از فایلهای محاسبه اندوکتانس این معادلات را تشکیل می‌دهد. در حالت شکستگی میله های روتور جهت محاسبه اندوکتانسهای استاتور و روتور موتور لزومی به استفاده از این فایلها نیست و در بدنه فایل

diffrential.M این اندوکتانسها محاسبه می‌شوند. اما برای محاسبه اندوکتانسهای متقابل بـین روتـور و اسـتاتور از فایلهـای محاسبه اندوکتانس استفاده شده است.

## 4-کدهای برنامه

**R_tf.M**

```
function y=r_tf(fi)

%parameters
n=40;     %number of rotor bars
gama=pi/86; %rotor bar view angle (rad)

%load data

alfa_r=2*pi/n;
if (fi<=0)
   y=0;
end;
if (fi>0)&(fi<=gama)
   y=(1/gama)*(fi);
end;
if (fi>gama)&(fi<=alfa_r)
   y=1;
end;
if (fi>alfa_r)&(fi<=gama+alfa_r)
   y=1-(1/gama)*(fi-alfa_r);
end;
if (fi>alfa_r+gama)
   y=0;
end;
```

**R_itf.M**

```
function y=r_itf(fi)
%parameters
n=40;        %number of rotor bars
gama=pi/86; %rotor bar view angle (rad)

%load data

alfa_r=2*pi/n;

if (fi<=0)
   y=0;
```

```
end;
if (fi>0)&(fi<=gama)
   y=(1/gama)*0.5*(fi)^2;
end;
if (fi>gama)&(fi<=alfa_r)
   y=r_itf(gama)+fi-gama;
end;
if (fi>alfa_r)&(fi<=gama+alfa_r)
   y=r_itf(alfa_r)+(fi-alfa_r)-(1/gama)*0.5*(fi-alfa_r)^2;
end;
if (fi>alfa_r+gama)
   y=r_itf(alfa_r+gama);
end;
```

## R_iitf.M

```
function y=r_iitf(fi)
%parameters
n=40;        %number of rotor bars
gama=pi/86; %rotor bar view angle (rad)

%load data

alfa_r=2*pi/n;

if (fi<=0)
   y=0;
end;
if (fi>0)&(fi<=gama)
   y=(1/gama)*0.5*(1/3)*(fi)^3;
end;
if (fi>gama)&(fi<=alfa_r)
   y=r_iitf(gama)+0.5*gama*(fi-gama)+0.5*(fi-gama)^2;
end;
if (fi>alfa_r)&(fi<=gama+alfa_r)
   y=r_iitf(alfa_r)+(alfa_r-0.5*gama)*(fi-alfa_r)+0.5*(fi-
alfa_r)^2-(1/gama)*0.5*(1/3)*(fi-alfa_r)^3;
end;
if (fi>alfa_r+gama)
   y=r_iitf(alfa_r+gama)+(alfa_r)*(fi-alfa_r-gama);
end;
```

## L_sra1.M

```
function y=lsr_a1(x)
teta=mod(x,pi);
%parameters
rotor_rad=0.0820;
stack_length=0.1100;
```

```
g=0.0008;
u=4*pi*1e-7;
Ns=56;

l0=u*rotor_rad*stack_length/g;

%load data

y1=r_iitf(pi/6-teta)-r_iitf(-teta);
y2=r_iitf(pi/2+pi/6-teta)-r_iitf(pi/2-teta);
y3=r_iitf(pi+pi/6-teta)-r_iitf(pi-teta);
y4=r_iitf(3*pi/2+pi/6-teta)-r_iitf(3*pi/2-teta);
y=(-y1+y2-y3+y4)*(12/pi)*l0*Ns-Ns*l0*r_itf(2*pi-teta);
```

**dL_sra1.M**

```
function y=dlsr_a1(x)
teta=mod(x,pi);
rotor_rad=0.0820;
stack_length=0.1100;
g=0.0008;
u=4*pi*1e-7;
Ns=56;

l0=u*rotor_rad*stack_length/g;

%load data
y1=r_itf(pi/6-teta)-r_itf(-teta);
y2=r_itf(pi/2+pi/6-teta)-r_itf(pi/2-teta);
y3=r_itf(pi+pi/6-teta)-r_itf(pi-teta);
y4=r_itf(3*pi/2+pi/6-teta)-r_itf(3*pi/2-teta);
y=-(-y1+y2-y3+y4)*(12/pi)*l0*Ns+Ns*l0*r_tf(2*pi-teta);
```

**Ls.M**

```
function y=Ls(teta,es,ed);
teta=mod(teta,2*pi);
%parameters
rotor_rad=0.0820;
stack_length=0.1100;
g=0.0008;
u=4*pi*1e-7;
Ns=56;
lls=12.1930e-3;

l0=u*rotor_rad*stack_length/g;

%load data
Ls=zeros(3,3);
```

```
%eccentricity parameters
e=sqrt(es^2+ed^2+2*es*ed*cos(teta));
if e==0
A=1;
B=0;
C=0;
end
if e>0
A=1/sqrt(1-e^2);
B=2*A*(1-1/A)/e;
C=2*A*((1-1/A)/e)^2;
end
if es==0
alfa=teta;
end
if ed==0
alfa=0;
end
if (es>0)&(ed>0)
alfa=atan(ed*sin(teta)/(es+ed*cos(teta)));              end

for i=1:3
for j=1:3
M_P=A;
M_stf_i=A*Ns+6*Ns*C*inv(pi^2)*cos(2*alfa-2*(2*i-1)*pi/3);
M_stf_j=A*Ns+6*Ns*C*inv(pi^2)*cos(2*alfa-2*(2*j-1)*pi/3);
M_sstf_ii=16*(Ns^2)*inv(9)*A+12*C*inv(pi^2)*(Ns^2)*cos(2*alfa
-2*((2*i-1)*pi/3));
M_smtf_ij=2*(Ns^2)*A/3-6*(Ns^2)*C*inv(pi^2)*cos(2*alfa-
2*(i+j-1)*pi/3);

if (i==j)
y(i,j)=2*pi*l0*(M_sstf_ii-inv(M_P)*M_stf_i*M_stf_i)+lls; end
if not(i==j)
y(i,j)=2*pi*l0*(M_smtf_ij-inv(M_P)*M_stf_i*M_stf_j);     end
end
end
```

## Lr.M

```
function y=Lr(teta,es,ed);
teta=mod(teta,2*pi);
%parameters
rotor_rad=0.0820;
stack_length=0.1100;
g=0.0008;
u=4*pi*1e-7;
```

```
l_bar=95e-9;
l_end=18e-9;
n=40;
gama=pi/86;

%load data

l0=u*rotor_rad*stack_length/g;

y=zeros(n,n);

%eccentricity parameters
e=sqrt(es^2+ed^2+2*es*ed*cos(teta));
if e==0
A=1;
B=0;
C=0;
end
if e>0
A=1/sqrt(1-e^2);
B=2*A*(1-1/A)/e;
C=2*A*((1-1/A)/e)^2;
end
if es==0
alfa=teta;
end
if ed==0
alfa=0;
end
if (es>0)&(ed>0)
alfa=atan(ed*sin(teta)/(es+ed*cos(teta)));
end
alfa=alfa-teta;

for i=1:n
for j=1:n
M_P=A;
RT0i=1/n;
RT1i=2*inv(pi*gama)*sin(gama/2)*sin(pi/n)*cos(alfa-(i-
1)*2*pi/n-pi/n-gama/2);
RT2i=0.5*inv(pi*gama)*sin(gama)*sin(2*pi/n)*cos(2*(alfa-(i-
1)*2*pi/n-pi/n-gama/2));
M_rtf_i=A*RT0i+B*RT1i+C*RT2i;;
RT0j=1/n;
RT1j=2*inv(pi*gama)*sin(gama/2)*sin(pi/n)*cos(alfa-(j-
1)*2*pi/n-pi/n-gama/2);
RT2j=0.5*inv(pi*gama)*sin(gama)*sin(2*pi/n)*cos(2*(alfa-(j-
1)*2*pi/n-pi/n-gama/2));
```

```
M_rtf_j=A*RT0j+B*RT1j+C*RT2j;;
RS0i=1/n-gama/(6*pi);
RS1i=2*inv(pi*gama)*(cos(pi/n-gama/2)-
(2/gama)*cos(pi/n)*sin(gama/2))*cos(alfa-(i-1)*2*pi/n-pi/n-
gama/2);
RS2i=0.5*inv(pi*gama)*(cos(2*pi/n-gama)-
(1/gama)*cos(2*pi/n)*sin(gama))*cos(2*(alfa-(i-1)*2*pi/n-
pi/n-gama/2));
M_rstf_ii=A*RS0i+B*RS1i+C*RS2i;
RP0i=gama/(12*pi);
RP1i=inv(pi*gama)*((2/gama)*sin(gama/2)-
cos(gama/2))*cos(alfa-(i-1)*2*pi/n-gama/2);
RP2i=inv(4*pi*gama)*((1/gama)*sin(gama)-
cos(gama))*cos(2*(alfa-(i-1)*2*pi/n-gama/2));
M_rmtf_i=A*RP0i+B*RP1i+C*RP2i;
RP0j=gama/(12*pi);
RP1j=inv(pi*gama)*((2/gama)*sin(gama/2)-
cos(gama/2))*cos(alfa-(j-1)*2*pi/n-gama/2);
RP2j=inv(4*pi*gama)*((1/gama)*sin(gama)-
cos(gama))*cos(2*(alfa-(j-1)*2*pi/n-gama/2));
M_rmtf_j = A*RP0j+B*RP1j+C*RP2j;

if (i==j)
y(i,j)=2*pi*l0*(M_rstf_ii-
inv(M_P)*M_rtf_i*M_rtf_i)+2*l_bar+2*l_end;
end
if (j==i-1)
y(i,j)=2*pi*l0*(M_rmtf_i-inv(M_P)*M_rtf_i*M_rtf_j)-l_bar;
end
if (i==j-1)
y(i,j)=2*pi*l0*(M_rmtf_j-inv(M_P)*M_rtf_j*M_rtf_i)-l_bar;
end
if ((i>j+1)|(i<j-1))
y(i,j)=-2*pi*l0*inv(M_P)*M_rtf_j*M_rtf_i;
end
if (i==1)&(j==n)
y(i,j)=2*pi*l0*(M_rmtf_i-inv(M_P)*M_rtf_i*M_rtf_j)-l_bar;
end
if (i==n)&(j==1)
y(i,j)=2*pi*l0*(M_rmtf_j-inv(M_P)*M_rtf_j*M_rtf_i)-l_bar;
end
end
end
```

## Lsr.M

```
function y=Lsr(teta,es,ed);
teta=mod(teta,2*pi);
```

```
%parameters
rotor_rad=0.0820;
stack_length=0.1100;
g=0.0008;
u=4*pi*1e-7;
Ns=56;
gama=pi/86;
N=40;

%load data

l0=u*rotor_rad*stack_length/g;
y=zeros(3,n);

%eccentricity parameters
e=sqrt(es^2+ed^2+2*es*ed*cos(teta));
if e==0
A=1;
B=0;
C=0;
end
if e>0
A=1/sqrt(1-e^2);
B=2*A*(1-1/A)/e;
C=2*A*((1-1/A)/e)^2;
end
if es==0
alfa=teta;
end
if ed==0
alfa=0;
end
if (es>0)&(ed>0)
alfa=atan(ed*sin(teta)/(es+ed*cos(teta)));                end

for i=1:3
for j=1:n
M_P=A;
M_stf_i=A*Ns+6*Ns*C*inv(pi^2)*cos(2*alfa-2*(2*i-1)*pi/3);
RT0j=1/n;
RT1j=2*inv(pi*gama)*sin(gama/2)*sin(pi/n)*cos(alfa-teta-(j-
1)*2*pi/n-pi/n-gama/2);
RT2j=0.5*inv(pi*gama)*sin(gama)*sin(2*pi/n)*cos(2*(alfa-teta-
(j-1)*2*pi/n-pi/n-gama/2));
M_rtf_j=A*RT0j+B*RT1j+C*RT2j;
```

```
M_rs_ij=(l0*2*pi*Ns/n+Lsr_a1(teta-(i-1)*2*pi/3+(j-
1)*2*pi/n))*(A+B*cos(teta+(j-1)*2*pi/n-
alfa)+C*cos(2*teta+2*(j-1)*2*pi/n-2*alfa));
y(i,j)=M_rs_ij-2*pi*l0*M_stf_i*M_rtf_j/M_P;
end
end
```

## dLs.M

```
function y=dLs(teta,es,ed);
teta=mod(teta,2*pi);
%parameters
rotor_rad=0.0820;
stack_length=0.1100;
g=0.0008;
u=4*pi*1e-7;
Ns=56;
lls=12.1930e-3;

%load data

l0=u*rotor_rad*stack_length/g;

%eccentricity parameters
e=sqrt(es^2+ed^2+2*es*ed*cos(teta));
if e==0
A=1;
B=0;
C=0;
end
if e>0
A=1/sqrt(1-e^2);
B=2*A*(1-1/A)/e;
C=2*A*((1-1/A)/e)^2;
end
if es==0
alfa=teta;
end
if ed==0
alfa=0;
end
if (es>0)&(ed>0)
alfa=atan((ed*sin(teta))/(es+ed*cos(teta)));
end

if not(e==0)
de=-es*ed*sin(teta)/e;
end
```

```
if e==0
de=0;
end
if e==0
dA=0;dB=0;dC=0;
end
if (es==0)|(ed==0)
dA=0;dB=0;dC=0;
end
if (es>0)&(ed>0)
dA=e*de/((1-e^2)^(1.5));
dB=de*(-2*inv(e^2)*A+2*inv(e^2))+2*inv(e)*dA;
dC=2*inv(e)*dB-2*de*inv(e^2)*B-2*dA;
end
if es==0 dalfa=1; end
if ed==0 dalfa=0; end
if (es>0)&(ed>0) dalfa=(es*ed*cos(teta)+ed^2)/(e^2); end

for i=1:3
for j=1:3
M_P=A;
dM_P=dA;
M_stf_i=A*Ns+6*Ns*C*inv(pi^2)*cos(2*alfa-2*(2*i-1)*pi/3);
dM_stf_i=dA*Ns+6*Ns*dC*inv(pi^2)*cos(2*alfa-2*(2*i-1)*pi/3)-
6*Ns*C*inv(pi^2)*2*dalfa*sin(2*alfa-2*(2*i-1)*pi/3);
M_stf_j=A*Ns+6*Ns*C*inv(pi^2)*cos(2*alfa-2*(2*j-1)*pi/3);
dM_stf_j=dA*Ns+6*Ns*dC*inv(pi^2)*cos(2*alfa-2*(2*j-1)*pi/3)-
6*Ns*C*inv(pi^2)*2*dalfa*sin(2*alfa-2*(2*j-1)*pi/3);
M_sstf_ii=16*(Ns^2)*inv(9)*A
+12*C*inv(pi^2)*(Ns^2)*cos(2*alfa-2*((2*i-1)*pi/3));
dM_sstf_ii=16*(Ns^2)*inv(9)*dA+12*dC*inv(pi^2)*(Ns^2)*cos(2*a
lfa-2*((2*i-1)*pi/3))-
12*C*inv(pi^2)*(Ns^2)*2*dalfa*sin(2*alfa-2*((2*i-1)*pi/3));
M_smtf_ij=2*(Ns^2)*A/3-6*(Ns^2)*C*inv(pi^2)*cos(2*alfa-
2*(i+j-1)*pi/3);
dM_smtf_ij=2*(Ns^2)*dA/3-6*(Ns^2)*dC*inv(pi^2)*cos(2*alfa-
2*(i+j-1)*pi/3)+6*(Ns^2)*C*inv(pi^2)*2*dalfa*sin(2*alfa-
2*(i+j-1)*pi/3);

if (i==j)
y(i,i)=2*pi*l0*(dM_sstf_ii-(inv(M_P)*2*(dM_stf_i*M_stf_i)-
M_stf_i*M_stf_i*dM_P/(M_P^2)));
end
if not(i==j)
y(i,j)=2*pi*l0*(dM_smtf_ij-
(inv(M_P)*(dM_stf_i*M_stf_j+M_stf_i*dM_stf_j)-
M_stf_i*M_stf_j*dM_P/(M_P^2)));
end
```

```
end
end
```

**dLr.M**

```
function y=dLr(teta,es,ed);
teta=mod(teta,2*pi);
%parameters
rotor_rad=0.0820;
stack_length=0.1100;
g=0.0008;
u=4*pi*1e-7;
l_bar=95e-9;
l_end=18e-9;
n=40;
gama=pi/86;

%load data

l0=u*rotor_rad*stack_length/g;

y=zeros(n,n);

%eccentricity parameters
e=sqrt(es^2+ed^2+2*es*ed*cos(teta));
if e==0
A=1;
B=0;
C=0;
end
if e>0
A=1/sqrt(1-e^2);
B=2*A*(1-1/A)/e;
C=2*A*((1-1/A)/e)^2;
end
if es==0
alfa=teta;
end
if ed==0
alfa=0;
end
if (es>0)&(ed>0)
alfa=atan((ed*sin(teta))/(es+ed*cos(teta)));
end

if not(e==0)
de=-es*ed*sin(teta)/e;
end
```

```
if e==0
de=0;
end
if e==0
dA=0;dB=0;dC=0;
end
if (es==0)|(ed==0)
dA=0;dB=0;dC=0;
end
if (es>0)&(ed>0)
dA=e*de/((1-e^2)^(1.5));
dB=de*(-2*inv(e^2)*A+2*inv(e^2))+2*inv(e)*dA;
dC=2*inv(e)*dB-2*de*inv(e^2)*B-2*dA;
end
if es==0 dalfa=1; end
if ed==0 dalfa=0; end
if (es>0)&(ed>0) dalfa=(es*ed*cos(teta)+ed^2)/(e^2);
end

for i=1:n
for j=1:n
M_P=A;
dM_P=dA;
RT0i=1/n;
RT1i=2*inv(pi*gama)*sin(gama/2)*sin(pi/n)*cos(alfa-(i-
1)*2*pi/n-pi/n-gama/2);
dRT1i=-2*inv(pi*gama)*sin(gama/2)*sin(pi/n)*dalfa*sin(alfa
-(i-1)*2*pi/n-pi/n-gama/2);
RT2i =0.5*inv(pi*gama)*sin(gama)*sin(2*pi/n)*cos(2*(alfa-(i-
1)*2*pi/n-pi/n-gama/2));
dRT2i=-
0.5*inv(pi*gama)*sin(gama)*sin(2*pi/n)*2*dalfa*sin(2*(alfa-
(i-1)*2*pi/n-pi/n-gama/2));
M_rtf_i=A*RT0i+B*RT1i+C*RT2i;
dM_rtf_i=dA*RT0i+dB*RT1i+dC*RT2i+B*dRT1i+C*dRT2i;

RT0j=1/n;
RT1j=2*inv(pi*gama)*sin(gama/2)*sin(pi/n)*cos(alfa-(j-
1)*2*pi/n-pi/n-gama/2);
dRT1j=-2*inv(pi*gama)*sin(gama/2)*sin(pi/n)*dalfa*sin(alfa-
(j-1)*2*pi/n-pi/n-gama/2);
RT2j=0.5*inv(pi*gama)*sin(gama)*sin(2*pi/n)*cos(2*(alfa-(j-
1)*2*pi/n-pi/n-gama/2));
dRT2j=-
0.5*inv(pi*gama)*sin(gama)*sin(2*pi/n)*2*dalfa*sin(2*(alfa-
(j-1)*2*pi/n-pi/n-gama/2));
M_rtf_j=A*RT0j+B*RT1j+C*RT2j;
dM_rtf_j=dA*RT0j+dB*RT1j+dC*RT2j+B*dRT1j+C*dRT2j;
```

```
RS0i=1/n-gama/(6*pi);
RS1i=2*inv(pi*gama)*(cos(pi/n-gama/2)-
(2/gama)*cos(pi/n)*sin(gama/2))*cos(alfa-(i-1)*2*pi/n-pi/n-
gama/2);
dRS1i=-2*inv(pi*gama)*(cos(pi/n-gama/2)-
(2/gama)*cos(pi/n)*sin(gama/2))*dalfa*sin(alfa-(i-1)*2*pi/n-
pi/n-gama/2);
RS2i=0.5*inv(pi*gama)*(cos(2*pi/n-gama)-
(1/gama)*cos(2*pi/n)*sin(gama))*cos(2*(alfa-(i-1)*2*pi/n-
pi/n-gama/2));
dRS2i=-0.5*inv(pi*gama)*(cos(2*pi/n-gama)-
(1/gama)*cos(2*pi/n)*sin(gama))*2*dalfa*sin(2*(alfa-(i-
1)*2*pi/n-pi/n-gama/2));
M_rstf_ii=A*RS0i+B*RS1i+C*RS2i;
dM_rstf_ii=dA*RS0i+dB*RS1i+dC*RS2i+B*dRS1i+C*dRS2i;
RP0i=gama/(12*pi);
RP1i=inv(pi*gama)*((2/gama)*sin(gama/2)-
cos(gama/2))*cos(alfa-(i-1)*2*pi/n-gama/2);
dRP1i=-inv(pi*gama)*((2/gama)*sin(gama/2)-
cos(gama/2))*dalfa*sin(alfa-(i-1)*2*pi/n-gama/2);
RP2i=inv(4*pi*gama)*((1/gama)*sin(gama)-
cos(gama))*cos(2*(alfa-(i-1)*2*pi/n-gama/2));
dRP2i=-inv(4*pi*gama)*((1/gama)*sin(gama)-
cos(gama))*2*dalfa*sin(2*(alfa-(i-1)*2*pi/n-gama/2));
M_rmtf_i=A*RP0i+B*RP1i+C*RP2i;
dM_rmtf_i=dA*RP0i+dB*RP1i+dC*RP2i+B*dRP1i+C*dRP2i;
RP0j=gama/(12*pi);
RP1j=inv(pi*gama)*((2/gama)*sin(gama/2)-
cos(gama/2))*cos(alfa-(j-1)*2*pi/n-gama/2);
dRP1j=-inv(pi*gama)*((2/gama)*sin(gama/2)-
cos(gama/2))*dalfa*sin(alfa-(j-1)*2*pi/n-gama/2);
RP2j=inv(4*pi*gama)*((1/gama)*sin(gama)-
cos(gama))*cos(2*(alfa-(j-1)*2*pi/n-gama/2));
dRP2j=-inv(4*pi*gama)*((1/gama)*sin(gama)-
cos(gama))*2*dalfa*sin(2*(alfa-(j-1)*2*pi/n-gama/2));
M_rmtf_j=A*RP0j+B*RP1j+C*RP2j;
dM_rmtf_j=dA*RP0j+dB*RP1j+dC*RP2j+B*dRP1j+C*dRP2j;

if (i==j)
y(i,j)=2*pi*l0*(dM_rstf_ii-
inv(M_P)*2*dM_rtf_i*M_rtf_i+M_rtf_i*M_rtf_i*dM_P/(M_P^2))
end
if (j==i-1)
y(i,j)=2*pi*l0*(dM_rmtf_i-
inv(M_P)*(dM_rtf_i*M_rtf_j+M_rtf_i*dM_rtf_j)+M_rtf_i*M_rtf_j*
dM_P/(M_P^2));
end
if (i==j-1)
```

```
y(i,j)=2*pi*l0*(dM_rmtf_j-
inv(M_P)*(dM_rtf_i*M_rtf_j+M_rtf_i*dM_rtf_j)+M_rtf_i*M_rtf_j*
dM_P/(M_P^2));
end
if ((i>j+1)|(i<j-1)) y(i,j)=-
2*pi*l0*((inv(M_P)*(dM_rtf_i*M_rtf_j+M_rtf_i*dM_rtf_j)-
M_rtf_i*M_rtf_j*dM_P/(M_P^2)));
end
if (i==1)&(j==n)    y(i,j)=2*pi*l0*(dM_rmtf_i-
inv(M_P)*(dM_rtf_i*M_rtf_j+M_rtf_i*dM_rtf_j)+M_rtf_i*M_rtf_j*
dM_P/(M_P^2));
end
if (i==n)&(j==1)
y(i,j)=2*pi*l0*(dM_rmtf_j-
inv(M_P)*(dM_rtf_i*M_rtf_j+M_rtf_i*dM_rtf_j)+M_rtf_i*M_rtf_j*
dM_P/(M_P^2));
end
end
end
```

**dLsr.M**

```
function y=dLsr(teta,es,ed);
teta=mod(teta,2*pi);

%parameters
rotor_rad=0.0820;
stack_length=0.1100;
g=0.0008;
u=4*pi*1e-7;
Ns=56;
gama=pi/86;
N=40;

%load data

l0=u*rotor_rad*stack_length/g;
y=zeros(3,n);


%eccentricity parameters
e=sqrt(es^2+ed^2+2*es*ed*cos(teta));
if e==0
A=1;
B=0;
C=0;
end
if e>0
```

```
A=1/sqrt(1-e^2);
B=2*A*(1-1/A)/e;
C=2*A*((1-1/A)/e)^2;
end
if es==0
alfa=teta;
end
if ed==0
alfa=0;
end
if (es>0)&(ed>0)
alfa=atan((ed*sin(teta))/(es+ed*cos(teta)));
end

if not(e==0)
de=-es*ed*sin(teta)/e;
end
if e==0
de=0;
end
if e==0
dA=0;dB=0;dC=0;
end
if (es==0)|(ed==0)
dA=0;dB=0;dC=0;
end
if (es>0)&(ed>0)
dA=e*de/((1-e^2)^(1.5));
dB=de*(-2*inv(e^2)*A+2*inv(e^2))+2*inv(e)*dA;
dC=2*inv(e)*dB-2*de*inv(e^2)*B-2*dA;
end
if es==0 dalfa=1; end
if ed==0 dalfa=0; end
if (es>0)&(ed>0) dalfa=(es*ed*cos(teta)+ed^2)/(e^2); end

for i=1:3
for j=1:n
M_P=A;
dM_P=dA;
M_stf_i=A*Ns+6*Ns*C*inv(pi^2)*cos(2*alfa-2*(2*i-1)*pi/3);
dM_stf_i=dA*Ns+6*Ns*dC*inv(pi^2)*cos(2*alfa-2*(2*i-1)*pi/3)-
6*Ns*C*inv(pi^2)*2*dalfa*sin(2*alfa-2*(2*i-1)*pi/3);
RT0j=1/n;
RT1j=2*inv(pi*gama)*sin(gama/2)*sin(pi/n)*cos(alfa-teta-(j-
1)*2*pi/n-pi/n-gama/2);
dRT1j=-2*inv(pi*gama)*sin(gama/2)*sin(pi/n)*(dalfa-
1)*sin(alfa-teta-(j-1)*2*pi/n-pi/n-gama/2);
```

```
RT2j=0.5*inv(pi*gama)*sin(gama)*sin(2*pi/n)*cos(2*(alfa-teta-
(j-1)*2*pi/n-pi/n-gama/2));
dRT2j=-0.5*inv(pi*gama)*sin(gama)*sin(2*pi/n)*2*(dalfa-
1)*sin(2*(alfa-teta-(j-1)*2*pi/n-pi/n-gama/2));
M_rtf_j=A*RT0j+B*RT1j+C*RT2j;
dM_rtf_j=dA*RT0j+dB*RT1j+dC*RT2j+B*dRT1j+C*dRT2j;
M_rs_ij=(l0*2*pi*Ns/n+Lsr_a1(teta-(i-1)*2*pi/3+(j-
1)*2*pi/n))*(A+B*cos(teta+(j-1)*2*pi/n-
alfa)+C*cos(2*teta+2*(j-1)*2*pi/n-2*alfa));
dM_rs_ij=(dLsr_a1(teta-(i-1)*2*pi/3+(j-
1)*2*pi/n))*(A+B*cos(teta+(j-1)*2*pi/n-
alfa)+C*cos(2*teta+2*(j-1)*2*pi/n-
2*alfa))+(l0*2*pi*Ns/n+Lsr_a1(teta-(i-1)*2*pi/3+(j-
1)*2*pi/n))*(dA+dB*cos(teta+(j-1)*2*pi/n-alfa)-B*(1-
dalfa)*sin(teta+(j-1)*2*pi/n-alfa)+dC*cos(2*teta+2*(j-
1)*2*pi/n-2*alfa)-C*(2-2*dalfa)*sin(2*teta+2*(j-1)*2*pi/n-
2*alfa));
y(i,j)=dM_rs_ij-
2*pi*l0*(dM_stf_i*M_rtf_j+M_stf_i*dM_rtf_j)/M_P+2*pi*l0*M_stf
_i*M_rtf_j*dM_P/M_P^2;
end
end
```

## diffrential.M (eccentricity)

```
function dy=diffrential(t,y)
es=0.25;ed=0; %eccentricity parameters
%parameters
n = 40;
m = 3;
P = 2;
J=0.05;
Tl=1;
rst = 1.75;
r_bar=31e-6;
r_end=2.2e-6;
rotor_rad=0.0820;
stack_length=0.1100;
g=0.0008;
u=4*pi*1e-7;
Vs=380;
Ws=314;

% load data

l0=u*rotor_rad*stack_length/g;
```

```
Rst=diag([rst rst rst]);

%RESISTANSE MATRIX OF ROTOR  (Rrt)
Rrt=zeros(n,n);
for i=1:n
for j=1:n
if i==j
Rrt(i,j) = 2*(r_bar+r_end);
end
if (i==j+1)|(i==j-1)
Rrt(i,j) = -r_bar;
end
if (i==1)&(j==n)
Rrt(i,j) = -r_bar;
end
if (i==n)&(j==1)
Rrt(i,j) = -r_bar;
end
end
end

%MAIN INDUCTANCE OF MOTOR
M =  Lsr(y(m+n+2),es,ed);
dM= dLsr(y(m+n+2),es,ed);
Lrt=  Lr(y(m+n+2),es,ed);
dLrt=dLr(y(m+n+2),es,ed);
Lst=  Ls(y(m+n+2),es,ed);
dLst=dLs(y(m+n+2),es,ed);

L =[Lst  M ;M'  Lrt];
dL=[dLst dM;dM' dLrt];
D = [Rst+y(m+n+1)*dLst y(m+n+1)*dM;y(m+n+1)*dM'
Rrt+y(m+n+1)*dLrt];

home, t

Ds=2*pi/m;
V_line=zeros(m+n,1);
V_line(1,1)=Vs*sin(Ws*t);
V_line(2,1)=Vs*sin(Ws*t+Ds);
V_line(3,1)=Vs*sin(Ws*t-Ds);

dy=zeros(m+n+2,1);
y(1)=-y(2)-y(3);
dy(1:m+n) =inv(L)*(-D*y(1:m+n)+V_line);
dy(m+n+1) =(1/J)*(0.5*(y(1:m+n))'*dL*y(1:m+n)-Tl);
dy(m+n+2) =y(m+n+1);
```

**diffrential.M (broken bar:bar resistance)**

```
function dy=diffrential(t,y)
%parameters
n = 40;
m = 3;
P = 2;
J=0.05;
Tl=1;
rst = 1.75;
r_bar=31e-6;
r_end=2.2e-6;
rotor_rad=0.0820;
stack_length=0.1100;
g=0.0008;
u=4*pi*1e-7;
Vs=380;
Ws=314;

% load data

l0=u*rotor_rad*stack_length/g;

%RESISTANSE MATRIX OF STATOR (Rst)
Rst = diag([rst rst rst]);

%RESISTANSE MATRIX OF ROTOR  (Rrt)
Rrt   = zeros(n,n);
for i=1:n
for j=1:n
if i==j
Rrt(i,j) = 2*(r_bar+r_end);
end
if (i==j+1)|(i==j-1)
Rrt(i,j) = -r_bar;
end
if (i==1)&(j==n)
Rrt(i,j) = -r_bar;
end
if (i==n)&(j==1)
Rrt(i,j) = -r_bar;
end
end
end


%INDUCTANCE MATRIX OF STATOR (Lst)
lsst=(14*pi/9)*(Ns^2)*l0+llst;
```

```
lmst=(-2*pi/3)*(Ns^2)*l0;
Lst= diag([lsst-lmst lsst-lmst lsst-lmst]);


%INDUCTANCE MATRIX OF ROTOR  (Lrt)
lsrt=(1-1/n)*(2*pi/n)*l0-(1/3)*bar_ang*l0;
lsrt_b=(1-2/n)*(4*pi/n)*l0-(1/3)*bar_ang*l0;
lmrt1=-(2*pi/n^2)*l0+(bar_ang/6)*l0;
lmrt2=-(2*pi/n^2)*l0;
Lrt=ones(n,n)*lmrt2;
for i=1:n
for j=1:n
if i==j
Lrt(i,j) = 2*(l_bar+l_end)+lsrt;
end
if (i==j+1)|(i==j-1)
Lrt(i,j) = -l_bar+lmrt1;          end
if (i==1)&(j==n)
Lrt(i,j) = -l_bar+lmrt1;
end
if (i==n)&(j==1)
Lrt(i,j) = -l_bar+lmrt1;          end
end
end
Ds=2*pi/m;
V_line=zeros(m+n,1);
V_line(1,1)=Vs*cos(Ws*t);
V_line(2,1)=Vs*cos(Ws*t-Ds);
V_line(3,1)=Vs*cos(Ws*t+Ds);

M=Msr(y(m+n+2));
dM=dMsr(y(m+n+2));
[t y(m+n+1) y(m+n+2)]
L = [Lst M;M' Lrt];
D = [Rst y(m+n+1)*dM;y(m+n+1)*dM' Rrt];
dy=zeros(m+n+2,1);
dy(1:m+n) =inv(L)*(-D*y(1:m+n)+V_line);
dy(m+n+1) =(1/J)*((y(1:m))'*dM*y(m+1:m+n)-Tl);
dy(m+n+2) =y(m+n+1);
```

**diffrential.M (broken bar:loop omitting)**

```
function dy=diffrential(t,y)
%parameters
n = 40;
m = 3;
P = 2;
J=0.05;
```

```
Tl=1;
rst = 1.75;
r_bar=31e-6;
r_end=2.2e-6;
rotor_rad=0.0820;
stack_length=0.1100;
g=0.0008;
u=4*pi*1e-7;
Vs=380;
Ws=314;

% load data

l0=u*rotor_rad*stack_length/g;

br=0;
nn=n-br;

%RESISTANSE MATRIX OF STATOR (Rst)
Rst = diag([rst rst rst]);


%RESISTANSE MATRIX OF ROTOR  (Rrt)
Rrt   = zeros(n,n);
for i=1:n
for j=1:n
if i==j
Rrt(i,j) = 2*(r_bar+r_end);
end
if (i==j+1)|(i==j-1)
Rrt(i,j) = -r_bar;
end
if (i==1)&(j==n)
Rrt(i,j) = -r_bar;
end
if (i==n)&(j==1)
Rrt(i,j) = -r_bar;
end
end
end
Rrt(1,1)=2*r_bar+2*(br+1)*r_end;

%INDUCTANCE MATRIX OF STATOR (Lst)
llst=0.009;
Ns=56;
lsst=(14*pi/9)*(Ns^2)*l0+llst;
lmst=(-2*pi/3)*(Ns^2)*l0;
Lst=diag([lsst-lmst lsst-lmst lsst-lmst]);
```

```
% STEP 5: INDUCTANCE MATRIX OF ROTOR  (Lrt)
lsrt=(1-1/n)*(2*pi/n)*l0-(1/3)*bar_ang*l0;
lsrt_b=(1-2/n)*(4*pi/n)*l0-(1/3)*bar_ang*l0;
lmrt1=-(2*pi/n^2)*l0+(bar_ang/6)*l0;
lmrt2=-(2*pi/n^2)*l0;
Lrt=ones(n,n)*lmrt2;
for i=1:n
for j=1:n
if i==j
Lrt(i,j) = 2*(l_bar+l_end)+lsrt;
end
if (i==j+1)|(i==j-1)
Lrt(i,j) = -l_bar+lmrt1;          end
if (i==1)&(j==n)
Lrt(i,j) = -l_bar+lmrt1;
end
if (i==n)&(j==1)
Lrt(i,j) = -l_bar+lmrt1;          end
end
end
Lrt(1,1)=(1-(br+1)/n)*(2*(br+1)*pi/n)*l0-
(1/3)*bar_ang*l0+2*(l_bar+(br+1)*l_end);
Lrt(1,2)=-(2*pi*(br+1)/n^2)*l0+(bar_ang/6)*l0-l_bar;
Lrt(2,1)=-(2*pi*(br+1)/n^2)*l0+(bar_ang/6)*l0-l_bar;
Lrt(1,nn)=-(2*pi*(br+1)/n^2)*l0+(bar_ang/6)*l0-l_bar;
Lrt(nn,1)=-(2*pi*(br+1)/n^2)*l0+(bar_ang/6)*l0-l_bar;
Lrt(1,br+1:nn-1)=((br+1)*lmrt2)*ones(1,nn-br-1);
Lrt(br+1:nn-1,1)=((br+1)*lmrt2)*ones(nn-br-1,1);

Ds=2*pi/m;
V_line=zeros(m+nn,1);
V_line(1,1)=Vs*cos(Ws*t);
V_line(2,1)=Vs*cos(Ws*t-Ds);
V_line(3,1)=Vs*cos(Ws*t+Ds);

M=Msr(y(m+nn+2));
dM=dMsr(y(m+nn+2));

MM=zeros(m,nn);
dMM=zeros(m,nn);

for i=1:br+1
MM(:,1)=MM(:,1)+M(:,i);
dMM(:,1)=dMM(:,1)+dM(:,i);
end
dMM(:,2:nn)=dM(:,br+2:40);
MM(:,2:nn)=M(:,br+2:40);
```

```
clear dM M;
M=MM;
dM=dMM;

[t y(m+nn+1)]
L = [Lst M;M' Lrt];
D = [Rst y(m+nn+1)*dM;y(m+nn+1)*dM' Rrt];

dy=zeros(m+nn+2,1);
dy(1:m+nn)   =inv(L)*(-D*y(1:m+nn)+V_line);
dy(m+nn+1)   =(1/J)*((y(1:m))'*dM*y(m+1:m+nn)-Tl);
dy(m+nn+2)   =y(m+nn+1);
```

# ضمیمه سه

# مقالات استخراج شده از پایان نامه

**A Precise Evaluation of Inductances of a Squirrel-Cage Induction Motor under Mixed Eccentric Conditions**

J. Faiz, *Senior Member, IEEE,* I. Tabatabaei, H. A. Toliyat, *Senior Member, IEEE*

*Abstract-* **This paper presents a more precise model for computation of three-phase squirrel cage induction machine inductances under different eccentric conditions. Generally, available techniques are based on the winding function theory and simplification and geometrical approximation of unsymmetrical models of the motor under mixed eccentricities. This paper determines a precise geometrical model under the mixed eccentricity conditions and evaluates the inductances. Meanwhile, the evaluated inductances are compared to those calculated using different approximate geometrical models and the best approximation is recommended for a geometrical modeling of induction motor under eccentricity conditions.**

# Extension of Winding Function Theory for Non-uniform Air-gap of Electric Machinery

J. Faiz, *Senior Member, IEEE,* I. Tabatabaei

**Abstract-This paper extends the winding function theory for non-uniform air-gap of rotating electric machinery. It is shown that the winding function for non-uniform air-gap differs with that used in the symmetrical case, while several papers employ the uniform air-gap winding function in order to study the electric motors performance in the fault conditions. It is particularly the case in the study of squirrel-cage induction motor in the non-uniform air-gap such as eccentricity of the rotor and stator.**

محاسبه تاثیر شرایط مختلف ناهم محوری روتور و استاتور

# بر روی اندوکتانسهای موتور القایی قفس سنجابی

جواد فیض** – ایمان طباطبایی*و** – عماد شریفی قزوینی*

پژوهشگاه نیرو* – دانشگاه تهران**

چکیده : در این مقاله پس از تشریح عملکرد موتورهای القایی قفس سنجابی در شرایط مختلف ناهم محوری روتور و استاتور به بررسی اثر این شرایط بر روی اندوکتانسهای موتور پرداخته می‌شود، سپس با در نظر گرفتن شرط ناهم‌محوری مرکب یک روش شبه تحلیلی جهت محاسبه اندوکتانسهای موتور القایی قفس سنجابی ارائه می‌شود. از آنجا که مدلهای مختلف موتور القایی در شرایط ناهم محوری به دلیل عدم تقارن هندسی بسیار پیچیده است، تحلیل کامپیوتری آن طولانی است و استفاده از روش شبه تحلیلی این مقاله باعث تسریع تحلیل مدل کامپیوتری موتور القایی در شرایط ناهم محوری خواهدشد. در انتها با استفاده از روش پیشنهاد شده به تحلیل کامپیوتری یک موتور القایی در حالت ناهم محوری مرکب پرداخته می‌شود.

# بررسی عملکرد موتور القائی سه فاز روتور قفسه‌ای تحت شرایط ناهم محوری مرکب

جواد فیض[**] – ایمان طباطبایی[*و**] – عماد شریفی قزوینی[*] – همایون مشکین کلک[*]

پژوهشگاه نیرو[*] – دانشگاه تهران[**]

**چکیده:** این مقاله یک روش جدید به منظور تحلیل عملکرد موتور القائی در حالت ناهم‌محوری مرکب میان روتور و استاتور ارائه می‌دهد، سپس به بررسی تاثیر این حالت بر روی طیف فرکانسی جریان سیم‌پیچی استاتور می‌پردازد. روش مذکور نسبت به سایر روشها دقیق‌تر و حجم کمتری از حافظه کامپیوتر را اشغال می‌کند. این روش قادر به در نظر گرفتن بسیاری از هارمونیکهای فضایی ناشی از شکل فاصله هوائی و سیمبندی فازهای استاتور می‌باشد که این امر در تحلیل عملکرد موتور در حالتهای مختلف ناهم محوری میان روتور و استاتور بسیار حائز اهمیت است. در این بررسی نشان داده می‌شود که این نوع ناهم محوری علاوه بر ایجاد دو هارمونیک اضافی در اطراف فرکانس اصلی موجب ایجاد دو هارمونیک اضافی در اطراف هارمونیک شیارگذاری روتور نیز می‌شود. لازم به ذکر است که در بررسی حاضر، شیارهای مورب روتور نیز در نظر گرفته شده است.